\documentclass{conm-p-l}
\usepackage{%
amssymb,%
righttag,%
pb-diagram,%
lamsarrow,%
pb-lams%
}

\usepackage[mathcal]{euscript}

\newcommand{\OO}        {{\mathcal{O}}}
\newcommand{\OW}        {{\mathcal{O}_W}}
\newcommand{\OX}        {{\mathcal{O}_X}}
\newcommand{\OXv}       {{\mathcal{O}_X^\vee}}
\newcommand{\OY}        {{\mathcal{O}_Y}}
\newcommand{\OZ}        {{\mathcal{O}_Z}}
\newcommand{\OG}        {{\mathcal{O}_G}}
\newcommand{\OGt}       {{\mathcal{O}_G^\times}}
\newcommand{\OC}        {{\mathcal{O}_C}}
\newcommand{\OCt}       {{\mathcal{O}_C^\times}}

\newcommand{\CC}        {{\mathcal{C}}}
\newcommand{\CD}        {{\mathcal{D}}}
\newcommand{\CE}        {{\mathcal{E}}}
\newcommand{\CF}        {{\mathcal{F}}}
\newcommand{\CI}        {{\mathcal{I}}}
\newcommand{\CJ}        {{\mathcal{J}}}
\newcommand{\CK}        {{\mathcal{K}}}
\newcommand{\CM}        {{\mathcal{M}}}
\newcommand{\CT}        {{\mathcal{T}}}
\newcommand{\CX}        {{\mathcal{X}}}
\newcommand{\FX}        {{\widehat{\mathcal{X}}}}

\newcommand{\BD}        {{\mathbb{D}}}
\newcommand{\BL}        {{\mathbb{L}}}
\newcommand{\BP}        {{\mathbb{P}}}

\newcommand{\era}       {\twoheadrightarrow}
\newcommand{\mra}       {\rightarrowtail}

\newcommand{\al}        {\alpha}
\newcommand{\bt}        {\beta}
\newcommand{\gm}        {\gamma}
\newcommand{\dl}        {\delta}
\newcommand{\ep}        {\epsilon}
\newcommand{\kp}        {\kappa}
\newcommand{\lm}        {\lambda}
\newcommand{\psib}      {\overline{\psi}}
\newcommand{\sg}        {\sigma}
\newcommand{\tht}       {\theta}
\newcommand{\om}        {\omega}

\newcommand{\Gm}        {\Gamma}
\newcommand{\Sg}        {\Sigma}
\newcommand{\Sgi}       {\Sigma^\infty}
\newcommand{\Lm}        {\Lambda}
\newcommand{\Om}        {\Omega}
\newcommand{\Omi}       {\Omega^\infty}

\newcommand{\rat}       {{\mathbb{Q}}}
\newcommand{\Zh}        {{\mathbb{Z}}}

\newcommand{\nat}       {{\mathbb{N}}}              
\newcommand{\Fp}        {{\mathbb{F}_p}}
\newcommand{\cplx}      {{\mathbb{C}}}

\newcommand{\aff}       {{\mathbb{A}}}
\newcommand{\haf}       {{\widehat{\mathbb{A}}}}
\newcommand{\ot}        {\otimes}
\newcommand{\hot}       {\widehat{\otimes}}
\newcommand{\st}        {\;|\;}
\newcommand{\pri}       {{\mathfrak{p}}}
\newcommand{\qri}       {{\mathfrak{q}}}
\newcommand{\rri}       {{\mathfrak{r}}}

\newcommand{\un}[1]     {\underline{#1}}
\newcommand{\ov}[1]     {\overline{#1}}

\newcommand{\colim}{\operatornamewithlimits{\underset{\longrightarrow}{lim}}}

\newcommand{\holim}{\operatornamewithlimits{\underset{\longleftarrow}{holim}}}
\newcommand{\invlim}{\operatornamewithlimits{\underset{\longleftarrow}{lim}}}

\newcommand{\fpl}{[\![} 
\newcommand{\fpr}{]\!]} 
\newcommand{\fps}[2]{{#1 \fpl #2 \fpr}} 
\newcommand{\bcf}[2]    
 {\left(\begin{array}{c}{#1}\\{#2}\end{array}\right)}

\newcommand{\Ab}        {\operatorname{Ab}}
\newcommand{\Alg}       {\operatorname{Alg}}
\newcommand{\Based}     {\operatorname{Based}}
\newcommand{\Bases}     {\operatorname{Bases}}
\newcommand{\Coord}     {\operatorname{Coord}}
\newcommand{\CW}        {{\operatorname{CW}}}
\newcommand{\Div}       {\operatorname{Div}}
\newcommand{\End}       {\operatorname{End}}
\newcommand{\Exp}       {\operatorname{Exp}}
\newcommand{\FGL}       {{\operatorname{FGL}}}
\newcommand{\FI}        {{\operatorname{FI}}}
\newcommand{\GL}        {\operatorname{GL}}
\newcommand{\Hom}       {\operatorname{Hom}}
\newcommand{\Ind}       {\operatorname{Ind}}
\newcommand{\IPS}       {{\operatorname{IPS}}}
\newcommand{\Iso}       {\operatorname{Iso}}
\newcommand{\Map}       {\operatorname{Map}}
\newcommand{\Mer}       {\operatorname{Mer}}
\newcommand{\Mod}       {\operatorname{Mod}}
\newcommand{\Mon}       {\operatorname{Mon}}
\newcommand{\Mult}      {\operatorname{Mult}}
\newcommand{\Nil}       {\operatorname{Nil}}
\newcommand{\Points}    {\operatorname{Points}}
\newcommand{\Prim}      {\operatorname{Prim}}
\newcommand{\Pro}       {\operatorname{Pro}}
\newcommand{\Sym}       {\operatorname{Sym}}
\newcommand{\Tor}       {\operatorname{Tor}}
\newcommand{\EPR}       {\operatorname{EPR}}
\newcommand{\OFG}       {\operatorname{OFG}}
\newcommand{\LEPR}      {\operatorname{LEPR}}
\newcommand{\LOFG}      {\operatorname{LOFG}}

\newcommand{\ann}       {{\operatorname{ann}}}
\newcommand{\ev}        {{\operatorname{ev}}}
\newcommand{\image}     {{\operatorname{image}}}
\newcommand{\odd}       {{\operatorname{odd}}}
\newcommand{\op}        {{\operatorname{op}}}
\newcommand{\red}       {{\operatorname{red}}}
\newcommand{\res}       {{\operatorname{res}}}
\newcommand{\zar}       {{\operatorname{zar}}}
\newcommand{\sol}       {{\operatorname{sol}}}
\newcommand{\spec}      {{\operatorname{spec}}}
\newcommand{\sch}       {{\operatorname{sch}}}
\newcommand{\spf}       {{\operatorname{spf}}}
\newcommand{\trace}     {{\operatorname{trace}}}

\newcommand{\Rings}     {{\operatorname{Rings}}}
\newcommand{\FRings}    {{\operatorname{FRings}}}
\newcommand{\LRings}    {{\operatorname{LRings}}}
\newcommand{\Sets}      {{\operatorname{Sets}}}
\newcommand{\Sheaves}   {{\operatorname{Sheaves}}}

\newcommand{\tX}        {\times_X}
\newcommand{\tY}        {\times_Y}
\newcommand{\tZ}        {\times_Z}

\newcommand{\tN}        {\tilde{N}}
\newcommand{\hE}        {\widehat{E}}
\newcommand{\sse}       {\subseteq}
\newcommand{\Smash}     {\wedge}

\newcommand{\bigWedge}  {\bigvee}

\newcommand{\cp}        {\mathbb{C}P}
\newcommand{\cpi}       {{\mathbb{C}P^\infty}}
\newcommand{\half}      {{\textstyle\frac{1}{2}}}
\newcommand{\hR}        {{\widehat{R}}}
\newcommand{\arw}       {\arrow}
\newcommand{\xla}       {\xleftarrow}
\newcommand{\xra}       {\xrightarrow}
\newcommand{\iffa}      {\Leftrightarrow}

\newcommand{\tm}        {\times}
\newcommand{\MG}        {\mathbb{G}_m}
\newcommand{\AFG}       {\widehat{\mathbb{G}}_a}
\newcommand{\ImJ}       {\operatorname{ImJ}}
\newcommand{\ZtBU}      {\Zh\tm BU}
\renewcommand{\div}     {\text{div}}

\renewcommand{\:}{\colon}

\newcommand{\dfn}[1]{{\emph{#1}}\index{#1}}
\newcommand{\idx}[1]{#1\index{#1}}

\newenvironment{diag}{
 \renewcommand{\typeout}[1]{}
 \begin{displaymath}
 \begin{diagram}}{
 \end{diagram}
 \end{displaymath}} 

\makeatletter
\@namedef{dgo@=>}{\let\dg@VECTOR=\dg@forkvector}%

\def\dg@forkvector(#1,#2)#3{%
   \begingroup
   \dg@XTEMP=#1\relax \dg@YTEMP=#2\relax
   \dg@XEND=1
   \multiply\dg@XEND\dg@XTEMP 
   \multiply\dg@XEND  50 \relax 
   \dg@WEND=1
   \multiply\dg@WEND\dg@YTEMP 
   \multiply\dg@WEND -50 \relax
   \begin{picture}(0,0)
      \put(\dg@WEND,\dg@XEND){\lamsvector(\dg@XTEMP,\dg@YTEMP){#3}}
      \multiply\dg@XEND -1 \relax
      \multiply\dg@WEND -1 \relax
      \put(\dg@WEND,\dg@XEND){\lamsvector(\dg@XTEMP,\dg@YTEMP){#3}}
   \end{picture}%
   \endgroup}%
\makeatother

\theoremstyle{definition}
\newtheorem{theorem}{Theorem}[section]

\newtheorem{corollary}  [theorem]{Corollary}
\newtheorem{lemma}      [theorem]{Lemma}

\newtheorem{proposition}[theorem]{Proposition}
\newtheorem{definition} [theorem]{Definition}
\newtheorem{remark}     [theorem]{Remark}
\newtheorem{example}    [theorem]{Example}

\makeindex
\begin{document}
\title{Formal Schemes and Formal Groups}
\author{Neil P.~Strickland}
\date{\today}
\bibliographystyle{abbrv}

\maketitle 
\tableofcontents

\section{Introduction}
\label{sec-intro}

In this paper we set up a framework for using algebraic geometry to
study the generalised cohomology rings that occur in algebraic
topology.  This idea was probably first introduced by
Quillen~\cite{qu:sec} and it implicitly or explicitly underlies much
of our understanding of complex oriented cohomology theories,
exemplified by the work of Morava.  Most of the results presented here
have close and well-known analogues in the algebro-geometric
literature, but with different definitions or technical assumptions
that are often inconvenient for topological applications.  Our aim
here is merely to put everything together in a systematic way that
naturally incorporates the phenomena that we see in topology while
discarding complications that never arise there.  In more detail, in
the classical situation one is often content to deal with finite
dimensional, Noetherian schemes.  Nilpotents are seen as a somewhat
peripheral phenomenon, and formal schemes are only introduced at a
late stage in the exposition.  Schemes are defined as spaces with
extra structure.  The idea of a scheme as a functor occurs in advanced
work (a nice example is~\cite{kama:ame}) but is usually absent from
introductory treatments.  For us, however, it is definitely most
natural to think of schemes as functors.  Our schemes are very often
not Noetherian or finite dimensional, and nilpotents are of crucial
importance.  We make heavy use of formal schemes, and we need to
define these in a more general way than is traditional.  On the other
hand, we can get a long way using only affine schemes, whereas the
usual treatment devotes a great deal of attention to the non-affine
case.

Section~\ref{sec-schemes} is an exposition of the basic facts of
algebraic geometry that is well adapted to the viewpoint discussed
above, together with a number of useful examples.

In Section~\ref{sec-non-affine}, we give a basic account of non-affine
schemes from our point of view.

In Section~\ref{sec-formal-schemes}, we give a very general definition
of formal schemes which follows naturally from our description of
ordinary (or ``informal'') schemes.  We then work out the basic
properties of the category of formal schemes, such as the existence of
limits and colimits and the behaviour of regular monomorphisms (or
``closed inclusions'').

In Section~\ref{sec-formal-groups}, we discuss the Abelian monoid and
group objects in the category of formal schemes.  We then specialise
in Section~\ref{sec-ordinary-groups} to the case of smooth,
commutative, one-dimensional formal groups, which we call ``ordinary
formal groups''.  

Finally, in Section~\ref{sec-topology}, we construct functors from the
homotopy category of spaces (or suitable subcategories) to the
category of formal schemes.  We use the work of Ravenel, Wilson and
Yagita~\cite{rawiya:bpc} to show that spaces whose Morava $K$-theory
is concentrated in even degrees give formal schemes with good
technical properties.  We also discuss what happens to a number of
popular spaces under our functors.  Further applications of this point
of view appear in~\cite{st:fsf,st:msg,grst:vlc,hoanst:esw} and a
number of other papers in preparation.

\subsection{Notation and conventions}
\label{subsec-notation}

We write $\Rings$\index{Rings@$\Rings$} for the category of rings (by
which we always mean commutative unital rings) and
$\Sets$\index{Sets@$\Sets$} for the category of sets.  For any ring
$R$, we write $\Mod_R$ for the category of $R$-modules, and $\Alg_R$
for the category of $R$-algebras.  Given a category $\CC$, we usually
write $\CC(X,Y)$ for the set of $\CC$-morphisms from $X$ to $Y$.  We
write $\CC_X$ for the category of objects of $\CC$ over $X$.  More
precisely, on object of $\CC_X$ is a pair $(Y,u)$ where $u\:Y\xra{}Z$,
and $\CC_X((Y,u),(Z,v))$ is the set of maps $f\:Y\xra{}Z$ in $\CC$
such that $vf=u$.

We write $\CF$\index{F@$\CF$} for the category of all functors
$\Rings\xra{}\Sets$.  

\subsection{Even periodic ring spectra}
\label{subsec-even}

We now give a basic topological definition, as background for some
motivating remarks to be made in subsequent sections.  Details of
topological applications will appear in Section~\ref{sec-topology}.
The definition below will be slightly generalised there, to deal with
unpleasantness at the prime $2$.

\begin{definition}\label{defn-intro-even-periodic}
 An \dfn{even periodic ring spectrum} is a commutative and associative
 ring spectrum $E$ such that
 \begin{enumerate}
 \item $\pi_1E=0$
 \item $\pi_2E$ contains a unit.
 \end{enumerate}
\end{definition}

The example to bear in mind is the complex $K$-theory spectrum $KU$.
Suitable versions of Morava $E$-theory and $K$-theory are also
examples, as are periodised versions of $MU$ and $H$; we write $MP$
and $HP$ for these.  See Section~\ref{sec-topology} for more details.

\section{Schemes}
\label{sec-schemes}

In this section we set up the basic categorical apparatus of schemes.
We then discuss limits and colimits of schemes, and various kinds of
subschemes.  We compare our functorial approach with more classical
accounts by discussing the Zariski space of a scheme.  We then discuss
various issues about nilpotent and idempotent functions.  We define
sheaves over functors, and show that our definition works as expected
for schemes.  We then define flatness and faithful flatness for maps
of schemes, and prove descent theorems for schemes and sheaves over
faithfully flat maps.  Finally, we address the question of defining a
``scheme of maps'' $\Map(X,Y)$ between two given schemes $X$ and $Y$.

\begin{definition}
 An \emph{affine scheme}\index{scheme!affine} is a covariant
 representable functor 
 \[ X\:\Rings\xra{}\Sets. \]
 We make little use of non-affine schemes, so we shall generally omit
 the word ``affine''.  A map of schemes is just a natural
 transformation.  We write $\CX$\index{X@$\CX$} for the category of
 schemes, which is a full subcategory of $\CF$.  We write
 $\spec(A)$\index{spec@$\spec$} for the functor represented by $A$, so
 $\spec(A)(R)=\Rings(A,R)$ and $\spec(A)$ is a scheme.
\end{definition}
\begin{remark}
 If $E$ is an even periodic ring spectrum and $Z$ is a finite spectrum
 we define $Z_E=\spec(E^0Z)$.  This gives a covariant functor
 $Z\mapsto Z_E$ from finite complexes to schemes.  We also write
 $S_E=\spec(E^0)$. 
\end{remark}

\begin{definition}
 We write $\aff^1$\index{AA1@$\aff^1$} for the forgetful functor
 $\Rings\xra{}\Sets$.  This is isomorphic to $\spec(\Zh[t])$ and thus
 is a scheme.  Given any functor $X\in\CF$, we write
 $\OX$\index{OX@$\OX$} for the set of natural maps $X\xra{}\aff^1$.
 (This can actually be a proper class for general $X$, but it will
 always be a set in the cases that we consider.)  Note that $\OX$ is a
 ring under pointwise operations.
\end{definition}

Our category of schemes is equivalent to the algebraic geometer's
category of affine schemes, which in turn is equivalent (by Yoneda's
lemma) to the opposite of the category of rings.  

We now describe the duality between schemes and rings in more detail.
The Yoneda lemma tells us that $\OO_{\spec(A)}$ is naturally
isomorphic to $A$.  For any functor $X\in\CF$ we have a tautological
map $\kp\:X\xra{}\spec(\OX)$\index{k@$\kappa$}.  To define $\kp$
explicitly, suppose we have a ring $R$ and an element $x\in X(R)$; we
need to produce a map $\kp_R(x)\:\OX\xra{}R$.  An element $f\in\OX$ is
a natural map $f\:X\xra{}\aff^1$, so it has a component
$f_R\:X(R)\xra{}R$, and we can define $\kp_R(x)(f)=f_R(x)$.  If
$X=\spec(A)$ then $\kp$ is easily seen to be bijective.  As schemes
are by definition representable, any scheme $X$ is equivalent to
$\spec(A)$ for some $A$, so we see that the map $X\xra{}\spec(\OX)$ is
always an isomorphism.  Thus, the functor $X\xra{}\OX$ is inverse to
the functor $\spec\:\Rings^{\op}\xra{}\CX$.

We next give some examples of schemes.
\begin{example}\label{eg-mult-group}
 A basic example is the ``multiplicative group'' \index{multiplicative
   group} $\MG$\index{Gm@$\MG$}, which is defined by
 \[ \MG(R) = R^{\tm} = \text{ the group of units of } R. \]
 This is a scheme because it is represented by $\Zh[x^{\pm 1}]$.
\end{example}

\begin{example}\label{eg-affine-space}
 The affine $n$-space $\aff^n$\index{AAn@$\aff^n$} is defined by
 $\aff^n(R)=R^n$.  This is a scheme because it is represented by
 $\Zh[x_1,\ldots,x_n]$.  If $f_1,\ldots,f_m$ are polynomials in $n$
 variables over $\Zh$ then there is an obvious natural map
 $R^m\xra{}R^n$ for all rings $R$, which sends
 $\un{a}=(a_1,\ldots,a_m)$ to $(f_1(\un{a}),\ldots,f_n(\un{a}))$.
 Thus, this gives a map $\aff^m\xra{}\aff^n$ of schemes.  These are in
 fact all the maps between these schemes.  The key point is of course
 that the set of ring maps
 $\Zh[y_1,\ldots,y_m]\xla{}\Zh[x_1,\ldots,x_n]$ bijects naturally with
 the set of such tuples $(f_1,\ldots,f_m)$.  It is a good exercise to
 work out all of the identifications going on here.
\end{example}

We next define the scheme $\FGL$\index{FGL@$\FGL$} of formal group
laws, which will play a central r\^{o}le in the applications of
schemes to algebraic topology.
\begin{example}\label{eg-FGL}
 A \dfn{formal group law} over a ring $R$ is a formal power series
 \[ F(x,y) = \sum_{k,l\ge 0} a_{kl} x^k y^l \in\fps{R}{x,y} \]
 satisfying
 \begin{align*}
  F(x,0)        & =  x                         \\
  F(x,y)        & =  F(y,x)                    \\
  F(F(x,y),z)   & =  F(x,F(y,z)).
 \end{align*}
 We can define a scheme $\FGL$ as follows:
 \[ \FGL(R) = \{ \text{ formal group laws over } R \}. \]
 To see that $\FGL$ is a scheme, we consider the ring
 $L_0=\Zh[a_{kl}\st k,l> 0 ]$ and the formal power series
 $F_0(x,y)=x+y+\sum a_{kl}x^k y^l\in\fps{L_0}{x,y}$.  We then let $I$
 be the ideal in $L_0$ generated by the coefficients of the power
 series $F_0(x,y)-F_0(y,x)$ and $F_0(F_0(x,y),z)-F_0(x,F_0(y,z))$.
 Finally, set $L=L_0/I$.  It is easy to see that $\FGL=\spec(L)$.
 The ring $L$\index{L@$L$} is called the \dfn{Lazard ring}.  It is a
 polynomial ring in countably many variables; there is a nice
 exposition of the proof in~\cite[Part II]{ad:shg}.  Recall that $MP$
 denotes the $2$-periodic version of $MU$; a fundamental theorem of
 Quillen~\cite{qu:fgl,qu:epc} (also proved in~\cite{ad:shg})
 identifies the scheme $S_{MP}:=\spec(MP^0)$ with $\FGL$. 
\end{example}

\begin{example}\label{eg-limits}
 Given any diagram of schemes $\{X_i\}$, we claim that the functor
 $X=\invlim_iX_i$ (which is defined by
 $(\invlim_iX_i)(R)=\invlim_i(X_i(R))$) is also a scheme.  Indeed,
 suppose that $X_i=\spec(A_i)$.  As $\spec\:\Rings^\op\xra{}\CX$ is
 an equivalence, we get a diagram of rings $A_i$ with arrows
 reversed.  It is well-known that the category of rings has colimits,
 and it is clear that $X=\spec(\colim_iA_i)$.

 In particular, if $X$ and $Y$ are schemes, we have a scheme $X\tm Y$
 with $(X\tm Y)(R)=X(R)\tm Y(R)$ and $\OO_{X\tm Y}=\OX\ot\OY$
 (because coproducts of rings are tensor products).  Similarly, if we
 have maps $X\xra{f}Z\xla{g}Y$ then we can form the pullback 
 \[ (X\tZ Y)(R)=X(R)\tm_{Z(R)}Y(R)=
     \{(x,y)\in X(R)\tm Y(R)\st f(x)=g(y)\}.
 \]
 This is represented by the tensor product $\OX\ot_{\OZ}\OY$.

 We write $1$ for any one-point set, and also for the constant functor
 $1(R)=1$.  Thus $1=\spec(\Zh)$, and this is the terminal object in
 $\CX$ or $\CF$.
\end{example}
\begin{example}
 Let $Z$ and $W$ be finite CW complexes, and let $E$ be an even
 periodic ring spectrum.  There is a natural map
 $(Z\tm W)_E\xra{}Z_E\tm_{S_E}W_E$.  This will be an isomorphism if
 $E^1Z=0=E^1W$ and we have a K\"unneth isomorphism
 $E^*(Z\tm W)=E^*(Z)\ot_{E^*}E^*(W)$.  This holds in particular if
 $H_*Z$ is a free Abelian group, concentrated in even degrees.
\end{example}
\begin{example}\label{eg-IPS}
 An \dfn{invertible power series} over a ring $R$ is a formal power
 series $f\in\fps{R}{x}$ such that $f(x)=wx+O(x^2)$ for some
 $w\in R^\tm$.  This implies, of course, that $f$ has a
 composition-inverse $g=f^{-1}$, so that $f(g(x))=x=g(f(x))$.  We
 write $\IPS(S)$\index{IPS@$\IPS$} for the set of such $f$, which is
 easily seen to be a scheme.  It is actually a group scheme, in that
 $\IPS(R)$ is a group (under composition), functorially in $R$.

 The group $\IPS$ acts on $\FGL$\index{FGL@$\FGL$} by
 \[ (f,F) \mapsto F_f \qquad\qquad F_f(x,y)=f(F(f^{-1}x,f^{-1}y)).
 \]

 An \dfn{isomorphism} between formal group laws $F$ and $G$ is an
 invertible series $f$ such that $f(F(a,b))=G(f(a),f(b))$.
 Let $\FI$\index{FI@$\FI$} be the following scheme:
 \[ \FI(R) = 
      \{ (F,f,G) \st F,G\in\FGL(R) \text{ and }
         f\: F\xra{}G\text{ is an isomorphism }\}.
 \]
 There is an evident composition map
 \[ \FI \tm_{\FGL} \FI \xra{} \FI \qquad\qquad
    ((F,f,G),(G,g,H)) \mapsto (F,gf,H).
 \]
 Moreover, there is an isomorphism
 \[ \IPS\tm\FGL \xra{} \FI \qquad\qquad (F,f) \mapsto (F,f,F_f).
 \]
 One can describe these maps by giving implicit formulae in the
 representing rings $\OO_\IPS$, $\OO_\FGL$ an $\OO_\FI$, but this
 should be avoided where possible.  Note that for each $R$ we can
 regard $\FGL(R)$ as the set of objects of a groupoid, whose morphism
 set is $\FI(R)$.  In other words, the schemes $\FGL$ and $\FI$ define
 a groupoid scheme.  It is known that $\FI=\spec(MP_0MP)$ (this
 follows easily from the description of $MU_*MU$ in~\cite{ad:shg}).
\end{example}

\begin{example}\label{eg-divisors}
 We now give an example for which representability is less obvious.
 We say that \emph{an effective divisor of degree $n$}\index{divisor}
 on $\aff^1$ over a scheme $Y$ is a subscheme
 $D\sse Y\tm\aff^1=\spec(\OY[x])$ such that $\OO_D$ is a quotient of
 $\OY[x]$ and is free of rank $n$ over $\OY$.  We let
 $X(R)=\Div_n^+(\aff^1)(R)$\index{Divpn@$\Div_n^+(C)$} denote the set
 of such divisors over $\spec(R)$, and we claim that
 $X=\Div_n^+(\aff^1)$ is a scheme.  Firstly, it is a functor of $R$:
 given a ring map $u\:R\xra{}R'$ and a divisor $D$ over $R$ we get a
 divisor $uD=\spec(R'\ot_R\OO_D)=\spec(R')\tm_{\spec(R)}D$ over $R'$.
 Next, given a divisor $D$ as above and an element $y\in R[x]$, we let
 $\lm(y)$ be the map $u\mapsto uy$, which is an $R$-linear
 endomorphism of the module $\OO_D\simeq R^n$.  The map $\lm(x)$ thus
 has a \idx{characteristic polynomial}
 $f_D(t)=\sum_{i=0}^na_i(D)t^{n-i}\in R[t]$.  One checks that the map
 $a_i\:X\xra{}\aff^1$ is natural, so we get an element $a_i$ of $\OX$.
 As $f_D(t)$ is monic, we have $a_0=1$.  The remaining $a_i$'s give us
 a map $X\xra{}\aff^n$.

 The Cayley-Hamilton theorem tells us that $f_D(\lm(x))=0$, but it is
 clear that $f_D(\lm(x))=\lm(f_D(x))$ and $f_D(x)=\lm(f_D(x))(1)$, so
 we find that $f_D(x)=0$ in $\OO_D$ and thus that $\OO_D$ is a
 quotient of $R[x]/f_D(x)$.  On the other hand, it is clear that
 $R[x]/f_D(x)$ is also free over $R$ of rank $n$, and it follows that
 $\OO_D=R[x]/f_D(x)$.  Given this, we see that $D$ is freely and
 uniquely determined by the coefficients $a_1,\ldots,a_n$, so that our
 map $X\xra{}\aff^n$ is an isomorphism.  This shows in particular that
 $X$ is a scheme.  (I learned this argument from~\cite{dega:ga}.)
\end{example}

\subsection{Points and sections}
\label{subsec-points}

Let $X$ be a scheme.  A \emph{point of $X$} means an element 
$x\in X(R)$ for some ring $R$.  We write $\OO_x$ for $R$, which
conveniently allows us to mention $x$ before giving $R$ a name.  
Recall that points $x\in X(R)$ biject with maps $\spec(R)\xra{}X$.
We say that $x$ is defined over $R$, or over $\spec(R)$.

We can also think of an element of $R$ as a point of the scheme
$\aff^1$ over $R$.  If $f\in\OX$ then $f$ is a natural map
$X(S)\xra{}S$ for all rings $S$, so in particular we have a map
$X(R)\xra{}R$.  We thus have $f(x)\in \OO_x=R$.

\begin{example}\label{eg-three-series}
 Let $F$ be a point of $\FGL$\index{FGL@$\FGL$}, in other words a
 formal group law over some ring $R$.  We can write
 \[ [3](x) = F(x,F(x,x)) = 3x + u(F) x^2 + v(F) x^3 + O(x^4) \]
 for certain scalars $u(F)$ and $v(F)$.  This construction associates
 to each point $F\in\FGL$ a point $v(F)\in\aff^1$ in a natural way,
 thus giving an element $v\in\OO_\FGL$.  Of course, we know that
 $\OO_\FGL$ is the \idx{Lazard ring} $L$, which is generated by the
 coefficients $a_{kl}$ of the universal formal group law
 \[ F_{\text{univ}}(x,y) = \sum_{k,l} a_{kl} x^k y^l \]
 Using this formal group law, we find that
 \[ [3](x) = 3 x + 3 a_{11} x^2 + (a_{11}^2 + 8 a_{12}) x^3 + O(x^4)\]
 This means that 
 \[ v(F_{\text{univ}}) = a_{11}^2 + 8 a_{12} \]
 It follows that for any $F$ over any ring $R$, the element $v(F)$ is
 the image of $a_{11}^2 + 8 a_{12}$ under the map $L\xra{}R$
 classifying $F$.
\end{example}

\begin{example}\label{eg-Ha}
 For any scalar $a$, we have a formal group law 
 \[ H_a(x,y) = x + y + a x y. \]
 The construction $a\mapsto H_a$ gives a natural transformation
 $h\:\aff^1(R) \xra{} \FGL(R)$\index{FGL@$\FGL$}, in other words a map
 of schemes $h\:\aff^1 \xra{} \FGL$.  This can be thought of as a
 family of formal group laws, parametrised by $a\in \aff^1$.  It can
 also be thought of as a single formal group law over
 $\Zh[a]=\OO_{\aff^1}$.
\end{example}

\begin{example}\label{eg-schiz}
 The point of view described above allows for some slightly
 schizo\-phrenic constructions, such as regarding the two projections
 $\pi_0,\pi_1\: X\tm X \xra{} X$ as two points of $X$ over $X^2$.
 Indeed, this is the universal example of a scheme $Y$ equipped with
 two points of $X$ defined over $Y$.  Similarly, we can think of the
 identity map $X\xra{}X$ as the universal example of a point of $X$.
 This is analogous to thinking of the identity map of $K(\Zh,n)$ as a
 cohomology class $u\in H^n K(\Zh,n)$; this is of course the universal
 example of a space with a given $n$-dimensional cohomology class.
\end{example}

\begin{definition}\label{defn-Points}
 For any functor $X\:\Rings\xra{}\Sets$, we define a category
 $\Points(X)$\index{PointsX@$\Points(X)$}, whose objects are pairs
 $(R,x)$ with $x\in X(R)$.  The maps $(R,x)\xra{}(S,y)$ are ring maps
 $f\:R\xra{}S$ such that $X(f)(x)=y$.
\end{definition}

\begin{remark}
 Let $X$ be a scheme.  The following categories are equivalent:
 \begin{itemize}
  \item[(a)] The category $\CX_X$ of schemes $Y$ equipped with a map
   $u\:Y\xra{}X$.  \index{XX@$\CX_X$}
  \item[(b)] The category of representable functors
   $Y'\:\Points(X)\xra{}\Sets$.  \index{PointsX@$\Points(X)$}
  \item[(c)] The category of representable functors
   $Y''\:\CX_X^\op\xra{}\Sets$. 
  \item[(d)] The category $\Alg^\op_\OX$ of algebras $R$ over $\OX$.
  \item[(e)] The category $\Points(X)^\op$ of pairs $(R,x)$ with
   $x\in X(R)$.
 \end{itemize}
 By Yoneda, an element $x\in X(R)$ corresponds to a map
 $x'\:\spec(R)\xra{}X$.  Similarly, a map $v\:Z\xra{}X$ gives a map
 $v^*\:\OX\xra{}\OZ$, making $\OZ$ into an $\OX$-algebra.  This can
 also be regarded as an element of $\spec(\OX)(\OZ)=X(\OZ)$.  With
 this notation, the equivalence is as follows.
 \begin{align*}
  Y(S)    &= \coprod_{z\in X(S)} Y'(S,z)                        \\
  Y'(S,z) &=
    \text{ preimage of $z\in X(S)$ under $u\:Y(S)\xra{}X(S)$ }  \\
          &= Y''(\spec(S)\xra{z'}X)                             \\
  Y''(Z\xra{v}X) &= Y'(\OZ,v^*)                                 \\
  R       &= \OY                                                \\
  Y       &= \spec(R).
 \end{align*}
 For us, the most important part of this will be the equivalence
 (a)$\iffa$(b).  
\end{remark}
\begin{remark}
 If $E$ is an even periodic ring spectrum and $S_E=\spec(E^0)$ then we
 can regard the construction $Z\mapsto Z_E=\spec(E^0Z)$ as a functor
 from finite complexes to $\CX_{S_E}$.
\end{remark}

\begin{definition}\label{defn-fibre}
 If $X$ is a scheme over another scheme $Y$, and $y\in Y(R)$ is a
 point of $Y$, we write $X_y=\spec(R)\tY X$, which is a scheme over
 $\spec(R)$.  Here we have used the map $\spec(R)\xra{}Y$
 corresponding to the point $y\in Y(R)$ to form the pullback
 $\spec(R)\tY X$.  We call $X_y$ the \dfn{fibre} of $X$ over the
 point $y$.
\end{definition}

\subsection{Colimits of schemes}
\label{subsec-colim-schemes}

The category of rings has limits for small diagrams, and the category
of schemes is dual to that of rings, so it has colimits for small
diagrams.  However, it seems that these colimits only interact well
with our geometric point of view if they have some additional
properties (this is also the reason for Mumford's geometric invariant
theory, which is much more subtle than anything that we consider
here.)  One good property that often occurs (with $\CC=\CX$ or
$\CC=\CX_Y$) is as follows.

\begin{definition}\label{defn-strong-colim}
 Let $\CC$ be a category with finite products, and let $\{X_i\}$ be a
 diagram in $\CC$.  We say that an object $X$ with a compatible system
 of maps $X_i\xra{}X$ is a \dfn{strong colimit} of the diagram if
 $W\tm X$ is the colimit of $\{W\tm X_i\}$ for each $W\in\CC$.  We
 define strong coproducts and strong coequalisers as special cases of
 this, in the obvious way.
\end{definition}

\begin{example}\label{eg-strong-coprod}
 The categories $\CX$ and $\CX_Y$ have strong finite coproducts, and
 $\OO_{\coprod_iX_i}=\prod_i\OO_{X_i}$.  Indeed, by the usual duality
 $\Rings^\op=\CX$, we see that the coproduct exists and has
 $\OO_{\coprod_iX_i}=\prod_i\OO_{X_i}$.  Thus, we need only check that
 $Z\tY \coprod_iX_i=\coprod_iZ\tY X_i$, or equivalently that
 $\OZ\ot_\OY\prod_i\OO_{X_i}=\prod_i\OZ\ot_\OY\OO_{X_i}$, which is
 clear because the indexing set is finite.  Note that when $Y=1$ is
 the terminal object, we have $\CX_Y=\CX$, so we have covered that
 case as well.
\end{example}

As a special case of the above, we can make the following definition.
\begin{definition}\label{defn-constant-scheme}
 Given a finite set $A$, we can define an associated \emph{constant
   scheme} \index{scheme!constant} $\un{A}$ by
 \[ \un{A}=\coprod_{a\in A}1 \]
 (where $1$ is the terminal object in $\CX$).  This has the property
 that $X\tm\un{A}=\coprod_{a\in A}X$ for all $X$.  We also have
 $\OO_{\un{A}}=F(A,\Zh)$, which denotes the ring of functions from the
 set $A$ to $\Zh$; this is a ring under pointwise operations.
\end{definition}

\begin{remark}\label{rem-points-coprod}
 It is not the case that $(X\amalg Y)(R)=X(R)\amalg Y(R)$ (unlike the
 case of products and pullbacks).  Instead, we have
 \[ (X\amalg Y)(R) = 
   \{(S,T,x,y) \st S,T \le R \;,\; R = S\tm T
       \;,\; x \in X(S) \;,\; y \in Y(T) \}. 
 \]
 To explain this, note that an element of $(X\amalg Y)(R)$ is (by
 Yoneda) a map $\spec(R)\xra{}X\amalg Y$.  This will be given by a
 decomposition $\spec(R)=\spec(S)\amalg\spec(T)$ and maps
 $\spec(S)\xra{}X$ and $\spec(T)\xra{}Y$.  Clearly, if $R$ does not
 split nontrivially as a product of smaller rings then we have the
 naive rule $(X\amalg Y)(R)=X(R)\amalg Y(R)$.

 Similarly, the initial scheme $\emptyset=\spec(0)$ has
 $\emptyset(R)=\emptyset$ unless $R=0$ in which case $\emptyset(R)$
 has a single element.
\end{remark}

\begin{example}\label{eg-strong-sym-prod}
 Let $f\:X\xra{}Y$ be a map of schemes.  Let $X^n_Y$ denote the fibre
 product of $n$ copies of $X$ over $Y$, so that the symmetric group
 $\Sg_n$ acts on $X^n_Y$, covering the trivial action on $Y$.  Suppose
 that the resulting map $f^*\:\OY\xra{}\OX$ makes $\OX$ into a free
 module over $\OY$.  We then claim that there is a
 \idx{strong colimit} for the action of $\Sg_n$ on $X^n_Y$.  To see
 this, write $A=\OX$ and $B=\OY$ and $C=A^{\ot_B n}$, so that
 $X^n_Y=\spec(C)$.  Our claim reduces easily to the statement that
 $B'\ot_B(C^{\Sg_n})=(B'\ot_B C)^{\Sg_n}$ for every algebra $B'$ over
 $B$.  To see that this holds, choose a basis for $A$ over $B$.  This
 gives an evident basis for $C$ over $B$, which is permuted by the
 action of $\Sg_n$.  Clearly $C^{\Sg_n}$ is a free module over $B$,
 with one generator for each $\Sg_n$-orbit in our basis for $C$.
 There is a similar description for $(B'\ot_B C)^{\Sg_n}$, which
 quickly implies our claim.
\end{example}

Some more circumstances in which colimits have unexpectedly good
behaviour are discussed in~\cite{grst:vlc}, which mostly follows ideas
of Quillen~\cite{qu:sec}.

\subsection{Subschemes}
\label{subsec-subschemes}

Recall that an element of $\OX$ is a natural map $X\xra{}\aff^1$.
Thus, if $x$ is a point of $X$ then $f(x)$ is a scalar (more
precisely, if $x\in X(R)$ then $f(x)\in R$) and we can ask whether
$f(x)=0$, or whether $f(x)$ is invertible.

\begin{definition}\label{defn-subschemes}
 Given a scheme $X$ and an ideal $I\leq\OX$, we define a scheme
 $V(I)$\index{VI@$V(I)$} by 
 \[ V(I)(R)=\{x\in X(R)\st f(x)=0\text{ for all } f\in I\}. \]
 One checks that $V(I)=\spec(\OX/I)$, so this really is a scheme.
 Schemes of this form are called \emph{closed subschemes}
 \index{subscheme!closed} of $X$.

 Given an element $f\in\OX$, we define a scheme
 $D(f)$\index{Df@$D(f)$} by
 \[ D(f)(R)=\{x\in X(R)\st f(x) \in R^\tm\}. \]
 One checks that $D(f)=\spec(\OX[1/f])$, so this really is a scheme.
 Schemes of this form are called \emph{basic open subschemes}
 \index{subscheme!open} of $X$.

 A \emph{locally closed subscheme} \index{subscheme!locally closed} is
 a basic open subscheme of a closed subscheme.  Such a thing has the
 form $D(f)\cap V(I)=\spec(\OX[1/f]/I)$.
\end{definition}

\begin{remark}\label{rem-regular-mono}
 Recall that a map $f\:R\xra{}S$ of rings is said to be a
 \dfn{regular epimorphism} if and only if it is the coequaliser of
 some pair of maps $T\arw{e,=>}R$, which happens if and only if it is
 the coequaliser of the obvious maps $R\tm_SR\arw{e,=>}R$.  It is easy
 to check that this holds if and only if $f$ is surjective.  Given
 this, we see that the \idx{regular monomorphism}s of schemes are
 precisely the closed inclusions\index{closed inclusion}, and that
 composites and pushouts of regular monomorphisms are regular
 monomorphisms.
\end{remark}

\begin{example}\label{eg-open-closed}
 The map $h$ in Example~\ref{eg-Ha} gives an isomorphism between
 $\aff^1$ and the closed subscheme $V((a_{ij}\st i+j>2))$ of $\FGL$.
 The multiplicative group $\MG$ is an open subscheme of $\aff^1$.
\end{example}

\begin{example}\label{eg-clopen}
 If $X$ is a scheme and $e\in\OX$ satisfies $e^2=e$ then it is easy to
 check that $D(e)=V(1-e)$, so this subscheme is both open and closed.
 Moreover, we have $X=D(e)\amalg D(1-e)$.  More generally, if we have
 idempotents $e_1,\ldots,e_m\in\OX$ with $\sum_ie_i=1$ and
 $e_ie_j=\dl_{ij}e_i$ then $X=\coprod_iD(e_i)$, and every splitting of
 $X$ as a finite coproduct occurs in this way.
\end{example}

\begin{example}\label{eg-affine-line}
 Suppose $X=\spec(k[x])$ is the affine line over a field $k$, and
 $\lm,\mu\in k$.  The closed subscheme
 $V(x-\lm)=\spec(k[x]/(x-\lm))\simeq \spec(k)$ corresponds to the
 point $\lm$ of the affine line; it is natural to refer to it as
 $\{\lm\}$.  The closed subscheme $V((x-\lm)(x-\mu))$ corresponds to
 the pair of points $\{\lm,\mu\}$.  If $\lm=\mu$, this is to be
 thought of as the point $\lm$ with multiplicity two, or as an
 infinitesimal thickening of the point $\lm$.
\end{example}

We can easily form the intersection of locally closed subschemes:
\[ D(a)\cap V(I) \cap D(b)\cap V(J) = D(ab) \cap V(I+J). \]
We cannot usually form the union of basic open subschemes and still
have an affine scheme.  Again, it would be easy enough to consider
non-affine schemes, but it rarely seems to be necessary.  Moreover, a
closed subscheme $V(a)$ determines the complementary open subscheme
$D(a)$ but not conversely; $D(a)=D(a^2)$ but $V(a)\neq V(a^2)$ in
general.

We say that a scheme $X$ is \dfn{reduced} if $\OX$ has no nonzero
nilpotents, and write
$X_\red=\spec(\OX/\sqrt{0})$\index{Xred@$X_\red$}, which is the
largest reduced closed subscheme of $X$.  Moreover, if $Y\sse X$ is
closed then $Y_\red=X_\red$ if and only if $X(k)=Y(k)$ for every
field $k$ (we leave the proof as an exercise).

We define the union of closed subschemes \index{subscheme!union} by
$V(I)\cup~V(J)=V(I\cap~J)$.  We also define the \emph{schematic union}
by $V(I)+V(J)=V(IJ)$.  This is a sort of ``union with multiplicity''
--- in particular, $V(I)+V(I)\neq V(I)$ in general.  In the previous
example, we have
\[ \{\lm\} \cup \{\lm\} = V((x-\lm)^2) \]
which is a thickening of $\{\lm\}$.  Note that 
$V(IJ)_\red=V(I\cap J)_\red$, because 
$(I\cap J)^2\leq IJ\leq I\cap J$.

We shall say that $X$ is \dfn{connected}\index{components} if it
cannot be split nontrivially as $Y\amalg Z$, if and only if there are
no idempotents in $\OX$ other than 0 and 1.  

We shall say that a scheme $X$ is \dfn{integral} if and only if $\OX$
is an integral domain, and that $X$ is \dfn{irreducible} if and only
if $X_\red$ is integral.  We also say that $X$ is \dfn{Noetherian} if
and only if the ring $\OX$ is Noetherian.  If so, then $X_\red$ can be
written in a unique way as a finite union $\bigcup_i Y_i$ with $Y_i$
an integral closed subscheme.  The schemes $Y_i$ are called the
\emph{irreducible components} of $X_\red$; they are precisely the
schemes $V(\pri_i)$ for $\pri_i$ a \idx{minimal prime} ideal of $\OX$.
See~\cite[section 6]{ma:crt} for this material.

Suppose that $X$ is Noetherian and reduced, say
$X=\bigcup_{i\in~S}Y_i$ as above for some finite set $S$.  Suppose
that $S=S'\amalg S''$.  Write $X'=\bigcup_{S'}Y_i=V(I')$, where
$I'=\bigcap_{S'}\pri_i$, and similarly for $X''$ and $I''$.  If we
then write
\[ \Gamma(I')=\{a\in\OX\st a (I')^N = 0 \text{ for } N\gg 0\},\]
we find that $\Gamma(I')=I''$ and thus $V(\Gamma(I'))=X''$.  

\begin{example}\label{eg-primary-dec}
 Take $Z=\spec(k[x,y]/(xy^2))$ and set 
 \begin{align*}
   X  &= V(y)   = \spec(k[x])            \\
   X' &= V(y^2) = \spec(k[x,y]/(y^2))    \\
   Y  &= V(x)   = \spec(k[y])           
 \end{align*}
 Then $X$ is the $x$-axis, $Y$ is the $y$-axis and $X'$ is an
 infinitesimal thickening of $X$.  The schemes $X$ and $Y$ are
 integral, and $X'$ is irreducible because $X'_\red=X$.  The scheme
 $Z$ is reducible, and its irreducible components are $X$ and $Y$.
\end{example}

\subsection{Zariski spectra and geometric points}
\label{subsec-zariski}

If $R$ is a ring, we define the 
\dfn{Zariski space}\index{zarR@$\zar(R)$}\index{Xzar@$X_\zar$} to be
\[ \zar(R) = \{ \text{ prime ideals } \pri<R \;\}. \]
If $X$ is a scheme, we write $X_\zar=\zar(\OX)$.  Note that
\begin{align*}
 V(I)_\zar &= \zar(\OX/I) = \{\pri\in X_\zar\st I\le\pri\} \\
 D(f)_\zar &= \zar(\OX[1/f]) =
    \{\pri\in X_\zar\st f\not\in\pri\} \\
 (X \amalg Y)_\zar &= X_\zar \amalg Y_\zar
\end{align*}
There is a map
\[ (X\tm Y)_\zar \xra{} X_\zar \tm Y_\zar, \]
but it is almost never a bijection.

Suppose that $Y,Z\le X$ are locally closed; then
\[ (Y\cap Z)_\zar = Y_\zar \cap Z_\zar. \]
If $Y$ and $Z$ are closed then
\[ (Y\cup Z)_\zar = (Y+Z)_\zar = Y_\zar \cup Z_\zar. \]

We give $X_\zar$ the topology with closed sets $V(I)_\zar$.  A map
of schemes $X\xra{}Y$ then induces a continuous map
$X_\zar\xra{}Y_\zar$.

Suppose that $R$ is an integral domain, and that $x\in X(R)$.  Then
$x$ gives a map $x^*\:\OX\xra{}R$, whose kernel $\pri_x$ is prime.  We
thus have a map $X(R)\xra{}X_\zar$, which is natural for monomorphisms
of $R$ and arbitrary morphisms of $X$.  

A \dfn{geometric point} of $X$ is an element of $X(k)$, for some
algebraically closed field $k$.  Suppose that either $\OX$ is a
$\rat$-algebra, or that some prime $p$ is nilpotent in $\OX$.  Let
$k$ be an algebraically closed field of the appropriate
characteristic, with transcendence degree at least the cardinality of
$\OX$.  Then it is easy to see that $X(k)\xra{}X_\zar$ is epi.

A useful feature of the Zariski space is that it behaves quite well
under colimits~\cite{qu:sec,grst:vlc}.  The following proposition is
another example of this.
\begin{proposition}\label{prop-zariski-quot}
 Suppose that a finite group $G$ acts on a scheme $X$.  Then
 $(X/G)_\zar=X_\zar/G$.
\end{proposition}
\begin{proof}
 Write $S=\OX$ and $R=S^G=\OO_{X/G}$.  Given a prime
 $\pri\in\zar(R)=(X/G)_\zar$, the fibre $F$ over $\pri$ in
 $\zar(S)=X_\zar$ is just $\zar(S_\pri/\pri S\pri)$ (see~\cite[Section
 7]{ma:crt}).  We need to prove that $F$ is nonempty, and that $G$
 acts transitively on $F$.  

 As localisation is exact, we have $(S_\pri)^G=R_\pri$, so we can
 replace $R$ by $R_\pri$ and thus assume that $R$ is local at $\pri$.
 With this assumption, we have $F=\zar(S/\pri S)$.  For $a\in S$ we
 write $f_a(t)=\prod_{g\in G}(t-ga)\in S[t]^G=R[t]$, so that $f_a$ is
 a monic polynomial with $f_a(a)=0$.  This shows that $S$ is an
 integral extension over $R$, so $F\neq\emptyset$ and there are no
 inclusions between the elements of $F$~\cite[Theorem 9.3]{ma:crt}.

 Let $\qri$ and $\rri$ be two points of $F$, so they are prime ideals
 in $S$ with $\qri\cap R=\qri^G=\pri$ and $\rri\cap R=\rri^G=\pri$.
 Write $I=\bigcap_{g\in G}g.\qri\leq S$.  If $a\in I$ then
 $g.a\in\qri$ for all $g$ so $f_a(t)\in t^{|G|}+\qri[t]$ but also
 $f_a(t)$ is $G$-invariant so
 $f_a(t)\in t^{|G|}+\qri^G[t]\sse t^{|G|}+\rri[t]$.  As $f_a(a)=0$ we
 conclude that $a$ is nilpotent mod $\rri$ but $\rri$ is prime so
 $a\in\rri$.  Thus $\bigcap_{g\in G}g.\qri\leq\rri$.  As $\rri$ is
 prime, we deduce that $g.\qri\leq\rri$ for some $g\in G$.  As there
 are no inclusions between the elements of $F$, we conclude that
 $g.\qri=\rri$.  Thus $G$ acts transitively on $F$, which proves that
 $(X/G)_\zar=X_\zar/G$.
\end{proof}

A number of interesting things can be detected by looking at Zariski
spaces.  For example, $X_\zar$ splits as a disjoint union if and only
if $X$ does --- see Corollary~\ref{cor-coproducts}.

We also use the space $X_\zar$ to define the \dfn{Krull dimension} of
$X$\index{dimX@$\dim(X)$}.
\begin{definition}\label{defn-dim}
 If there is a chain $\pri_0<\ldots<\pri_n$ in $X_\zar$, but no longer
 chain, then we say that $\dim(X)=n$.  If there are arbitrarily long
 chains then $\dim(X)=\infty$.
\end{definition}
\begin{example}\label{eg-dim}
 The terminal object $1$ has dimension one (because there are chains
 $(0)<(p)$ of prime ideals in $\Zh$).  If $\OX$ is a field then
 $\dim(X)=0$.  If $\OX$ is Noetherian then $\dim(\MG\tm X)=1+\dim(X)$
 and $\dim(\aff^n\tm X)=n+\dim(X)$~\cite[Section 15]{ma:crt}.  In
 particular, we have $\dim(\aff^n)=\dim(1\tm\aff^n)=n+1$.
\end{example}
\begin{example}\label{eg-infinite-dim}
 The schemes $\FGL$, $\IPS$ and $\FI$ all have infinite dimension.
\end{example}

\subsection{Nilpotents, idempotents and connectivity}
\label{subsec-nilpotents}

\begin{proposition}
 Suppose that $e\in R$ is \idx{idempotent}, and $f=1-e$.  Then
 \[ eR = R/f = R[e^{-1}] = \{a\in R\st fa=0\}. \]
 Moreover, this is a ring with unit $e$, and we have $R=eR\tm fR$ as
 rings. \qed
\end{proposition}

\begin{proposition}
 If $X$ is a scheme, then splittings $X=\coprod_{i=1}^nX_i$ biject
 with systems of idempotents $\{e_1,\ldots,e_n\}$ with $\sum_ie_i=1$
 and $e_ie_j=\dl_{ij}e_j$. \qed
\end{proposition}

\begin{example}\label{eg-Mult-n}
 Let $\Mult(n)$\index{Multn@$\Mult(n)$} be the scheme of polynomials
 $\phi(u)$ of degree at most $n$ such that $\phi(1)=1$ and
 $\phi(uv)=\phi(u)\phi(v)$.  Such a series can be written as
 $\phi(u)=\sum_{i=0}^ne_iu^i$, and the conditions on $\phi$ are
 equivalent to $\sum_ie_i=1$ and $e_ie_j=\dl_{ij}e_j$.  In other
 words, the elements $e_i$ are orthogonal idempotents.  Using this, we
 see easily that $\Mult(n)=\coprod_{i=0}^n 1$.
\end{example}

\begin{example}\label{eg-E-n}
 Now let $E(n)$\index{En@$E(n)$} be the scheme of $n\tm n$ matrices
 $A$ over $R$ such that $A^2=A$.  Define
 $\al_A(u)=uA+(1-A)=(u-1)A+1\in M_n(R[u])$ and
 $\phi_A(u)=\det(\al_A(u))\in R[u]$.  We find easily that $\al_A(1)=1$
 and $\al_A(uv)=\al_A(u)\al_A(v)$, so $\phi_A(u)\in\Mult(n)(R)$.  This
 construction gives a map $E(n)\xra{}\Mult(n)=\coprod_{i=0}^n1$, which
 gives a splitting $E(n)=\coprod_{i=0}^nE(n,i)$, where $E(n,i)$ is the
 scheme of $n\tm n$ matrices $A$ such that $A^2=A$ and
 $\phi_A(u)=u^i$.

 Note that the function $A\mapsto\trace(A)$ lies in $\OO_{E(n)}$ and
 that $E(n,i)$ is contained in the closed subscheme
 $E'(n,i)=\{A\st\trace(A)=i\}$.  However, if $n>0$ but $n=0$ in $R$
 then $E'(n,0)(R)$ and $E'(n,n)(R)$ are not disjoint, which shows that
 $E'(n,i)\neq E(n,i)$ in general.
\end{example}

For any ring $R$, we let $\Nil(R)$\index{NilR@$\Nil(R)$} denote the
set of nilpotents\index{nilpotent} in $R$.
\begin{proposition}\label{prop-nilpotents}
 $\Nil(R)$ is the intersection of all prime ideals in $R$.
\end{proposition}
\begin{proof}
 \cite[Section 1]{ma:crt}
\end{proof}

\begin{proposition}[Idempotent Lifting]\label{prop-lift}
 Suppose that $e\in R/\Nil(R)$ is \idx{idempotent}.  Then there is a
 unique idempotent $\tilde{e}\in R$ lifting $e$.
\end{proposition}
\begin{proof}
 Choose a (not necessarily idempotent) lift of $e$ to $R$, call it
 $e$, and write $f=1-e$.  We know that $ef$ is nilpotent, say
 $e^nf^n=0$.  Define
 \[ c = e^n + f^n - 1 = e^n + f^n - (e + f)^n \]
 This is visibly divisible by $ef$, hence nilpotent; thus
 $e^n+f^n=1+c$ is invertible.  Define
 \[ \tilde{e} = e^n/(1+c) \qquad\qquad 
    \tilde{f} = f^n/(1+c) = 1-\tilde{e}
 \]
 Then $\tilde{e}$ is an idempotent lifting $e$.  If $\tilde{e}_1$ is
 another such then $\tilde{e}_1\tilde{f}$ is idempotent.  It lifts
 $ef=0$, so it is also nilpotent.  It follows that
 $\tilde{e}_1\tilde{f}=0$ and $\tilde{e}_1=\tilde{e}\tilde{e}_1$.
 Similarly, $\tilde{e}=\tilde{e}\tilde{e}_1$, so $\tilde{e}=\tilde{e}_1$.
\end{proof}

\begin{theorem}[Chinese Remainder Theorem]\label{thm-chinese}
 \index{Chinese remainder theorem}
 Suppose that $\{I_\al\}$ is a finite family of ideals in $R$, which
 are pairwise coprime (i.e.\ $I_\al+I_\bt=R$ when $\al\neq\bt$).  Then
 \[ R/\bigcap_\al I_\al = \prod_\al R/I_\al \]
\end{theorem}
\begin{proof}
 \cite[Theorems 1.3,1.4]{ma:crt}
\end{proof}

\begin{corollary}\label{cor-coproducts}
 Suppose that $\zar(R)=\coprod_\al\zar(R/I_\al)$ (a finite coproduct).
 Then there are unique ideals $J_\al\le I_\al\le \sqrt{J_\al}$ such
 that $R\simeq \prod_\al R/J_\al$.
\end{corollary}
\begin{proof}
 Proposition~\ref{prop-nilpotents} implies that $\bigcap_\al I_\al$ is
 nilpotent.  If $\al\neq\bt$ then no prime ideal contains
 $I_\al+I_\bt$, so $I_\al+I_\bt=R$.  Now use the Chinese remainder
 theorem, followed by idempotent lifting.
\end{proof}
\begin{remark}
 There are nice topological applications of these ideas
 in~\cite{kasttu:mkh,grst:vlc}, for example.
\end{remark}

\subsection{Sheaves, modules and vector bundles}
\label{subsec-sheaves}

The simplest definition of a sheaf over a scheme $X$ is just as a
module over the ring $\OX$.  (It would be more accurate to refer to
this as a quasi-coherent sheaf of $\OO$-modules over $X$, but we shall
just call it a sheaf.)  However, we shall give a different (but
equivalent) definition which fits more neatly with our emphasis on
schemes as functors, and which generalises more easily to formal
schemes. 

\begin{definition}\label{defn-sheaf}
 A \dfn{sheaf} over a functor $X\in\CF$ consists of the following
 data:
 \begin{itemize}
 \item[(a)] For each $(R,x)\in\Points(X)$, a module $M_x$ over $R$.
 \item[(b)] For each map $f\:(R,x)\xra{}(S,y)$ in $\Points(X)$, an
  isomorphism $\tht(f)=\tht(f,x)\:S\ot_RM_x\xra{}M_y$ of $S$-modules.
 \end{itemize}
 The maps $\tht(f,x)$ are required to satisfy the functorality
 conditions 
 \begin{itemize}
 \item[(i)]  In the case $f=1\:(R,x)\xra{}(R,x)$ we have
  $\tht(1,x)=1\:M_x\xra{}M_x$. 
 \item[(ii)] Given maps $(R,x)\xra{f}(S,y)\xra{g}(T,z)$, the map
  $\tht(gf,x)$ is just the composite
  \[ T\ot_RM_x=T\ot_SS\ot_RM_x \xra{1\ot\tht(f,x)}
     T\ot_SM_y \xra{\tht(g,y)} M_z.
  \]
 \end{itemize}
 We write $\Sheaves_X$\index{SheavesX@$\Sheaves_X$} for the category
 of sheaves over $X$.  This has direct sums (with
 $(M\oplus N)_x=M_x\oplus N_x$) and tensor products (with
 $(M\ot N)_x=M_x\ot_RN_x$ when $x\in X(R)$).  The unit for the tensor
 product is the sheaf $\OO$, which is defined by $\OO_x=R$ for all
 $x\in X(R)$.
\end{definition}
\begin{remark}
 If $M$ and $N$ are sheaves over a sufficiently bad functor $X$, it
 can happen that $\Sheaves_X(M,N)$ is a proper class.  This will not
 be the case if $X$ is a scheme or a formal scheme, however.
\end{remark}

\begin{example}\label{eg-origin}
 Let $x$ be a point of $\aff^1(R)$, or in other words an element of
 $R$.  Define $M_x=R/x$; this gives a sheaf over $\aff^1$.  Note that
 $M_x=0$ if $x$ is invertible, but $M_x=R$ if $x=0$.  Thus, $M$ is
 concentrated at the origin of $\aff^1$.
\end{example}

\begin{definition}\label{defn-Gamma}
 \begin{enumerate}
 \item Let $X$ be a functor in $\CF$.  If $N$ is a module over the
  ring $\OX=\CF(X,\aff^1)$, we define a sheaf $\tN$\index{N@$\tN$}
  over $X$ by $\tN_x=R\ot_\OX N$, where we use $x$ to make $R$ into an
  algebra over $\OX$.
 \item If $M$ is a sheaf over $X$ and $R$ is a ring, we write
  $\aff(M)(R)=\coprod_{x\in X(R)}M_x$.  \index{AAM@$\aff(M)$} If
  $f\:R\xra{}S$ is a homomorphism, we define a map
  $\aff(M)(R)\xra{}\aff(M)(S)$, which sends $M_x$ to $M_{f(x)}$ by
  $m\mapsto\tht(f,x)(1\ot m)$.  This gives a functor
  $\aff(M)\in\CF_X$.
 \item If $M$ is a sheaf over $X$, we define
  $\Gm(X,M)=\CF_X(X,\aff(M))$.  \index{GXM@$\Gm(X,M)$} Thus, an
  element $u\in\Gm(X,M)$ is a system of elements $u_x\in M_x$ for all
  rings $R$ and points $x\in X(R)$, which behave in the obvious way
  under maps of rings.  If $M=\OO$ then $\aff(\OO)=\aff^1\tm X$ and
  $\Gm(X,\OO)=\OX$.  It follows that $\Gm(X,M)$ is a module over $\OX$
  for all $M$.
 \item If $Y$ is a scheme over $X$, we also define
  $\Gm(Y,M)=\CF_X(Y,\aff(M))$.
 \end{enumerate}
\end{definition}

\begin{proposition}\label{prop-Gamma-adjoint}
 For any functor $X\in\CF$, the functor
 $\Gm(X,-)\:\Sheaves_X\xra{}\Mod_\OX$ is right adjoint to the functor
 $N\mapsto\tN$. \index{N@$\tN$}\index{GXM@$\Gm(X,M)$}
\end{proposition}
\begin{proof}
 For typographical convenience, we will write $TN$ for $\tN$ and $GM$
 for $\Gm(X,M)$.  We define maps $\eta\:N\xra{}GTN$ and
 $\ep\:TGM\xra{}M$ as follows.  Let $n$ be an element of $N$; for each
 point $x\in X(R)$, we define $\eta(n)_x=1\ot n\in R\ot_\OX N=(TN)_x$,
 giving a map $\eta$ as required.  Next, we define
 $\ep_x\:(TGM)_x=R\ot_\OX\Gm(X,M)\xra{}M_x$ by $\ep_x(a\ot u)=au_x$.
 We leave it to the reader to check the triangular identities
 $(\ep_T)(T\eta)=1_T$ and $(G\ep)(\eta_G)=1_G$, which show that we
 have an adjunction.
\end{proof}

\begin{proposition}\label{prop-Gamma-equiv}
 Let $X$ be a scheme, and let $x_0\in X(\OX)$ be the tautological
 point, which corresponds to the identity map of $\OX$ under the
 isomorphism $X=\spec(\OX)$.  Then there is a natural isomorphism
 $\Gm(X,M)=M_{x_0}$, and $\Gm(X,-)\:\Sheaves_X\xra{}\Mod_\OX$ is an
 equivalence of categories.\index{GXM@$\Gm(X,M)$}
\end{proposition}
\begin{proof}
 First, we define a map $\al\:\Gm(X,M)\xra{}M_{x_0}$ by
 $u\mapsto u_{x_0}$.  Next, suppose that $m\in M_{x_0}$.  If
 $x\in X(R)$ for some ring $R$ then we have a corresponding ring map
 $\hat{x}\:f\mapsto f(x)$ from $(\OX,x_0)$ to $(R,x)$.  We define
 $\bt(m)_x=\tht(\hat{x},x_0)(m)\in M_x$.  One can check that this
 gives an element $\bt(m)\in\Gm(X,M)$, and that
 $\bt\:M_{x_0}\xra{}\Gm(X,M)$ is inverse to $\al$.  It follows that
 $\Gm(X,\tN)=\tN_{x_0}$, which is easily seen to be the same as $N$.
 Also, if $N=M_{x_0}$ then $\tN_x=R\ot_\OX M_{x_0}$, and
 $\tht(\hat{x},x_0)$ gives an isomorphism of this with $M_x$, so
 $\tN=M$.  It follows that the functor $N\mapsto\tN$ is inverse to
 $\Gm(X,-)$. 
\end{proof}

It follows that when $X$ is a scheme, the category
$\Sheaves_X$\index{SheavesX@$\Sheaves_X$} is
Abelian.  Because tensor products preserve colimits and finite
products, we see that the functors $M\mapsto M_x$ preserve colimits
and finite products.

We next need some recollections about finitely generated projective
modules.\index{projective module} If $M$ is such a module over a ring
$R$ and $\pri\in\zar(R)$ then $M_\pri$ is a finitely generated module
over the local ring $R_\pri$ and thus is free~\cite[Theorem
2.5]{ma:crt}, of \idx{rank} $r_\pri(M)$ say.  Note that $r_\pri(M)$ is
also the dimension of $\kp(\pri)\ot_RM$ over the field
$\kp(\pri)=R_\pri/\pri R_\pri$.  If this is independent of $\pri$ then
we call it $r(M)$ and say that $M$ has constant rank.  Clearly, if any
two of $M$, $N$ and $M\oplus N$ have constant rank then so does the
third and $r(M\oplus N)=r(M)+r(N)$.  Also, if $r(M)=0$ then $M=0$.

\begin{definition}\label{defn-vec-bundle}
 Let $M$ be a sheaf over a functor $X$.  If $M_x$ is a finitely
 generated \idx{projective module} over $\OO_x$ for all $x\in X$, we
 say that $M$ is a \dfn{vector bundle} or \emph{locally free sheaf}
 \index{sheaf!locally free} over $X$.  If in addition $M_x$ has rank
 one for all $x$, we say that $M$ is a \dfn{line bundle} or
 \emph{invertible sheaf}\index{sheaf!invertible}.
\end{definition}
If $X$ is a scheme, a sheaf $M$ is a vector bundle if and only if
$\Gm(X,M)$ is a finitely generated projective module over $\OX$.
However, this does not generalise easily to formal schemes, so we do
not take it as the definition.  It is not hard to check that $M_x$ has
constant rank $r$ for all $R$ and all $x\in X(R)$ if and only if $M_x$
has dimension $r$ over $K$ for all algebraically closed fields $K$ and
all $x\in X(K)$.

\begin{remark}\label{rem-usually-free}
 In algebraic topology, it is very common that the naturally occurring
 projective modules are free, and thus that the corresponding vector
 bundles and line bundles are trivialisable.  However, they are
 typically not equivariantly trivial for important groups of
 automorphisms, so it is conceptually convenient to avoid choosing
 bases.  The main example is that if $Z$ is a finite complex and $V$
 is a complex vector bundle over $Z$ with Thom complex $Z^V$ then
 $\widetilde{E}^0Z^V$ gives a line bundle over $Z_E$.  A choice of
 complex orientation on $E$ gives a Thom class and thus a
 trivialisation, but this is not invariant under automorphisms of
 $E$. 
\end{remark}

\begin{example}\label{eg-projector}
 Recall the scheme $E(n)=\coprod_{i=0}^nE(n,i)$ \index{En@$E(n)$} of
 Example~\ref{eg-E-n}.  A point of $E(n)(R)$ is an $n\tm n$ matrix $A$
 over $R$ with $A^2=A$.  This means that $M_A=A.R^n$ is a finitely
 generated projective $R$-module, so this construction defines a
 vector bundle $M$ over $E(n)$.  If $A$ is a point of $E(n,i)$ (so
 that $\det((u-1)A+1)=u^i\in R[u]$) and $R$ is an algebraically closed
 field, then elementary linear algebra shows that $A$ has rank $i$.
 It follows that the restriction of $M$ to $E(n,i)$ has rank $i$.

 Let $N$ be a \idx{vector bundle} over an arbitrary scheme $X$.  The
 associated projective $\OX$-module \index{projective module} is then
 a retract of a finitely generated free module, so there is a matrix
 $A\in E(n)(\OX)$ such that $\Gm(X,N)=A.\OO_X^n$ for some $n$.  The
 point $A\in E(n)(\OX)$ corresponds to a map $\al\:X\xra{}E(n)$, and
 we find that $\al^*M=N$.  If $X_i$ denotes the preimage of $E(n,i)$
 under $\al$, then $X=\coprod_iX_i$ and the restriction of $N$ to
 $X_i$ has rank $i$.
\end{example}

Let $X$ be a scheme.  Using equivalence $\Sheaves_X\simeq\Mod_\OX$
\index{SheavesX@$\Sheaves_X$} again, we see that there are sheaves
$\Hom(M,N)$ such that
\[ \Sheaves_X(L,\Hom(M,N))=\Sheaves_X(L\ot M,N). \]
In particular, we define $M^\vee=\Hom(M,\OO)$.  If $M$ is a vector
bundle then we have $\Hom(M,N)_x=\Hom_R(M_x,N_x)$ and thus
$(M^\vee)_x=\Hom(M_x,R)$.  In that case $M^\vee$ is again a vector
bundle and $M^{\vee\vee}=M$.  If $M$ is a line bundle then we also
have $M\ot M^\vee=\OO$.

\begin{example}\label{eg-ideal-sheaf}
 Let $Y$ be a closed subscheme of $X$, with inclusion map
 $j\:Y\xra{}X$.  Then
 $I_Y=\{f\in\OX\st f(y)=0\text{ for all points } y\in Y\}$ is an
 ideal in $\OX$ and $\OY=\OX/I_Y$.  We define $j_*\OO$ to be the
 sheaf over $X$ corresponding to the $\OX$-module $\OY$.  More
 explicitly, we have
 \[ (j_*\OO)_x = \OO_x/(f(x)\st f\in J_Y\sse\OX). \]
 We also let $\CI_Y$ be the sheaf associated to the $\OX$-module
 $I_Y$, so that $(\CI_Y)_x=\OO_x\ot_\OX I_Y$ for all points $x$ of
 $X$.  Note that the sequence $\CI_Y\mra\OO\era j_*\OO$ is short exact
 in $\Sheaves_X$, even though the sequences
 $(\CI_Y)_x\xra{}\OX\era(j_*\OO)_x$ need only be right exact.
\end{example}

\begin{example}\label{eg-sheaf-pullback}
 Given a sheaf $N$ over a functor $Y$ and a map $f\:X\xra{}Y$, we can
 define a sheaf $f^*N$ over $X$ by $(f^*N)_x=N_{f(x)}$.  The functor
 $f^*\:\Sheaves_Y\xra{}\Sheaves_X$ clearly preserves colimits and tensor
 products.  If $N$ is a vector bundle then so is $f^*N$ and we have
 $f^*\Hom(N,M)=\Hom(f^*N,f^*M)$ for all $M$.  If $X$ and $Y$ are
 schemes, we find that $\Gm(X,f^*N)=\OX\ot_\OY\Gm(Y,N)$.
\end{example}
\begin{example}
 If the functor $f^*$ defined above has a right adjoint, we call it
 $f_*$.  If $X$ and $Y$ are schemes then we know from
 Proposition~\ref{prop-Gamma-equiv} that there is an essentially
 unique functor $f_*\:\Sheaves_X\xra{}\Sheaves_Y$ such that
 $\Gm(Y,f_*M)=\Gm(X,M)$ (where the right hand side is regarded as an
 $\OY$-module using the map $\OX\xra{}\OY$ induced by $f$).  Using the
 fact that $\Gm(X,f^*N)=\OX\ot_\OY\Gm(Y,N)$ one checks that $f_*$ is
 right adjoint to $f^*$ as required.
\end{example}

\begin{proposition}\label{prop-aff-scheme}
 If $M$ is a vector bundle over a scheme $X$, then $\aff(M)$
 \index{AAM@$\aff(M)$} is a scheme.
\end{proposition}
\begin{proof}
 Write $N=\Mod_\OX(\Gm(X,M),\OX)$.  Then for any map
 $(x\:\OX\xra{}R)\in X(R)$ we have $M_x=\Mod_\OX(N,R)$, where $R$ is
 considered as an $\OX$-module via $x$.  If we let $S$ be the
 symmetric algebra $\Sym_\OX[N]$ then we have $M_x=\Alg_\OX(S,R)$.  It
 follows easily that
 $\Rings(S,R)=\coprod_x\Alg_{\OX,x}(S,R)=\coprod_xM_x=\aff(M)(R)$, so
 $\aff(M)$ is representable as required.
\end{proof}
\begin{definition}\label{defn-gen-L}
 Given a line bundle $L$ over a functor $X$, we define a functor
 $\aff(L)^\tm$ \index{AALt@$\aff(L)^\tm$} over $X$ by
 \[ \aff(L)^\tm(R) =
     \coprod_{x\in X(R)}
      \{\text{ isomorphisms $u\:R\xra{}L_x$ of $R$-modules }\}.
 \]
 If $X$ is a scheme, an argument similar to the one for $\aff(M)$
 shows that $\aff(L)^\tm=\spec(\bigoplus_{n\in\Zh}N^{\ot n})$, where
 $N=\Mod_\OX(\Gm(X,L),\OX)$ and $N^{\ot(-n)}$ means the dual of
 $N^{\ot n}$.  In particular, $\aff(L)^\tm$ is a scheme in this case.
\end{definition}

\subsection{Faithful flatness and descent}
\label{subsec-descent}

\begin{definition}\label{defn-ff}
 Let $f\:X\xra{}Y$ be a map of schemes, and $f^*\:\CX_Y\xra{}\CX_X$
 the associated pullback functor.  We say that $f$ is \dfn{flat} if
 $f^*$ preserves finite colimits.  By Example~\ref{eg-strong-coprod},
 it is equivalent to say that $f^*$ preserves coequalisers.  We say
 that $f$ is \emph{faithfully flat}\index{flat!faithfully} if $f^*$
 preserves finite colimits and reflects isomorphisms.
\end{definition}

\begin{remark}\label{rem-refl-colim}
 Let $f\:X\xra{}Y$ be faithfully flat.  We claim that $f^*$ reflects
 finite colimits, so that $f^*Z=\colim_if^*Z_i$ if and only if
 $Z=\colim_iZ_i$.  More precisely, if $\{Z_i\}$ is a finite diagram in
 $\CX_Y$ and $\{Z_i\xra{}Z\}$ is a cone under the diagram, then
 $\{f^*Z_i\xra{}f^*Z\}$ is a colimit cone in $\CX_X$ if and only if
 $\{Z_i\xra{}Z\}$ is a colimit cone in $\CX_Y$.  The ``if'' part is
 clear.  For the ``only if'' part, write $Z'=\colim_iZ_i$, so we have
 a canonical map $u\:Z'\xra{}Z$.  As $f$ is flat we have
 $f^*Z'=\colim_i f^*Z_i=f^*Z$.  As $f^*$ reflects isomorphisms, we see
 that $u$ is an isomorphism if $f^*u$ is an isomorphism.  The claim
 follows.
\end{remark}

\begin{remark}\label{rem-ff-classical}
 Classically, a module $M$ over a ring $A$ is said to be flat if the
 functor $M\ot_A(-)$ is exact.  It is said to be faithfully flat if in
 addition, whenever $M\ot_AL=0$ we have $L=0$.  It turns out that $f$
 is (faithfully) flat \index{flat!faithfully} if and only if the
 associated ring map $\OY\xra{}\OX$ makes $\OX$ into a (faithfully)
 flat module over $\OY$.  We leave this as an exercise (consider
 schemes of the form $\spec(\OX\oplus L)$, where $L$ is an $\OX$
 module and the ring structure is such that $L.L=0$).
\end{remark}

\begin{remark}
 The idea of faithful flatness was probably first used in topology by
 Quillen~\cite{qu:sec}.  He observed that if $V$ is a complex vector
 bundle over a finite complex $Z$ and $F$ is the bundle of complete
 flags in $V$, then the projection map $F_E\xra{}Z_E$ is faithfully
 flat.  This idea was extended and used to great effect
 in~\cite{hokura:ggc}.  
\end{remark}

We next define some other useful properties of maps, which do not seem
to fit anywhere else.
\begin{definition}\label{defn-vflat}
 We say that a map $f\:X\xra{}Y$ is \emph{very flat}\index{flat!very}
 if it makes $\OX$ into a free module over $\OY$.  A very flat map is
 flat, and even faithfully flat provided that $X\neq\emptyset$.
\end{definition}
\begin{definition}\label{defn-finite}
 We say that a map $f\:X\xra{}Y$ is \emph{finite}\index{finite map} if
 it makes $\OX$ into a finitely generated module over $\OY$.
\end{definition}

\begin{remark}\label{rem-ff-zariski}
 A flat map $f\:X\xra{}Y$ is faithfully flat\index{flat!faithfully} if
 and only if the resulting map $f_\zar\:X_\zar\xra{}Y_\zar$ is
 surjective~\cite[Theorem 7.3]{ma:crt}.
\end{remark}

\begin{example}\label{eg-open-ff}
 An open inclusion $D(a)\xra{}X$ (where $a\in\OX$) is always flat.  If
 $a_1,\ldots,a_m\in\OX$ generate the unit ideal then
 $\coprod_kD(a_k)\xra{}X$ is faithfully flat.
\end{example}

\begin{example}\label{eg-divisor-ff}
 If $D$ is a \idx{divisor} on $\aff^1$ over $Y$ (as in
 Example~\ref{eg-divisors}) then $D\xra{}Y$ is very flat and thus
 faithfully flat.
\end{example}

\begin{definition}\label{defn-Omega}
 Given a ring $R$ and an $R$-algebra $S$, we write $I$ for the kernel
 of the multiplication map $S\ot_RS\xra{}S$, and $\Om^1_{S/R}=I/I^2$,
 \index{O1@$\Om^1$} which is a module over $S$.  Given a map of
 schemes $X\xra{}Y$, we define $\Om^1_{X/Y}=\Om^1_{\OX/\OY}$, which we
 think of as a sheaf over $X$.  We say that $X$ is \dfn{smooth} over
 $Y$ of relative dimension $n$ if the map $X\xra{}Y$ is flat and
 $\Om^1_{X/Y}$ is a vector bundle of rank $n$ over $X$ (we allow the
 case $n=\infty$).  In that case, we write $\Om^k_{X/Y}$ for the
 $k$'th exterior power of $\Om^1_{X/Y}$ over $\OX$, which is a vector
 bundle over $X$ of rank $\bcf{n}{k}$.
\end{definition}

\begin{remark}
 If $X$ and $Y$ are reduced affine algebraic varieties over $\cplx$,
 and $X$ is \idx{smooth} over $Y$ then the preimage of each point
 $y\in Y$ is a smooth variety of dimension independent of $y$.  The
 converse is probably not true but at least that is roughly the right
 idea.  It has nothing to do with the question of whether the map
 $X\xra{}Y$ is a smooth map of manifolds.  The latter only makes sense
 if $X$ and $Y$ are both smooth varieties (in other words, smooth over
 $\spec(\cplx)$), and in that case it holds automatically for any
 algebraic map $X\xra{}Y$.
\end{remark}

The following two propositions summarise the basic properties of
(faithfully) flat maps.
\begin{proposition}\label{prop-ff-composites}
 Let $X\xra{f}Y\xra{g}Z$ be maps of schemes.  Then:
 \begin{itemize}
 \item[(a)] If $f$ and $g$ are flat then $gf$ is flat.
 \item[(b)] If $f$ and $g$ are faithfully flat then $gf$ is faithfully
  flat.
 \item[(c)] If $f$ is faithfully flat and $gf$ is flat then $g$ is
  flat.
 \item[(d)] If $f$ and $gf$ are faithfully flat then $g$ is faithfully
  flat.
 \end{itemize}\index{flat}\index{flat!faithfully}
\end{proposition}
\begin{proof}
 All this follows easily from the definitions.
\end{proof}

\begin{proposition}\label{prop-ff-pullbacks}
 Suppose we have a pullback diagram of schemes
 \begin{diag}
  \node{W} \arw{s,l}{f} \arw{e,t}{r} \node{X} \arw{s,r}{g} \\
  \node{Y} \arw{e,b}{s} \node{Z.}
 \end{diag}
 Then:
 \begin{itemize}
 \item[(a)] If $s$ is flat then $r$ is flat.
 \item[(b)] If $s$ is faithfully flat then $r$ is faithfully flat.
 \item[(c)] If $g$ is faithfully flat and $r$ is flat then $s$ is
  flat.
 \item[(d)] If $g$ and $r$ are faithfully flat so $s$ is faithfully
  flat. 
 \end{itemize}\index{flat}\index{flat!faithfully}
\end{proposition}
\begin{proof}
 Consider the functor $f_*\:\CX_W\xra{}\CX_Y$, which sends a scheme
 $U\xra{u}W$ over $W$ to the scheme $U\xra{fu}Y$ over $Y$.  Colimits
 in $\CX_W$ are constructed by forming the colimit in $\CX$ and
 equipping it with the obvious map to $W$.  This means that $f_*$
 preserves and reflects colimits, as does $g_*$.  For any scheme $V$
 over $X$, we have $W\tX V=(Y\tZ X)\tX V=Y\tZ V$, or in other words
 $f_*r^*V=s^*g_*V$ in $\CX_Y$.  It follows that if $s^*$ preserves or
 reflects finite colimits then so does $r^*$, which gives~(a) and~(b).

 For part~(c), suppose that $g$ is faithfully flat and $r$ is flat.
 This implies that $sf=gr$ is flat.  Also, part~(b) says that any
 pullback of a faithfully flat map is faithfully flat, and $f$ is a
 pullback of $g$ so $f$ is faithfully flat.  As $sf$ is flat, part~(c)
 of the previous proposition tells us that $s$ is flat, as required.
 A similar argument proves~(d).
\end{proof}

\begin{proposition}\label{prop-ff-strong}
 Let $f\:X\xra{}Y$ be a faithfully flat map\index{flat!faithfully},
 and let $\{V_i\}$ be a finite diagram in $\CX_Y$.  If $\{f^*V_i\}$
 has a \idx{strong colimit} in $\CX_X$, then $\{V_i\}$ has a strong
 colimit in $\CX_X$.  In other words, $f^*$ reflects strong finite
 colimits.
\end{proposition}
\begin{proof}
 Write $V=\colim_iV_i$.  Given a map $g\:X'\xra{}X$, we need to show
 that $g^*V=\colim_ig^*V_i$.  To see this, form the pullback square
 \begin{diag}
  \node{Y'} \arrow{e,t}{f'} \arrow{s,l}{g'} \node{X'} \arrow{s,r}{g}\\
  \node{Y}  \arrow{e,b}{f}  \node{X.}
 \end{diag}
 We know from Proposition~\ref{prop-ff-pullbacks} that $f'$ is
 faithfully flat.  Because $f$ is flat, we have
 $f^*V=\colim_if^*V_i$.  By hypothesis, this colimit is strong, so
 $(g')^*f^*V=\colim_i(g')^*f^*V_i$.  As $gf'=fg'$, we have
 $(f')^*g^*V=\colim_i(f')^*g^*V_i$.  As $f'$ is faithfully flat, the
 functor $(f')^*$ reflects colimits, so $g^*V=\colim_ig^*V_i$ as
 required. 
\end{proof}

\begin{proposition}\label{prop-ff-reg-epi}
 If $f\:X\xra{}Y$ is faithfully flat \index{flat!faithfully} and
 $Y\xra{}Z$ is arbitrary then the diagram
 \[ X \tY  X \arw{e,=>} X \xra{f} Y \]
 is a \idx{strong coequaliser} in $\CX_Z$.
\end{proposition}
\begin{proof}
 As $f^*\:\CX_Y\xra{}\CX_X$ reflects strong coequalisers, it is enough
 to show that the above diagram becomes a strong coequaliser after
 applying $f^*$.  Explicitly, we need to show that the following is a
 strong coequaliser:
 \[ X \tY  X \tY  X \arw{e,tb,=>}{d_0}{d_1}
    X \tY  X \xra{d} X,
 \]
 where
 \begin{align*}
  d_0(a,b,c) &= (b,c)           \\
  d_1(a,b,c) &= (a,c)           \\
  d(a,b)     &= b.
 \end{align*}
 In fact, one can check that this is a split coequaliser, with
 splitting given by the maps 
 \[ X \tY  X \tY  X \xla{s} X \tY  X \xla{t} X,
 \]
 where
 \begin{align*}
  s(a,b) &= (a,b,b)     \\
  t(a)   &= (a,a).
 \end{align*}
 As split coequalisers are preserved by all functors, they are
 certainly strong coequalisers.
\end{proof}

Now suppose that $f\:X\xra{}Y$ is faithfully flat, and that $U$ is a
scheme over $X$.  We will need to know when $U$ descends to $Y$, in
other words when there is a scheme $V$ over $Y$ such that $U=V\tY X$.
Given a point $a\in X(R)$, we regard $a$ as a map $\spec(R)\xra{}X$
and write $U_a$ for the pullback of $U$ along this map, which is a
scheme over $\spec(R)$.

\begin{definition}\label{defn-descent-datum}
 Let $f\:X\xra{}Y$ be a map of schemes, and let $U$ be a scheme over
 $X$.  A \emph{system of descent data}\index{descent data} for $U$
 consists of a collection of maps $\tht_{a,b}\:U_a\xra{}U_b$ of
 schemes over $\spec(R)$, for any ring $R$ and any pair of points
 $a,b\in X(R)$ with $f(a)=f(b)$.  These maps are required to be
 natural in $(a,b)$, and to satisfy the cocycle conditions
 $\tht_{a,a}=1$ and $\tht_{a,c}=\tht_{b,c}\circ\tht_{a,b}$.

 We write $\CX_f$\index{Xf@$\CX_f$} for the category of pairs
 $(U,\tht)$, where $U$ is a scheme over $X$ and $\tht$ is a system of
 descent data.
\end{definition}
\begin{remark}\label{rem-descent-univ}
 The naturality condition for the maps $\tht_{a,b}$ just means that
 they give rise to a map $\pi_0^*U\xra{}\pi_1^*U$ of schemes over
 $X\tY X$.  
\end{remark}
\begin{remark}\label{rem-theta-iso}
 Note also that the cocycle conditions imply that
 $\tht_{a,b}\circ\tht_{b,a}=1$, so $\tht_{a,b}$ is an isomorphism.
\end{remark}

\begin{definition}\label{defn-descent-effective}
 If $V$ is a scheme over $Y$ and $f\:X\xra{}Y$ then there is an
 obvious system of descent data for $U=f^*V$, in which $\tht_{a,b}$ is
 the identity map of $U_a=V_{f(a)}=V_{f(b)}=U_b$.  We can thus
 consider $f^*$ as a functor $\CX_Y\xra{}\CX_f$.  We say that a system
 of descent data $\tht$ on $U$ is \emph{effective}\index{descent
   data!effective} if $(U,\tht)$ is equivalent to an object in the
 image of $f^*$.  It is equivalent to say that there is a scheme $V$
 over $Y$ and an isomorphism $\phi\:U\simeq f^*V$ such that
 \[ \tht_{a,b}= (U_a\xra{\phi}V_{f(a)} = V_{f(b)}\xra{\phi^{-1}}U_b)
 \]
 for all $(a,b)$.
\end{definition}

\begin{definition}\label{defn-descent-functor}
 Given a map $f\:X\xra{}Y$, a scheme $U\xra{g}X$ over $X$, and a
 system of \idx{descent data} $\tht$ for $U$, we define $U\xra{q}QU$
 to be the coequaliser of the maps $d_0,d_1\:U\tY X\xra{}U$ defined by
 \begin{align*}
  d_0(u,a) &= u                         \\
  d_1(u,a) &= \tht_{g(u),a}(u).
 \end{align*}
 We note that $d_0$ and $d_1$ have a common splitting
 $s\:u\mapsto(u,g(u))$, so we have a reflexive coequaliser.  We also
 note that there is a unique map $r\:QU\xra{}Y$ such that $rq=fg$, so
 we can think of $QU$ as a scheme over $Y$.
\end{definition}

\begin{proposition}[Faithfully flat descent]\label{prop-descent}
 If $f\:X\xra{}Y$ is faithfully flat\index{flat!faithfully}, then the
 functor $f^*\:\CX_Y\xra{}\CX_f$ is an equivalence, with inverse given
 by $Q$.  Moreover, the coequaliser in $\CX_Y$ that defines $QU$ is a
 strong coequaliser.
\end{proposition}
\begin{proof}
 Firstly, it is entirely formal to check that $Q$ is left adjoint to
 $f^*$.  Next, we claim that $Qf^*=1$, or in other words that the
 projection map $f^*V=V\tY X\xra{}V$ is a coequaliser of the maps
 $d_0,d_1\:V\tY X\tY X\xra{}V\tY X$.  Explicitly, we need to show
 that $(v,a)\mapsto v$ is the coequaliser of $(v,b,a)\mapsto(v,b)$ and
 $(v,b,a)\mapsto(v,a)$.  This is just the same as
 Proposition~\ref{prop-ff-reg-epi}.  Thus $Qf^*=1$ as claimed.

 We now show that $f^*QU=U$.  As $f^*$ preserves
 coequalisers, it will be enough to show that the projection
 $f^*U=U\tY X\xra{}U$ is the coequaliser of the fork
 $U\tY X\tY X\arw{e,tb,=>}{f^*d_0}{f^*d_1}U\tY X$.  More
 explicitly, we need to show that the map $(u,a)\mapsto u$ is the
 coequaliser of the maps $(u,a,b)\mapsto(u,b)$ and
 $(u,a,b)\mapsto(\tht_{g(u),a}(u),b)$.  In fact, it is a split
 coequaliser, with splitting given by the maps $u\mapsto(u,g(u))$ and
 $(u,a)\mapsto(u,a,a)$.  Thus, $f^*Q=1$ as claimed.  We also see that
 the coequaliser defining $QU$ becomes split and thus strong after
 applying $f^*$.  It follows from Proposition~\ref{prop-ff-strong}
 that it was a strong coequaliser in the first place.
\end{proof}
\begin{corollary}
 If $f\:X\xra{}Y$ is faithfully flat, then the functor
 $f^*\:\CX_Y\xra{}\CX_X$ is faithful. \qed
\end{corollary}

We also have a similar result for sheaves.
\begin{definition}
 Let $f\:X\xra{}Y$ be a map of schemes, and let $M$ be a sheaf over
 $X$.  A \emph{system of decent data}\index{descent data} for $M$
 consists of a collection of maps $\tht_{a,b}\:M_a\xra{}M_b$ of
 $R$-modules, for every ring $R$ and every pair of points
 $a,b\in X(R)$ with $f(a)=f(b)$.  These are supposed to be natural in
 $(a,b)$ and to satisfy the conditions $\tht_{a,a}=1$ and
 $\tht_{b,c}\circ\tht_{a,b}=\tht_{a,c}$.  We write
 $\Sheaves_f$\index{Sheavesf@$\Sheaves_f$} for the category of sheaves
 over $X$ equipped with descent data.  The pullback functor $f^*$ can
 be regarded as a functor from $\Sheaves_Y$ to $\Sheaves_f$.
\end{definition}

\begin{proposition}\label{prop-descent-sheaves}
 If $f$ is faithfully flat, then the functor
 $f^*\:\Sheaves_X\xra{}\Sheaves_f$ is an equivalence of categories.
\end{proposition}
The proof is similar to that of Proposition~\ref{prop-descent}, and is
omitted.

We shall say that a statement holds \emph{locally in the flat
  topology} \index{flat!topology} or \emph{fpqc locally} \index{fpqc}
if it is true after pulling back along a faithfully flat map. (fpqc
stands for fid\`{e}lement plat et quasi-compact; the compactness
condition is automatic for affine schemes).  Suppose that a certain
statement $S$ is true whenever it holds fpqc-locally.  We then say
that $S$ is an fpqc-local statement.

\begin{remark}\label{rem-topologies}
 Let $X$ be a topological space.  We say that a statement $S$ holds
 locally on $X$ if and only if there is an open covering
 $X=\bigcup_i U_i$ such that $S$ holds on each $U_i$.  Write
 $Y=\coprod_i U_i$, so $Y \xra{} X$ is a coproduct of open inclusions
 and is surjective.  We could call such a map an ``disjoint covering
 map''.  We would then say that $S$ holds locally if and only if it
 holds after pulling back along a disjoint covering map.  One can get
 many analogous concepts varying the class of maps in question.  For
 example, we could use covering maps in the ordinary sense.  In the
 category of compact smooth manifolds, we could use submersions.  This
 is the conceptual framework in which the above definition is supposed
 to fit.
\end{remark}

\begin{example}
 Suppose that $N$ is a sheaf on $Y$ which vanishes fpqc-locally.  This
 means that there is a faithfully flat map $f\:X\xra{}Y$ such that
 $\Gm(X,f^*N)= \OX\ot_\OY\Gm(Y,N)=0$.  By the classical definition of
 faithful flatness, this implies that $N=0$.  In other words, the
 vanishing of $N$ is an fpqc-local condition.  
\end{example}
\begin{example}
 Let $N$ be a sheaf over $Y$, and let $n$ be an element of $\Gm(Y,N)$
 that vanishes fpqc-locally.  This means that there is a faithfully
 flat map $f\:X\xra{}Y$ such that the image of $n$ in 
 $\Gm(X,f^*N)=\OX\ot_\OY\Gm(Y,N)$ is zero.  Let $g$ be the projection
 $X\tY X\xra{}Y$.  One can show that the diagram
 \[ \Gm(Y,N) \xra{f^*} \Gm(X,f^*N) \arw{e,=>} \Gm(X\tY X,g^*N) \]
 is an equaliser.  Indeed, it becomes split after tensoring with $\OX$
 over $\OY$, and that functor reflects equalisers by the classical
 definition of faithful flatness.  In particular, the map marked $f^*$
 is injective, so $n=0$.  Thus, the vanishing of $n$ is an fpqc-local
 condition.
\end{example}
\begin{example}\label{eg-locfree}
 Suppose that $M$ is a \idx{vector bundle} of rank $r$ over a scheme
 $X$.  We claim that $M$ is fpqc-locally free of rank $r$, in other
 words that there is a faithfully flat map $f\:W\xra{}X$ such that
 $f^*M\simeq\OO^r$.  To prove this, choose a matrix $A\in M_n(\OX)$
 such that $\Gm(X,M)=A.\OO_X^n$.  If $R$ is a ring and $x\in X(R)$
 then $A(x)\in M_n(R)$ and $M_x=A(x).R^n$.  Let $W(R)$ be the set of
 triples $(x,P,Q)$ such that $x\in X(R)$ and $P$ and $Q$ are matrices
 over $R$ of shape $r\tm n$ and $n\tm r$ such that $\det(PA(x)Q)$ is
 invertible.  This is easily seen to be a scheme over $X$.  In fact,
 it is an open subscheme of the scheme of all triples $(x,P,Q)$, which
 can be identified with $\aff^{2nr}\tm X$.  It follows that $W$ is
 flat over $X$.  Moreover, if $R$ is a field then elementary linear
 algebra tells us that the map $W(R)\xra{}X(R)$ is surjective, so that
 $W$ is faithfully flat over $R$.  If $(x,P,Q)$ is a point of $W$ then
 $A(x)Q\:R^r\xra{}M_x$ is a split monomorphism.  By comparison of
 ranks, it is an isomorphism.  It follows that $M$ becomes free after
 pulling back to $W$.
\end{example}
\begin{example}\label{eg-ff-ff-local}
 Proposition~\ref{prop-ff-pullbacks} tells us that flatness and
 faithful flatness are themselves fpqc-local properties.
\end{example}
\begin{example}\label{eg-bases-M}
 Let $M$ be a vector bundle of rank $r$ over a scheme $X$, as in
 Example~\ref{eg-locfree}.  Let $\Bases(M)$\index{BasesM@$\Bases(M)$}
 be the functor of pairs $(x,B)$ where $x$ is a point of $X$ and
 $B\:\OO_x^r\xra{}M_x$ is an isomorphism.  Note that $\Bases(M)(R)$
 can be identified with the set of tuples
 $(x,b_1,\ldots,b_r,\bt_1,\ldots,\bt_r)$ such that $b_i\in M_x$ and
 $\bt_j\in M_x^\vee$ and $\bt_j(b_i)=\dl_{ij}$, so $\Bases(M)$ is a
 closed subscheme of $\aff(M)^r_X\tX \aff(M^\vee)^r_X$.

 It is clear that $M$ becomes free after pulling back along the
 projection 
 \[ f\:\Bases(M)\xra{}X. \]
 If $M=\OO^r$ is free, then $\Bases(M)$ is just the scheme $\GL_r\tm
 X$, where $\GL_r$ is the scheme of invertible $r\tm r$ matrices.
 It's not hard to see that
 $\OO_{\GL_r}=\Zh[x_{i,j}\st 0\le i,j<r][\det(x_{ij})^{-1}]$ is
 torsion-free, and clearly $\GL_r(k)\neq\emptyset$ for all fields $k$,
 and one can conclude that the map $\GL_r\xra{}1=\spec(\Zh)$ is
 faithfully flat.  It follows that $\Bases(M)$ is faithfully flat over
 $X$ when $M$ is free.  Even if $M$ is not free, we see from
 Example~\ref{eg-locfree} that it is fpqc-locally free, so the map
 $\Bases(M)\xra{}X$ is fpqc-locally faithfully flat.  As remarked in
 Example~\ref{eg-ff-ff-local}, faithful flatness is itself a local
 condition, so $\Bases(M)\xra{}X$ is faithfully flat.
\end{example}
\begin{example}
 Any monic polynomial $f\in R[x]$ can be factored as a product of
 linear terms, locally in the flat topology.  Indeed, suppose
 \[ f=\sum_{0}^m(-1)^{m-k}a_{m-k} x^k \]
 with $a_0=1$.  It is well known that $S=\Zh[x_1,\ldots x_m]$ is free
 of rank $m!$ over $T=S^{\Sg_m}=\Zh[\sg_1,\ldots\sg_m]$, where $\sg_k$
 is the $k$'th elementary symmetric function in the $x$'s.  A basis is
 given by the monomials $x^\al=\prod x_k^{\al_k}$ for which $\al_k<k$.
 We can map $T$ to $R$ by sending $\sg_k$ to $a_k$, and then observe
 that $U=S\ot_TR$ is free and thus faithfully flat over $R$.  Clearly
 $f(x)=\prod_k (x-x_k)$ in $U[x]$, as required.
\end{example}

We conclude this section with some remarks about open mappings.  We
have to make a slightly twisted definition, because in our affine
context we do not have enough open subschemes.  Suppose that
$f\:X\xra{}Y$ is a map of spaces, and that $W\sse X$ is closed.  We
can then define $W'=\{y\in Y\st f^{-1}y\sse W\}=f(W^c)^c$.  Clearly $f$
is open iff ($W$ closed implies $W'$ closed).  We will define openness
for maps of schemes by analogy with this.  

\begin{definition}\label{defn-open-map}
 Let $f\:X\xra{}Y$ be a map of schemes.  For any closed subscheme
 $W\sse X$, we define a subfunctor $W'$ of $Y$ by
 \[ W'(R) = \{ y\in Y(R) \st W_y = X_y \}. \]
 We say that $f$ is \emph{open}\index{open map} if for every $W$, the
 corresponding subfunctor $W'\sse Y$ is actually a closed subscheme.
\end{definition}

\begin{proposition}\label{prop-free-open}
 A very flat map is open.\index{flat!very}\index{open map}
\end{proposition}
\begin{proof}
 Let $f\:X\xra{}Y$ be very flat.  Write $A=\OX$ and $B=\OY$, and
 choose a basis $A=B\{e_\al\}$.  Suppose that $W=V(I)$ is a closed
 subscheme of $X$.  Let $\{g_\bt\}$ be a system of generators of $I$,
 so we can write $g_\bt=\sum_\al g_{\bt\al}e_\al$ for suitable
 elements $g_{\al\bt}\in A$.  Consider a point $y\in Y(R)$,
 corresponding to a map $y^*\:B\xra{}R$.  This will lie in $W'(R)$ iff
 $R\ot_B~A=R\ot_B(A/I)$, iff the image of $I$ in $R\ot_BA=R\{e_\al\}$
 is zero.  This image is generated by the elements
 $h_\bt=\sum_\al~y^*(g_{\bt\al})e_\al$.  Thus, it vanishes iff
 $y^*(g_{\bt\al})=0$ for all $\al$ and $\bt$.  This shows that
 $W'=V(I')$, where $I'=(g_{\bt\al})$, so $W'$ is a closed subscheme as
 required.
\end{proof}

\subsection{Schemes of maps}
\label{subsec-map-schemes}

\begin{definition}\label{defn-maps-informal}
 Let $Z$ be a functor $\Rings\xra{}\Sets$, and let $X$ and $Y$ be
 functors over $Z$.  For any ring $R$, we let $\Map_Z(X,Y)(R)$
 \index{MapZXY@$\Map_Z(X,Y)$} be the class of pairs $(z,u)$, where
 $z\in Z(R)$ and $u\:X_z\xra{}Y_z$ is a map of functors over
 $\spec(R)$.  If this is a set (rather than a proper class) for all
 $R$, then we get a functor $\Map_Z(X,Y)\in\CF$.  This is clearly the
 case whenever $X$, $Y$ and $Z$ are all schemes.  However, the functor
 $\Map_Z(X,Y)$ need not itself be a scheme.

 When $Z=1$ is the terminal scheme we will usually write $\Map(X,Y)$
 rather than $\Map_1(X,Y)$.
\end{definition}
\begin{remark}
 It is formal to check that 
 \[ \CF_Z(W,\Map_Z(X,Y)) =
    \CF_Z(W\tZ X,Y) = 
    \CF_W(W\tZ X,W\tZ Y).
 \]
 In particular, if $X$, $Y$, $Z$ and $\Map_Z(X,Y)$ are all schemes
 then we have
 \[ \CX_Z(W,\Map_Z(X,Y)) =
    \CX_Z(W\tZ X,Y) = 
    \CX_W(W\tZ X,W\tZ Y).
 \]
\end{remark}

\begin{example}\label{eg-map-aff}
 It is not hard to see that maps
 $\aff^n\tm\spec(R)\xra{}\aff^m\tm\spec(R)$ over $\spec(R)$ biject
 with $m$-tuples of polynomials over $R$ in $n$ variables, so
 $\Map(\aff^n,\aff^m)(R)=R[x_1,\ldots,x_n]^m$, which is isomorphic to
 $\bigoplus_{n\in\nat}R$ (naturally in $R$).  This functor is not a
 representable (it does not preserve infinite products, for example)
 so $\Map(\aff^n,\aff^m)$ is not a scheme.  It is a formal scheme,
 however.  
\end{example}
\begin{example}\label{eg-Map-Dn}
 Write $D(n)(R)=\{a\in R\st a^{n+1}=0\}$, so
 \[ D(n)=\spec(\Zh[x]/x^{n+1}) \]
 is a scheme.  We find that
 $\Map(D(n),\aff^1)(R)=R[x]/x^{n+1}\simeq\prod_{i=0}^n R$, so that
 $\Map(D(n),\aff^1)\simeq\aff^{n+1}$ is a scheme.
\end{example}
\begin{example}
 Let $E$ be an even periodic ring spectrum.  As $\Om U(n)$ is a
 commutative $H$-space, we see that $E_0(\Om U(n))$ is a ring, so we
 can define a scheme $\spec(E_0(\Om U(n)))$.  We will see later that
 there is a canonical isomorphism
 \[ \spec(E_0(\Om U(n))) \simeq \Map_{S_E}((\cp^{n-1})_E,\MG). \]
\end{example}

We now give a proposition which generalises the last two examples.
\begin{proposition}\label{prop-map-informal}
 Let $Z$ be a scheme and let $X$ and $Y$ be schemes over $Z$, and
 suppose that $X$ is finite and very flat over $Z$.  Then
 $\Map_Z(X,Y)$ is a scheme. \index{finite map}\index{flat!very}
\end{proposition}
\begin{proof}
 Let $R$ be a ring, and $z$ a point of $Z(R)$, giving a map
 $\hat{z}\:\OZ\xra{}R$.  We need to produce an algebra $B$ over $\OZ$
 such that the maps $B\xra{}R$ of $\OZ$-algebras biject with maps
 $X_z\xra{}Y_z$ of schemes over $\spec(R)$, or equivalently with maps
 $R\ot_\OZ\OY\xra{}R\ot_\OZ\OX$ of $R$-algebras, or equivalently with
 maps $\OY\xra{}R\ot_\OZ\OX$ of $\OZ$-algebras.

 Write $\OXv=\Hom_\OZ(\OX,\OZ)$ and
 $A=\Sym_{\OZ}[\OXv\otimes_\OZ\OY]$.  Then
 \[ \Alg_\OZ(A,R) = \Hom_\OZ(\OXv\ot_\OZ\OY,R) =
     \Hom_\OZ(\OY,R\ot_\OZ\OX).
 \]
 A suitable quotient $B$ of $A$ will pick out the algebra maps from
 $\OY$ to $\OW\ot_\OZ\OX$.  To be more explicit, let
 $\{e_1,\ldots,e_n\}$ be a basis for $\OX$ over $\OZ$, with
 $1=\sum_id_ie_i$ and $e_ie_j=\sum_kc_{ijk}e_k$.  Let $\{\ep_i\}$ be
 the dual basis for $\OXv$.  Then $B$ is $A$ mod the relations
 \begin{align*}
  \ep_k \ot ab &= \sum_{i,j} c_{ijk}(\ep_i\ot a)(\ep_j\ot b)   \\
  \ep_i \ot 1  &= d_i.
 \end{align*}
 More abstractly, if we write $q$ for the projection $A\xra{}B$ and
 $j$ for the inclusion $\OXv\ot\OY\xra{}A$, then $B$ is the
 largest quotient of $A$ such that the following diagrams commute:
 \[
 \renewcommand{\typeout}[1]{}
 \begin{diagram}
  \node{\OY\ot\OY\ot\OXv}
  \arw{e,t}{\mu_Y\ot 1}
  \arw{s,l}{1\ot\mu_X^\vee}
  \node{\OY\ot\OXv}
  \arw[3]{s,r}{qj}\\
  \node{\OY\ot\OY\ot\OXv\ot\OXv}
  \arw{s,l}{\text{twist}}\\
  \node{\OY\ot\OXv\ot\OY\ot\OXv}
  \arw{s,l}{qj\ot qj}\\
  \node{B \ot B}
  \arw{e,b}{\mu_B}
  \node{B}
 \end{diagram}\hspace{8em}
 \begin{diagram}
  \node{\OXv}
  \arw{e,t}{\eta^\vee}
  \arw{s,l}{1\ot\eta}
  \node{\OZ}
  \arw{s,r}{\eta}\\
  \node{\OXv\ot\OY}
  \arw{e,b}{qj}
  \node{B}
 \end{diagram}
 \] 
 We conclude that $\spec(B)$ has the defining property of
 $\Map_Z(X,Y)$.
\end{proof}

\subsection{Gradings}
\label{subsec-gradings}

In this section, we show that graded rings are essentially the same as
schemes with an action of the multiplicative group $\MG$.
\begin{definition}
 A \dfn{grading} of a ring $R$ is a system of additive subgroups
 $R_k\leq R$ for $k\in\Zh$ such that $R=\bigoplus_kR_k$ and $1\in R_0$
 and $R_jR_k\sse R_{j+k}$ for all $j,k$.  We say that a  map
 $g\:R\xra{}S$ between graded rings is \emph{homogeneous} if
 $g(R_k)\sse S_k$ for all $k$.
\end{definition}

\begin{proposition}\label{prop-grading-action}
 Let $X$ be a scheme.  Then gradings of $\OX$ biject with actions of
 the group scheme $\MG$ on $X$.  Given such actions on $X$ and $Y$, a
 map $f\:X\xra{}Y$ is $\MG$-equivariant if and only if the
 corresponding map $\OY\xra{}\OX$ is homogeneous.
\end{proposition}
\begin{proof}
 Given an action of $\MG$ on $X$, we define $(\OX)_k$ to be the set
 of maps $f\:X\xra{}\aff^1$ such that $f(u.x)=u^kf(x)$ for all rings
 $R$ and points $u\in\MG(R)$, $x\in X(R)$.  It is clear that
 $1\in(\OX)_0$ and that $(\OX)_j(\OX)_k\sse(\OX)_{j+k}$.  We need to
 check that $\OX=\bigoplus_k(\OX)_k$.  For this, we consider the map
 $\al^*\:\OX\xra{}\OO_{\MG\tm X}=\OX[u^{\pm 1}]$.  If
 $\al^*(f)=\sum_ku^kf_k$ (so $f_k=0$ for almost all $k$), then we find
 that the $f_k$ are the unique functions $X\xra{}\aff^1$ such that
 $f(u.x)=\sum_ku^kf_k(x)$ for all $u$ and $x$.  By taking $u=1$, we
 see that $f=\sum_kf_k$.  We also find that 
 \[ \sum_k u^k v^k f_k(x) = f((uv).x) = 
    f(u.(v.x)) = \sum_{j,k} u^j v^k f_{kj}(x).
 \]
 By working in the universal case $R=\OX[u^{\pm 1},v^{\pm 1}]$ and
 comparing coefficients, we see that $f_{kj}=\dl_{jk}f_k$ so that
 $f_k\in(\OX)_k$.  It follows easily that the addition map
 $\bigoplus_k(\OX)_k\xra{}\OX$ is an isomorphism, with inverse
 $f\mapsto(f_k)_{k\in\Zh}$.  Thus, we have a grading of $\OX$.

 Conversely, suppose we have a grading $(\OX)_*$.  We can then write any
 element $f\in\OX$ as $\sum_kf_k$ with $f_k\in(\OX)_k$ and $f_k=0$ for
 almost all $k$.  We define $\al^*(f)=\sum_ku^kf_k$, and check that
 this gives a ring map $\OX\xra{}\OX[u^{\pm 1}]$.  One can also check
 that $\al=\spec(\al^*)\:\MG\tm X\xra{}X$ is an action, and that this
 construction is inverse to the previous one.
\end{proof}

\begin{example}\label{eg-grading-FGL}
 Recall the scheme $\FGL$ from Example~\ref{eg-FGL}\index{FGL@$\FGL$}.
 We can let $\MG$ act on $\FGL$ by $(u.F)(x,y)=uF(x/u,y/u)$; this
 gives a grading of $\OO_\FGL$.  Write
 $F(x,y)=\sum_{i,j}a_{ij}(F)x^iy^j$, and recall that the elements
 $a_{ij}$ generate $\OO_\FGL$.  It is clear that
 $(u.F)(x,y)=\sum_{i,j}u^{1-i-j}a_{ij}(F)x^iy^j$, so that
 $a_{ij}(u.F)=u^{1-i-j}a_{ij}(F)$, so $a_{ij}$ is homogeneous of
 degree $1-i-j$.  This is of course the same as the grading coming
 from the isomorphisms $\OO_\FGL=\pi_0MP=\pi_*MU$, except that all
 degrees are halved.
\end{example}

\section{Non-affine schemes}
\label{sec-non-affine}

Let $\CE$\index{E@$\CE$} be the category of (not necessarily affine)
schemes in the classical sense\index{scheme!non-affine}, as discussed
in~\cite{ha:ag} for example.  In this section we show that $\CE$ can
be embedded as a full subcategory of $\CF$, containing our category
$\CX$ of affine schemes.  We show that our definition of sheaves over
functors gives the right answer for functors coming from non-affine
schemes, and we investigate the schemes $\BP^n$ from this point of
view.  This theory is useful in topology when one wants to study
elliptic cohomology, for example~\cite{hoanst:esw}.  The results here
are surely known to algebraic geometers, but I do not know a
reference. 

Given a ring $A$, we write $\zar(A)$\index{zarA@$\zar(A)$} for the
Zariski spectrum of $A$, considered as an object of $\CE$ in the usual
way.  The results of this section will allow us to identify $\zar(A)$
with $\spec(A)$.  Of course, in most treatments, $\spec(A)$ is defined
to be what we call $\zar(A)$.

\begin{definition}
 Given a scheme $X\in\CE$, we define a functor $FX\in\CF$ by 
 \[ FX(R) = \CE(\zar(R),X). \]
 It is well-known that
 \[ \CE(\zar(R),\zar(A))=\Rings(A,R), \]
 so that $F(\zar(A))=\spec(A)$.
\end{definition}

\begin{proposition}\label{prop-non-affine}
 The functor $F\:\CE\xra{}\CF$ is full and faithful.
\end{proposition}
\begin{proof}
 Let $X,Y\in\CE$ be schemes; we need to show that the map
 $F\:\CE(X,Y)\xra{}\CF(FX,FY)$ is an isomorphism.  First suppose that
 $X$ is affine, say $X=\zar(A)$.  Then the Yoneda lemma tells us that 
 \[ \CF(FX,FY)=\CF(\spec(A),FY)=FY(A)=\CE(\zar(A),Y)=\CE(X,Y) \]
 as required.

 Now let $X$ be an arbitrary scheme.  We can cover $X$ by open affine
 subschemes $X_i$, and for each $i$ and $j$ we can cover $X_i\cap X_j$
 by open affine subschemes $X_{ijk}$.  This gives rise to a diagram as
 follows.
 \begin{diag}
  \node{\CE(X,Y)}   
   \arrow{e,V}
   \arrow{s,l}{F}
  \node{\textstyle{\prod_i\CE(X_i,Y)}}
   \arrow{e,=>}
   \arrow{s,lr}{F}{\simeq}
  \node{\textstyle{\prod_{ijk}\CE(X_{ijk},Y)}}
   \arrow{s,lr}{F}{\simeq} \\
  \node{\CF(FX,FY)}
   \arrow{e,b}{J} 
  \node{\textstyle{\prod_i\CF(FX_i,FY)}}
   \arrow{e,=>}
  \node{\textstyle{\prod_{ijk}\CF(FX_{ijk},FY).}}
 \end{diag}
 Standard facts about the category $\CE$ show that the top line is an
 equaliser.  The affine case of our proposition shows that the middle
 and right-hand vertical arrows are isomorphisms.  If we can prove
 that the map $J$ is injective, then a diagram chase will show that
 the left-hand vertical map is an isomorphism, as required.

 Suppose we have two maps $f,g\:FX\xra{}FY$ and that $Jf=Jg$, or in
 other words $f|_{FX_i}=g|_{FX_i}$ for all $i$.  We need to show that
 $f=g$.  Consider a ring $R$ and a point $x\in FX(R)$, or equivalently
 a map $W=\zar(R)\xra{x}X$.  We need to show that $f(x)=g(x)$ as maps
 from $W$ to $Y$.  We can cover $W$ by open affine subschemes $W_s$
 such that $x\:W_s\xra{}X$ factors through $X_i$ for some $i$.  As
 $f|_{FX_i}=g|_{FX_i}$, we see that $f(x)\circ j_s=g(x)\circ j_s$,
 where $j_s\:W_s\xra{}W$ is the inclusion.  As the schemes $W_s$ cover
 $W$, we see that $f(x)=g(x)$ as required.
\end{proof}

\begin{proposition}\label{prop-non-affine-sheaves}
 Let $X\in\CE$ be a scheme.  Then the category of quasicoherent
 sheaves of $\OO$-modules over $X$ is equivalent to the category of
 sheaves over $FX$.\index{sheaf}
\end{proposition}
\begin{proof}
 Let $M$ be a quasicoherent sheaf of $\OO$-modules over $X$.  Consider
 a ring $R$ and a point $x\in FX(R)$, corresponding to a map
 $x\:\zar(R)\xra{}X$.  We can pull $M$ back along this map to get a
 quasicoherent sheaf of $\OO$-modules over $\zar(R)$, whose global
 sections form a module $G(M)_x=\Gm(\zar(R),x^*M)$ over $R$.  It is
 not hard to see that this construction gives a sheaf $GM$ over the
 functor $FX$.  If $X$ is affine then we know from
 Proposition~\ref{prop-Gamma-equiv} that sheaves over $FX$ are the
 same as modules over $\OX$, and it is classical that these are the
 same as quasicoherent sheaves of $\OO$-modules over $X$, so the
 functor $G$ is an equivalence in this case.  

 Now let $X\in\CE$ be an arbitrary scheme, and let $N$ be a sheaf over
 $FX$.  We can cover $X$ by open affine subschemes $X_i$, and we can
 cover $X_i\cap X_j$ by open affine subschemes $X_{ijk}$.  By the
 affine case of the proposition, we can identify $N_i=N|_{FX_i}$ with
 a quasicoherent sheaf $M_i$ of $\OO$-modules over $X_i$.  The obvious
 isomorphism $N_i|_{FX_{ijk}}=N_j|_{FX_{ijk}}$ gives an isomorphism
 $M_i|_{X_{ijk}}=M_j|_{X_{ijk}}$ (because our functor $G$ is an
 equivalence for the affine scheme $X_{ijk}$).  One checks that these
 isomorphisms satisfy the relevant cocycle condition, so we can glue
 together the sheaves $M_i$ to get a quasicoherent sheaf $M$ over $X$.
 One can also check that this construction is inverse to our previous
 one, which implies that $G$ is an equivalence of categories.
\end{proof}

{}From now on we will not usually distinguish between $X$ and $FX$.

We next examine how projective spaces\index{projective space} fit into
our framework.  Let $\BP^n$\index{Pn@$\BP^n$} be the scheme obtained
by gluing together $n+1$ copies of $\aff^n$ in the usual way.  In more
detail, we consider the scheme $\aff^{n+1}=\prod_{i=0}^n\aff^1$, and
let $U_i$ be the closed subscheme where $x_i=1$, so $U_i\simeq\aff^n$.
If $j\neq i$ we let $V_{ij}$ be the open subscheme of $U_i$ where
$x_j$ is invertible.  We define $\phi_{ij}\:V_{ij}\xra{}V_{ji}$ by
\[ \phi_{ij}(x_0,\ldots,x_n) = (x_0,\ldots,x_n)/x_j. \]
We use these maps to glue the $U_i$'s together to get a scheme
$\BP^n$.  

We define a sheaf $L_i$ over $U_i$ by 
$L_{i,\un{a}}=R\un{a}\leq R^{n+1}$ for $\un{a}\in U_i(R)$.  Note that
if $\pi_i\:R^{n+1}\xra{}R$ is the $i$'th projection then $\pi_i$
induces an isomorphism $L_{i,\un{a}}\xra{}R$, so $L_{i,\un{a}}$ is a
line bundle over $U_i$.  If $\un{a}\in V_{ij}(R)$ then it is clear
that $L_{i,\un{a}}=L_{j,\phi_{ij}(\un{a})}$.  It follows that the
bundles $L_i$ glue together to give a line bundle $L$ over $\BP^n$.
{}From the construction, we see that there is a short exact sequence
$L\mra\OO^{n+1}\era V$, in which $V$ is a vector bundle of rank $n$.
We also write $\OO(k)$ for the $(-k)$'th tensor power of $L$, which is
again a line bundle over $\BP^n$.

\begin{proposition}\label{prop-proj-scheme}
 For any ring $R$, we can identify $\BP^n(R)=\CE(\zar(R),\BP^n)$ with
 the set of submodules $M\leq R^{n+1}$ such that $M$ is a summand and
 has rank one.  
\end{proposition}
This will be proved after a lemma.

\begin{definition}
 Write $Q^n(R)$ for the set of submodules $M\leq R^{n+1}$ such that
 $L$ is a rank-one projective module and a summand, or equivalently
 $R^{n+1}/M$ is a projective module of rank $n$.  Given a map
 $R\xra{}R'$ we have a map $Q^n(R)\xra{}Q^n(R')$ sending $M$ to
 $R'\ot_RM$, which makes $Q^n$ into a functor.  

 We now define a map $\gm\:\BP^n\xra{}Q^n$, which will turn out to be
 an isomorphism.  Consider a ring $R$ and a point $x\in\BP^n(R)$,
 corresponding to a map $x\:\spec(R)\xra{}\BP^n$.  By pulling back the
 sequence $L\mra\OO^{n+1}\era V$ and identifying sheaves over
 $\spec(R)$ with $R$-modules, we get a short exact sequence
 $x^*L\mra R^{n+1}\era x^*V$.  Here $x^*L$ and $x^*V$ are projective,
 with ranks one and $n$ respectively, so $x^*L\in Q^n(R)$.  We define
 $\gm(x)=x^*L$.
\end{definition}

\begin{lemma}\label{lem-Qn-equaliser}
 Let $W$ be an affine scheme, and let $W_1,\ldots,W_m$ be a finite
 cover of $W$ by basic affine open subschemes $W_i=D(a_i)$.  Then
 there is an equaliser diagram
 \[ \CF(W,Q^n) \xra{} \prod_i \CF(W_i,Q^n)
       \arrow{e,=>} \prod_{ij} \CF(W_i\cap W_j,Q^n).
 \]
\end{lemma}
\begin{proof}
 Write $W'=\coprod_iW_i$ and $W''=\coprod_{ij}W_i\cap W_j$, so that
 the evident map $f\:W'\xra{}W$ is faithfully flat and
 $W''=W'\tm_WW'$.  We can thus use
 Proposition~\ref{prop-descent-sheaves} to identify $\Sheaves_W$ with
 the category $\Sheaves_f$ of sheaves on $W'$ equipped with descent
 data.  It follows that for any sheaf $F$ on $W$, the subsheaves of
 $F$ biject with subsheaves $K\leq f^*F$ that are preserved by the
 descent data for $f^*F$.  This condition is equivalent to the
 condition $\pi_0^*K=\pi_1^*K\leq(f\pi_0)^*F=(f\pi_1)^*F$.  Now take
 $F=\OO^{n+1}$, and the lemma follows easily.
\end{proof}

\begin{proof}[Proof of Proposition~\ref{prop-proj-scheme}]
 Suppose we have two points $x\in U_i(R)\subset\BP^n(R)$ and
 $y\in U_j(R)\subset\BP^n(R)$, and that $\gm(x)=\gm(y)$.  It then
 follows easily from the definitions that $x=y$.

 Now suppose we have two points $x,y\in\BP^n(R)$ such that
 $\gm(x)=\gm(y)$.  We write $W=\spec(R)$, so $x\:W\xra{}\BP^n$.  We
 can cover $W$ by basic affine open subsets $W_1,\ldots,W_m$ with the
 property that each $x(W_k)$ is contained in some $U_i$, and each
 $y(W_k)$ is contained in some $U_j$.  This implies (by the previous
 paragraph) that $x=y$ as maps $W_k\xra{}\BP^n$.  We can now deduce
 from Lemma~\ref{lem-Qn-equaliser} that $x=y$.  Thus,
 $\gm\:\BP^n(R)\xra{}Q^n(R)$ is always injective.

 Now consider a point $M\in Q^n(R)$, so $M$ is a sheaf over
 $W=\spec(R)$.  We claim that we can cover $W$ by basic open
 subschemes $V$ such that $M|_V$ lies in the image of
 $\gm\:\CF(V,\BP^n)\xra{}\CF(V,Q^n)$.  Indeed, as $M$ is projective,
 we can start by covering $W$ with basic open subschemes on which $M$
 is free.  It is easy to see that over such a subscheme, there exist
 maps $\OO\xra{u}\OO^{n+1}\xra{v}\OO$ such that the image of $u$ is
 $M$ and $vu=1$.  If we write $u$ and $v$ in terms of bases in the
 obvious way then $\sum_i u_iv_i=1$, so the elements $u_i$ generate
 the unit ideal, so the basic open subschemes $D(u_i)$ form a
 covering.  On $D(u_i)$ we can define
 $x=(u_0,\ldots,u_n)/u_i\in U_i$, and it is clear that $\gm(x)=M$.  

 We can thus choose a basic open covering $W=W_1\cup\ldots\cup W_m$
 and maps $x_k\:W_k\xra{}\BP^n$ such that $\gm(x_k)=M|_{W_k}$.  Let
 $x_{jk}$ be the restriction of $x_j$ to $W_{jk}=W_j\cap W_k$.  We
 then have $\gm(x_{jk})=M|_{W_{jk}}=\gm(x_{kj})$ and $\gm$ is
 injective so $x_{jk}=x_{kj}$.  We also have a diagram 
 \begin{diag}
  \node{\CF(W,\BP^n)}             
   \arrow{e,V}
   \arrow{s,l}{\gm}
  \node{\textstyle{\prod_i\CF(W_i,\BP^n)}}
   \arrow{e,=>}
   \arrow{s,l}{\gm}
  \node{\textstyle{\prod_{ij}\CF(W_{ij},\BP^n)}}
   \arrow{s,r}{\gm} \\
  \node{\CF(W,Q^n)}   
   \arrow{e}
  \node{\textstyle{\prod_i\CF(W_i,Q^n)}}
   \arrow{e,=>}
  \node{\textstyle{\prod_{ij}\CF(W_{ij},Q^n).}}
 \end{diag}
 The top row is unchanged if we replace $\CF$ by $\CE$, and this makes
 it clear that it is an equaliser diagram.  The bottom row is an
 equaliser diagram by Lemma~\ref{lem-Qn-equaliser}.  We have already
 seen that the vertical maps are injective.  The elements $x_i$ give
 an element of $\prod_i\CF(W_i,\BP^n)$, whose image in
 $\prod_i\CF(W_i,Q^n)$ is the same as that of $M\in\CF(W,Q^n)$.  We
 conclude by diagram chasing that there is an element
 $x\in\CF(W,\BP^n)$ such that $\gm(x)=M$.  Thus $\gm$ is also
 surjective, as required.
\end{proof}

\begin{definition}
 Suppose that we have elements $a_0,\ldots,a_n\in R$, which generate
 the unit ideal, say $\sum_i b_i a_i=1$.  Let $M$ be the submodule of
 $R^{n+1}$ generated by $\un{a}=(a_0,\ldots,a_n)$.  The elements $b_j$
 define a map $R^{n+1}\xra{}R$ which carries $L$ isomorphically to
 $R$.  It follows that $M\in Q^n(R)$; the submodules $M$ that occur in
 this way are precisely those that are free over $R$.  We write
 $[a_0:\ldots:a_n]$ for the corresponding point of $\BP^n(R)$.  Most
 of the time, when working with points of $\BP^n$\index{Pn@$\BP^n$},
 we can assume that they have this form, and handle the general case
 by localising.
\end{definition}

We finish this section with a useful lemma.

\begin{lemma}
 We have $[a_0:\ldots:a_n]=[a'_0:\ldots:a'_n]$ if and only if there is
 a unit $u\in R^\tm$ such that $ua'_j=a_j$ for all $j$, if and only if
 $a_i a'_j=a_j a'_i$ for all $i$ and $j$.
\end{lemma}
\begin{proof}
 The first equivalence is clear if we think in terms of $Q^n(R)$.  For
 the second, suppose that $ua'_j=a_j$ for all $j$.  Then
 $a_i a'_j=u^{-1}a_i a_j=a'_i a_j$ as required.  Conversely, suppose
 that $a_i a'_j=a_j a'_i$ for all $i$ and $j$.  We can choose
 sequences $b_0,\ldots,b_n$ and $b'_0,\ldots,b'_n$ such that
 $\sum_ia_ib_i=1$ and $\sum_ia'_ib'_i=1$.  Now define
 $u=\sum_ia_ib'_i$ and $v=\sum_ja'_jb'_j$.  Then
 \[ u a'_j = \sum_i b'_i a_i a'_j = \sum_i b'_i a'_i a_j = a_j. \]
 Moreover, we have
 \[ u \sum_j b_j a'_j = \sum_j b_j a_j = 1, \]
 so $u$ is a unit as required.
\end{proof}

\section{Formal schemes}
\label{sec-formal-schemes}

In this section we define formal schemes, and set up an extensive
categorical apparatus for dealing with them, and generalise our
results for schemes to formal schemes as far as possible.  We define
the subcategory of solid formal schemes, which is convenient for some
purposes.  We also define functors from various categories of
coalgebras to the category of formal schemes, which are useful
technical tools.  Finally, we study the question of when $\Map_Z(X,Y)$
is a formal scheme.

\begin{definition}
 A \emph{formal scheme}\index{scheme!formal} is a functor
 $X\:\Rings\xra{}\Sets$ that is a small filtered colimit of schemes.
 More precisely, there must be a small filtered category $\CI$ and a
 functor $i\mapsto X_i$ from $\CI$ to $\CX\sse\CF=[\Rings,\Sets]$ such
 that $X=\colim_iX_i$ in $\CF$, or equivalently $X(R)=\colim_iX_i(R)$
 for all $R$.  We call such a diagram $\{X_i\}$ a \dfn{presentation}
 of $X$.  We write $\FX$\index{Xh@$\FX$} for the category of formal
 schemes.
\end{definition}

\begin{example}
 The most basic example is the functor $\haf^1$\index{AA1h@$\haf^1$}
 defined by $\haf(R)=\Nil(R)$.  This is clearly the colimit over $N$
 of the functors $D(N)=\spec(\Zh[x]/x^{N+1})$.  We also define
 $\haf^n(R)=\Nil(R)^n$.\index{AAnh@$\haf^n$}
\end{example}
\begin{example}
 More generally, given a scheme $X$ and a closed subscheme $Y=V(I)$,
 we define a formal scheme
 $X^\wedge_Y=\colim_NV(I^N)$.\index{XY@$X^\wedge_Y$} 
\end{example}
\begin{example}\label{eg-haf-infty}
 For a common example not of the above type, consider the functor
 $\haf^{(\infty)}(R)=\bigoplus_{n\in\nat}\Nil(R)$, so
 $X=\colim_n\haf^n$, which is again a formal scheme.
\end{example}
\begin{example}\label{eg-ZE-formal}
 If $Z$ is an infinite CW complex and $\{Z_\al\}$ is the collection of
 finite subcomplexes and $E$ is an even periodic ring spectrum, we
 define $Z_E=\colim_\al(Z_\al)_E$.  This is clearly a formal scheme.
\end{example}

We can connect this with the framework of~\cite[Section
8]{grve:sgaivei} by taking $\CC$ to be the category $\Rings^\op$.
{}From this point of view, a formal scheme is an ind-representable
contravariant functor from $\Rings^\op$ to $\Sets$.  We shall omit any
mention of universes here, leaving the set-theoretically cautious
reader to lift the necessary details
from~\cite[Appendice]{grve:sgaivei}, or to avoid the problem in some
other way.

Given two filtered diagrams $X\:\CI\xra{}\CX$ and $Y\:\CJ\xra{}\CX$ we
know from~\cite[8.2.5.1]{grve:sgaivei} that 
\[ \FX(\colim_iX_i,\colim_jY_j) = \invlim_i\colim_j\CX(X_i,Y_j).\]
It follows that $\FX$ is equivalent to the category whose objects are
pairs $(\CI,X)$ and whose morphisms are given by the above formula.
We will feel free to use either model for $\FX$ where convenient.

\begin{proposition}\label{prop-ind-cart}
 A functor $X\:\Rings\xra{}\Sets$ is a formal scheme
 \index{scheme!formal} if and only if
 \begin{itemize}
 \item[(a)] $X$ preserves finite limits, and
 \item[(b)] There is a set of schemes $X_i$ and natural maps
  $X_i\xra{}X$ such that the resulting map $\coprod_iX_i(R)\xra{}X(R)$
  is surjective for all $R$.
 \end{itemize}
\end{proposition}
\begin{proof}
 This is essentially~\cite[Th\'eor\`eme 8.3.3]{grve:sgaivei}.  To see
 this, let $\CD$ be the category of schemes over $X$.  A map
 $\spec(R)\xra{}X$ is the same (by Yoneda) as an element of $X(R)$, so
 $\CD^\op$ is equivalent to the category $\Points(X)$.  This category
 corresponds to the category $\CC_{/F}$ of the cited theorem.  Thus,
 by the equivalence (i)$\iffa$(iii) of that theorem, we see that $X$
 is a formal scheme if and only if $X$ preserves finite limits, and
 $\CD$ has a small cofinal subcategory.  (Grothendieck actually talks
 about finite colimits, but in our case that implicitly refers to
 colimits in $\Rings^\op$ and thus limits in $\Rings$.)  It is shown
 in the proof of the theorem that if $X$ preserves finite limits, then
 $\CD$ is a filtered category, so we can
 use~\cite[Proposition~8.1.3(c)]{grve:sgaivei} to recognise cofinal
 subcategories.  This means that a small collection $\{X_i\}$ of
 schemes over $X$ gives a cofinal subcategory if and only if each map
 from a scheme $Y$ to $X$ factors through some $X_i$.  By writing
 $Y=\spec(R)$ and using the Yoneda lemma, it is equivalent to say that
 the map $\coprod_iX_i(R)\xra{}X(R)$ is surjective for all $R$.
\end{proof}

\subsection{(Co)limits of formal schemes}
\label{subsec-formal-colim}

\begin{proposition}\label{prop-formal-colim}
 The category $\FX$ has all small colimits.  The inclusion
 $\CX\xra{}\FX$ preserves finite colimits, and the inclusion
 $\FX\xra{}\CF=[\Rings,\Sets]$ preserves filtered colimits.  Moreover,
 if $X\in\CX$ then the functor $\FX(X,-)\:\FX\xra{}\Sets$ also
 preserves colimits.
\end{proposition}
\begin{proof}
 Apart from the last sentence, the proof is the same as that
 of~\cite[Theorem VI.1.6]{jo:ss}.  Johnstone assumes that $\CC$ (which
 is $\Rings^\op$ in our case) is small, but he does not really use
 this.  The last sentence is~\cite[Lemma VI.1.8]{jo:ss}.
\end{proof}
\begin{example}
 It is not hard to see that the functor $Z\mapsto Z_E$ of
 example~\ref{eg-ZE-formal} converts filtered homotopy colimits to
 colimits of formal schemes.
\end{example}

Suppose we have a diagram of formal schemes $X\:\CI\xra{}\FX$.  For
each $i\in\CI$ we then have a filtered category $\CJ(i)$ and a functor
$X(i,-)\:\CJ(i)\xra{}\CX$ such that $X(i)=\colim_{\CJ(i)}X(i,j)$.  For
many purposes, it is convenient if we can take all the categories
$\CJ(i)$ to be the same.  This motivates the following definition.

\begin{definition}\label{defn-rectifiable}
 A category $\CI$ is \dfn{rectifiable} if for every functor
 $X\:\CI\xra{}\FX$ there is a filtered category $\CJ$ and a functor
 $Y\:\CI\tm\CJ\xra{}\CX$ such that $X(i)=\colim_\CJ Y(i,j)$ as
 functors of $i$.
\end{definition}

\begin{proposition}\label{prop-rect-finite}
 If $\CI$ is a finite category such that $\CI(i,i)=\{1\}$ for all
 $i\in\CI$, then $\CI$ is \idx{rectifiable}.
\end{proposition}
\begin{proof}
 See~\cite[Proposition 8.8.5]{grve:sgaivei}.
\end{proof}

\begin{proposition}\label{prop-rect-discrete}
 If $\CI$ is a discrete small category (in other words, a set), then
 $\CI$ is \idx{rectifiable}.
\end{proposition}
\begin{proof}
 As $X(i)$ is a formal scheme, there is a filtered category $\CJ(i)$
 and a functor $Z(i,-)\:\CJ(i)\xra{}\CX$ such that
 $X(i)=\colim_{\CJ(i)}Z(i,j)$.  Write $\CJ=\prod_i\CJ(i)$, let
 $\pi_i\:\CJ\xra{}\CJ(i)$ be the projection, and let $Y(i,-)$ be the
 composite functor $\CJ\xra{\pi_i}\CJ(i)\xra{Z(i,-)}\CX$.  It is easy
 to check that $\CJ$ is filtered and that $\pi_i$ is cofinal, so
 $X(i)=\colim_\CJ Y(i,j)$, as required.
\end{proof}

\begin{proposition}\label{prop-formal-lim}
 The category $\FX$ has finite limits, and the inclusions
 $\CX\xra{}\FX\xra{}\CF$ preserve all limits that exist.  Moreover,
 finite limits in $\FX$ commute with filtered colimits.
\end{proposition}
\begin{proof}
 First consider a diagram $X\:\CI\xra{}\FX$ indexed by a finite
 rectifiable category.  We define $U(R)=\invlim_\CI X(i)(R)$, which
 gives a functor $\Rings\xra{}\Sets$.  It is well-known that this is
 the inverse limit of the diagram $X$ in the functor category $\CF$,
 so it will suffice to show that $U$ is a formal scheme.  As $\CI$ is
 rectifiable, we can choose a diagram $Y\:\CI\tm\CJ\xra{}\CX$ as in
 Definition~\ref{defn-rectifiable}.  As $\CX$ has limits, we can
 define $Z(j)=\invlim_iY(i,j)\in\CX$, and then define
 $W=\colim_jZ(j)\in\FX$.  Then $W(R)=\colim_j\invlim_iY(i,j)(R)$.  As
 filtered colimits commute with finite limits in the category of sets,
 this is the same as
 $\invlim_i\colim_jY(i,j)(R)=\invlim_iX(i)(R)=V(R)$.  Thus $V=W$ is a
 formal scheme, as required.

 Both finite products and equalisers can be considered as limits
 indexed by rectifiable categories, and we can write any finite limit
 as the equaliser of two maps between finite products.  This shows
 that $\FX$ has finite limits.  

 Now let $\{X_i\}$ be a diagram of formal schemes, let $X$ be a formal
 scheme, and let $\{f_i\:X\xra{}X_i\}$ be a cone.  If this is a limit
 cone in $\FX$ then we must have
 $X(R)=\FX(\spec(R),X)=\invlim_i\FX(\spec(R),X_i)=\invlim_iX_i(R)$,
 which means that it is a limit cone in $\CF$ (because limits in
 functor categories are computed pointwise).  The converse is equally
 easy, so the inclusion $\FX\xra{}\CF$ preserves and reflects limits.
 Similarly, the inclusion $\CX\xra{}\CF$ preserves and reflects
 limits, and it follows that the same is true of the inclusion
 $\CX\xra{}\FX$.
\end{proof}

\subsection{Solid formal schemes}
\label{subsec-solid}

\begin{definition}\label{defn-spf}
 A \dfn{linear topology} on a ring $R$ is a topology such that the
 cosets of open ideals are open and form a basis of open sets.  One
 can check that such a topology makes $R$ into a topological ring.  We
 write $\LRings$\index{LRings@$\LRings$} for the category of rings
 with a given linear topology, and continuous homomorphisms.  For any
 ring $S$, the discrete topology is a linear topology on $S$, so we
 can think of $\Rings$ as a full subcategory of $\LRings$.  Given a
 linearly topologised ring $R$, we define $\spf(R)\:\Rings\xra{}\Sets$
 \index{spf@$\spf$} by
 \[ \spf(R)(S) = \LRings(R,S) = \colim_J\Rings(R/J,S), \]
 where $J$ runs over the directed set of open ideals.  Clearly this
 defines a functor $\spf\:\LRings^\op\xra{}\FX$.
\end{definition}

\begin{definition}\label{defn-completion}
 Let $R$ be a linearly topologised ring.  The \dfn{completion} of $R$
 is the ring $\hR=\invlim_IR/I$\index{Rh@$\hR$}, where $I$ runs over
 the open ideals in $R$.  There is an evident map $R\xra{}\hR$, and
 the composite $R\xra{}\hR\xra{}R/I$ is surjective so we have
 $R/I=\hR/\ov{I}$ for some ideal $\ov{I}\leq\hR$.  These ideals form a
 filtered system, so we can give $\hR$ the linear topology for which
 they are a base of neighbourhoods of zero.  It is easy to check that
 $\widehat{\hR}=\hR$ and that $\spf(\hR)=\spf(R)$.  We say that $R$ is
 \emph{complete}, or that it is a \dfn{formal ring}, if $R=\hR$.  Thus
 $\hR$ is always a formal ring.  We write
 $\FRings$\index{FRings@$\FRings$} for the category of formal rings.
\end{definition}

\begin{definition}\label{defn-OX-formal}
 Given a formal scheme $X$, we recall that $\OX=\FX(X,\aff^1)$.  This
 is again a ring under pointwise operations.  If $\{X_i\}$ is a
 presentation of $X$ then $\OX=\invlim_i\OO_{X_i}$.

 For any point $x$ of $X$ we define
 $I_x=\{f\in\OX\st f(x)=0\in\OO_x\}$.  {}From a slightly different point
 of view, we can think of $x$ as a map $Y=\spec(\OO_x)\xra{}X$ and
 $I_x$ as the kernel of the resulting map $\OX\xra{}\OY$.  As the
 informal schemes over $X$ form a filtered category, we see that the
 ideals $I_x$ form a directed system.  Thus, there is a unique linear
 topology on $\OX$, such that the ideals $I_x$ form a base of
 neighbourhoods of zero.  With this topology, if $\{X_i\}$ is a
 presentation of $X$, then $\OX=\invlim_i\OO_{X_i}$ as topological
 rings.

 Note that 
 \[ \FX(X,\spf(R))=\invlim_i\FX(X_i,\spf(R))= 
    \invlim_i\LRings(R,\OO_{X_i})=\LRings(R,\OX),
 \]
 so that $\OO\:\FX\xra{}\LRings^\op$ is left adjoint to
 $\spf\:\LRings^\op\xra{}\FX$.  In particular, we have a unit map
 $X\xra{}\spf(\OX)$ in $\FX$, and a counit map $R\xra{}\OO_{\spf(R)}$
 in $\LRings$.  The latter is just the completion map $R\xra{}\hR$. 
\end{definition}

\begin{definition}\label{defn-solid}
 We say that a formal scheme $X$ is \dfn{solid} if it is isomorphic to
 $\spf(R)$ for some linearly topologised ring $R$.  We write
 $\FX_\sol$\index{Xsol@$\FX_\sol$} for the category of solid formal
 schemes.
\end{definition}

In the earlier incarnation of this paper~\cite{st:fpf} we defined
formal schemes to be what we now call solid formal schemes.  While
only solid formal schemes seem to occur in the cases of interest, the
category of all formal schemes has rather better categorical
properties, so we use it instead.
\begin{example}\label{eg-informal-solid}
 Any informal scheme $X$ is a solid formal scheme (because the zero
 ideal is open).
\end{example}
\begin{example}\label{eg-haf-solid}
 The formal scheme $\haf^n$ is solid.  To see this, consider the
 formal power series ring $R=\fps{\Zh}{x_1,\ldots,x_n}$, with the
 usual linear topology defined by the ideals $I^k$, where
 $I=(x_1,\ldots,x_k)$.  This is clearly a formal ring, and
 $\haf^n=\spf(R)$.
\end{example}
\begin{example}\label{eg-Noetherian-local}
 If $R$ is a complete Noetherian semilocal ring with Jacobson radical
 $I$ (for example, a complete Noetherian local ring with maximal ideal
 $I$) then it is natural to give $R$ the linear topology defined by
 the ideals $I^k$, and to define $\spf(R)$ using this.  With this
 convention, the set $\FX(\spf(R),\spf(S))$ (where $S$ is another ring
 of the same type) is just the set of local homomorphisms $S\xra{}R$.
 Thus, the categories of formal schemes used in~\cite{st:fsf}
 and~\cite{grst:vlc} embed as full subcategories of our category
 $\FX$.  
\end{example}
\begin{example}
 Let $Z$ be an infinite CW complex with finite subcomplexes
 $\{Z_\al\}$, and let $E$ be an even periodic ring spectrum.  Let
 $J_\al$ be the kernel of the map $E^0Z\xra{}E^0Z_\al$.  These ideals
 define a linear topology on $E^0Z$.  In good cases $E^0Z$ will be
 complete and we will have $Z_E=\spf(E^0Z)$, so this is a solid formal
 scheme.  See Section~\ref{sec-topology} for technical results that
 guarantee this.
\end{example}

\begin{proposition}\label{prop-solid}\ \\ \vspace{-2ex}
 \begin{itemize}
 \item[(a)] If $X$ is a \idx{solid} formal scheme then $\OX$ is a
  formal ring.
 \item[(b)] A formal scheme $X$ is solid if and only if it is
  isomorphic to $\spf(R)$ for some formal ring $R$, if and only if the
  natural map $X\xra{}\spf(\OX)$ is an isomorphism.
 \item[(c)] The functor $X\mapsto X_\sol=\spf(\OX)$ is left adjoint to
  the inclusion of $\FX_\sol$ in $\FX$.
 \item[(d)] The functor $R\mapsto\hR$ is left adjoint to the inclusion
  of $\FRings$ in $\LRings$.
 \item[(e)] The functors $R\mapsto\spf(R)$ and $X\mapsto\OX$ give an
  equivalence between $\FX_\sol$ and $\FRings^\op$.
 \end{itemize}
\end{proposition}
\begin{proof}
 (a): If $X$ is solid then $X=\spf(R)$ for some linearly topologised
 ring $R$, so $\OX=\OO_{\spf(R)}=\hR$ which is a formal ring.

 (b): If $X$ is solid then $X=\spf(R)$ as above, but
 $\spf(R)=\spf(\hR)$ so we may assume that $R$ is formal.  We find as
 in~(a) that $\OX=R$ and thus that the map $X\xra{}\spf(\OX)=\spf(R)$
 is an isomorphism.  The converse is easy.

 (c): Let $T$ denote the functor $X\mapsto X_\sol$.  This arises from
 an adjunction, so it is a monad.  On the other hand, if $R=\OX$ then
 $R$ is formal by~(a), so $R=\OO_{\spf(R)}=\OO_{X_\sol}$.  By applying
 $\spf(-)$, we see that $(X_\sol)_\sol=X_\sol$, so $T^2=T$ and $T$ is
 an idempotent monad.  Moreover, $\FX_\sol$ is the subcategory of
 formal schemes for which the unit map $\eta_X\:X\xra{}TX$ is an
 isomorphism.  It is well-known that this is automatically a
 reflective subcategory.  In outline, if $Y$ is solid and $X$ is
 arbitrary and $f\:X\xra{}Y$, then 
 $f'=\eta_Y^{-1}\circ Tf\:X_\sol\xra{}Y$ is the unique map such that
 $f'\circ\eta_X=f$.

 (d): The proof is similar.

 (e): If $R$ is formal then $\spf(R)$ is solid and
 $\OO_{\spf(R)}=\hR=R$.  If $X$ is solid then $\OX$ is formal (by~(a))
 and $X=\spf(\OX)$ (by~(b)).  
\end{proof}

\begin{definition}\label{defn-hot}
 Let $R$, $S$ and $T$ be linearly topologised rings, and let
 $R\xra{}S$ and $R\xra{}T$ be continuous homomorphisms.  We then give
 $S\ot_RT$ the linear topology defined by the ideals $I\ot T+S\ot J$,
 where $I$ runs over open ideals in $S$ and $J$ runs over open ideals
 in $T$.  This is easily seen to be the pushout of $S$ and $T$ under
 $R$ in $\LRings$.  We also define $S\hot_RT$ to be the completion of
 $S\ot_RT$.  If $R$, $S$ and $T$ are formal then $S\hot_RT$ is the
 pushout in $\FRings$ (because completion is left adjoint to the
 inclusion $\FRings\xra{}\LRings$).
\end{definition}

\begin{proposition}\label{prop-solid-products}
 The subcategory $\FX_\sol\sse\FX$\index{Xsol@$\FX_\sol$} is closed
 under finite products and arbitrary coproducts.  It also has its own
 colimits for arbitrary diagrams, which need not be preserved by the
 inclusion $\FX_\sol\xra{}\FX$.
\end{proposition}
\begin{proof}
 One can check that $\spf(R\ot S)=\spf(R\hot S)=\spf(R)\tm\spf(S)$,
 which gives finite products.  Let $\{R_i\st i\in\CI\}$ be a family of
 formal rings, and write $R=\prod_iR_i$.  We give this ring the
 product topology, which is the same as the linear topology defined by
 the ideals of the form $\prod_iJ_i$, where $J_i$ is open in $R_i$ and
 $J_i=R_i$ for almost all $i$.  We claim that
 $\spf(R)=\coprod_i\spf(R_i)$.  

 To see this, let $\CJ$ denote the set ideals $J=\prod_iJ_i$ as above.
 This is easily seen to be a directed set.  For $J\in\CJ$ we see that
 $R/J=\prod_iR_i/J_i$, where almost all terms in the product are
 zero.  Thus $\spec(R/J)=\coprod_{i\in\CI}\spec(R_i/J_i)$, where
 almost all terms in the coproduct are empty.  As colimits commute
 with coproducts, we see that
 $\spf(R)=\coprod_\CI\colim_\CJ\spec(R_i/J_i)$.  As the projection
 from $\CJ$ to the set of open ideals in $R_i$ is cofinal, we see that
 $\colim_\CJ\spec(R_i/J_i)=\spf(R_i)$, so that
 $\spf(R)=\coprod_\CI\spf(R_i)$ as claimed.

 Now let $\{X_i\}$ be an arbitrary diagram of solid formal schemes,
 and let $X$ be its colimit in $\FX$.  As the functor
 $Y\mapsto Y_\sol$ is left adjoint to the inclusion
 $\FX_\sol\xra{}\FX$, we see that $X_\sol$ is the colimit of our
 diagram in $\FX_\sol$.
\end{proof}
\begin{remark}
 We will see in Corollary~\ref{cor-solid-limits} that $\FX_\sol$ is
 actually closed under finite limits.
\end{remark}

\begin{example}\label{eg-const-formal}
 As a special case of the preceeding proposition, consider an infinite
 set $A$.  Let $R$ be the ring of functions $u\:A\xra{}\Zh$ with the
 product topology, so that $\un{A}=\spf(R)=\coprod_{a\in A}1$.  We
 call formal schemes of this type \emph{constant formal schemes}.
 \index{scheme!constant} More generally, given a formal scheme $X$ we
 write $\un{A}_X=\coprod_{a\in A}X$.  If $X$ is solid then
 $\un{A}_X=\spf(C(A,\OX))$, where $C(A,\OX)$ is the ring of functions
 $A\xra{}\OX$, under the evident product topology.  Clearly, if $E$ is
 an even periodic ring spectrum and we regard $A$ as a discrete space
 then $A_E=\un{A}\tm S_E$.
\end{example}

\subsection{Formal schemes over a given base}
\label{subsec-formal-based}

Let $X$ be a formal scheme.  Write $\FX_X$\index{XXh@$\FX_X$} for the
category of formal schemes over $X$, and $\CX_X$\index{XX@$\CX_X$} for
the full subcategory of informal schemes over $X$.  We also write
$\Points(X)$\index{PointsX@$\Points(X)$} for the category of pairs
$(R,x)$, where $R$ is a ring and $x\in X(R)$; the maps are as in
Definition~\ref{defn-Points}.  Again, the Yoneda isomorphism
$X(R)=\FX(\spec(R),X)$ gives an equivalence $\Points(X)=\CX_X^\op$.
Moreover, formal schemes $Y$ over $X$ biject with ind-representable
functors $Y'\:\Points(X)\xra{}\Sets$ by the rules
\begin{align*}
 Y'(R,x) &= \text{ preimage of $x$ under the map $Y(R)\xra{}X(R)$} \\
 Y(R)    &= \coprod_{x\in X(R)} Y'(R,x).
\end{align*}

Now consider a formal scheme $X$ with presentation $\{X_i\}$, indexed
by a filtered category $\CI$.  We next investigate the relationship
between the categories $\FX_X$ and $\FX_{X_i}$, which we now define.
\begin{definition}\label{defn-filtered-base}
 Given a diagram $\{X_i\}$ as above, we write $\CD_{\{X_i\}}$
 \index{DXi@$\CD_{\{X_i\}}$} for the
 category of diagrams $\{Y_i\}\:\CI\xra{}\FX$ equipped with a map of
 diagrams $\{Y_i\}\xra{}\{X_i\}$.  For any such diagram $\{Y_i\}$ and
 any map $u\:i\xra{}j$ in $\CI$, we have a commutative square
 \begin{diag}
  \node{Y_i} \arw{s} \arw{e,t}{Y_u} \node{Y_j} \arw{s} \\
  \node{X_i}           \arw{e,b}{X_u} \node{X_j.}
 \end{diag}
 We write $\FX_{\{X_i\}}$ \index{XXhi@$\FX_{\{X_i\}}$} for the full
 subcategory of $\CD_{\{X_i\}}$ consisting of diagrams $\{Y_i\}$ for
 which all such squares are pullbacks. 

 We define functors $F\:\CD_{\{X_i\}}\xra{}\FX_X$ and
 $G\:\FX_X\xra{}\CD_{\{X_i\}}$ by 
 \begin{align*}
  F\{Y_i\} &= \colim_i Y_i              \\
  GY       &= \{Y\tX X_i\}.
 \end{align*}
\end{definition}
\begin{proposition}\label{prop-filtered-base}
 The functor $F$ is left adjoint to $G$, and it preserves finite
 limits.  The functor $G$ is full and faithful, and its image is
 $\FX_{\{X_i\}}$.  The functors $F$ and $G$ give an equivalence
 between $\FX_X$ and $\FX_{\{X_i\}}$.

 Moreover, if $W$ is an informal scheme over $X$ and
 $\{Y_i\}\in\FX_{\{X_i\}}$, then any factorisation $W\xra{}X_i\xra{}X$
 of the given map $W\xra{}X$ gives an isomorphism
 $W\tX F\{Y_i\}=W\tm_{X_i}Y_i$. 
\end{proposition}
\begin{proof}
 A map $F\{Y_i\}\xra{}Z$ is the same as a compatible system of maps
 $Y_i\xra{}Z$ over $X$.  As the map $Y_i\xra{}X$ has a given
 factorisation through $X_i$, this is the same as a compatible system
 of maps $Y_i\xra{}Z\tX X_i=G(Z)_i$ over $X_i$, or in other words a
 map $\{Y_i\}\xra{}G(Z)$.  Thus $F$ is left adjoint to $G$.

 As filtered colimits commute with finite limits, we see that
 $FG(Y)=\colim_i(Y\tX X_i)=Y\tX \colim_iX_i=Y$.  This means that 
 \[ \CD_{\{X_i\}}(GY,GZ) = \FX_X(Y,FGZ) = \FX_X(Y,Z), \]
 so $G$ is full and faithful.  This means that $G$ is an equivalence
 of $\FX_X$ with its image, and it is clear that the image is
 contained in $\FX_{\{X_i\}}$.  The commutation of finite limits and
 filtered colimits also implies that $F$ preserves finite limits.

 We now prove the last part of the proposition; afterwards we will
 deduce that the image of $G$ is precisely $\FX_{\{x_i\}}$.  Consider
 an informal scheme $W$ and a map $f\:W\xra{}X$, and an object
 $\{Y_i\}$ of $\FX_{\{X_i\}}$.  Let $\CJ$ be the category of pairs
 $(i,g)$, where $i\in\CI$ and $g\:W\xra{}X_i$ and the composite
 $W\xra{g}X_i\xra{}X$ is the same as $f$.  It is not hard to check
 that $\CJ$ is filtered and that the projection functor $\CJ\xra{}\CI$
 is cofinal.  For each $(i,g)\in\CJ$ we have a pullback diagram
 \begin{diag}
  \node{W\tm_{X_i}Y_i} \arw{s} \arw{e} \node{Y_i} \arw{s} \\
  \node{W} \arw{e} \node{X_i.}
 \end{diag}
 By taking the colimit over $\CJ$ we get a pullback diagram
 \begin{diag}
  \node{\colim W\tm_{X_i}Y_i} \arw{s} \arw{e}
  \node{F\{Y_i\}} \arw{s} \\
  \node{W} \arw{e} \node{X.}
 \end{diag}
 On the other hand, for each map $u\:(i,g)\xra{}(j,h)$ in $\CJ$ we
 have $Y_i=X_i\tm_{X_j}Y_j$ (by the definition of $\FX_{\{X_i\}}$) and
 thus $W\tm_{X_i}Y_i=W\tm_{X_j}Y_j$.  It follows easily that for each
 $(i,g)$ the map $W\tm_{X_i}Y_i\xra{}\colim W\tm_{X_j}Y_j$ is an
 isomorphism, and thus (by the diagram) that
 $W\tX F\{Y_i\}=W\tm_{X_i}Y_i$.  

 Now take $W=X_i$ and $g=1$ in the above.  We find that
 $X_i\tX F\{Y_i\}=Y_i$, and thus that $FG\{Y_i\}=\{Y_i\}$, and thus
 that $\{Y_i\}$ is in the image of $G$.  This shows that the image of
 $G$ is precisely $\FX_{\{X_i\}}$, as required.
\end{proof}

\begin{definition}
 Let $Y$ be a formal scheme over a formal scheme $X$.  We say that $Y$
 is \dfn{relatively informal} over $X$ if for all informal schemes
 $X'$ over $X$, the pullback $Y\tX X'$ is informal.
\end{definition}

\begin{proposition}\label{prop-rel-inf-lim}
 The category of relatively informal schemes over $X$ has limits,
 which are preserved by the inclusion into $\FX_X$.
\end{proposition}
\begin{proof}
 We can write $X$ as the colimit of a filtered diagram of informal
 schemes $X_i$.  It is clear that the category of relatively informal
 schemes is equivalent to the subcategory $\CC$ of $\FX_{\{X_i\}}$
 consisting of systems $\{Y_i\}$ of informal schemes.  As the category
 of informal schemes has limits, we see that the category of informal
 schemes over $X_i$ has limits.  Moreover, for each map
 $X_i\xra{}X_j$, the functor
 $X_i\tm_{X_j}(-)\:\CX_{X_j}\xra{}\CX_{X_i}$ preserves limits.  Given
 this, it is easy to check that $\CC$ has limits, as required.  As the
 inclusion $\CX\xra{}\FX$ preserves limits, one can check that the
 same is true of the inclusions $\CX_{X_i}\xra{}\FX_{X_i}$ and
 $\CC\xra{}\FX_{\{X_i\}}=\FX_X$.
\end{proof}

\subsection{Formal subschemes}
\label{subsec-formal-subschemes}

\begin{definition}\label{defn-closed-formal}
 We say that a map $f\:X\xra{}Y$ of formal schemes is a
 \dfn{closed inclusion} if it is a \idx{regular monomorphism} in $\FX$.
 (This means that it is the equaliser of some pair of arrows
 $Y\arw{e,=>}Z$, or equivalently that it is the equaliser of the pair
 $Y\arw{e,=>}Y\amalg_XY$.)  A \emph{closed formal subscheme}
 \index{subscheme!closed} of a formal scheme $Y$ is a subfunctor $X$
 of $Y$ such that $X$ is a formal scheme and the inclusion $X\xra{}Y$
 is a closed inclusion. 
\end{definition}
\begin{remark}
 The functor $Z\mapsto Z(R)$ is representable (by $\spec(R)$).  It
 follows that if $f\:V\xra{}W$ is a monomorphism in $\FX$ then
 $V(R)\xra{}W(R)$ is injective for all $R$, so $V$ is isomorphic to a
 subfunctor of $W$.  If $f$ is a regular monomorphism, then the
 corresponding subfunctor is a closed subscheme.
\end{remark}

\begin{example}\label{eg-formal-VJ}
 Let $J$ be an ideal in $\OX$, generated by elements
 $\{f_i\st i\in I\}$ say.  We define
 \[ V(J)(R)=\{x\in X(R)\st f(x)=0\text{ for all } f\in J\} 
     = \{x\st f_i(x)=0\text{ for all }i\}.
 \]
 Define a scheme $\aff^I$ by $\aff^I(R)=\prod_{i\in I}R$ (this is
 represented by the polynomial algebra $\Zh[x_i\st i\in I]$).  This is
 just the product $\prod_{i\in I}\aff^1$; by
 Proposition~\ref{prop-formal-lim}, it does not matter whether we
 interpret this in $\CX$ or $\FX$.  It follows that there is a map
 $f\:X\xra{}\aff^I$ with components $f_i$, and another map
 $g\:X\xra{}\aff^I$ with components $0$.  Clearly $V(J)$ is the
 equaliser of $f$ and $g$, and thus it is a closed formal subscheme of
 $X$.  There is a natural map $\OX/J\xra{}\OO_{V(J)}$ which is an
 isomorphism in most cases of interest, but I suspect that this is not
 true in general (compare Remark~\ref{rem-lim-1}).
\end{example}
\begin{example}\label{eg-include-completion}
 If $X$ is an informal scheme and $Y$ is a closed informal subscheme
 of $X$ then the evident map $X^\wedge_Y\xra{}X$ is a closed
 inclusion.\index{XY@$X^\wedge_Y$}
\end{example}

\begin{proposition}\label{prop-closed-same}
 A map $f\:X\xra{}Y$ of informal schemes is a \idx{closed inclusion}
 in $\FX$ if and only if it is a closed inclusion in $\CX$.
\end{proposition}
\begin{proof}
 It follows from Proposition~\ref{prop-formal-colim} that the pushout
 $Y\amalg_XY$ is the same whether constructed in $\CX$ or $\FX$.  It
 follows in turn from Proposition~\ref{prop-formal-lim} that the
 equaliser of the two maps $Y\arw{e,=>}Y\amalg_XY$ is the same
 whether constructed in $\CX$ or $\FX$.  The map $f$ is a closed
 inclusion if and only if $X$ maps isomorphically to this equaliser,
 so the proposition follows. 
\end{proof}

\begin{proposition}\label{prop-closed-colim}
 If $X\in\FX$ and $Y\in\CX$, then a map $f\:X\xra{}Y$ is a
 \idx{closed inclusion} if and only if there is a directed set of
 closed informal subschemes $Y_i$ of $Y$ such that $X=\colim_i Y_i$.
\end{proposition}
\begin{proof}
 First suppose that $f$ is a closed inclusion.  We can write $X$ as a
 colimit of informal schemes, say $X=\colim_{i\in\CI}X_i$.  Write
 $Z_i=Y\amalg_{X_i}Y$.  One checks that these schemes give a
 functor $\CI\xra{}\CX$, and that $\colim_i Z_i=Y\amalg_XY$.  Let
 $Y_i$ be the equaliser of the two maps $X\arw{e,=>}Z_i$, so
 that $Y_i$ is a closed informal subscheme of $X$, and again the
 schemes $Y_i$ give a functor $\CI\xra{}\CX$.  As finite limits
 commute with filtered colimits in $\FX$, we see that
 $\colim_i Y_i$ is the equaliser of the maps
 $Y\arw{e,=>}\colim_i Z_i=Y\amalg_XY$.  This is just $X$,
 because $f$ is assumed to be a regular monomorphism.

 Conversely, suppose that $\{Y_i\}$ is a directed family of closed
 subschemes of an informal scheme $Y$.  Write $Z_i=Y\amalg_{Y_i}Y$
 and $Z=\colim_i Z_i$.  By much the same logic as above, we see
 that there is a pair of maps $Y\arw{e,=>}Z$ whose equaliser in
 $X=\colim_i Y_i$, so that $X$ is a closed formal subscheme of
 $Y$. 
\end{proof}

\begin{proposition}\label{prop-pb-refl}
 A map $f\:X\xra{}Y$ in $\FX$ is a \idx{closed inclusion} if and only
 if for all informal schemes $Y'$ and all maps $Y'\xra{}Y$, the
 pulled-back map $f'\:X'\xra{}Y'$ is a closed inclusion.
\end{proposition}
\begin{proof}
 It is clear that the condition is necessary, because in any category
 a pullback of a regular monomorphism is a regular monomorphism.  For
 sufficiency, suppose that $f\:X\xra{}Y$ is such that all maps of the
 form $f'\:X'\xra{}Y'$ are closed inclusions.  Write $Y$ as a colimit
 of informal schemes $Y_i$ in the usual way, and let
 $f_i\:X_i\xra{}Y_i$ be the pullback of $f$ along the map
 $Y_i\xra{}Y$.  As finite limits in $\FX$ commute with filtered
 colimits, we see that $X=\colim_i X_i$.  By assumption, $f_i$
 is a closed inclusion.  Write $Z_i=Y_i\amalg_{X_i}Y_i$, so
 $X_i$ is the equaliser of the fork $Y_i\arw{e,=>}Z_i$.  Write
 $Z=\colim_i Z_i$.  As finite limits in $\FX$ commute with
 filtered colimits, we see that $X$ is the equaliser of the maps
 $Y\arw{e,=>}Z$, and thus that $f$ is a closed inclusion.
\end{proof}

\begin{proposition}\label{prop-closed-comp}
 Let $X\xra{f}Y\xra{g}Z$ be maps of formal schemes.  If $f$ and $g$
 are closed inclusions, then so is $gf$.  Conversely, if $gf$ is a
 \idx{closed inclusion} and $g$ is a monomorphism then $f$ is a closed
 inclusion.  
\end{proposition}
\begin{proof}
 The second part is a formal statement which holds in any category: if
 we have maps $X\xra{f}Y\xra{g}Z$ such that $gf$ is the equaliser of a
 pair $Z\arw{e,tb,=>}{p}{q}W$, then a diagram chase shows that
 $X\xra{f}Y$ is the equaliser of $pg$ and $qg$ and thus is a regular
 monomorphism.  

 For the first part, we can assume by Proposition~\ref{prop-pb-refl}
 that $Z$ is an informal scheme.  We then know from
 Proposition~\ref{prop-closed-colim} that there is a filtered system
 of closed subschemes $Z_i$ of $Z$ such that $Y$ is the colimit of
 the $Z_i$.  The maps $Y\xra{}Z$ and $Z_i\xra{}Y\xra{}Z$ are
 closed inclusions, so the second part tells us that $Y_i\xra{}Y$ is
 a closed monomorphism.  Let $X_i$ be the preimage of $Z_i\sse Y$
 under the map $f\:X\xra{}Y$.  The maps $X_i\xra{}Z_i$ and
 $Z_i\xra{}Z$ are closed inclusions of informal schemes, so the
 composite $X_i\xra{}Z$ is easily seen to be a closed inclusion
 (because closed inclusions in the informal category are just dual to
 surjections of rings).  As filtered colimits commute with pullbacks,
 we see that $X=\colim_i X_i$.  It follows from
 Proposition~\ref{prop-closed-colim} that $X\xra{}Z$ is a closed
 inclusion. 
\end{proof}

\begin{proposition}
 Any closed formal subscheme\index{subscheme!closed} of a \idx{solid}
 formal scheme is again solid.
\end{proposition}
\begin{proof}
 Let $W\xra{f}X\arw{e,tb,=>}{g}{h}Y$ be an equaliser diagram, and
 suppose that $X$ is solid.  We need to show that $W$ is solid.
 Choose a presentation $Y=\colim_{i\in\CI}Y_i$ for $Y$.  Let $\CJ$ be
 the set of tuples $j=(J,i,g',h')$, where $J$ is an open ideal in
 $\OX$ and $i\in\CI$ and $g',h'\:V(J)\xra{}Y_i$ and the following
 diagram commutes.
 \begin{diag}
  \node{V(J)} \arrow{e,tb,=>}{g'}{h'} \arrow{s,V} 
  \node{Y_i}                          \arrow{s}         \\
  \node{X}    \arrow{e,tb,=>}{g}{h}   \node{Y}
 \end{diag}
 One can make $\CJ$ into a filtered category so that $j\mapsto J$ is
 a cofinal functor to the directed set of open ideals of $\OX$, and
 $j\mapsto i$ is a cofinal functor to $\CI$ (see the proof
 of~\cite[Proposition 8.8.5]{grve:sgaivei}).
 The equaliser of $g'$ and $h'$ is a closed subscheme of $V(J)$, so it
 has the form $V(I_j)$ for some ideal $I_j\ge J$.  As equalisers
 commute with filtered colimits, we see that $W=\colim_\CJ V(I_j)$.
 Let $\CK$ be the set of ideals of the form $I_j$ for some $j$.  The
 functor $j\mapsto I_j$ from $\CJ$ to $\CK$ is cofinal, so we have
 $W=\colim_{I\in\CK}V(I)$.  We can define a new linear topology on
 $R=\OX$ by letting the ideals $I\in\CK$ be a base of neighbourhoods
 of zero, and we conclude that $W=\spf(R)$.  Thus, $W$ is solid.
\end{proof}
\begin{remark}\label{rem-lim-1}
 In the above proof, suppose that $Y$ is also solid, and let $K$ be
 the ideal in $\OX$ generated by elements of the form $g^*u-h^*u$ with
 $u\in\OY$.  One can then check that $\OW=\invlim_J\OX/(K+J)$, where
 $J$ runs over the open ideals in $\OX$.  The kernel of the map
 $\pi\:\OX\xra{}\OW$ is $\bigcap_J(J+K)$, which is just the closure of
 $K$.  One would like to say that $\pi$ was surjective, but in fact
 its cokernel is $\invlim^1_J(J+K)$, which can presumably be nonzero.
\end{remark}
\begin{corollary}\label{cor-solid-limits}
 The subcategory $\FX_\sol\sse\FX$ of \idx{solid} formal schemes is
 closed under finite limits.
\end{corollary}
\begin{proof}
 We know from Proposition~\ref{prop-solid-products} that a finite
 product of solid schemes is solid, and a finite limit is a closed
 formal subscheme of a finite product.
\end{proof}

\subsection{Idempotents and formal schemes}
\label{subsec-formal-idempotents}

\begin{proposition}\label{prop-formal-idempotents}
 Let $X$ be a formal scheme.  Then systems of formal subschemes $X_i$
 such that $X=\coprod_iX_i$ biject with systems of \idx{idempotent}s
 $e_i\in\OX$ such that $e_ie_j=\dl_{ij}e_i$ and $\sum_ie_i$ converges
 to $1$ in the natural topology in $\OX$.  More explicitly, we require
 that for every open ideal $J\leq\OX$ the set
 $S=\{i\st e_i\not\in J\}$ is finite, and $\sum_Se_i=1\pmod{J}$.
\end{proposition}
\begin{proof}
 Suppose that $X=\coprod_{i\in\CI}X_i$.  Then
 $\OX=\FX(X,\aff^1)=\prod_i\FX(X_i,\aff^1)=\prod_i\OO_{X_i}$ as
 rings.  If $K$ is a finite subset of $\CI$, we write
 $X_K=\coprod_{i\in K}X_i$.  We then have $X=\colim_KX_K$, and this is
 a filtered colimit, so $X(R)=\colim_KX_K(R)$ for all $R$.  Using
 this, it is not hard to check that $\OX=\prod_i\OO_{X_i}$ as
 topological rings, where the right hand side is given the product
 topology.  Note that the product topology is defined by the ideals of
 the form $\prod_iJ_i$, where $J_i$ is an open ideal in $\OO_{X_i}$
 and $J_i=\OO_{X_i}$ for almost all $i$.

 For each $i$ there is an evident idempotent $e_i$ in
 $\OX=\prod_i\OO_{X_i}$, whose $j$'th component is $\dl_{ij}$.  This
 gives a system of idempotents as described in the proposition.

 Conversely, suppose we start with such a system of idempotents.  For
 any idempotent $e\in\OX$ it is easy to check that $D(e)=V(1-e)$, so
 we can define $X_i=D(e_i)=V(1-e_i)$.  We need to check that
 $X=\coprod_iX_i$.  We can write $X=\colim_\CJ Y_j$ for some filtered
 system of informal schemes $Y_j$.  Let $e_{ij}$ be the image of $e_i$
 in $\OO_{Y_j}$ and write $Z_{ij}=D(e_{ij})=V(1-e_{ij})\sse Y_j$.  As
 $Y_j$ is informal we know that the kernel of the map
 $\OX\xra{}\OO_{Y_j}$ is open and thus that $e_{ij}=0$ for almost
 all $i$.  We thus have a decomposition $Y_j=\coprod_iZ_{ij}$, in
 which only finitely many factors are nonempty.  If we fix $i$, it is
 easy to check that the schemes $Z_{ij}$ are functors of $j$, and that
 $\colim_jZ_{ij}=X_i$.  As colimits commute with coproducts, we find
 that $X=\coprod_iX_i$ as claimed.
\end{proof}

\begin{corollary}\label{cor-strong-coprod}
 Coproducts in $\FX$ or $\FX_X$ are strong. \index{strong coproduct}
\end{corollary}
\begin{proof}
 Let $\{Y_i\}$ be a family of schemes over $X$, and write
 $Y=\coprod_iY_i$.  Let $Z$ be another scheme over $X$, and write
 $Z_i=Z\tX Y_i$.  We need to show that $Z\tX Y=\coprod_iZ_i$.  To see
 this, take idempotents $e_i\in\OY$ as in the proposition, so that
 $Y_i=D(e_i)=V(1-e_i)$.  Let $e'_i$ be the image of $e_i$ under the
 evident map $\OZ\xra{}\OY$; it is easy to check that $Z_i=D(e'_i)$.
 As the idempotents $e_i$ are orthogonal and sum to $1$ and the map
 $\OO_{Z\tX Y}\xra{}\OY$ is a continuous map of topological rings, we
 see that the $e'_i$ are also orthogonal idempotents whose sum is $1$.
 This shows that $Z\tX Y=\coprod_iZ_i$ as claimed.
\end{proof}

\subsection{Sheaves over formal schemes}
\label{subsec-formal-sheaves}

In Section~\ref{subsec-sheaves}, we defined sheaves and vector bundles
over all functors, and in particular over formal schemes.

\begin{remark}\label{rem-aff-formal}
 If $M$ is a \idx{vector bundle} and $L$ is a \idx{line bundle} over a
 formal scheme $X$, we can define functors $\aff(M)(R)$
 \index{AAM@$\aff(M)$} and $\aff(L)^\tm(R)$ \index{AALt@$\aff(L)^\tm$}
 just as in Definitions~\ref{defn-Gamma} and~\ref{defn-gen-L}.  We
 claim that these are formal schemes.  Given a map $f\:W\xra{}X$, it
 is easy to check that $f^*\aff(M)=\aff(f^*M)$ (where the pullback on
 the left hand side is computed in the functor category $\CF$).  In
 particular, if $W$ is informal then Proposition~\ref{prop-aff-scheme}
 shows that $f^*\aff(M)$ is a scheme.  Now write $X=\colim_iX_i$ in
 the usual way, and let $M_i$ be the pullback of $M$ over $X_i$.  We
 find easily that $\aff(M)=\colim_i\aff(M_i)$, so $\aff(M)$ is a
 formal scheme.  Similarly, $\aff(L)^\tm$ is a formal scheme.
\end{remark}
\begin{remark}\label{rem-aff-infinite}
 If $M$ is a sheaf such that $M_x$ is an infinitely generated free
 module for all $x$, we find that $\aff(M)$\index{AAM@$\aff(M)$} is a
 formal scheme over $X$.  Unlike the case of a vector bundle, it is
 not relatively informal over $X$.  We leave the proof as an exercise.
\end{remark}

\begin{remark}\label{rem-complete-modules}
 Let $\{X_i\}$ be a presentation of a formal scheme $X$.  If $M$ is a
 sheaf over $X$ then one can check that
 $\Gm(X,M)=\invlim_i\Gm(X_i,M)$.  In particular, if $X$ is solid and
 $M_J=\Gm(V(J),M)$ for all open ideals $J\leq\OX$ we find that
 $\Gm(X,M)=\invlim_JM_J$.  Moreover, if $J\leq K$ we find that
 $M_K=M_J/KM_J$.  

 In particular, if $N$ is an $\OX$-module we find that
 $\Gm(X,\tN)=\invlim_JN/JN$.  We say that $N$ is \emph{complete} if
 $N=\invlim_JN/JN$.  It follows that the functor $N\mapsto\tN$ embeds
 the category of complete modules as a full subcategory of
 $\Sheaves_X$.  {\bf Warning:} it seems that the functor
 $N\mapsto\invlim_JN/JN$ need not be idempotent in bad cases, so
 $\invlim_JN/JN$ need not be complete.
\end{remark}

We next consider the problem of constructing sheaves over filtered
colimits.
\begin{definition}
 Let $\{X_i\}$ be a filtered diagram of functors, with colimit $X$.
 Let $\Sheaves_{\{X_i\}}$ \index{SheavesXi@$\Sheaves_{\{X_i\}}$}
 denote the category of systems $(\{M_i\},\phi)$ of the following
 type:
 \begin{itemize}
 \item[(a)] For each $i$ we have a sheaf $M_i$ over $X_i$.
 \item[(b)] For each $u\:i\xra{}j$ (with associated map
  $X_u\:X_i\xra{}X_j$) we have an isomorphism
  $\phi(u)\:M_i\simeq X_u^*M_j$.
 \item[(c)] In the case $u=1\:i\xra{}i$ we have $\phi(1)=1$.
 \item[(d)] Given $i\xra{u}j\xra{v}k$ we have
  $\phi(vu)=(X_u^*\phi(v))\circ\phi(u)$. 
 \end{itemize}
\end{definition}

\begin{proposition}
 Let $\{X_i\st i\in\CI\}$ be a filtered diagram of functors, with
 colimit $X$.  The category $\Sheaves_{\{X_i\}}$ is equivalent to
 $\Sheaves_X$.
\end{proposition}
\begin{proof}
 Given a sheaf $M$ over $X$, we define a system of sheaves
 $M_i=v_i^*M$, where $v_i\:X_i\xra{}X$ is the given map.  If
 $u\:i\xra{}j$ then $v_j\circ X_u=v_i$ so we have a canonical
 identification $M_i=X_u^*M_j$, which we take as $\phi(u)$.  This
 gives an object of $\Sheaves_{\{X_i\}}$.

 On the other hand, suppose we start with an object $\{M_i\}$ of
 $\Sheaves_{\{X_i\}}$, and we want to construct a sheaf $M$ over $X$.
 Given a ring $R$ and a point $x\in X(R)$, we need to define a module
 $M_x$ over $R$.  As $X=\colim_iX_i(R)$, we can choose $i\in\CI$ and
 $y\in X_i(R)$ such that $v_i(y)=x$.  We would like to define
 $M_x=M_{i,y}$, but we need to check that this is canonically
 independent of the choices made.  We thus let $\CJ$ be the category
 of all such pairs $(i,y)$.  Because $X(R)=\colim_iX_i(R)$, we see
 that $\CJ$ is filtered.  For each $(i,y)\in\CJ$ we have an $R$-module
 $M_{i,y}$, and the maps $\phi(u)$ make this a functor
 $\CJ\xra{}\Mod_R$.  We define $M_x=\colim_\CJ M_{i,y}$.  Because this
 is a filtered diagram of isomorphisms, each of the canonical maps
 $M_{i,y}\xra{}M_x$ is an isomorphism.  We leave it to the reader to
 check that this construction produces a sheaf, and that it is inverse
 to our previous construction.
\end{proof}

\begin{corollary}\label{cor-fibrewise-sheaf}
 Let $X\xra{}Y$ be a map of formal schemes.  To construct a sheaf over
 $X$, it suffices to construct sheaves over $W\tY X$ in a sufficiently
 natural way, for all informal schemes $W$ over $Y$.  It also suffices
 to construct sheaves over $X_y$ in a sufficiently natural way, for
 all points $y$ of $Y$.
\end{corollary}
\begin{proof}
 The two claims are really the same, as points of $Y$ biject with
 informal schemes over $Y$ by sending a point $y\in Y(R)$ to the usual
 map $\spec(R)\xra{y}Y$.

 For the first claim, we choose a presentation $Y=\colim_iY_i$ and
 write $X_i=Y_i\tY X$, and note that $X=\colim_iX_i$.  By assumption,
 we have sheaves $M_i$ over $X_i$.  ``Sufficiently natural'' means
 that we have maps $\phi(u)$ making $\{M_i\}$ into an object of
 $\Sheaves_{\{X_i\}}$, so the proposition gives us a sheaf over $X$.
\end{proof}

\subsection{Formal faithful flatness}
\label{subsec-formal-flat}

\begin{definition}\label{defn-ff-formal}
 Let $f\:X\xra{}Y$ be a map of formal schemes.  We say that $f$ is
 \dfn{flat} if the pullback functor $f^*\:\FX_Y\xra{}\FX_X$ preserves
 finite colimits.  We say that $f$ is \emph{faithfully flat}
 \index{flat!faithfully} if $f^*$ preserves and reflects finite
 colimits.
\end{definition}

\begin{remark}
 For any map $f\:X\xra{}Y$ of formal schemes, we know that $f^*$
 preserves all small coproducts.  Thus $f$ is flat if and only if
 $f^*$ preserves coequalisers, if and only if $f^*$ preserves all
 small colimits.
\end{remark}

Definition~\ref{defn-ff-formal} could in principle conflict with
Definition~\ref{defn-ff}; the following proposition shows that this is
not the case.
\begin{proposition}\label{prop-ff-consistent}
 A map $f\:X\xra{}Y$ of informal schemes is flat (resp.\ faithfully
 flat) as a map of informal schemes if and only if it is flat (resp.\ 
 faithfully flat) as a map of formal schemes.
\end{proposition}
\begin{proof}
 Recall that the inclusion $\CX\xra{}\FX$ preserves finite colimits.
 Given this, we see easily that a map that is formally flat (resp.\ 
 faithfully flat) flat is also informally flat (resp.\ faithfully
 flat).

 Now suppose that $f$ is informally flat.  Let
 $U\arw{e,=>}V\xra{}W$ be a coequaliser in $\FX_Y$.  By
 Proposition~\ref{prop-rect-finite}, we can find a filtered system of
 diagrams $U_i\arw{e,=>}V_i$ (with $U_i$ and $V_i$ in $\CX$) whose
 colimit is the diagram $U\arw{e,=>}V$.  We define $W_i$ to be the
 coequaliser of $U_i\arw{e,=>}V_i$.  As colimits commute, we have
 $W=\colim_iW_i$.  Clearly all this can be thought of as happening
 over $W$ and thus over $Y$.  By assumption, the diagram
 $f^*U_i\arw{e,=>}f^*V_i\xra{}f^*W_i$ is a coequaliser.  We now take
 the colimit over $i$, noting that $f^*$ commutes with filtered
 colimits and that colimits of coequalisers are coequalisers.  This
 shows that $f^*U\arw{e,=>}f^*V\xra{}f^*W$ is a coequaliser.  Thus,
 $f$ is flat.

 Now suppose that $f$ is informally faithfully flat, and let
 $u\:U\xra{}V$ be a map of formal schemes over $Y$ such that $f^*u$ is
 an isomorphism.  Choose a presentation $V=\colim_iV_i$ and write
 $U_i=U\tm_VV_i$, so that $U=\colim_iU_i$.  As $f^*$ preserves
 pullbacks, we see that the map $f^*U_i\xra{}f^*V_i$ is the pullback
 of the isomorphism $f^*U\xra{}f^*V$ along the map $f^*V_i\xra{}f^*V$,
 and thus that the map $f^*U_i\xra{}f^*V_i$ is itself an isomorphism.
 As $f$ is informally faithfully flat, we conclude that
 $U_i\simeq V_i$.  By passing to colimits, we see that $U\simeq V$ as
 claimed. 
\end{proof}

\begin{remark}
 Propositions~\ref{prop-ff-composites}, \ref{prop-ff-pullbacks},
 \ref{prop-ff-reg-epi} and \ref{prop-descent} are general nonsense,
 valid in any category with finite limits and colimits.  They
 therefore carry over directly to formal schemes.
 \index{flat} \index{flat!faithfully}
\end{remark}

\begin{lemma}\label{lem-flat-inf-test}
 Let $f\:X\xra{}Y$ be a map of formal schemes.  Let $\CX_Y$ be the
 category of informal schemes with a map to $Y$, and let
 $f^*_0\:\CX_Y\xra{}\FX_X$ be the restriction of $f^*$ to $\CX_Y$.  If
 $f^*_0$ preserves coequalisers, then $f$ is flat.
 \index{flat} \index{flat!faithfully}
\end{lemma}
\begin{proof}
 Suppose that $f_0^*$ preserves coequalisers.  Let
 $U\arw{e,=>}V\xra{}W$ be a coequaliser in $\FX_Y$.  By
 Proposition~\ref{prop-rect-finite}, we can find a filtered system of
 diagrams $U_i\arw{e,=>}V_i$ (with $U_i$ and $V_i$ in $\CX$) whose
 colimit is the diagram $U\arw{e,=>}V$.  We define $W_i$ to be the
 coequaliser of $U_i\arw{e,=>}V_i$.  As colimits commute, we have
 $W=\colim_iW_i$.  Clearly all this can be thought of as happening
 over $W$ and thus over $Y$.  By assumption, the diagram
 $f^*U_i\arw{e,=>}f^*V_i\xra{}f^*W_i$ is a coequaliser.  We now take
 the colimit over $i$, noting that $f^*$ commutes with filtered
 colimits and that colimits of coequalisers are coequalisers.  This
 shows that $f^*U\arw{e,=>}f^*V\xra{}f^*W$ is a coequaliser.  Thus,
 $f$ is flat.
\end{proof}

\begin{proposition}\label{prop-ff-filtered-base}
 Let $f\:X\xra{}Y$ be a map of formal schemes.  Suppose that $Y$ has a
 presentation $Y=\colim_iY_i$ for which the maps
 $f_i\:X_i=f^*Y_i\xra{}Y_i$ are (faithfully) flat.  Then $f$ is
 (faithfully) flat.
 \index{flat} \index{flat!faithfully}
\end{proposition}
\begin{proof}
 First suppose that each $f_i$ is flat.  Let $U\arw{e,=>}V\xra{}W$ be
 a coequaliser of informal schemes over $Y$.  By
 Lemma~\ref{lem-flat-inf-test}, it is enough to check that
 $f^*U\arw{e,=>}f^*V\xra{}f^*W$ is a coequaliser.  We know from
 Proposition~\ref{prop-formal-colim} that
 $\FX(W,Y)=\colim_i\FX(W,Y_i)$, so we can choose a factorisation
 $W\xra{}Y_i\xra{}Y$ of the given map $W\xra{}Y$, for some $i$.  We
 then have $f^*W=W\tY X=W\tm_{Y_i}Y_i\tY X=W\tm_{Y_i}X_i=f_i^*W$.
 Similarly, we have $f^*V=f_i^*V$ and $f^*U=f_i^*U$.  As $f_i$ is
 flat, we see that $f^*U\arw{e,=>}f^*V\xra{}f^*W$ is a coequaliser, as
 required.

 Now suppose that each $f_i$ is faithfully flat.  Let $s\:U\xra{}V$ be
 a morphism in $\FX_Y$ such that $f^*s$ is an isomorphism.  We need to
 show that $s$ is an isomorphism.  We have a pullback square of the
 following form.
 \begin{diag}
  \node{X_i} \arrow{s,l}{u_i} \arrow{e,t}{f_i}
  \node{Y_i} \arrow{s,r}{v_i} \\
  \node{X}   \arrow{e,b}{f}   \node{Y.}
 \end{diag}
 As $f^*s$ is an isomorphism, we see that $f_i^*v_i^*s=u_i^*f^*s$ is
 an isomorphism.  As $f_i$ is faithfully flat, we conclude that
 $v_i^*s\:v_i^*U\xra{}v_i^*V$ is an isomorphism for all $i$.  We also
 know that $U=\colim_iv_i^*U$ and $V=\colim_iv_i^*V$, and it follows
 easily that $s$ is an isomorphism.
\end{proof}

\begin{proposition}\label{prop-locfree-formal}
 Let $M$ be a \idx{vector bundle} of rank $r$ over a formal scheme
 $X$.  Then there is a faithfully flat map $f\:\Bases(M)\xra{}X$
 \index{BasesM@$\Bases(M)$} such that $f^*M\simeq\OO^r$.
 \index{flat!faithfully}
\end{proposition}
\begin{proof}
 Let $\Bases(M)(R)$ be the set of pairs $(x,B)$, where $x\in X(R)$ and
 $B\:R^r\xra{}M_x$ is an isomorphism.  Define $f\:\Bases(M)\xra{}X$ by
 $f(x,B)=x$.  As in the informal case (Example~\ref{eg-bases-M}) we see
 that $\Bases(M)$ is a formal scheme over $X$, and that
 $f^*M\simeq\OO^r$.  If $X_i$ is an informal scheme and
 $u\:X_i\xra{}X$ then one checks that $u^*\Bases(M)=\Bases(u^*M)$,
 which is faithfully flat over $X_i$ by Example~\ref{eg-bases-M}.  It
 follows from Proposition~\ref{prop-ff-filtered-base} that $\Bases(M)$
 is faithfully flat over $X$.
\end{proof}

\begin{definition}\label{defn-vflat-formal}
 A map $f\:X\xra{}Y$ of formal schemes is \emph{very flat}
 \index{flat!very} if for all informal schemes $Y'$ over $Y$, the
 scheme $X'=f^*Y'$ is informal and the map $X'\xra{}Y'$ is very flat
 (in other words, $\OO_{X'}$ is a free module over $\OO_{Y'}$).
 Similarly, we say that $f$ is \emph{finite} \index{finite map} if for
 all such $Y'$, the scheme $X'$ is informal and the map $X'\xra{}Y'$
 is finite.
\end{definition}

\subsection{Coalgebraic formal schemes}
\label{subsec-coalgebraic}

Fix a scheme $Z$, and write $R=\OZ$.  We next study the category
$\CC_Z$ \index{CZ@$\CC_Z$} of coalgebras over $R$, and a certain full
subcategory $\CC'_Z$ \index{CZp@$\CC'_Z$}.  It turns out that there is
a full and faithful embedding $\CC'_Z\xra{}\FX_Z$, and that the
categorical properties of $\CC_Z$ are in some respects superior to
those of $\FX_Z$.  Because of this, the categories $\CC_Z$ and
$\CC'_Z$ are often useful tools for constructing objects of $\FX_Z$
with specified properties.  Our use of coalgebras was inspired by
their appearance in~\cite{de:lpd}, although it is assumed there that
$R$ is a field, which removes many technicalities.

We will use $R$ and $Z$ as interchangeable subscripts, so
\[ \FX_R=\FX_Z=\{\text{formal schemes over }Z\}, \]
for example.  Write $\CM_R=\CM_Z$ \index{MZ@$\CM_Z$} and $\CC_R=\CC_Z$
for the categories of modules and coalgebras over $R$.  (All
coalgebras will be assumed to be cocommutative and counital.)  It is
natural to think of $\CC_Z$ as a ``geometric'' category, and we choose
our notation to reflect this point of view.  In particular, we shall
see shortly that $\CC_Z$ has finite products; we shall write them as
$U\tm V$, although they are actually given by the tensor product over
$R$.  We also write $1$ for the terminal object, which is the
coalgebra $R$ with $\psi_R=\ep_R=1_R$.

The following result is well-known when $R$ is a field, but we outline
a proof to show that nothing goes wrong for more general rings.
\begin{proposition}\label{prop-coalg-limits}
 The category $\CC_Z$ \index{CZ@$\CC_Z$} has finite products, and
 strong colimits for all small diagrams.  The forgetful functor to
 $\CM_Z$ creates colimits.
\end{proposition}
\begin{proof}
 Given two coalgebras $U,V$, we make $U\ot V$ into a coalgebra with
 counit $\ep_U\ot\ep_V\:U\ot V\xra{}R$ and coproduct
 \[ U\ot V \xra{\psi_U\ot\psi_V} 
    U\ot U \ot V\ot V \xra{1\ot\tau\ot 1} U\ot V\ot U\ot V.
 \]
 This is evidently functorial in $U$ and $V$.  There are two
 projections $\pi_U=1\ot\ep_V\:U\ot V\xra{}U$ and
 $\pi_V=\ep_U\ot 1\:U\ot V\xra{}V$, and one checks that these are
 coalgebra maps.  One also checks that a pair of maps $f\:W\xra{}U$
 and $W\xra{}V$ yield a coalgebra map
 $h=(f,g)=(f\ot g)\circ\psi_W\:W\xra{}U\ot V$, and that this is the
 unique map such that $\pi_U\circ h=f$ and $\pi_V\circ h=g$.  Thus,
 $U\ot V$ is the categorical product of $U$ and $V$.  Similarly, we
 can make $R$ into a coalgebra with $\psi_R=\ep_R=1_R$, and this makes
 it a terminal object in $\CC_Z$.

 Now suppose we have a diagram of coalgebras $U_i$, and let
 $U=\colim_iU_i$ denote the colimit in $\CM_Z$.  Because tensor
 products are right exact, we see that
 $U\ot U=\colim_{i,j}U_i\ot U_j$, so there is an obvious map
 $U_i\ot U_i\xra{}U\ot U$.  By composing with the coproduct on $U_i$,
 we get a map $U_i\xra{}U\ot U$.  These maps are compatible with the
 maps of the diagram, so we get a map $U=\colim_iU_i\xra{}U\ot U$.  We
 use this as the coproduct on $U$.  The counit maps $U_i\xra{}R$ also
 fit together to give a counit map $U\xra{}R$, and this makes $U$ into
 a coalgebra.  One can check that this gives a colimit in the category
 $\CC_Z$.  Thus, $\CC_Z$ has colimits and they are created in $\CM_Z$.
 It is clear from the construction that
 $V\tm\colim_iU_i=\colim_i(V\tm U_i)$, because tensoring with $V$ is
 right exact. 
\end{proof}

Let $f\:R\xra{}S=\OY$ be a map of rings, and let
$T_f\:\CM_Z\xra{}\CM_Y$ be the functor $M\mapsto S\ot_RM$.  This
clearly gives a functor $\CC_Z\xra{}\CC_Y$ which preserves finite
products and all colimits.

We now introduce a class of coalgebras with better than usual
behaviour under duality.
\begin{definition}\label{defn-good-basis}
 Let $U$ be a coalgebra over $R$, and suppose that $U$ is free as an
 $R$-module, say $U=R\{e_i\st i\in I\}$.  For any finite set $J$ of
 indices, we write $U_J=R\{e_i\st i\in J\}$; if this is a subcoalgebra
 of $U$, we call it a \dfn{standard subcoalgebra}.  We say that
 $\{e_i\}$ is a \dfn{good basis} if each finitely generated submodule
 of $U$ is contained in a standard subcoalgebra.  We write $\CC'_Z$
 \index{CZp@$\CC'_Z$} for the category of those coalgebras that admit
 a good basis.  It is easy to see that $\CC'_Z$ is closed under finite
 products.
\end{definition}

\begin{proposition}\label{prop-sch}
 There is a full and faithful functor
 $\sch=\sch_Z\:\CC'_Z\xra{}\FX_Z$, \index{schZ@$\sch_Z$} which
 preserves finite products and commutes with base change.  Moreover,
 $\sch(U)$ is always \idx{solid} and we have
 $\OO_{\sch(U)}=U^\vee:=\Hom_R(U,R)$.
\end{proposition}
\begin{proof}
 Let $U$ be a coalgebra in $\CC'_Z$.  For each subcoalgebra $V\leq U$
 such that $V$ is a finitely generated free module over $R$, we define
 $V^\vee=\Hom_R(V,R)$.  We can clearly make this into an $R$-algebra
 using the duals of the coproduct and counit maps, so we have a scheme
 $\spec(V^\vee)$ over $Z$.  We define
 $\sch(U)=\colim_V\spec(V^\vee)\in\FX_Z$.  If we choose a good basis
 $\{e_i\st i\in I\}$for $U$ then it is clear that the standard
 subcoalgebras form a cofinal family of $V$'s, so we have
 $\sch(U)=\colim_J\spec(U_J^\vee)$, where $J$ runs over the finite
 subsets of $I$ for which $U_J$ is a subcoalgebra.  This is clearly a
 directed, and thus filtered, colimit.  It follows that
 $\OO_{\sch(U)}=\invlim_JU_J^\vee=U^\vee$.  The resulting topology on
 $U^\vee=\Hom_R(U,R)$ is just the topology of pointwise convergence,
 where we give $R$ the discrete topology.  We can also think of this
 as $\prod_IR$, and the topology is just the product topology.  It is
 clear from this that $\sch(U)$ is solid.

 If $V$ is another coalgebra with good basis, then the obvious basis
 for $U\ot_RV$ is also good.  Moreover, if $U_J$ and $V_K$ are
 standard subcoalgebras of $U$ and $V$, then $U_J\ot_RV_J$ is a
 standard subcoalgebra of $U\ot_RV$, and the subcoalgebras of this
 form are cofinal among all standard subcoalgebras of $U\ot_RV$.  It
 follows easily that $\sch(U\tm V)=\sch(U\ot_RV)=
  \colim_{J,K}\spec(U_J^\vee)\tZ \spec(V_K^\vee)$.  As finite
 limits commute with filtered colimits in $\FX$, this is the same as
 $\sch(U)\tZ \sch(V)$.  

 Now consider a map $Y=\spec(S)\xra{}Z$ of schemes.  The claim is that
 the functors $\sch_Y$ and $\sch_Z$ commute with base change, in other
 words that $\sch_Y(S\ot_R U)=Y\tZ \sch_Z(U)$.  As pullbacks commute
 with filtered colimits, the right hand side is just
 $\colim_J\spec(S\ot_RU_J)$, which is the same as the left hand side.
\end{proof}

\begin{definition}\label{defn-coalgebraic}
 Let $Z$ be an informal scheme.  We write $\FX'_Z$ for the image of
 $\sch_Z$, which is a full subcategory of $\FX_Z$.  We say that a
 formal scheme $Y$ is \dfn{coalgebraic} over $Z$ if it lies in
 $\FX'_Z$.  We say that $Y$ is \emph{finitely coalgebraic} over $Z$
 \index{coalgebraic!finitely} if $\OY$ is a finitely generated free
 module over $\OZ$, or equivalently $Y$ is finite and very flat over
 $Z$; this easily implies that $Y$ is coalgebraic over $Z$.

 More generally, let $Z$ be a formal scheme, and $Y$ a formal scheme
 over $Z$.  We say that $Y$ is (finitely) coalgebraic over $Z$ if for
 all informal schemes $Z'$ over $Z$, the pullback $Z'\tZ Y$ is
 (finitely) coalgebraic over $Z'$.  We again write $\FX'_Z$ for the
 category of coalgebraic formal schemes over $Z$.
\end{definition}

\begin{example}
 Let $Z$ be a space such that $H_*(Z;\Zh)$ is a free Abelian group,
 concentrated in even degrees.  It is not hard to check that $E_0Z$ is
 a coalgebra over $E^0$ which admits a good basis, and that
 $Z_E=\sch_{E^0}(E_0Z)$.  Details are given in
 Section~\ref{sec-topology}. 
\end{example}
\begin{remark}\label{rem-sch-inverse}
 The functor $\sch_X\:\CC'_X\xra{}\FX'_X$ is an equivalence of
 categories, with inverse
 $Y\mapsto cY=\Hom^{\text{cts}}_{\OX}(\OY,\OX)$. \index{cY@$cY$}
\end{remark}

\begin{remark}\label{rem-group-like}
 For any coalgebra $U$, we say that an element $u\in U$ is
 \emph{group-like} if $\ep(u)=1$ and $\psi(u)=u\ot u$, or equivalently
 if the map $R\xra{}U$ defined by $r\mapsto ru$ is a coalgebra map.
 We write $\GL(U)=\CC_R(R,U)$ for the set of group-like elements.  If
 $U$ is a finitely generated free module over $R$, then it is easy to
 check that $\GL(U)=\Alg_R(U^\vee,R)$.  {}From this one can deduce that
 \[ \FX_Z(Y,\sch_Z(U)) = \GL(\OY\ot_R U), \]
 where we regard $\OY\ot_R U$ as a coalgebra over $\OY$.  This gives
 another useful characterisation of $\sch_Z(U)$.
\end{remark}

\begin{proposition}\label{prop-strong-colim}
 Let $\{U_i\}$ be a diagram in $\CC_Z$ with colimit $U$, and suppose
 that $U$ and $U_i$ actually lie in $\CC'_Z$.  Then $\sch(U)$ is the
 \idx{strong colimit} in $\FX_Z$ of the formal schemes $\sch(U_i)$.
\end{proposition}
\begin{proof}
 Note that $U=\colim_iU_i$ as $R$-modules (because colimits in $\CC_Z$
 are created in $\CM_Z$), and it follows immediately that
 $U^\vee=\invlim_iU_i^\vee$ as rings.  There are apparently two
 possible topologies on $U^\vee$.  The first is as in the definition
 of $\sch_R(U)$, where the basic neighbourhoods of zero are the
 submodules $\ann(M)$, where $M$ runs over finitely generated
 submodules of $U$.  The second is the inverse limit topology: for
 each index $i$ and each finitely generated submodule $N$ of $U_i$,
 the preimage of the annihilator of $N$ under the evident map
 $U^\vee\xra{}U_i^\vee$ is a neighbourhood of zero.  This is just the
 same as the annihilator of the image of $N$ in $U$, and
 neighbourhoods of this form give a basis for the inverse limit
 topology.  Given this, it is easy to see that the two topologies in
 question are the same.  We thus have an inverse limit of topological
 rings.  As the category of formal schemes is just dual to the
 category of formal rings, we have a colimit diagram of formal
 schemes, so $\sch(U)=\colim_i\sch(U_i)$.

 We need to show that the colimit is strong, in other words that for
 any formal scheme $T$ over $Z$ we have
 $T\tZ \sch_Z(U)=\colim_i(T\tZ \sch_Z(U_i))$.  First suppose that
 $T=\spec(B)$ is an informal scheme.  We then have
 $T\tZ \sch_Z(U)=\sch_T(B\ot_RU)$ and similarly for each $U_i$, and
 $B\ot_RU=\colim_iB\ot_RU_i$ because tensor products are right exact.
 By the first part of the proof (with $R$ replaced by $B$) we see that
 $T\tZ \sch_Z(U)=\colim_i(T\tZ \sch_Z(U_i))$ as required.

 If $T$ is a formal scheme, we write it as a strong filtered colimit
 of informal schemes $T_k$.  The colimit of the isomorphisms 
 $T_k\tZ \sch_Z(U)=\colim_i(T_k\tZ \sch_Z(U_i))$ is the required
 isomorphism $T\tZ\sch_Z(U)=\colim_i(T\tZ \sch_Z(U_i))$.
\end{proof}
\begin{example}\label{eg-strong-sym-formal}
 If $X$ is coalgebraic over $Y$ we claim that $X^n_Y/\Sg_n$ is a
 strong colimit for the action of $\Sg_n$ on $X^n_Y$.  To see this, we
 first suppose that $Y$ is informal and $X=\sch_Y(U)$ for
 some coalgebra $U$ that is free over $X$ with good basis
 $\{e_i\st i\in I\}$ say.  Then $X^n_Y=\sch_Y(U^{\ot n})$, and
 the set of terms $e_{\un{i}}=e_{i_1}\ot\ldots\ot e_{i_n}$ for
 $\un{i}=(i_1,\ldots,i_n)\in I^n$ is a good basis for $I^n$.  For each
 orbit $j\in I^n/\Sg_n$, we choose an element $\un{i}$ of the orbit
 and let $f_j$ be the image of $e_{\un{i}}$ in $U^{\ot n}/\Sg_n$.
 We find that the terms $f_j$ form a good basis for $U^{\ot n}/\Sg_n$,
 so this coalgebra lies in $\CC'_Y$.  It follows from
 Proposition~\ref{prop-strong-colim} that
 $X^n_Y=\sch_Y(U^{\ot n}/\Sg_n)$, and that this is a strong colimit.
 For a general base $Y$, we choose a presentation $Y=\colim_iY_i$ and
 write $X_i=X\tY Y_i$ and $Z_i=(X_i)^n_{Y_i}/\Sg_n$.  By what we have
 just proved, this is an object of $\FX_{\{Y_i\}}$, with
 $\colim_iZ_i=X^n_Y/\Sg_n$.  It is now easy to see that this is a
 strong colimit, using the ideas of
 Proposition~\ref{prop-filtered-base}.  
\end{example}

We conclude this section with a result about gradings.
\begin{proposition}\label{prop-graded-coalgebra}
 Let $Y$ be a coalgebraic formal scheme over an informal scheme $X$,
 and suppose that $X$ and $Y$ have compatible actions of $\MG$.  Then
 $cY$ has a natural structure as a graded coalgebra over
 $\OX$. \index{grading}
\end{proposition}
\begin{proof}
 Write $R=\OX$ and $U=cY$.  Proposition~\ref{prop-grading-action}
 makes $R$ into a graded ring.  Next, observe that $\OY=U^\vee$ and
 $\OO_{\MG\tm Y}=U^\vee\hot\Zh[t^{\pm 1}]$, which is the ring of
 doubly infinite Laurent series $\sum_{k\in\Zh}a_kt^k$ such that
 $a_k\in U^\vee$ and $a_k\xra{}0$ as $|k|\xra{}\infty$.  Thus, the
 action $\al\:\MG\tm Y\xra{}Y$ gives a continuous homomorphism
 $\al^*\:U^\vee\xra{}U^\vee\hot\Zh[t^{\pm 1}]$, say
 $\al^*(a)=\sum_ka_kt^k$.  The basic neighbourhoods of zero in
 $U^\vee$ are the kernels of the maps $U^\vee\xra{}W^\vee$, where $W$
 is a standard subcoalgebra of $U$.  Similarly, the basic
 neighbourhoods of zero in $U^\vee\hot\Zh[t^{\pm 1}]$ are the kernels
 of the maps to $V^\vee[t^{\pm 1}]$, where $V$ is a standard
 subcoalgebra.  Thus, continuity means that for every standard
 subcoalgebra $V\leq U$, there is a standard subcoalgebra $W$ such
 that whenever $a(W)=0$ we have $a_k(V)=0$ for all $k$.  In
 particular, it follows that the map $\pi_k\:a\xra{}a_k$ is
 continuous.  Just as in the proof of
 Proposition~\ref{prop-grading-action}, we see that $\sum_ka_k=a$ and
 that $\pi_j\pi_k=\dl_{jk}\pi_k$.  It follows that $U^\vee$ is a kind
 of completed direct sum of the subgroups $\image(\pi_k)$.  We would
 like to dualise this and thus split $U$ as an honest direct sum.

 First, we need to show that the maps $\pi_i$ have a kind of
 $R$-linearity.  Let $r$ be an element of $R$, and let $r_i$ be the
 part in degree $i$, so that $r=\sum_ir_i$ and $r_i=0$ for almost all
 $i$.  Using the compatibility of the actions, we find that
 $(ra)_i=\sum_j r_ja_{i-j}$ (which is really a finite sum).  

 Suppose that $u\in U$.  Choose a standard subcoalgebra $V$ containing
 $u$, and let $W$ be a standard subcoalgebra such that whenever
 $a(W)=0$ we have $a_i(V)=0$ for all $i$.

 Suppose that $a\in U^\vee$.  It follows from our asymptotic condition
 on Laurent series that $a_i(u)=0$ when $|i|$ is large, so we can
 define $\chi_k(u)(a)=\sum_i a_i(u)_{i+k}\in R$.  We then have
 \begin{align*}
  \chi_k(u)(ra) &= \sum_i ((ra)_i(u))_{i+k}             \\
                &= \sum_{i,j} (r_j a_{i-j}(u))_{i+k}    \\
                &= \sum_{i,j} r_j (a_{i-j}(u))_{i+k-j}  \\
                &= \sum_{m,j} r_j a_m(u)_{m+k}          \\
                &= r \chi_k(u)(a).
 \end{align*}
 Thus, the map $\chi_j(u)\:U^\vee\xra{}R$ is $R$-linear.  Clearly, if
 $a(W)=0$ then $\chi_j(u)(a)=0$, so $\chi_j(u)$ can be regarded as an
 element of $(U^\vee/\ann(W))^\vee=W^{\vee\vee}=W$ (because $W$ is a
 finitely generated free module).  More precisely, there is a unique
 element $u_j\in U$ such that $\chi_j(u)(a)=a(u_j)$ for all $a$, and
 in fact $u_j\in W$.

 Next, we choose a finite set of elements in $U^\vee$ which project to
 a basis for $W^\vee$.  We can then choose a number $N$ such that
 $b_i(u)=0$ whenever $b$ lies in that set and $|i|>N$.  Because
 $a_j(u)=0$ for all $j$ whenever $a(W)=0$, we conclude that $a_i(u)=0$
 for all $a\in U^\vee$ and all $i$ such that $|i|>N$.  It follows that
 $u_i=0$ when $|i|>N$.  This justifies the following manipulation:
 $a(u)=\sum_{i,j}a_i(u)_j=\sum_ja(u_j)=a(\sum_ju_j)$.  We conclude
 that $u=\sum_ju_j$.  We define a map $\phi_i\:U\xra{}U$ by
 $\phi_i(u)=u_i$, and we define $U_i=\image(\phi_i)$.  We leave it to
 the reader to check that $\phi_i\phi_j=\dl_{ij}\phi_j$, so that
 $U=\bigoplus_iU_i$, and that this grading is compatible with the
 $R$-module structure and the coalgebra structure.
\end{proof}

\subsection{More mapping schemes}
\label{subsec-more-mappings}

Recall the functor $\Map_Z(X,Y)$ \index{MapZXY@$\Map_Z(X,Y)$}, given
in Definition~\ref{defn-maps-informal}.  We now prove some more
results which tell us when $\Map_Z(X,Y)$ is a scheme or a formal
scheme.

First, note that for any functor $W$ over $Z$, we have
\[ \CF_Z(W,\Map_Z(X,Y)) =
   \CF_Z(W\tZ X,Y) = 
   \CF_W(W\tZ X,W\tZ Y).
\]
Indeed, if $W$ is informal then this follows from the definitions and
the Yoneda lemma, by writing $W$ in the form $\spec(R)$.  The general
case follows from this by taking limits, because every functor is the
colimit of a (not necessarily small or filtered) diagram of
representable functors.  It is also not hard to give a direct proof.

Conversely, suppose we have a functor $M$ over $Z$ and a natural
isomorphism $\CF_Z(W,M)\simeq\CF_Z(W\tZ X,Y)$ for all informal schemes
$W$ over $Z$.  It is then easy to identify $M$ with $\Map_Z(X,Y)$.

\begin{lemma}\label{lem-MapZXY-exists}
 Let $X$ and $Y$ be functors over $Z$, and suppose that $X$ and $Z$
 are formal schemes.  Then $\Map_Z(X,Y)(R)$ is a set for all $R$, so
 the functor $\Map_Z(X,Y)$ exists.
\end{lemma}
\begin{proof}
 We have only a set of elements $z\in Z(R)$, so it suffices to check
 that for any such $z$ there is only a set of maps $X_z\xra{}Y_z$ of
 functors over $\spec(R)$.  Here $X_z$ is a formal scheme, with
 presentation $\{W_i\}$ say.  Clearly $\CF(W_i,Y_z)=Y_z(\OO_{W_i})$ is
 a set, and $\CF_{\spec(R)}(X_z,Y_z)$ is a subset of
 $\prod_i\CF(W_i,Y_z)$.  
\end{proof}

Recall also from Proposition~\ref{prop-map-informal} that
$\Map_Z(X,Y)$ is a scheme when $X$, $Y$ and $Z$ are all informal
schemes, and $X$ is finite and very flat over $Z$.
\begin{definition}
 We say that a formal scheme $Y$ over $Z$ is of finite presentation if
 there is an equaliser diagram in $\FX_Z$ of the form
 \[ Y \xra{} \aff^n\tm Z \arw{e,=>} \aff^m\tm Z. \]
\end{definition}

\begin{theorem}\label{thm-map-formal}
 Let $X$ and $Y$ be formal schemes over $Z$.  Then $\Map_Z(X,Y)$
 \index{MapZXY@$\Map_Z(X,Y)$} is a formal scheme if
 \begin{itemize}
  \item[(a)] $X$ is coalgebraic over $Z$ and $Y$ is relatively informal
   over $Z$, or
  \item[(b)] $X$ is finite and very flat over $Z$, or
  \item[(c)] $X$ is very flat over $Z$ and $Y$ is of finite
   presentation over $Z$.
 \end{itemize}
\end{theorem}
This will be proved at the end of the section, after some auxiliary
results.

\begin{lemma}\label{lem-map-pullback}
 If $Z'$ is a functor over $Z$ then
 $\Map_{Z'}(X\tZ Z',Y\tZ Z')=\Map_Z(X,Y)\tZ Z'$.
 \index{MapZXY@$\Map_Z(X,Y)$} 
\end{lemma}
\begin{proof}
 If $W$ is a scheme over $Z'$ then 
 \begin{align*}
     \CF_{Z'}(W,\Map_Z(X,Y)\tZ Z') 
  &= \CF_Z(W,\Map_Z(X,Y))                     \\
  &= \CF_Z(W\tZ X, Y)                         \\
  &= \CF_{Z'}(W\tZ X,Y\tZ Z')                 \\
  &= \CF_{Z'}(W\tm_{Z'}(X\tZ Z'),Y\tZ Z').
 \end{align*}
 Thus, $\Map_Z(X,Y)\tZ Z'$ has the required universal property.
\end{proof}

\begin{lemma}\label{lem-map-colim}
 If $X$ is a strong colimit of formal schemes $X_i$ and
 $\Map_Z(X_i,Y)$ is a formal scheme and is relatively informal over
 $Z$ for all $i$ then $\Map_Z(X,Y)$ is a formal scheme and is equal to
 $\invlim_i\Map_Z(X_i,Y)$ (where the inverse limit is computed in
 $\FX_Z$).  
\end{lemma}
Note that coproducts and filtered colimits are always strong, so the
lemma applies in those cases.
\begin{proof}
 Because $\Map_Z(X_i,Y)$ is relatively informal,
 Proposition~\ref{prop-rel-inf-lim} allows us to form the limit
 $\invlim_i\Map_Z(X_i,Y)$ in $\FX_Z$.  If $W$ is a formal
 scheme over $Z$ then we have
 \begin{align*}
  \FX_Z(W,\invlim_i\Map_Z(X_i,Y)) 
  &= \invlim_i\FX_Z(W,\Map_Z(X_i,Y))           \\
  &= \invlim_i\FX_Z(W\tZ X_i,Y)                \\
  &= \FX_Z(\colim_i W\tZ X_i,Y)                \\
  &= \FX_Z(W\tZ X,Y).
 \end{align*}
 This proves that $\invlim_i\Map_Z(X_i,Y)=\Map_Z(X,Y)$ as required.
\end{proof}

We leave the next lemma to the reader.
\begin{lemma}\label{lem-map-lim}
 Suppose that $Y$ is an inverse limit of a finite diagram of formal
 schemes $\{Y_i\}$ over $Z$.  Then
 $\Map_Z(X,Y)=\invlim_i\Map_Z(X,Y_i)$, where the limit is computed in
 $\CF_Z$.  Thus, if $\Map_Z(X,Y_i)$ is a formal scheme for all $i$,
 then $\Map_Z(X,Y)$ is a formal scheme.  \qed
\end{lemma}

\begin{lemma}\label{lem-map-filt-base}
 Let $\{Z_i\}$ be a filtered system of informal schemes with colimit
 $Z$.  Let $X$ and $Y$ be formal schemes over $Z$, with
 $X_i=X\tZ Z_i$ and $Y_i=Y\tZ Z_i$.  If $\Map_{Z_i}(X_i,Y_i)$ is a
 formal scheme for all $i$ then $\Map_Z(X,Y)$ is a formal scheme and
 is equal to $\colim_i\Map_{Z_i}(X_i,Y_i)$.
\end{lemma}
\begin{proof}
 Lemma~\ref{lem-map-pullback} tells us that the system of formal
 schemes 
 \[ M_i=\Map_{Z_i}(X_i,Y_i) \]
 defines an object of the category $\FX_{\{Z_i\}}$ of
 Proposition~\ref{prop-filtered-base}.  Thus, if we define
 $M=\colim_iM_i$ we find that $\FX_Z(W,M)$ is the set of maps of
 diagrams $\{W\tZ Z_i\}\xra{}\{M_i\}$ over $\{Z_i\}$.  This is the
 same as the set of maps of diagrams
 $\{W\tZ X_i\}=\{W\tZ Z_i\tm_{Z_i}X_i\}\xra{}\{Y_i\}$ over $\{Z_i\}$.
 By the adjunction in Proposition~\ref{prop-filtered-base}, this is
 the same as the set of maps $W\tZ X=\colim_iW\tZ X_i\xra{}Y$ over
 $Z$.  Thus, $M$ has the defining property of $\Map_Z(X,Y)$.
\end{proof}

\begin{lemma}\label{lem-map-filt-target}
 Let $X$ be relatively informal over $Z$, and let $\{Y_i\}$ be a
 filtered system of formal schemes over $Z$ with colimit $Y$.  If
 $\Map_Z(X,Y_i)$ is a formal scheme for all $i$, then $\Map_Z(X,Y)$
 is a formal scheme and is equal to $\colim_i\Map_Z(X,Y_i)$.
\end{lemma}
\begin{proof}
 Write $M=\colim_i\Map_Z(X,Y_i)$.  Let $W$ be an informal scheme over
 $Z$.  As $X$ is relatively informal, we see that $W\tZ X$ is
 informal.  It follows that the functors $\FX_Z(W,-)$ and
 $\FX_Z(W\tZ X,-)$ preserve filtered colimits.  We thus have
 \begin{align*}
  \FX_Z(W,M) 
  &= \colim_i\FX_Z(W,\Map_Z(X,Y_i)) \\
  &= \colim_i\FX_Z(W\tZ X,Y_i)  \\
  &= \FX_Z(W\tZ X,\colim_iY_i) \\
  &= \FX_Z(W\tZ X,Y),
 \end{align*}
 as required.
\end{proof}

\begin{lemma}\label{lem-map-aff}
 If $X$ and $Z$ are informal and $X$ is very flat over $Z$ then
 the functor $\Map_Z(X,\aff^1\tm Z)$ is a formal scheme.
\end{lemma}
\begin{proof}
 We can choose a basis for $\OX$ over $\OZ$ and thus write $\OX$ as a
 filtered colimit of finitely generated free modules $M_i$ over $\OZ$.
 {}From the definitions we see that $\Map_Z(X,\aff^1)(R)$ is the set of
 pairs $(x,u)$, where $x\in X(R)$ (making $R$ into an $\OX$-algebra)
 and $u$ is a map $R[t]\xra{}R\ot_\OZ\OX$ of $R$-algebras.  This is of
 course equivalent to an element of $R\ot_\OZ\OX=\colim_iR\ot_\OZ
 M_i$.  Thus, we see that $\Map_Z(X,\aff^1\tm Z)=\colim_i\aff(M_i)$,
 which is a formal scheme.
\end{proof}

\begin{proof}[Proof of Theorem~\ref{thm-map-formal}]
 We shall prove successively that $\Map_Z(X,Y)$ is a formal scheme
 under any of the following hypotheses.  Cases~(3), (5) and~(7) give
 the results claimed in the theorem.
 \begin{itemize}
 \item[(1)] $X$, $Y$ and $Z$ are informal, and $X$ is finite and very
  flat.  In this case $\Map_Z(X,Y)$ is informal.
 \item[(2)] $Y$ is informal, and $X$ is finite and very flat.  In this
  case $\Map_Z(X,Y)$ is relatively informal.
 \item[(3)] $X$ is finite and very flat.
 \item[(4)] $Y$ and $Z$ are informal, and $X$ is coalgebraic.  In this
  case, $\Map_Z(X,Y)$ is informal.
 \item[(5)] $Y$ is relatively informal, and $X$ is coalgebraic.  In
  this case, $\Map_Z(X,Y)$ is relatively informal.
 \item[(6)] $X$ and $Z$ are informal, $X$ is very flat, and $Y$ is of
  finite presentation.
 \item[(7)] $X$ is very flat and $Y$ is of finite presentation.
 \end{itemize}
 Proposition~\ref{prop-map-informal} gives case~(1).  For case~(2),
 write $Z=\colim_iZ_i$ in the usual way.  Then case~(1) tells us that
 $\Map_{Z_i}(X\tZ Z_i,Y\tZ Z_i)$ is an informal scheme.  Using this
 and Lemma~\ref{lem-map-filt-base}, we see that $\Map_Z(X,Y)$ is a
 formal scheme.  Using case~(1) and Lemma~\ref{lem-map-pullback} we
 see that $\Map_Z(X,Y)$ is relatively informal.  In case~(3), we write
 $Y$ as a filtered colimit of informal schemes $Y_j$.  Case~(2) tells
 us that $\Map_Z(X,Y_j)$ is a relatively informal scheme, so
 Lemma~\ref{lem-map-filt-target} tells us that $\Map_Z(X,Y)$ is a
 formal scheme.  In case~(4), it follows easily from the definitions
 that $X$ can be written as the filtered colimit of a system of
 finite, very flat schemes $X_i$.  It then follows from case~(1) that
 $\Map_Z(X_i,Y)$ is an informal scheme.  Using
 Lemma~\ref{lem-map-colim} we see that
 $\Map_Z(X,Y)=\invlim_i\Map_Z(X_i,Y)$.  This is an inverse limit of
 informal schemes, and thus is an informal scheme.  We deduce~(5)
 from~(4) in the same way that we deduced~(2) from~(1).  Case~(6)
 follows easily from Lemmas~\ref{lem-map-aff} and~\ref{lem-map-lim}.
 Again, the argument for (1)$\Rightarrow$(2) also gives
 (6)$\Rightarrow$(7).
\end{proof}

\section{Formal curves}
\label{sec-formal-curves}

In this section, we define formal curves.  We also study
\idx{divisor}s, differentials, and meromorphic functions on such
curves.

Let $X$ be a formal scheme, and let $C$ be a formal scheme over $X$.
We say that $C$ is a \dfn{formal curve} over $X$ if it is isomorphic
in $\FX_X$ to $\haf^1\tm X$.  (In some sense, it would be better to
allow formal schemes that are only isomorphic to $\haf^1\tm X$
fpqc-locally on $X$, but this seems unnecessary for the topological
applications so we omit it.)  A \dfn{coordinate} on $C$ is a map
$x\:C\xra{}\haf^1$ giving rise to an isomorphism $C\simeq\haf^1\tm X$.
\begin{example}
 If $E$ is an even periodic ring spectrum then $(\cp^\infty)_E$ and
 $(\mathbb{H}P^\infty)_E$ are formal curves over $S_E$.
\end{example}

\subsection{Divisors on formal curves}
\label{subsec-divisors-formal-curves}

Let $C$ be a formal curve over $X$, and let $D$ be a closed subscheme
of $X$.  If $X$ is informal, we say that $D$ is a \emph{effective
  divisor of degree $n$ on $C$} \index{divisor} if $D$ is informal,
and $\OO_D$ is a free module of rank $n$ over $\OX$.  If $X$ is a
general formal scheme, we say that $D$ is a divisor if $D\tX X'$ is a
divisor on $C\tX X'$, for all informal schemes $X'$ over $X$.  If $Y$
is a formal scheme over $X$, we refer to divisors on $C\tX Y$ as
divisors on $C$ over $Y$.

\begin{proposition}\label{prop-divisors-formal}
 There is a formal scheme $\Div_n^+(C)$ \index{Divpn@$\Div^+_n(C)$}
 over $X$ such that maps $Y\xra{}\Div_n^+(C)$ over $X$ biject with
 effective \idx{divisor}s of degree $n$ on $C$ over $Y$.  Moreover, a
 choice of coordinate on $C$ gives rise to an isomorphism
 $\Div_n^+(C)\simeq\haf^n\tm X$.
\end{proposition}
\begin{proof}
 This is much the same as Example~\ref{eg-divisors}.  We define 
 \begin{multline*}
  \Div_n^+(C)(R) = \\
    \{ (a,D)\st a\in X(R)
    \text{ and $D$ is an effective divisor of degree $n$ on $C_a$ } 
    \}.
 \end{multline*}
 We make this a functor by pullback, just as in
 Example~\ref{eg-divisors}.  To see that $\Div^+_n(C)$ is a formal
 scheme, choose a coordinate $x$ on $C$.  Given a point $(a,D)$ as
 above, we find that $C_a=C\tX\spec(R)=\spf(\fps{R}{x})$, where the
 topology on $\fps{R}{x}$ is defined by the ideals $(x^k)$.  We know
 that $D$ is a closed subscheme of $C_a$, and that $D$ is informal.  It
 follows that $D=\spec(\fps{R}{x}/J)$ for some ideal $J$ such that
 $x^k\in J$ for some $k$.  Let $\lm(x)$ be the endomorphism of $\OO_D$
 given by multiplication by $x$, and let
 $f_D(t)=\sum_{i=0}^na_i(D)t^{n-i}$ be the characteristic polynomial
 of $\lm(x)$.  As $x^k\in J$, we see that $\lm(x)^k=0$.  If $R$ is a
 field, then we deduce by elementary linear algebra that $f_D(t)=t^n$.
 If $\pri$ is a prime ideal in $R$ then by considering the divisor
 $\spec(\kp(\pri))\tm_{\spec(R)}D$, we conclude that
 $f_D(t)=t^n\pmod{\pri[t]}$.  Using Proposition~\ref{prop-nilpotents},
 we deduce that $a_i(D)\in\Nil(R)$ for $i>0$.  Thus, the $a_i$'s give
 a map $\Div^+_n(C)\xra{}\haf^n\tm X$.  As in
 Example~\ref{eg-divisors}, the Cayley-Hamilton theorem tells us that
 $f_D(x)\in J$ and thus that $\OO_D=R[x]/f_D(x)=\fps{R}{x}/f_D(x)$.

 Conversely, suppose we have elements $b_0,\ldots,b_n$ with $b_0=1$
 and $b_i\in\Nil(R)$ for $i>0$ and we define $g(t)=\sum_ib_it^{n-i}$
 and $D=\spf(\fps{R}{x}/g(x))$.  In $\OO_D$ we have
 $x^n=-\sum_{i>0}b_ix^{n-i}$, which is nilpotent, so $x$ is nilpotent,
 so $(g(x))$ is open in $\fps{R}{x}$.  This means that $D$ is informal
 and that $\OO_D=R[x]/g(x)$, which is easily seen to be a free module
 of rank $n$ over $R$.  Thus, $D$ is an effective divisor of rank $n$
 on $C_a$.  We conclude that $\Div_n^+(C)$ is isomorphic to $\haf^n$,
 and in particular is a formal scheme.

 If $Y$ is an arbitrary formal scheme over $X$, we can choose a
 presentation $Y=\colim_iY_i$, so $Y_i$ is an informal scheme over
 $X$.  The above tells us that maps $Y_i\xra{}\Div_n^+(C)$ over $X$
 biject with effective divisors of degree $n$ on $C$ over $Y_i$.
 Thus, maps $Y\xra{}\Div_n^+(C)$ over $X$ biject with systems of
 divisors $D_i$ over $Y_i$, such that for each map $Y_i\xra{}Y_j$ we
 have $D_i=D_j\tm_{Y_j}Y_i$.  Using
 Proposition~\ref{prop-filtered-base}, we see that these biject with
 effective divisors of degree $n$ on $C$ over $Y$.
\end{proof}
\begin{example}\label{eg-BUn-div}
 It is essentially well-known that $BU(n)_E=\Div_n^+(G_E)$, where
 $G_E=(\cp^\infty)_E$.  A proof will be given in
 Section~\ref{sec-topology}. 
\end{example}
\begin{remark}\label{rem-Div-pullback}
 It is not hard to check that for any map $Y\xra{}X$ of formal schemes
 and any formal curve $C$ over $X$ we have
 $\Div_n^+(C\tX Y)=\Div_n^+(C)\tX Y$ (because both sides represent
 the same functor $\FX_Y\xra{}\Sets$). \index{divisor}
 \index{Divpn@$\Div^+_n(C)$} 
\end{remark}

\begin{definition}\label{defn-JD}
 Let $D$ be an effective \idx{divisor} on a curve $C$ over $X$.  We
 shall define an associated \idx{line bundle} $J(D)$ \index{JD@$J(D)$}
 over $C$.  By Corollary~\ref{cor-fibrewise-sheaf}, it is enough to do
 this in a sufficiently natural way when $X$ is an informal scheme.
 In that case we have $\OO_D=\OC/J(D)$ for some ideal $J(D)$ in $\OC$.
 In terms of a coordinate $x$, we see from the proof of
 Proposition~\ref{prop-divisors-formal} that $J(D)$ is generated by a
 monic polynomial $f(x)$ whose lower coefficients are nilpotent.  Thus
 $f(x)=x^n-g(x)$ where $g(x)^k=0$ say.  If $fh=0$ then
 $x^{nk}h=g^kh=0$ so $h=0$, so $f$ is not a zero-divisor and $J(D)$ is
 free of rank one over $\OC$.  Thus, $J(D)$ can be regarded as a line
 bundle over $C$ as required (using
 Remark~\ref{rem-complete-modules}).
\end{definition}

\begin{proposition}\label{prop-divisor-sum}
 There is a natural commutative and associative addition
 $\sg\:\Div_j^+(C)\tX \Div_k^+(C)\xra{}\Div_{j+k}^+(C)$, such that
 $J(D+E)=J(D)\ot J(E)$. \index{JD@$J(D)$} \index{Divpn@$\Div^+_n(C)$}
\end{proposition}
\begin{proof}
 Let $a\:\spec(R)\xra{}X$ be an element of $X(R)$, and let $D$ and $E$
 be effective divisors of degrees $j$ and $k$ on $C_a$ over $\spec(R)$.
 We then have $D=V(J(D))$ and $E=V(J(E))$ where $J(D)$ and $J(E)$ are
 ideals in $\OO_{C_a}$.  We define $F=V(J(D)J(E))$.  If we choose a
 coordinate $x$ on $C$ we find (as in the proof of
 Proposition~\ref{prop-divisors-formal}) that $J(D)=(f_D(x))$ and
 $J(E)=(f_E(x))$, where $f_D$ and $f_E$ are monic polynomials whose
 lower coefficients are nilpotent.  This means that $g=f_Df_E$ is a
 polynomial of the same type, and it follows that $F=V(g)$ is a
 divisor of degree $j+k$ as required.  We define $\sg(D,E)=D+E=F$.  It
 is clear from the construction that $J(D+E)=J(D)\ot J(E)$.
\end{proof}

\begin{proposition}\label{prop-divisor-sym}
 Let $C$ be a \idx{formal curve} over a formal scheme $X$.  Then
 $\Div_n^+(C)=C^n_X/\Sg_n$,\index{Divpn@$\Div^+_n(C)$} and this is a
 \idx{strong colimit}.  Moreover, the quotient map
 $C^n_X\xra{}C^n_X/\Sg_n$ is faithfully flat.\index{flat!faithfully}
\end{proposition}
\begin{proof}
 First consider the case $n=1$.  Fix a ring $R$ and a point
 $a\in X(R)$, and write $C_a=C\tX \spec(R)$, which is a formal curve
 over $Y=\spec(R)$.  A point $c\in C$ lying over $a$ is the same as a
 section of the projection $C_a\xra{}Y$.  Such a section is a split
 monomorphism, and thus a closed inclusion; we write $[c]$ for its
 image, which is a closed formal subscheme of $C_a$.  The projection
 $C_a\xra{}Y$ carries $[c]$ isomorphically to $Y$, which shows that
 $[c]$ is an effective divisor of degree $1$ on $C$ over $Y$.  Thus,
 this construction gives a map $C\xra{}\Div_1^+(C)$.  If $x$ is a
 coordinate on $C$ then it is easy to see that $x(c)\in\Nil(R)$ and
 $[c]=\spf(\fps{R}{x}/(x-x(c)))$.  Using this, we see easily that our
 map is an isomorphism, giving the case $n=1$ of the Proposition.

 We now use the iterated addition map
 $C^n_X=\Div_1^+(C)^n_X\xra{}\Div_n^+(C)$ to get a map
 $C^n_X/\Sg_n\xra{}\Div_n^+(C)$.  

 Next, because $C\simeq\haf^1\tm X$, it is easy to see that $C$ is
 coalgebraic over $X$ and thus (by Example~\ref{eg-strong-sym-formal})
 that $C^n_X/\Sg_n$ is a strong colimit.  Given this, we can reduce
 easily to the case where $X$ is informal, say $X=\spec(R)$.  Choose a
 coordinate $x$ on $C$.  This gives isomorphisms
 $\OO_{\Div_n^+(C)}=\fps{R}{a_1,\ldots,a_n}=S$ and
 $\OO_{C^n_X}=\fps{R}{x_1,\ldots,x_n}=T$ and
 $\OO_{C^n_X/\Sg_n}=T^{\Sg_n}$.  The claim is thus that the map
 $S\xra{}T^{\Sg_n}$ is an isomorphism, and that $T$ is faithfully flat
 over $T^{\Sg_n}$.  The map $S\xra{}T^{\Sg_n}$ sends $a_i$ to the
 coefficient of $x^{n-i}$ in $\prod_j(x-x_j)$, which is (up to sign)
 the $i$'th elementary symmetric function of the variables $x_j$.  It
 is thus a well-known theorem of Newton that $S=T^{\Sg_n}$.  It is
 also well-known that the elements of the form
 $\prod_{j=1}^nx_j^{d_j}$ with $0\le d_j<j$ form a basis for $T$ over
 $T^{\Sg_n}$, so that $T$ is indeed faithfully flat over $T^{\Sg_n}$.
\end{proof}

We next consider pointed curves,\index{formal curve!pointed} in other
words curves $C$ equipped with a specified ``zero-section''
$0\:X\xra{}C$ such that the composite $X\xra{0}C\xra{}X$ is the
identity.  If $C$ is such a curve and $x$ is a coordinate on $C$, we
say that $x$ is \emph{normalised} \index{coordinate!normalised} if
$x(0)=0$.  If $y$ is an unnormalised coordinate then $x=y-y(0)$ is a
normalised one, so normalised coordinates always exist.

\begin{definition}\label{defn-Div}
 Let $C$ be a pointed \idx{formal curve} over $X$.  Define
 \[ f\:\Div_n^+(C)\xra{}\Div_{n+1}^+(C) \]
 by $f(D)=D+[0]$.  For $n\in\Zh$ with $n<0$ we write
 $\Div_n^+(C)=\emptyset$.  Define
 \begin{align*}
  \Div^+(C)     &= \coprod_{n\ge 0} \Div^+_n(C)                 \\
  \Div_n(C)     &=
   \colim(\Div_n^+(C) \xra{f} \Div_{n+1}^+(C) \xra{f} \ldots)   \\
  \Div(C)       &= \coprod_{n\in\Zh} \Div_n(C)                  \\
  &= \colim(\Div^+(C) \xra{f} \Div^+(C) \xra{f} \ldots).
 \end{align*}
 \index{Diva@$\Div^+(C)$}\index{Div@$\Div(C)$}
\end{definition}

It is not hard to see that $f^k$ induces an isomorphism
$\Div_n(C)\simeq\Div_{n+k}(C)$, so $\Div(C)$ can be identified with
$\coprod_n\Div_0(C)=\un{\Zh}\tm\Div_0(C)$.

A choice of normalised coordinate on $C$ gives an isomorphism
$\Div_n^+(C)\simeq X\tm\haf^n$.  Under this identification, $f$
becomes the map 
\[ (x,a_1,\ldots,a_n)\mapsto(x,a_1,\ldots,a_n,0). \]  
We thus have an isomorphism $\Div_0(C)=\haf^{(\infty)}$ (using the
notation of Example~\ref{eg-haf-infty}) and thus
$\Div=\un{\Zh}\tm\haf^{(\infty)}$.

\begin{definition}\label{defn-thom-sheaf}
 Given a divisor $D$ on a pointed curve $C$ over $X$, we define the
 \dfn{Thom sheaf} of $D$ to be the \idx{line bundle} $L(D)=0^*J(D)$
 over $X$.  It is clear that $L(D+E)=L(D)\ot L(E)$.  Note that a
 coordinate on $C$ gives a generator $f_D(x)$ for $J(D)$ and thus a
 generator $u_D$ for $L(D)$, which we call the \dfn{Thom class}.  This
 is natural for maps of $X$, and satisfies $u_{D+E}=u_D\ot u_E$.
\end{definition}

\begin{definition}\label{defn-coord}
 If $C$ is a pointed formal curve over $X$, we define a functor
 $\Coord(C)\in\CF_X$ \index{CoordC@$\Coord(C)$} by
 \[ \Coord(C)(R) =
     \{(a,x)\st a\in X(R)
            \text{ and $x$ is a normalised coordinate on $C_a$ }\}.
 \]    
 \index{coordinate!normalised}
\end{definition}
\begin{proposition}\label{prop-coord-scheme}
 The functor $\Coord(C)$ \index{CoordC@$\Coord(C)$} is a formal scheme
 over $X$, and is unnaturally isomorphic to $\MG\tm\aff^\infty\tm X$.
\end{proposition}
\begin{proof}
 Choose a normalised coordinate $x$ on $C$, and suppose that
 $a\in X(R)$.  Then any normalised function $y\:C_a\xra{}\haf^1$ has
 the form 
 \[ y(c)=f(x(c))=\sum_{k>0}u_kx(c)^k \]
 for a uniquely determined sequence of coefficients $u_k$.  Moreover,
 $y$ is a coordinate if and only if
 $f\:\haf^1\tm\spec(R)\xra{}\haf^1\tm\spec(R)$ is an isomorphism, if
 and only if there is a power series $g$ with $g(f(t))=t=f(g(t))$.  It
 is well-known that this happens if and only if $u_1$ is invertible.
 Thus, the set of coordinates on $C_a$ bijects naturally with
 $(\MG\tm\aff^\infty)(R)$, and $\Coord(C)\simeq \MG\tm\aff^\infty\tm
 X$ is a formal scheme, as required.
\end{proof}
\begin{remark}
 We will see later that when $E$ is an even periodic ring spectrum and
 $G_E=(\cp^\infty)_E$ we have $\Coord(G_E)=\spec(E_0MP)$.
\end{remark}

\subsection{Weierstrass preparation}
\label{subsec-weierstrass}

\begin{definition}\label{defn-weierstrass}
 A \dfn{Weierstrass series} over a ring $R$ is a formal power series
 $g(x)=\sum_ka_kx^k\in\fps{R}{x}$ such that there exists an integer
 $n$ such that $a_k$ is nilpotent for $k<n$, and $a_n$ is a unit.  The
 integer $n$ is called the \dfn{Weierstrass degree} of $g(x)$.  (It is
 clearly well-defined unless $R=0$).  A \dfn{Weierstrass polynomial}
 over a ring $R$ is a monic polynomial $h(x)=\sum_{k=0}^n b_kx^k$ such
 that $b_k$ is nilpotent for $k<n$.
\end{definition}

The following result is a version of the Weierstrass Preparation
Theorem; see~\cite[Theorem 3]{fr:fg} (for example) for a more
classical version.
\begin{lemma}\label{lem-weierstrass}
 Let $R$ be a ring, and let $g(x)$ be a \idx{Weierstrass series} over
 $R$, of \idx{Weierstrass degree} $n$.  Then there is a unique ring
 map $\al\:\fps{R}{y}\xra{}\fps{R}{x}$ sending $y$ to $g(x)$, and this
 makes $\fps{R}{x}$ into a free module over $\fps{R}{y}$ with basis
 $\{1,x,\ldots,x^{n-1}\}$.
\end{lemma}
\begin{proof}
 We can easily reduce to the case where $a_n=1$.  It is also easy to
 check that there is a unique map $\al$ sending $y$ to $g(x)$, and
 that it sends any series $\sum_j b_jy^j$ to the sum
 $\sum_jb_jg(x)^j$, which is $x$-adically convergent.

 For any $j\ge 0$ and $0\leq k<n$ we define $z_{jk}=g(x)^jx^k$.  Given
 any $m\ge 0$ we can write $m=nj+k$ for some $j\ge 0$ and $0\le k<n$,
 and we put $w_m=z_{jk}$.  For any $R$-module $M$, we define a map 
 \[ \bt_M \: \prod_m M \xra{} \fps{M}{x} \]
 by $\bt_M(b)=\sum_m b_m w_m$.  It is easy to check that this
 sum is again $x$-adically convergent.  The claim in the lemma is
 equivalent to the statement that $\bt_R$ is an isomorphism.

 Write $I=(a_0,\ldots,a_{n-1})$.  This is finitely generated, so the
 same is true of $I^r$ for all $r$, and it follows that
 $I^r\prod_mM=\prod_mI^rM$ and so on.  We also see that
 $w_m=x^m\pmod{I,x^{m+1}}$.  

 Now consider a module $M$ with $IM=0$, so that
 $bw_m=bx^m\pmod{x^{m+1}}$ for $b\in M$.  Given any series
 $c(x)=\sum_mc_mx^m\in\fps{M}{x}$, we see by induction on $m$ that
 there is a unique sequence $(b_j)$ such that
 $\sum_{j<m}b_mw_m=c(x)\pmod{x^m}$ for all $m$.  It follows that
 $\bt_M$ is an isomorphism whenever $IM=0$.  Next, whenever we have a
 short exact sequence $L\mra M\era N$ we have short exact sequences
 $\prod_mL\mra\prod_mM\era\prod_mN$ and
 $\fps{L}{x}\mra\fps{M}{x}\era\fps{N}{x}$, and we can use the
 five-lemma to see that $\bt_M$ is iso if $\bt_L$ and $\bt_N$ are.
 Using this we see by induction that $\bt_{R/I^r}$ is iso for all
 $r$.  On the other hand, when $R$ is large we have $I^r=0$ and so
 $\bt_R$ is an isomorphism.
\end{proof}

\begin{corollary}\label{cor-weierstrass-quotient}
 In the situation of the lemma, the quotient ring $\fps{R}{x}/g(x)$ is
 a free module of rank $n$ over $R$, with basis
 $\{1,\ldots,x^{n-1}\}$. \qed
\end{corollary}

\begin{corollary}\label{cor-weierstrass-factorisation}
 If $g(x)$ is a \idx{Weierstrass series} over a ring $R$ then there is
 a unique factorisation $g(x)=h(x)u(x)$, where $h(x)$ is a
 \idx{Weierstrass polynomial}, and $u(x)$ is invertible.
\end{corollary}
\begin{proof}
 By the previous corollary, we have
 $-x^n=\sum_{j=0}^{n-1}b_jx^j\pmod{g(x)}$ for some unique sequence
 $b_0,\ldots,b_{n-1}\in R$.  Put $h(x)=x^n+\sum_jb_jx^j$, so $h$ is a
 monic polynomial of degree $n$ with $h(x)=0\pmod{g(x)}$, say
 $h(x)=g(x)v(x)$.  Now write $g(x)$ in the form $\sum_ka_kx^k$ and put
 $I=(a_0,\ldots,a_{n-1})$, so $I$ is a nilpotent ideal.  Modulo $I$ we
 find that $g(x)$ is a unit multiple of $x^n$, so
 $x^n=0\pmod{I,g(x)}$.  The uniqueness argument applied mod $I$ now
 tells us that $h(x)=x^n\pmod{I}$, so $h(x)$ is a Weierstrass
 polynomial.  It is also clear that $v(x)$ becomes a unit mod
 $\fps{I}{x}$, but $\fps{I}{x}$ is nilpotent so $v(x)$ is a unit.  We
 can thus take $u(x)=1/v(x)$ to get the required factorisation.
\end{proof}

We now give a more geometric restatement of the above results.
\begin{definition}
 Let $C\xra{q}X$ and $D\xra{r}X$ be formal curves over a formal scheme
 $X$, and let $f\:C\xra{}D$ be a map over $X$.  We then have a curve
 $r^*C=C\tm_XD$ over $D$, with projection map $s\:(c,d)\mapsto d$.  We
 also have a map $f'\:C\xra{}r^*C$ of formal schemes over $D$, given
 by $f'(c)=(c,f(c))$.  We say that $f$ is an \dfn{isogeny} if the map
 $f'$ makes $C$ into a divisor on $r^*C$ over $D$.  This implies in
 particular that $f$ is finite and very flat.
\end{definition}

\begin{lemma}\label{lem-weierstrass-isogeny}
 Let $X$ be an informal scheme, and let $f\:C\xra{}D$ be a map of
 formal curves over $X$.  Let $x$ and $y$ be coordinates on $C$ and
 $D$ respectively, and suppose that $f^*y=g(x)$ for some
 \idx{Weierstrass series} $g(x)$.  Then $f$ is an \idx{isogeny}.
\end{lemma}
\begin{proof}
 Write $R=\OX$, and let $n$ be the Weierstrass degree of $g(x)$.  We
 then have $C=\spf(\fps{R}{x})$ and $D=\spf(\fps{R}{y})$ and
 $r^*C=\spf(\fps{R}{x,y})$.  In this last case we think of $x$ as the
 coordinate on $r^*C$ and $y$ as a parameter on the base.  The map
 $f'$ corresponds to the map $\al\:\fps{R}{x,y}\xra{}\fps{R}{x}$ that
 sends $x$ to $x$ and $y$ to $g(x)$.  We thus need to show that $\al$
 is surjective (making $f'$ a closed inclusion) and that it makes
 $\fps{A}{x}$ into a free module of rank $n$ over $\fps{A}{y}$.  The
 surjectivity is clear, and the freeness follows from
 Lemma~\ref{lem-weierstrass}.
\end{proof}
\begin{example}
 One can check that the evident map
 $\cp^\infty\xra{}\mathbb{H}P^\infty$ gives an isogeny
 $(\cp^\infty)_E\xra{}(\mathbb{H}P^\infty)_E$ of formal curves.
\end{example}

\begin{definition}
 Let $X$ be an informal scheme, and $C$ a formal curve over $X$.  We
 then let $\CM_{C/X}$ \index{MCX@$\CM_{C/X}$} be the ring obtained
 from $\OC$ by inverting all coordinates on $C$.  We refer to this as
 the ring of \emph{meromorphic functions} \index{meromorphic} on $C$.
\end{definition}

\begin{lemma}\label{lem-invert-coords}
 Let $X$ be an informal scheme, and $C$ a formal curve over $X$, and
 $x$ a coordinate on $C$.  Then $\CM_{C/X}=\OC[1/x]$, which is the
 ring of series $\sum_{k\in\Zh} a_kx^k$ such that $a_k\in\OX$ and
 $a_k=0$ for $k\ll 0$.\index{MCX@$\CM_{C/X}$}
\end{lemma}
\begin{proof}
 Let $y$ be another coordinate on $C$; it will suffice to check that
 $y$ becomes invertible in $\OC[1/x]$.  As $x$ and $y$ are both
 coordinates, we find that $y=\sum_{k\ge 0}a_kx^k$ for some series
 such that $a_0$ is nilpotent and $a_1$ is a unit.  In other words, we
 have $y=b + x c(x)$, where $b$ is nilpotent and $c(0)$ is invertible
 in $\OX$, so $c(x)$ is invertible in $\OC$.  It is thus clear that
 $y-b$ has inverse $x^{-1}c(x)^{-1}$ in $\OC[1/x]$.  The sum of a unit
 and a nilpotent element is always invertible, so $y$ is a unit as
 required. 
\end{proof}
\begin{remark}\label{rem-poles}
 The elements of $\CM_{C/X}$ should be thought of as Laurent
 expansions of functions whose poles are all very close to the origin,
 the expansion being valid outside a small disc containing all the
 poles.
\end{remark}

\begin{lemma}\label{lem-mero-units}
 Let $x$ be a coordinate on $C$, and let $f(x)=\sum_{k\in\Zh}a_kx^k$
 be an element of $\CM_{C/X}$.\index{MCX@$\CM_{C/X}$} Then $f$ is
 invertible in $\CM_{C/X}$ if and only if $X$ can be written as a
 coproduct $X=\coprod_{k\in\Zh}X_k$, where $X_k=\emptyset$ for almost
 all $k$, such that $a_j$ is nilpotent on $X_k$ for $j<k$, and $a_k$
 is invertible on $X_k$.
\end{lemma}
\begin{proof}
 First suppose that $f(x)$ is invertible, say $f(x)g(x)=1$ with
 $g(x)=\sum_{j\in\Zh}b_jx^j$.  Write $I_j=(a_jb_{-j})$ and
 $J_j=\sum_{k\neq j}I_j$.  Because $f(x)g(x)=1$ it is clear that
 $\sum_jI_j=\OX$ and thus $J_i+J_j=\OX$ when $i\neq j$.  There exists
 $K$ such that $a_{-j}=b_{-j}=0$ when $j>K$.  It follows that $I_j=0$
 and $J_j=\OX$ when $|j|>K$.  Next, let $\pri$ be a prime ideal in
 $\OX$.  As $\OX/\pri$ is an integral domain, it is clear that modulo
 $\pri$ we must have $f(x)=a_kx^k+\ldots$ and
 $g(x)=b_{-k}x^{-k}+\ldots$ for some $k$.  This implies that
 $a_ib_j\in\pri$ whenever $i+j<0$.  As the intersection of all prime
 ideals is the set of nilpotents, the elements $a_ib_j$ must be
 nilpotent when $i+j<0$.  If $i\neq j$ then either $i-j$ or $j-i$ is
 negative, so $a_ib_{-i}a_jb_{-j}$ is nilpotent.  It follows that
 $I_iI_j$ is nilpotent when $i\neq j$, and thus that $\bigcap_jJ_j$ is
 nilpotent.  It follows from the results of
 Section~\ref{subsec-nilpotents} that there are unique ideals $J'_j$
 such that $J_j\leq J'_j\leq \sqrt{J_j}$ and $\OX=\prod_j\OX/J'_j$.
 We take $X_j=\spec(\OX/J'_j)$; one can check that this has the
 claimed properties.

 Conversely, suppose that $X$ has a decomposition of the type
 discussed.  We reduce easily to the case where $X=X_k$ for some $k$.
 After replacing $f$ by $x^{-k}f$, we may assume that $k=0$.  This
 means that $f(x)=\sum_{j\in\Zh}a_jx^j$, where $a_0$ is invertible and
 $a_j$ is nilpotent for $j<0$ and $a_j=0$ for $j\ll 0$. The
 invertibility or otherwise of $f$ is unaffected if we subtract off a
 nilpotent term, so we may assume that $a_j=0$ for $j<0$.  The
 resulting series is invertible in $\OC$ and thus certainly in
 $\CM_{C/X}$. 
\end{proof}

\begin{definition}\label{defn-mero-deg}
 Let $x$ be a coordinate on $C$, and let $f$ be an invertible element
 of $\CM_{C/X}$,\index{MCX@$\CM_{C/X}$} so we have a decomposition
 $X=\coprod_kX_k$ as above.  If $X=X_k$ then we say that $f$ has
 constant \idx{degree} $k$.  More generally, we let $\deg(f)$ be the
 map from $X$ to the constant scheme $\un{\Zh}$ that takes the value
 $k$ on $X_k$.  One can check that these definitions are independent
 of the choice of coordinate.
\end{definition}

\begin{lemma}\label{lem-mero-factor}
 Let $x$ be a coordinate on $C$, and let $f$ be an invertible element
 of $\CM_{C/X}$,\index{MCX@$\CM_{C/X}$} with constant degree $k$.
 Then there is a unique factorisation $f(x)=x^ku(x)g(x)$, where
 $u(x)\in\OCt$, and $g(x)=\sum_{j\ge 0} b_j x^{-j}$ where $b_0=1$
 and $b_j$ is nilpotent for $j>0$ and $b_j=0$ for $j\gg 0$.
\end{lemma}
\begin{proof}
 Clearly we have $h(x)=x^Nf(x)\in\OC$ for some $N>0$.  We see from
 Lemma~\ref{lem-mero-units} that $h(x)$ is a Weierstrass power series
 of Weierstrass degree $N+k$.  It follows from
 Corollary~\ref{cor-weierstrass-factorisation} that $h(x)$ has a
 unique factorisation of the form $h(x)=k(x)u(x)$, where $k(x)$ is a
 Weierstrass polynomial of degree $N+k$, and $u(x)\in\OCt$.  We
 write $g(x)=x^{-N-k}h(x)$; this clearly gives a factorisation of the
 required type, and one can check that it is unique.
\end{proof}

\begin{proposition}\label{prop-mero-ses}
 Let $C$ be a formal curve over a formal scheme $X$.  For any ring
 $R$, we define 
 \[ \Mer(C,\MG)(R) = 
    \{(u,f)\st u\in X(R) \;,\; f\in\CM_{C_u/\spec(R)}^\tm \}. 
 \]
 \index{MerCGm@$\Mer(C,\MG)$} Then $\Mer(C,\MG)$ is a formal scheme
 over $X$, and there is a short exact sequence of formal groups
 \[ \Map(C,\MG) \mra \Mer(C,\MG) \era \Div(C), \]
 which admits a non-canonical splitting.
\end{proposition}
\begin{proof}
 As $\Map(C,\MG)(R)=\{(u,f)\st u\in X(R)\;,\;f\in\OO_{C_u}^\tm\}$,
 there is an obvious inclusion $\Map(C,\MG)\xra{}\Mer(C,\MG)$ of
 group-valued functors.  Next, let $Y(R)$ be the set of series
 $g(x)=\sum_{j\ge 0}b_jx^{-j}$ such that $b_0=1$ and $b_j$ is
 nilpotent for $j>0$ and $b_j=0$ for $j\gg 0$.  Then
 $Y=\colim_k\prod_{0<j<k}\haf^1$ is a formal scheme, and
 Lemma~\ref{lem-mero-factor} gives an isomorphism
 $\Mer(C,\MG)\simeq\Map(C,\MG)\tm\un{\Zh}\tm Y$.  This shows that
 $\Mer(C,\MG)$ is a formal scheme.  We next define a map
 $\div\:\Mer(C,\MG)\xra{}\Div(C)$.  Suppose that $f\in\OO_{C_u}$ is
 such that $\OO_{C_u}/f$ is a free module of rank $n$ over $R$.  Then
 $D=\spf(\OO_{C_u}/f)\in\Div_n(G)(R)$and we define $\div(f)=D$.
 Given another such function $g\in\OO_{C_u}$, we define
 $\div(f/g)=\div(f)-\div(g)$.  This is well-defined, because if
 $f/g=f'/g'$ then $fg'=f'g$ (because series of this form are never
 zero-divisors) and thus $\div(f)+\div(g')=\div(f')+\div(g)$ and so
 $\div(f)-\div(g)=\div(f')-\div(g')$.  It is easy to see that we get a
 well-defined homomorphism $\div\:\Mer(C,\MG)\xra{}\Div(C)$, which is
 zero on $\Map(C,\MG)$.  Conversely, suppose that $\div(f/g)=0$, so
 that $\div(f)=\div(g)$.  Then $f$ and $g$ are non-zero-divisors and
 they generate the same ideal in $\OO_{C_u}$, so they are unit
 multiples of each other and thus $f/g\in\Map(C,\MG)(R)$.  Thus
 $\ker(\div)=\Map(C,\MG)$.  

 Now let $j\:C\xra{}\Div(C)$ be the evident inclusion.  Given a point
 $a\in C(R)$, we also define $\sg(a)=x-x(a)=x(1-x(a)/x)\in\Mer(C)(R)$.
 This gives a map $\sg\:C\xra{}\Mer(C,\MG)$, and it is clear that
 $\div\circ\sg=j$.  As $\Div(C)$ is the free Abelian formal group
 generated by $C$, we see that there is a unique homomorphism
 $\tau\:\Div(C)\xra{}\Mer(C,\MG)$ with $\tau\circ j=\sg$.  We thus
 have $\div\circ\tau\circ j=j$ and thus $\div\circ\tau=1$.  It follows
 that the sequence $\Map(C,\MG)\mra\Mer(C,\MG)\era\Div(C)$ is a split
 exact sequence.  The splitting depends on a choice of coordinate, but
 the other maps are canonical.
\end{proof}

\subsection{Formal differentials}
\label{subsec-formal-diffs}

We next generalise Definition~\ref{defn-Omega} to a certain (rather
small) class of formal schemes.  
\begin{definition}
 We say that a formal scheme $W$ over $X$ is \dfn{formally smooth} of
 dimension $n$ over $X$ if it is isomorphic in $\FX_X$ to
 $\haf^n\tm X$.  In particular, $W$ is formally smooth of dimension
 one if and only if it is a formal curve.
\end{definition}

\begin{definition}\label{defn-Omega-formal}
 Let $W$ be formally smooth of dimension $n$ over $X$; we shall define
 a vector bundle $\Om^1_{W/X}$\index{O1@$\Om^1$} of rank $n$ over $W$.
 By Corollary~\ref{cor-fibrewise-sheaf}, it suffices to do this in a
 sufficiently natural way whenever $X$ is an informal scheme.  In that
 case, we let $J$ be the kernel of the multiplication map
 $\OO_{W\tm_XW}=\OW\hot_\OX\OW\xra{}\OW$, so that $V(J)$ is the
 diagonal subscheme in $W\tm_XW$.  We then write $\Om^1_{W/X}=J/J^2$,
 which is a module over $\OO_{W\tm_XW}/J=\OW$.  If $f\in\OW$ then we
 write $d(f)=f\ot 1-1\ot f+J^2\in\Om^1_{W/X}$, and note that
 $d(fg)=fd(g)+gd(f)$ as usual.  As $W$ is formally smooth, we can
 choose $x_1,\ldots,x_n\in\OW$ giving an isomorphism $W\simeq\haf^n\tm
 X$ and thus $\OW\simeq\fps{\OX}{x_1,\ldots,x_n}$.  One checks that
 $\Om^1_{W/X}$ is freely generated over $\OW$ by
 $\{d(x_1),\ldots,d(x_n)\}$.  Thus, $\Om^1_{W/X}$ can be regarded as a
 vector bundle of rank $n$ over $W$, as required.  We write
 $\Om^k_{W/X}$ for the $k$'th exterior power of $\Om^1_{W/X}$.
\end{definition}

Any map $f\:V\xra{}W$ of formally smooth schemes over $X$ gives rise
to a map $f^*\:f^*\Om^1_{W/X}\xra{}\Om^1_{V/X}$ of vector bundles over
$V$.  If we have coordinates $x_i$ on $W$ and $y_j$ on $V$ then
$x_i\circ f=u_i(y_1,\ldots,y_d)$ for certain power series $u_i$ over
$\OX$, and we have 
$f^*(dx_i)=\sum_j(\partial u_i/\partial y_j) dy_j$.  Thus, $f^*$ is a
coordinate-free way of encoding the partial derivatives of the series
$u_i$.  

If $W$ is formally smooth over $X$ and $g\:Y\xra{}X$ then
$g^*W=Y\tm_XW$ is easily seen to be formally smooth over $Y$, with
$\Om^1_{g^*W/Y}=h^*\Om^1_{W/X}$, where $h\:g^*W\xra{}W$ is the evident
projection map.

\begin{definition}
 If $R$ is an $\Fp$-algebra, then we have a ring map $\phi_R$ from $R$
 to itself defined by $\phi_R(a)=a^p$.  We call this the
 \dfn{algebraic Frobenius map}.  Now let $X$ be a functor over
 $\spec(\Fp)$.  If $R$ is an $\Fp$-algebra, we define
 $(F_X)_R=X(\phi_R)\:X(R)\xra{}X(R)$.  If $R$ is not an $\Fp$-algebra
 then $\spec(\Fp)(R)=\emptyset$ and thus $X(R)=\emptyset$ and we
 define $(F_X)_R=1\:X(R)\xra{}X(R)$.\index{FX@$F_X$} This gives a map
 $F_X\:X\xra{}X$, which we call the \dfn{geometric Frobenius map}.
\end{definition}

\begin{remark}
 If $h\:X\xra{}Y$ is a map of functors over $\spec(\Fp)$ then one can
 check that $ F_Y\circ h=h\circ F_X$.  If $X$ is a scheme then $F_X$
 is characterised by the fact that $g( F_X(a))=g(a)^p$ for all rings
 $R$, points $a\in X(R)$, and functions $g\in\OX$.  If $X=\spec(A)$
 then $F_X=\spec(\phi_A)$.
\end{remark}

\begin{definition}
 Let $X$ be a functor over $\spec(\Fp)$, and let $W$ be functor over
 $X$, with given map $q\:W\xra{}X$.  We then have a functor $F_X^*W$
 over $X$ defined by 
 \[ (F_X^*W)(R) = \{(a,b)\in W(R)\tm X(R)\st q(a)=F_X(b)\}. \]
 We define a map $F_{W/X}\:W\xra{}F_X^*W$ \index{FWX@$F_{W/X}$} by
 $F_{W/X}(a)=(F_W(a),q(a))$.  Note that if $W$ is
 \idx{formally smooth} over $X$ then the same is true of $F_X^*W$.
 Moreover, if we have coordinates $x_i$ on $W$ and we use the obvious
 resulting coordinates $y_i$ on $F_X^*W$ then we have
 $y_i(F_{W/X}(a))=x_i(a)^p$.
\end{definition}

\begin{proposition}\label{prop-zero-differential}
 Let $f\:V\xra{}W$ be a map of formally smooth schemes over $X$, and
 suppose that $f^*=0\:\Om^1_{W/X}\xra{}\Om^1_{V/X}$.\index{O1@$\Om^1$}
 \begin{itemize}
  \item[(a)] If $X$ lies over $\spec(\rat)$ then there is a unique map
   $g\:X\xra{}W$ of schemes over $X$ such that $f$ is the composite
   $V\xra{}X\xra{g}W$.  In other words, $f$ is constant on the fibres
   of $V$.
  \item[(b)] If $X$ lies over $\spec(\Fp)$ for some prime $p$, then
   there is a unique map $f'\:F_X^*W\xra{}V$ such that
   $f=f'\circ F_{W/X}$. \index{FWX@$F_{W/X}$}
 \end{itemize}
\end{proposition}
\begin{proof}
 Choose coordinates $x_i$ on $W$ and $y_j$ on $V$, so
 $x_i\circ f=u_i(y_1,\ldots,y_d)$ for certain power series $u_i$ over
 $\OX$. We have $0=f^*(dx_i)=\sum_j(\partial u_i/\partial y_j)dy_j$,
 so $\partial u_i/\partial y_j=0$ for all $i$ and $j$.  In the
 rational case we conclude that the series $u_i$ are constant, and in
 the mod $p$ case we conclude that
 $u_i(y_1,\ldots,y_d)=v_i(y_1^p,\ldots,y_d^p)$ for some unique series
 $v_i$.  The conclusion follows easily.
\end{proof}

\subsection{Residues}
\label{subsec-residues}

We now describe an algebraic theory of residues, which is essentially
the same as that discussed in~\cite{he:acc} and presumably identical
to the unpublished definition by Cartier mentioned in~\cite{qu:fgl}.
\begin{definition}
 If $f(x)=\sum_{k\in\Zh}a_kx^k\in\fps{R}{x}[1/x]$, we define
 $\rho(f)=a_{-1}$.  
\end{definition}

\begin{remark}\label{rem-residues}
 Recall from Remark~\ref{rem-poles} that the elements of
 $\fps{R}{x}[1/x]$ should be compared with meromorphic functions on a
 neighbourhood of zero in $\cplx$ of moderate size, whose poles are
 concentrated very near the origin.  The expansion in terms of $x$
 should be thought of as a Laurent expansion that is valid outside a
 tiny disc containing all the poles.  Thus, the coefficient of $1/x$
 is the sum of the \idx{residue}s at all the poles, and not just the
 pole at the origin.  To justify this, note that if $a$ is nilpotent
 (say $a^N=0$) we have $1/(x-a)=\sum_{k=0}^{N-1}a^k/x^{k+1}$ so
 $\rho(1/(x-a))=1$.
\end{remark}

\begin{proposition}
 For any $f\in\fps{R}{x}[1/x]$ we have $\rho(f')=0$.  If $f$ is
 invertible we have $\rho(f'/f)=\deg(f)$, where $\deg(f)$ is as in
 Definition~\ref{defn-mero-deg}.  Moreover, we have $\rho(f^n.f')=0$
 for $n\neq -1$.
\end{proposition}
\begin{proof}
 It is immediate from the definitions that $\rho(f')=0$.  Now let $f$
 be invertible; we may assume that it has constant degree $d$ say.
 Lemma~\ref{lem-mero-factor} gives a factorisation $f(x)=x^du(x)g(x)$,
 where $u(x)\in\fps{R}{x}^\tm$, and $g(x)=\sum_{j\ge 0} b_j x^{-j}$
 where $b_0=1$ and $b_j$ is nilpotent for $j>0$ and $b_j=0$ for
 $j\gg 0$.  We then have $f'/f=d/x+u'/u+g'/g$.  It is clear that
 $u'/u\in\fps{R}{x}$ so $\rho(u'/u)=0$.  Similarly, we find that $g'$
 only involves powers $x^k$ with $k<-1$.  Moreover, if $h(x)=1-g(x)$
 then $h$ is a polynomial in $1/x$ and is nilpotent, and
 $1/g=\sum_{k=0}^Nh^k$ for some $N$ so $1/g$ is a polynomial in
 $1/x$.   It follows that $g'/g$ only involves powers $x^k$ with
 $k<-1$, so $\rho(g'/g)=0$.  Thus $\rho(f'/f)=d$ as claimed.

 Finally, suppose that $n\neq -1$.  Note that
 $(n+1)\rho(f^n.f')=\rho((f^{n+1})')=0$.  If $R$ is torsion-free we
 conclude that $\rho(f^n.f')=0$.  If $R$ is not torsion-free, we
 recall that $f(x)$ has the form $\sum_{i=m}^{\infty}a_ix^i$ for some
 $m$, where $a_i$ is nilpotent for $i<d$ and $a_d$ is invertible.
 Thus there is some $N>0$ such that $a_i^N=0$ for all $i<d$.  Define
 $R'=\Zh[b_i\st i\ge m][1/b_d]/(b_i^N\st m\leq i<d)$ and
 $g(x)=\sum_ib_ix^i\in\fps{R'}{x}[1/x]^\tm$.  It is clear that $R'$ is
 torsion-free and thus that $\rho(g^n.g')=0$.  There is an evident map
 $R'\xra{}R$ carying $g$ to $f$, so we deduce that $\rho(f^n.f')=0$ as
 claimed.
\end{proof}

\begin{corollary}
 If $g(x)\in\fps{R}{x}$ is a Weierstrass series of degree $d>0$ and
 $f(x)\in\fps{R}{x}[1/x]$ then $\rho(f(g(x))g'(x))=d\rho(f(x))$.
\end{corollary}
\begin{proof}
 Suppose that $f(x)=\sum_{k\ge m}a_kx^k$.  We first observe that the
 claim makes sense: as $g$ is a Weierstrass series of degree $d>0$ we
 know that $g(0)$ is nilpotent, so $g^N\in\fps{R}{x}x$ for some $N$,
 so $g^{Nk}\in\fps{R}{x}x^k$ for $k\ge 0$.  Moreover,
 Lemma~\ref{lem-mero-units} implies that $g$ is invertible in
 $\fps{R}{x}[1/x]$.  Thus, the terms in the sum
 $f(g(x))=\sum_{k\ge m}a_k g(x)^k$ are all defined, and the sum is
 convergent.  We thus have 
 \[ \rho(f(g(x)))=\sum_k a_k\rho(g^k.g')=d\, a_{-1}=d\rho(f) \]
 as required.
\end{proof}

\begin{definition}
 Let $C$ be a formal curve over an affine scheme $X$.  We write
 $\CM\Om^1_{C/X}$ for $\CM_{C/X}\ot_{\OC}\Om^1_{C/X}$, which is a free
 module of rank one over $\CM_{C/X}$.  It is easy to check that there
 is a unique map $d\:\CM_{C/X}\xra{}\CM\Om^1_{C/X}$ extending the
 usual map $d\:\OC\xra{}\Om^1_{C/X}$ and satisfying
 $d(fg)=f\,d(g)+g\,d(f)$. 
\end{definition}

\begin{corollary}\label{cor-invariant-residue}
 Let $C$ be a formal curve over an affine scheme $X$.  Then there is a
 natural \idx{residue} map $\res=\res_{C/X}\:\CM\Om^1_{C/X}\xra{}\OX$
 such that
 \begin{itemize}
  \item[(a)] $\res(df)=0$ for all $f\in\CM_{C/X}$.
  \item[(b)] $\res((df)/f)=\deg(f)$ for all $f\in\CM^\tm_{C/X}$.
  \item[(c)] If $q\:C\xra{}C'$ is an isogeny then
   $\res(q^*\al)=\deg(q)\res(\al)$ for all $\al$.
 \end{itemize}
\end{corollary}
\begin{proof}
 Choose a coordinate $x$ on $C$, so that any
 $\al\in\CM_{C/X}\ot_{\OC}\Om^1_{C/X}$ has a unique expression
 $\al=f(x)dx$ for some $f\in\fps{\OX}{t}[1/t]$.  Define
 $\res(\al)=\rho(f)$.  If $y$ is a different coordinate then $x=g(y)$
 for some Weirstrass series $g$ of degree $1$ and $dx=g'(y)dy$ so
 $\al=f(g(y))g'(y)dy$ and we know that $\rho(f(g(y))g'(y))=\rho(f)$ so
 our definition is independent of the choice of the coordinate.  The
 rest of the corollary is just a translation of the properties of
 $\rho$.  
\end{proof}
See Remark~\ref{rem-Lin-residue} for a topological application of
this. 

\section{Formal groups}
\label{sec-formal-groups}

A \dfn{formal group} over a formal scheme $X$ is just a group object
in the category $\FX_X$.  In this section, we will study formal groups
in general.  In the next, we will specialise to the case of
commutative formal groups $G$ over $X$ with the property that the
underlying scheme is a formal curve; we shall call these
\emph{ordinary formal groups}.  For technical reasons, it is
convenient to compare our formal groups with group objects in suitable
categories of coalgebras.  To combine these cases, we start with a
discussion of Abelian group objects in an arbitrary category with
finite products.  We then discuss the existence of free Abelian formal
groups, or of schemes of homomorphisms between formal groups.  As a
special case, we discuss the Cartier duality functor
$G\mapsto\Hom(G,\MG)$.  Finally, we define torsors over a commutative
formal group, and show that they form a strict Picard category.

\subsection{Group objects in general categories}
\label{subsec-group-objects}

Let $\CD$ be a category with finite products (including an empty
product, in other words a terminal object).  There is an evident
notion of an Abelian group object in $\CD$; we write
$\Ab\CD$\index{AbD@$\Ab\CD$} for the category of such objects.  We
also consider the category $\Mon\CD$\index{MonD@$\Mon\CD$} of
Abelian monoids in $\CD$.  A \dfn{basepoint} for an object $U$ of
$\CD$ is a map from the terminal object to $U$.  We write
$\Based\CD$\index{BasedD@$\Based\CD$} for the category of objects of
$\CD$ equipped with a specified basepoint.  There are evident
forgetful functors
\[ \Ab\CD \xra{} \Mon\CD \xra{} \Based\CD \xra{} \CD.\]

If $U\in\CD$ and $G\in\Ab\CD$ then the set $\CD(U,G)$ has a natural
Abelian group structure.  In fact, to give such a group structure is
equivalent to giving maps $1\xra{0}G\xla{\sg}G\tm G$ making it an
Abelian group object, as one sees easily from Yoneda's lemma.  

Let $\{G_i\}$ be a diagram in $\Ab\CD$, and suppose that the
underlying diagram in $\CD$ has a limit $G$.  Then
$\CD(U,G)=\invlim_i\CD(U,G_i)$ has a natural Abelian group structure.
It follows that there is a unique way to make $G$ into an Abelian
group object such that the maps $G\xra{}G_i$ become homomorphisms, and
with this structure $G$ is also the limit in $\Ab\CD$.  In other
words, the forgetful functor $\Ab\CD\xra{}\CD$ creates limits.
Similarly, we see that all the functors
$\Ab\CD\xra{}\Mon\CD\xra{}\Based\CD\xra{}\CD$ and their composites
create limits.

Suppose that $G$, $H$ and $K$ are Abelian group objects in $\CD$ and
that $f\:G\xra{}K$ and $g\:H\xra{}K$ are homomorphisms.  One checks
that the composite $G\tm H\xra{f\tm g}K\tm K\xra{\sg}K$ is also a
homomorphism, and that it is the unique homomorphism whose composites
with the inclusions $G\xra{}G\tm H$ and $H\xra{}G\tm H$ are $f$ and
$g$.  This means that $G\tm H$ is the coproduct of $G$ and $H$ in
$\Ab\CD$, as well as being their product.

We next investigate another type of colimit in $\Ab\CD$.
\begin{definition}\label{defn-refl-fork}
 A \dfn{reflexive fork} in any category $\CD$ is a pair of objects
 $U,V$, together with maps $d_0,d_1\:U\xra{}V$ and $s\:V\xra{}U$ such
 that $d_0s=1=d_1s$.  The \emph{coequaliser} of such a fork means the
 coequaliser of the maps $d_0$ and $d_1$.
\end{definition}

\begin{proposition}\label{prop-refl-fork}
 Let $\CD$ be a category with finite products.  Let
 \[ V\arw{e,t}{s}U\arw{e,tb,=>}{d_0}{d_1}V \]
 be a \idx{reflexive fork} in $\Mon\CD$\index{Mond@$\Mon\CD$}, and let
 $U\arw{e,tb,=>}{d_0}{d_1}V\arw{e,t}{e}W$ be a strong coequaliser in
 $\CD$.  Then there is a monoid structure on $W$ such that $e$ is a
 homomorphism, and this makes the above diagram into a coequaliser in
 $\Mon\CD$.
\end{proposition}
\begin{proof}
 Let $\sg_U\:U\tm U\xra{}U$ and $\sg_V\:V\tm V\xra{}V$ be the addition
 maps.  We have a commutative diagram as follows:
 \begin{diag}
  \node{V\tm U}
  \arw{s,lr,=>}{1\tm d_0}{1\tm d_1}
  \arw{e,t}{s\tm 1}
  \node{U\tm U}
  \arw{s,lr,=>}{d_0\tm d_0}{d_1\tm d_1}
  \arw{e,t}{\sg_U}
  \node{U}
  \arw{s,lr,=>}{d_0}{d_1} \\
  \node{V\tm V}
  \arw{e,=}
  \node{V\tm V}
  \arw{e,b}{\sg_V}
  \node{V}
 \end{diag}
 The right hand square commutes because $d_0$ and $d_1$ are
 homomorphisms, and the left hand one because $d_0s=d_1s=1$.  Using
 this, we see that $e\sg_V(1\tm d_0)=e\sg_V(1\tm d_1)$, and a similar
 proof shows that $e\sg_V(d_0\tm 1)=e\sg_V(d_1\tm 1)$.  In terms of
 elements, this just says that $e(d_0(u)+v)=e(d_1(u)+v)$.  As our
 coequaliser diagram was assumed to be strong, we see that the diagram
 \[ V\tm U \arw{e,tb,=>}{1\tm d_0}{1\tm d_1} V\tm V 
     \arw{e,t}{1\tm e} V\tm W
 \]
 is a coequaliser.  This implies that there is a unique map
 $\tau\:V\tm W\xra{}W$ with $\tau(1\tm e)=e\sg_V\:V\tm V\xra{}W$.
 Now consider the diagram 
 \begin{diag}
  \node{U\tm V}
  \arw{e,tb,=>}{d_0\tm 1}{d_1\tm 1}
  \arw{s,l}{1\tm e}
  \node{V\tm V}
  \arw{s,r}{1\tm e}
  \arw{e,t}{\sg_V}
  \node{V}
  \arw{s,r}{e} \\
  \node{U\tm W} 
  \arw{e,tb,=>}{d_0\tm 1}{d_1\tm 1}
  \node{V\tm W}
  \arw{e,b}{\tau}
  \node{W.}
 \end{diag}
 We have already seen that $e\sg_V(d_0\tm 1)=e\sg_V(d_1\tm 1)$, and it
 follows that $\tau(d_0\tm 1)(1\tm e)=\tau(d_1\tm 1)(1\tm e)$.  As the
 relevant coequaliser is preserved by the functor $U\tm(-)$, we see
 that $1_U\tm e$ is an epimorphism, so we can conclude that
 $\tau(d_0\tm 1)=\tau(d_1\tm 1)$.  As the functor $(-)\tm W$ preserves
 our coequaliser, this gives us a unique map $\sg_W\:W\tm W\xra{}W$ such
 that $\sg_W(e\tm 1)=\tau\:V\tm W\xra{}W$.  One checks that this makes
 $W$ into an Abelian group object, and that $e$ is a homomorphism.
 One can also check that this makes $W$ into a coequaliser in
 $\Ab\CD$.
\end{proof}

\begin{remark}\label{rem-refl-fork}
 The same result holds, with essentially the same proof, with
 $\Mon\CD$\index{MonD@$\Mon\CD$} replaced by
 $\Ab\CD$\index{AbD@$\Ab\CD$} or
 $\Based\CD$\index{BasedD@$\Based\CD$}.  The same methods also show
 that a reflexive fork in the category of $R$-algebras (for any ring
 $R$) has the same coequaliser when computed in the category of
 $R$-algebras, or of $R$-modules, or of sets.
\end{remark}

We next try to construct free Abelian groups or monoids on objects of
$\CD$ or $\Based\CD$.  If $U\in\CD$ and $V\in\Based\CD$, we ``define''
objects $M^+(U),N^+(V)\in\Mon\CD$ and $M(U),N(V)\in\Ab(\CD)$ by the
equations
\begin{align*}
 \Mon\CD(M^+(U),M) &= \CD(U,M)          \index{MUa@$M^+(U)$}    \\ 
 \Mon\CD(N^+(V),M) &= \Based\CD(V,M)    \index{NVa@$N^+(V)$}    \\
 \Ab\CD(M(U),G)    &= \CD(U,G)          \index{MU@$M(U)$}       \\
 \Ab\CD(N(V),G)    &= \Based\CD(V,G).   \index{NV@$N(V)$}
\end{align*}
More precisely, if there is an object $H\in\Ab\CD$ with a natural
isomorphism 
\[ \Ab\CD(H,G)=\Based\CD(V,G) \]
for all $G\in\Ab\CD$, then $H$ is unique up to canonical isomorphism,
and we write $N(V)$ for $H$.  Similar remarks apply to the other three
cases.  Given a monoid object $M$, we also ``define'' its group
completion $G(M)\in\Ab\CD$ by the equation
$\Ab\CD(G(M),H)=\Mon\CD(M,H)$.

There are fairly obvious ways to try to construct free group and
monoid objects, using a mixture of products and colimits.  However,
there are two technical points to address.  Firstly, it turns out that
we need our colimits to be strong colimits in the sense of
Definition~\ref{defn-strong-colim}.  Secondly, in some places we can
arrange to use reflexive coequalisers, which is technically
convenient.
  
\begin{proposition}\label{prop-exists-MpU}
 Let $U$ be an object of $\CD$.  For each $k\ge 0$, the symmetric
 group $\Sg_k$ acts on $U^k$.  Suppose that the quotient $U^k/\Sg_k$
 exists as a strong colimit and also that
 $L=\coprod_{k\ge 0}U^k/\Sg_k$ exists as a strong coproduct.  Then
 $L=M^+(U)$.\index{MUa@$M^+(U)$}
\end{proposition}
\begin{proof}
 Let $\CI$ be the category with object set $\nat$, and with morphisms
 \[ \CI(j,k) = \begin{cases}
      \emptyset & \text{ if } j\neq k \\
      \Sg_k     & \text{ if } j=k.
    \end{cases}
 \]
 It is easy to see that there is a functor $k\mapsto U^k$ from $\CI$
 to $\CD$, and that $L$ is a strong colimit of this functor.  It
 follows that $L\tm U^m$ is the colimit of the functor
 $k\mapsto U^k\tm U^m$, and thus (using the ``Fubini theorem'' for
 colimits) that 
 \[ L\tm L=\colim_{(k,m)\in\CI\tm\CI} U^k\tm U^m. \]
 Similarly, $L\tm L\tm L$ is the colimit of the functor
 $(k,m,n)\mapsto U^k\tm U^m\tm U^n$ from $\CI\tm\CI\tm\CI$ to $\CD$.

 Let $j_k\:U^k\xra{}L$ be the evident map.  We then have maps
 $U^k\tm U^m\simeq U^{k+m}\xra{j_{k+m}}L$, and these fit together to
 give a map $\sg\:L\tm L\xra{}L$.  We also have a zero map
 $0=j_0\:1=U^0\xra{}L$.  We claim that this makes $L$ into a
 commutative monoid object in $\CD$.  To check associativity, for
 example, we need to show that
 $\sg\circ(\sg\tm 1)=\sg\circ(1\tm\sg)\:L^3\xra{}L$.  By the above
 colimit description of $L^3$, it is enough to check this after
 composing with the map $j_k\tm j_m\tm j_n\:U^{k+m+n}\xra{}L^3$, and
 it is easy to check that both the resulting composites are just
 $j_{k+m+n}$.  We leave the rest to the reader.

 Now suppose we have a monoid $M\in\Mon\CD$ and a map $f\:U\xra{}M$ in
 $\CD$.  We then have maps $f_k=(U^k\xra{f^k}M^k\xra{\sg}M)$, which
 are easily seen to be invariant under the action of $\Sg_k$, so we
 get an induced map $f'\:L\xra{}M$ in $\CD$.  We claim that this is a
 homomorphism.  It is clear that $f'\circ 0=0$, so we need only check
 that $f\circ \sg=\sg\circ(f\tm f)\:L^2\xra{}M$.  Again, we need only
 check this after composing with the map
 $j_k\tm j_m\:U^{k+m}\xra{}L^2$, and it then becomes easy.  We also
 claim that $f'$ is the unique homomorphism $g\:L\xra{}M$ such that
 $g\circ j_1=f$.  Indeed, we have $j_k=(U^k\xra{j_1^k}L^k\xra{\sg}L)$,
 so for any such $g$ we have
 $g\circ j_k=(U^k\xra{f^k}M^k\xra{\sg}M)=f'\circ j_k$.  By our colimit
 description of $L$, we see that $g=f'$ as claimed.

 This shows that monoid maps $g\:L\xra{}M$ biject naturally with maps
 $f\:U\xra{}M$, by the correspondence $g\mapsto g\circ j_1$.  This
 means that $L=M^+(U)$ as claimed.
\end{proof}

\begin{proposition}\label{prop-exists-NpU}
 Let $V$ be an object of $\Based\CD$, and suppose that $V_k=V^k/\Sg_k$
 exists as a strong colimit for all $k\ge 0$.  The basepoint of $U$
 then induces maps $V_k\xra{}V_{k+1}$.  Suppose also that the sequence
 of $V_k$'s has a strong colimit $L$.  Then
 $L=N^+(V)$.\index{NVa@$N^+(V)$}
\end{proposition}
\begin{proof}
 This is essentially the same as the proof of
 Proposition~\ref{prop-exists-MpU}, and is left to the reader.
\end{proof}

We next try to construct group completions of monoid objects.  We
digress briefly to introduce some convenient language.  Let $M$ be a
monoid object, so that $\CD(U,M)$ is naturally a monoid for all $U$.
We thus have a map $f_U\:\CD(U,M^3)=\CD(U,M)^3\xra{}\CD(U,M^2)$
defined by $f(a,b,c)=(c+2a,3b+c)$ (for example).  This is natural in
$U$, so Yoneda's lemma gives us a map $f\:M^3\xra{}M^2$.  {}From now on,
we will allow ourselves to abbreviate this definition by saying ``let
$f\:M^3\xra{}M^2$ be the map $(a,b,c)\mapsto(c+2a,3b+c)$''.  This is
essentially the same as thinking of $\CD$ as a subcategory of
$[\CD^\op,\Sets]$, by the Yoneda embedding.

Given a monoid object $M$, we define maps $d_0,d_1\:M^3\xra{}M^2$ and
$s\:M^2\xra{}M$ by
\begin{align*}
 d_0(a,b,x) &= (a,b)            \\
 d_1(a,b,x) &= (a+x,b+x)        \\
 s(a,b)     &= (a,b,0).
\end{align*}
This is clearly a \idx{reflexive fork} in $\Mon\CD$.
\begin{proposition}
 If the above fork has a strong coequaliser $q\:M^2\xra{}H$ in $\CD$,
 then $H$ has a unique group structure making $q$ into a homomorphism
 of monoids, and with that group structure we have $H=G(M)$.
\end{proposition}
\begin{proof}
 Firstly, Proposition~\ref{prop-refl-fork} tells us that there is a
 unique monoid structure on $H$ making $q$ into a monoid map, and that
 this makes $H$ into the coequaliser in $\Mon\CD$.  We define a monoid
 map $\nu'\:M^2\xra{}H$ by $\nu'(a,b)=q(b,a)$.  Clearly
 $\nu'd_0(a,b,x)=qd_0(b,a,x)$ and $\nu'd_1(a,b,x)=qd_1(b,a,x)$ but
 $qd_0=qd_1$ so $\nu'd_0=\nu'd_1$, so there is a unique map
 $\nu\:H\xra{}H$ with $\nu'=\nu q$.  We then have
 \begin{align*}
  q(a,b)+\nu q(a,b)
  &= q(a,b) + q(b,a) \\
  &= q(a+b,a+b)      \\
  &= q d_1(0,0,a+b)  \\
  &= q d_0(0,0,a+b)  \\
  &= q(0,0) = 0.
 \end{align*}
 This shows that $(1+\nu)q=0$, but $q$ is an epimorphism so
 $1+\nu=0$.  This means that $\nu$ is a negation map for $H$, making
 it into a group object.  We let $j\:M\xra{}H$ be the map
 $a\mapsto q(a,0)$, which is clearly a homomorphism of monoids.
 Clearly $q(a,b)=q(a,0)+q(0,b)=j(a)+\nu j(b)=j(a)-j(b)$.

 Now let $K$ be another Abelian group object, and let $f\:M\xra{}K$ be
 a homomorphism of monoids.  We define $f'\:M^2\xra{}K$ by
 $f'(a,b)=f(a)-f(b)$.  It is clear that $f'd_0=f'd_1$, so we get a unique
 monoid map $f''\:H\xra{}K$ with $f''q=f'$.  In particular, we have
 $f''j(a)=f''q(a,0)=f'(a,0)=f(a)$, so that $f''j=f$.  If $g\:H\xra{}K$
 is another homomorphism with $gj=f$ then
 $gq(a,b)=g(j(a)-j(b))=f(a)-f(b)=f''q(a,b)$, and $q$ is an epimorphism
 so $g=f''$.

 This shows that group maps $H\xra{}K$ biject with monoid maps
 $M\xra{}K$ by the correspondence $g\mapsto gj$, which means that
 $H=G(M)$ as claimed.
\end{proof}

\subsection{Free formal groups}
\label{subsec-formal-groups}

We next discuss the existence of free Abelian formal groups.
\begin{proposition}\label{prop-mon-filtered-base}
 Let $Y$ be a formal scheme over a formal scheme $X$.  Write $X$ as a
 filtered colimit of informal schemes $X_i$, and put $Y_i=Y\tX X_i$.
 If $M^+(Y_i)$ exists in $\Mon\FX_{X_i}$ for all $i$, then $M^+(Y)$
 exists and is equal to $\colim_iM^+(Y_i)$.  Similar remarks apply to
 $M(Y)$ and (if $Y$ has a given section $0\:X\xra{}Y$) to $N^+(Y)$ and
 $N(Y)$.
\end{proposition}
\begin{proof}
 We use the notation of Definition~\ref{defn-filtered-base} and
 Proposition~\ref{prop-filtered-base}.  It is clear that
 $\{M^+(Y_i)\}$ is the free Abelian monoid object on $\{Y_i\}$ in the
 category $\CD_{\{X_i\}}$.  As the functor
 $F\:\CD_{\{X_i\}}\xra{}\FX_X$ preserves finite limits, we see that
 $L=\colim_iM^+(Y_i)=F\{M^+(Y_i)\}$ is an Abelian monoid object in
 $\FX_X$.  Using the fact that $F$ preserves finite products and is
 left adjoint to $G$, we see that
 \[ \FX_X(L^m,Z) = \CD_{\{X_i\}}(\{M^+(Y_i)^m_{X_i}\},\{Z\tX X_i\})
 \] 
 for all $Z\in\FX_Z$.  Using this, one can check that 
 \begin{align*}
  \Mon\FX(L,M) &= \Mon\CD_{\{X_i\}}(\{M^+(Y_i)\},\{M\tX X_i\}) \\
               &= \CD_{\{X_i\}}(\{Y_i\},\{M\tX X_i\}) = \FX_X(Y,M),
 \end{align*}
 as required.  We leave the case of $M(Y)$ and so on to the reader.
\end{proof}

\begin{proposition}\label{prop-exists-MpY}
 If $Y$ is a coalgebraic formal scheme over $X$, then the free Abelian
 monoid scheme $M^+(Y)$\index{MUa@$M^+(U)$} exists.  If $Y$ also has a
 specified section $0\:X\xra{}Y$ (making it an object of
 $\Based\FX_X$) then $N^+(Y)$\label{NVa@$N^+(V)$} exists.
\end{proposition}
\begin{proof}
 By the previous proposition, we may assume that $X$ is informal, and
 that $Y=\sch_X(U)$ for some coalgebra $U$ over $R=\OX$ with a good
 basis $\{e_i\st i\in I\}$.  We know from
 Example~\ref{eg-strong-sym-formal} that $Y^k_X/\Sg_k$ is a strong
 colimit for the action of $\Sg_k$ on $Y^k_X$.  Moreover,
 $\coprod_k Y^k/\Sg_k$ exists as a strong coproduct by
 Corollary~\ref{cor-strong-coprod}.  We conclude from
 Proposition~\ref{prop-exists-MpU} that $M^+(Y)=\coprod_kY^k/\Sg_k$.
 In the based case, we observe that the diagram $\{Y^k/\Sg_k\}$ is
 just indexed by $\nat$ and thus is filtered, and filtered colimits
 exists and are strong in $\FX_X$ by
 Proposition~\ref{prop-formal-lim}.  Given this,
 Proposition~\ref{prop-exists-NpU} completes the proof.
\end{proof}
\begin{remark}\label{rem-symmetric-algebra}
 If $X$ is informal we see that the coalgebra $cM^+(Y)$ is just the
 symmetric algebra generated by $cY$ over $\OO_X$.  In the based case,
 if $e_0\in cY$ is the basepoint then $cN^+(Y)=cM^+(Y)/(e_0-1)$.
\end{remark}

We next show that in certain cases of interest, the free Abelian
monoid $N^+(Y)$ constructed above is actually a group.
\begin{definition}
 A \emph{good filtration} of a coalgebra $U$ over a ring $R$ is a
 sequence of submodules $F_sU$ for $s\ge 0$ such that
 \begin{itemize}
 \item[(a)] $\ep\:F_0U\xra{}A$ is an isomorphism.
 \item[(b)] For $s>0$ the quotient $G_sU=F_sU/F_{s-1}U$ is a finitely
  generated free module over $R$.
 \item[(c)] $\bigcup_sF_sU=U$
 \item[(d)] $\psi(F_sU)\sse\sum_{s=t+u}F_sU\ot F_tU$.
 \end{itemize}
 We write $\CC''=\CC''_R=\CC''_Z$ for the category of coalgebras that
 admit a good filtration.  Given a good filtration, we write $\eta$
 for the composite $A\xra{\ep^{-1}}F_0U\mra U$.  One can check that
 this is a coalgebra map, so it makes $U$ into a based coalgebra.  A
 \emph{good basepoint} for $U$ is a basepoint which arises in this
 way.  We say that a \emph{very good basis} for $U$ is a basis
 $\{e_0,e_1,\ldots\}$ for $U$ over $R$ such that
 \begin{itemize}
 \item[(i)]    $e_0=\eta(1)$
 \item[(ii)]   $\ep(e_i)=0$ for $i>0$
 \item[(iii)]  There exist integers $N_s$ such that $\{e_i\st i<N_s\}$
  is a basis for $F_sU$.
 \end{itemize}
 One can check that very good bases exist, and that a very good basis
 is a good basis.
\end{definition}

\begin{proposition}\label{prop-filt-prod}
 If $U$ and $V$ lie in $\CC''_Z$ then so do $U\tm V$ and $U^k/\Sg_k$.
 If we choose a good basepoint for $U$ then we can define
 $N^+(U)$\index{NVa@$N^+(V)$}, and it again lies in $\CC''_Z$.
\end{proposition}
\begin{proof}
 Choose good filtrations on $U$ and $V$.  Define a filtration on
 $U\tm V=U\ot V$ by setting $F_s(U\ot V)=\sum_{s=t+u}F_tU\ot F_uV$.
 It is not hard to check that this is good.  Essentially the same
 procedure gives a filtration of $U^{\ot m}$.  This is invariant under
 the action of the symmetric group $\Sg_m$, so we get an induced
 filtration of the group of coinvariants $U^{\ot m}_{\Sg_m}$.  Our
 filtrations on these groups are compatible as $m$ varies, so we get
 an induced filtration of $N(U)=\colim_mU^{\ot m}_{\Sg_m}$.  Using a
 very good basis for $U$ and the associated monomial basis for $N(U)$,
 we can check that the filtration of $N(U)$ is good.
\end{proof}

\begin{proposition}\label{prop-filt-mon}
 Let $U$ be an Abelian monoid object in $\CC_Z$, with addition map
 $\sg\:U\tm U=U\ot U\xra{}U$.  If $U$ admits a good filtration such
 that the basepoint is good and $\sg(F_sU\ot F_tU)\sse F_{s+t}U$ for
 all $s,t\ge 0$, then $U$ is actually an Abelian group object.
\end{proposition}
\begin{proof}
 First note that we can use $\sg$ to make $U$ into a ring.  We need to
 construct a negation map (otherwise known as an antipode)
 $\chi\:U\xra{}U$, which must be a coalgebra map satisfying
 $\sg(1\ot\chi)\psi=\eta\ep$.  In terms of elements, if
 $\psi(a)=1\ot a + \sum a'\ot a''$ then the requirement is that
 $\chi(a)=\eta\ep(a)-\sum a'\chi(a'')$.  The idea is to use this
 formula to define $\chi$ on $F_sU$ by recursion on $s$.

 Write $\psib=\psi-\eta\ot 1\:U\xra{}U\ot U$.  Note that
 $\psib(F_sU)\sse\sum_{t=0}^sF_{s-t}U\ot F_tU$, and that
 $(\ep\ot 1)\psib=0$.  Choose a very good basis $\{e_i\}$ for $U$, and
 write $\psib(e_i)=\sum_{j,k}a_{ijk}e_j\ot e_k$.  Suppose that
 $N_{s-1}\le i<N_s$, so that $e_i\in F_sU\setminus F_{s-1}U$.  If
 $j>0$ and $k\ge N_{s-1}$ then $e_j\ot e_k\not\in F_s(U\ot U)$ so
 $a_{ijk}=0$.  On the other hand, the equation
 $(\ep\ot 1)\psib(e_i)=0$ gives $\sum_ma_{i0m}e_m=0$ for all $m$, so
 $a_{i0k}=0$, so $a_{ijk}=0$ for all $j$.  This applies for all
 $k\ge N_{s-1}$, and thus in particular for $k\ge i$.  

 We now define $\chi(e_i)$ recursively by $\chi(e_0)=e_0$ and
 \[ \chi(e_i)=-\sum_{0\le k<i}a_{ijk}e_j\chi(e_k) \]
 for $i>0$.  By the previous paragraph, we actually have
 $\chi(e_i)=-\sum_{k\ge 0}a_{ijk}e_j\chi(e_k)$, and it follows that
 $\sg(1\ot\chi)\psi=\eta\ep$ as required.  We still have to check that
 $\chi$ is a coalgebra map.  For the counit, it is clear that
 $\ep\chi(e_0)=\ep(e_0)$.  If we assume inductively that
 $\ep(\chi(e_k))=\ep(e_k)=0$ for $0<k<i$ then we find that
 \[ \ep\chi(e_i)=-\sum_{0\le k<i}a_{ijk}\ep(e_j)\ep\chi(e_k)
     = a_{i00}=(\ep\ot\ep)\psi(e_i)=\ep(e_i)=0.
 \]
 A similar, but slightly more complicated, induction shows that
 $\psi\chi=(\chi\ot\chi)\psi$, so $\chi$ is a coalgebra map as
 required. 
\end{proof}

\begin{proposition}\label{prop-Mplus-div}
 Let $C$ be a pointed formal curve over a formal scheme $X$.  Then
 there are natural isomorphisms 
 \begin{align*}
  M^+(C) &= \Div^+(C)                 \\
  N^+(C) &= N(C) = \Div_0(C)          \\
  M(C)   &= \Div(C).
 \end{align*}
 \index{MUa@$M^+(U)$} \index{NVa@$N^+(V)$}
 \index{MU@$M(U)$} \index{NV@$N(V)$}
 \index{Diva@$\Div^+(C)$} \index{Div@$\Div(C)$}
\end{proposition}
\begin{proof}
 This follows easily from the constructions in
 Section~\ref{subsec-group-objects} and the results above.
\end{proof}

\subsection{Schemes of homomorphisms}
\label{subsec-hom-schemes}

\begin{definition}
 Given formal groups $G$ and $H$ over $X$ and a ring $R$, we let
 $\Hom_X(G,H)(R)$ be the set of pairs $(x,u)$, where $x\in X(R)$ and
 $u\:G_x\xra{}H_x$ is a homomorphism of formal groups over
 $\spec(R)$.  This is a subfunctor of $\Map_X(G,H)$, so we have
 defined an object $\Hom_X(G,H)\in\CF$.  It is not hard to define an
 equaliser diagram
 \[ \Hom_X(G,H) \xra{} \Map_X(G,H) \arw{e,tb,=>}{d_0}{d_1}
     \Map_X(G\tm_XG,H).
 \]
 In more detail, note that a point of $\Map_X(G,H)$ is a map
 $x\:\spec(R)\xra{}X$ together with a map $f\:G_x\xra{}H_x$ of schemes
 over $\spec(R)$.  Given such a pair $(x,f)$, we define
 $g,h\:G_x\tm_{\spec(R)}G_x\xra{}H_x$ by $g(a,b)=f(a+b)$ and
 $h(a,b)=f(a)+f(b)$, and then we define $d_i$ by $d_0(f)=g$ and
 $d_1(f)=h$.  
\end{definition}

\begin{proposition}\label{prop-exists-Hom}
 Let $G$ and $H$ be formal groups over $X$.  If $G$ is finite and very
 flat over $X$, or if $G$ is coalgebraic and $H$ is relatively
 informal, or if $G$ is very flat and $H$ is of finite presentation,
 then $\Hom_X(G,H)$ is a formal scheme and there is a natural
 isomorphism 
 \[ \FX_X(Y,\Hom_X(G,H))= \Ab\FX_Y(G\tX Y,H\tX Y) \]
 for all $Y\in\FX_X$.
\end{proposition}
\begin{proof}
 Theorem~\ref{thm-map-formal} tells us that $\Map_X(G,H)$ and
 $\Map_X(G\tX G,H)$ are formal schemes, and $\FX_X$ is closed under
 finite limits in $\CF$, so $\Hom_X(G,H)$ is a formal scheme.  The
 natural isomorphism comes from the Yoneda lemma when $Y$ is informal,
 and follows in general by passage to colimits.
\end{proof}

\begin{example}
 Let $\AFG$ be the additive formal group (over the terminal scheme
 $1=\spec(\Zh)$) defined by $\AFG(R)=\Nil(R)$, with the usual
 addition.  Thus, the underlying scheme of $\AFG$ is just $\haf^1$.
 This is coalgebraic over $1$, so we see that
 $\End(\AFG)=\Hom_1(\AFG,\AFG)$\index{EndGa@$\End(\AFG)$} exists.  One
 checks that any map $\haf^1\tm Y\xra{}\haf^1\tm Y$ over $Y$ is given
 by a unique power series $f(x)\in\fps{\OY}{x}$ such that $f(0)$ is
 nilpotent.  It follows easily that $\End(\AFG)(R)$ is the set of
 power series $f\in\fps{R}{x}$ such that
 $f(x+y)=f(x)+f(y)\in\fps{R}{x,y}$.  If $R$ is an algebra over $\Fp$,
 then a well-known lemma says that $f(x+y)=f(x)+f(y)$ if and only if
 $f$ can be written in the form $f(x)=\sum_k a_k x^{p^k}$, for
 uniquely determined coefficients $a_k\in R$.  One can deduce that
 $\spec(\Fp)\tm\End(\AFG)=\spec(\Fp[a_k\st k\ge 0])$.
\end{example}
\begin{example}
 A similar analysis shows that $\End(\MG)(R)$\index{EndGm@$\End(\MG)$}
 is the set of Laurent polynomials $f\in R[u^{\pm 1}]$ such that
 $f(u)f(v)=f(uv)$ and $f(1)=1$.  If $f(u)=\sum_ke_ku^k$, we find that
 the elements $e_k\in R$ are orthogonal idempotents with $\sum_ke_k=1$.
 It follows that $\End(\MG)$ is the constant formal scheme $\un{\Zh}$,
 with the $n$'th piece in the coproduct corresponding to the
 endomorphism $u\mapsto u^n$.
\end{example}
\begin{example}
 We can also form the scheme $\Exp=\Hom(\AFG,\MG)$\index{Exp@$\Exp$}.
 In this case, $\Exp(R)$ is the set of power series
 $f(x)=\sum_ka^{[k]}x^k$ such that $f(0)=1$ and $f(x+y)=f(x)f(y)$, or
 equivalently $a^{[0]}=1$ and $a^{[i]}a^{[j]}=\bcf{i+j}{i}a^{[i+j]}$.
 In other words, a point of $\Exp(R)$ is an element $a=a^{[1]}$ of $R$
 together with a specified system of divided powers for $a$.  Clearly,
 if $R$ is a $\rat$-algebra then there is a unique possible system of
 divided powers, viz. $a^{[k]}=a^k/k!$, so
 $\spec(\rat)\tm\Exp\simeq\spec(\rat)\tm\aff^1$.

 Now let $R$ be an $\Fp$-algebra.  Given an element $a\in R$ with
 $a^p=0$, we define $T(a)(x)=\sum_{j=0}^{p-1}a^jx^j/j!$; it is not
 hard to see that $T(a)\in\Exp(R)$.  Given a sequence of such elements
 $\un{a}=(a_0,a_1,\ldots)$, we define
 $T(\un{a})(x)=\prod_iT(a_i)(x^{p^i})$; it is not hard to check that
 the product is convergent in the $x$-adic topology on $\fps{R}{x}$,
 and that $T(\un{a})\in\Exp(R)$.  Thus $T$ defines a map
 $\spec(\Fp)\tm D_p^\nat\xra{}\spec(\Fp)\tm\Exp$.  It can be shown
 that this is an isomorphism.

 More generally, we have $\Exp=\spec(D_\Zh[a])$, where $D_\Zh[a]$ is
 the divided-power algebra on one generator $a$ over $\Zh$.  The
 previous paragraph is equivalent to the fact that
 $D_\Fp[a]=D_\Zh[a]/p=\Fp[a_k\st k\ge 0]/(a_k^p)$, where
 $a_k=a^{[p^k]}$.
\end{example}

\subsection{Cartier duality}
\label{subsec-cartier}

Let $G$ be a coalgebraic commutative formal group over a formal scheme
$X$.  By Proposition~\ref{prop-exists-Hom}, we can define the group
scheme $DG=\Hom_X(G,\MG\tm X)$\index{DG@$DG$}.  We call this the
\dfn{Cartier dual} of $G$.  Note also that the product structure on
$G$ makes $cG$ into commutative group in the category of coalgebras,
in other words a Hopf algebra, and in particular an algebra over
$\OX$.  We can thus define $H=\spec(cG)$, which is an informal scheme
over $X$.  The coproduct on $cG$ gives a product on $H$, making it
into a group scheme over $X$.  Moreover, we know that $cG$ is a free
module over $\OX$, so that $H$ is very flat over $\OX$.  Thus, by
Proposition~\ref{prop-exists-Hom}, we can define a formal group scheme
$DH=\Hom_X(H,\MG\tm X)$.  We again call this the Cartier dual of $H$.
These definitions appear in various levels of generality in many
places in the literature; the treatment in~\cite{de:lpd} is similar in
spirit to ours, although restricted to the case where $\OX$ is a
field. 

\begin{proposition}\label{prop-cartier}
 If $G$ and $H$ are as above, then $DG=H$ and $DH=G$.
\end{proposition}
\begin{proof}
 First suppose that $X=\spec(R)$ is informal.  We shall analyse the
 set $\FX_X(X,DG)$ of sections of the map $DG\xra{}X$.  {}From the
 definitions, we see that a section of the map $DG\xra{}X$ is the same
 as a map $G\xra{}\MG\tm X$ of formal groups over $X$, or equivalently
 a map of Hopf algebras $\OO_{\MG\tm X}\xra{}\OG$.  As
 $\OO_{\MG\tm X}=R[u^{\pm 1}]$ with $\ep(u)=1$ and $\psi(u)=u\ot u$,
 such a map is equivalent to an element $v\in\OGt$ with
 $\ep(v)=1$ and $\psi(v)=v\ot v$.  In fact, if $v$ is any element with
 $\ep(v)=1$ and $\psi(v)=v\ot v$ then the Hopf algebra axioms imply
 that $v\chi(v)=1$ so we do not need to require separately that $v$ be
 invertible.  As $G$ is coalgebraic we have $\OG=\Hom_R(cG,R)$, so
 we can regard $v$ as a map $cG\xra{}R$ of $R$-modules.  The
 conditions $\ep(v)=1$ and $\psi(v)=v\ot v$ then become $v(1)=1$ and
 $v(ab)=v(a)v(b)$, so the set of such $v$'s is just
 $\Alg_R(cG,R)=\FX_X(X,H)$.

 Now let $X$ be arbitrary.  The above (together with the commutation
 of various constructions with pullbacks, which we leave to the
 reader) shows that for any informal scheme $W$ over $X$ we have
 $\FX_X(W,DG)=\FX_W(W,D(G\tm_XW))=\FX_W(W,H\tm_XW)=\FX_X(W,H)$.  It
 follows that $DG=H$ as claimed.

 We now show that $DH=G$.  Just as previously, we may assume that
 $X=\spec(R)$ is informal, and it is enough to show that $DH$ and $G$
 have the same sections.  Again, the sections of $DH$ are just the
 elements $v\in\OO_H=cG$ with $\ep(v)=1$ and $\psi(v)=v\ot v$.  In
 this case, we identify $cG$ with the continuous dual of $\OG$, so
 $v$ is a continuous map $\OG\xra{}R$ of $R$-algebras, and thus a
 section of $\spf(\OG)=G$ as required.
\end{proof}

\subsection{Torsors}
\label{subsec-torsors}

Let $G$ be a formal group over a formal scheme $X$.  Let $T$ be a
formal scheme over $X$ with an action of $G$.  More explicitly, we
have an action map $\al\:G\tX T\xra{}T$, so whenever $g$ and $t$ are
points of $G$ and $T$ with the same image in $X$, we can define
$g.t=\al(g,t)$.  This is required to satisfy $1.t=t$ and
$g.(h.t)=(gh).t$ (whenever $g$, $h$ and $t$ all have the same image in
$X$).  We write $G\FX_X$\index{GXX@$G\FX_X$} for the category of such
$T$.  Note that $G$ itself can be regarded as an object of $G\FX_X$.

If $Y$ is a scheme with a specified map $p\:Y\xra{}X$ we shall allow
ourselves to write $G\FX_Y$ instead of $(p^*G)\FX_Y$.  It is easy to
see that $p^*$ gives a functor $G\FX_X\xra{}G\FX_Y$.

\begin{definition}\label{defn-torsor}
 Let $G$ be a formal group over a formal scheme $X$, and let $T$ be a
 formal scheme over $X$ with an action of $G$.  We say that $T$ is a
 \emph{$G$-torsor over $X$}\index{torsor} if there exists a faithfully
 flat map $p\:Y\xra{}X$ such that $p^*T\simeq p^*G$ in $G\FX_Y$.  We
 write $G\CT_X$\index{GTX@$G\CT_X$} for the category of $G$-torsors
 over $X$.
\end{definition}

\begin{example}\label{eg-bases-torsor}
 Let $M$ be a \idx{vector bundle} over $X$ of rank $d$, and let
 $\Bases(M)$\index{BasesM@$\Bases(M)$} be as in
 Example~\ref{eg-bases-M}.  Let $\GL_d$ be the group scheme of
 invertible $d\tm d$ matrices.  Then $\GL_d\tm X$ acts on $\Bases(M)$,
 and if $M$ is free then $\Bases(M)\simeq\GL_d\tm X$.  As we can
 always pull back along a faithfully flat map $p\:Y\xra{}X$ to make
 $M$ free, and $\Bases(p^*M)=p^*\Bases(M)$, we find that $\Bases(M)$
 is a torsor for $\GL_d\tm X$.
\end{example}
\begin{example}\label{eg-coord-torsor}
 Let $C$ be a pointed formal curve over $X$, let
 $\Coord(C)$\index{CoordC@$\Coord(C)$} be as in
 Definition~\ref{defn-coord}, and let $\IPS$ be as in
 Example~\ref{eg-IPS}.  Then $\Coord(C)$ is a torsor for group scheme
 $\IPS\tm X$.  In fact, this torsor is trivialisable (i.e. isomorphic
 to $\IPS\tm X$ even without pulling back) but not canonically so.
\end{example}

\begin{proposition}\label{prop-torsor-groupoid}
 Every morphism in $G\CT_X$\index{GTX@$G\CT_X$} is an isomorphism, so
 $G\CT_X$ is a groupoid.
\end{proposition}
\begin{proof}
 First, let $u\:G\xra{}G$ be a map of $G$-torsors.  As $u$ is
 $G$-equivariant we have $u(g)=g.u(1)$, so $h\mapsto h.u(1)^{-1}$ is
 an inverse for $u$.  Now let $u\:S\xra{}T$ be an arbitrary map of
 $G$-torsors.  Then there is a faithfully flat map $p\:Y\xra{}X$ such
 that $p^*S\simeq p^*T\simeq p^*G$, so the first case tells us that
 $p^*u$ is an isomorphism.  As $p$ is faithfully flat, we see that
 $p^*$ reflects isomorphisms, so $u$ is an isomorphism.
\end{proof}

\begin{proposition}\label{prop-torsor-functors}
 Every homomorphism $\phi\:G\xra{}H$ of formal groups over $X$ gives
 rise to functors
 \begin{align*}
  \phi^\bullet &\: H\FX_X\xra{}G\FX_X   \\
  \phi_\bullet &\: G\CT_X\xra{}H\CT_X,
 \end{align*}
 such that
 \[ H\FX_X(\phi_\bullet T,U) = G\FX_X(T,\phi^\bullet U) \]
 for all $U\in H\FX_X$.
\end{proposition}
\begin{proof}
 The functor $\phi^\bullet$ is just $\phi^\bullet U=U$, with
 $G$-action $g.u:=\phi(g).u$.  Let $\phi_\bullet T$ be the coequaliser
 of the maps $(h,g,t)\mapsto(h\phi(g),t)$ and $(h,g,t)\mapsto(h,g.t)$
 from $H\tX G\tX T$ to $H\tX T$.  Note that these maps have a
 common splitting $(h,t)\mapsto (h,1,t)$, so we have a reflexive
 fork.  In the case $T=G$, the coequaliser is just the map
 $H\tX G\xra{}H$ given by $(h,g)\mapsto h\phi(g)$.  In fact, this
 coequaliser is split by the maps $h\mapsto (h,1)$ and
 $(h,g)\mapsto(h,g,1)$, so it is a strong coequaliser.  

 Now consider a general $G$-torsor $T$.  We claim that the coequaliser
 that defines $\phi_\bullet T$ is strong.  By
 proposition~\ref{prop-ff-strong}, we can check this after pulling
 back along a faithfully flat map $p\:Y\xra{}X$.  We can choose $p$ so
 that $p^*T\simeq p^*G$, and then the claim follows from the previous
 paragraph.  

 We can let $H$ act on the left on $H\tX G\tX T$ and $H\tX G$, and
 then the maps whose coequaliser defines $\phi_\bullet T$ are both
 $H$-equivariant.  The reader can easily check that if a fork in
 $H\FX_X$ has a strong coequaliser in $\FX_X$ then the coequaliser has
 a unique $H$-action making it the coequaliser in $H\FX_X$.  This
 implies that $\phi_\bullet T$ is the coequaliser of our fork in
 $H\FX_X$, and one can deduce that
 \[ H\FX_X(\phi_\bullet T,U) = G\FX_X(T,\phi^\bullet U) \]
 for all $U\in H\FX_X$.

 All that is left is to check that $\phi_\bullet T$ is a torsor.  For
 this, we just choose a faithfully flat map $p$ such that
 $p^*T\simeq p^*G$, and observe that
 $p^*\phi_\bullet T=\phi_\bullet p^*T\simeq p^*H$.
\end{proof}

\begin{proposition}\label{prop-torsor-tensor}
 If $G$ is an Abelian formal group over $X$, then there is a functor
 $\ot\:G\CT_X\tm G\CT_X\xra{}G\CT_X$ which makes $G\CT_X$
 \index{GTX@$G\CT_X$} into a symmetric monoidal category with unit
 $G$.  Moreover, the twist map of $T\ot T$ is always the identity, and
 every object is invertible under $\ot$, so that $G\CT_X$ is a strict
 Picard category.
\end{proposition}
\begin{proof}
 If $S$ and $T$ are $G$-torsors over $X$, then it is easy to see that
 $S\tX T$ has a natural structure as a $G\tX G$-torsor.  As $G$ is
 Abelian, the multiplication map $\mu\:G\tX G\xra{}G$ is a
 homomorphism, so we can define $S\ot T=\mu_\bullet(S\tX T)$.  We
 leave it to the reader to check that this gives a symmetric monoidal
 structure with unit $G$.  If we let $\chi\:G\xra{}G$ denote the map
 $g\mapsto g^{-1}$ then $\chi$ is also a homomorphism, so we can
 define $T^{-1}=\chi_\bullet T$.  We then have
 $T\ot T^{-1}=(\mu(1\tm\chi))_\bullet(T\tX T)=0_\bullet(T\tX T)=G$,
 so $T^{-1}$ is an inverse for $T$.  Finally, we need to show that the
 twist map $\tau\:T\ot T\xra{}T\ot T$ is the identity.  As the map
 $q\:T\tX T\xra{}T\ot T$ is a regular epimorphism, it suffices to
 show that $\tau q=q$, and clearly $\tau q(a,b)=q(b,a)$ so we need to
 show that $q(a,b)=q(b,a)$.  In the case $T=G$ we have $T\ot T=G$ and
 the map $q$ is just $q(a,b)=ab$, so the claim holds.  For general
 $T$, we just pull back along a faithfully flat map $p$ such that
 $p^*T\simeq p^*G$ and use the fact that $p^*$ is faithful.
\end{proof}

\begin{proposition}
 Let $\MG$ denote the multiplicative group, which is defined by
 $\MG(R)=R^\tm$.  Then the functor $L\mapsto\aff(L)^\tm$
 \index{AALt@$\aff(L)^\tm$} (as in Definition~\ref{defn-gen-L} and
 Remark~\ref{rem-aff-formal}) is an equivalence from the category of
 \idx{line bundle}s over $X$ and isomorphisms, to the category of
 $\MG$-torsors over $X$.  Moreover, this equivalence respects tensor
 products.
\end{proposition}
\begin{proof}
 Let $L$ be a line bundle over $X$.  For any $x\in X(R)$, we have a
 rank one projective module $L_x$ over $R$, and clearly $R^\tm=\MG(R)$
 acts on the set of bases for $L_x$ (even though this set may be
 empty).  If $L$ is free then it is clear that
 $\aff(L)^\tm\simeq\aff(\OO)^\tm=\MG\tm X$, and thus that
 $\aff(L)^\tm$ is a torsor.  In general, we know from
 Proposition~\ref{prop-locfree-formal} that $L$ is fpqc-locally
 isomorphic to $\OO$, so $\aff(L)^\tm$ is fpqc-locally isomorphic to
 $\aff(\OO)^\tm=\MG\tm X$, and thus is a torsor.

 In the opposite direction, let $T$ be a $\MG$-torsor over $X$.
 Define a formal scheme $A$ over $X$ by the coequaliser
 \[ \aff^1\tm\MG\tm T\arw{e,tb,=>}{\lm}{\rho}\aff^1\tm T\xra{}A, \]
 where $\lm(a,u,t)=(au,t)$ and $\rho(a,u,t)=(a,ut)$.  Locally in the
 flat topology we may assume that $T=\MG\tm X$, and it is easy to
 check that $\aff^1\tm X$ is the split coequaliser of the fork.  Thus
 Proposition~\ref{prop-ff-strong} tells us that $A$ is the strong
 coequaliser of the original fork.  Also, we can make
 $\aff^1\tm\MG\tm T$ and $\aff^1\tm T$ into modules over the ring
 scheme $\aff^1$.  As the functor $\aff^1\tm(-)$ preserves our
 coequaliser, the formal scheme $A$ is also a module over $\aff^1$.
 This means that if we define $L_x$ to be the preimage of
 $x\in X(R)$ under the map $A(R)\xra{}X(R)$, then $L_x$ is an
 $R$-module.  Locally on $X$ we have $T\simeq\MG\tm X$ and thus
 $A\simeq\aff^1$ and thus $L_x\simeq R$.  One can deduce that $L$ is a
 line bundle over $X$, with $\aff(L)=A$ and thus $\aff(L)^\tm=T$.  

 We leave it to the reader to check that this gives an equivalence of
 categories, which preserves tensor products.
\end{proof}

\section{Ordinary formal groups}
\label{sec-ordinary-groups}

Recall that an \dfn{ordinary formal group} over a scheme $X$ is a
formal group $G$ over $X$ that is isomorphic to $X\tm\haf^1$ as a
formal scheme over $X$.  In particular, $G$ is a pointed formal curve
\index{formal curve!pointed} over $X$, so we can choose a normalised
coordinate $x$ on $G$ giving an isomorphism $G\simeq\haf^1\tm X$ in
$\Based\FX_X$.  However, for the usual reasons it is best to proceed
as far as possible in a coordinate-free way.  Lazard's
book~\cite{la:cfg} gives an account in this spirit, but in a somewhat
different framework.

If we do choose a coordinate $x$ on $G$ then we have a function
$(g,h)\mapsto x(g+h)$ from $G\tX G$ to $\haf^1$.  
As $G\tX G\simeq\haf^2\tm X$, we see that this can be written uniquely
in the form $x(g+h)=\sum_{i,j}a_{ij}x(g)^ix(h)^j=F_x(x(g),x(h))$ for
some power series $F_x(s,t)\in\fps{\OX}{s,t}$.  It is easy to see that
this is a formal group law (Example~\ref{eg-FGL}), so we get a map
$X\xra{}\FGL$\index{FGL@$\FGL$}.  This construction gives a canonical
map $\Coord(G)\xra{}\FGL$\index{CoordC@$\Coord(C)$}.  We can let the
group scheme $\IPS$\index{IPS@$\IPS$} act on $\FGL$ as in
Example~\ref{eg-IPS}, and on $\Coord(G)$ by $f.x=f(x)$.  It is easy to
see that the map $\Coord(G)\xra{}\FGL$ is $\IPS$-equivariant.

\begin{definition}\label{defn-inv-diff}
 Let $G$ be a formal group over an affine scheme $X$.  Let $I$ be the
 ideal in $\OX$ of functions $g\:X\xra{}\aff^1$ such that $g(0)=0$.

 Define $\om_G=\om_{G/X}=I/I^2$\index{oGX@$\om_{G/X}$}, and let $d_0(g)$
 denote the image of $g$ in $\om_{G/X}$.  We also define
 \[ \Prim(\Om^1_{G/X})=
     \{\al\in\Om^1_{G/X}\st
       \sg^*\al=\pi_0^*\al+\pi_1^*\al\in\Om_{G\tm_XG/X}\}.
 \]
 \index{O1@$\Om^1$} \index{PrimOGX@$\Prim(\Om^1_{G/X})$}
 Here $\pi_0,\pi_1\:G\tm_XG\xra{}G$ are the two projections, and
 $\sg\:G\tm_XG\xra{}G$ is the addition map.
\end{definition}

We now give a formal version of the fact that left-invariant
differential forms on a Lie group biject with elements of the
cotangent space at the identity element.
\begin{proposition}\label{prop-inv-diff}
 $\om_{G/X}$\index{oGX@$\om_{G/X}$} is a free module on one generator
 over $\OX$.  Moreover, there are natural isomorphisms
 $\om_{G/X}\simeq\Prim(\Om^1_{G/X})$ and
 $\Om^1_{G/X}=\OG\ot_{\OX}\om_{G/X}$. \index{O1@$\Om^1$}
 \index{PrimOGX@$\Prim(\Om^1_{G/X})$}
\end{proposition}
\begin{proof}
 Let $x$ be a normalised coordinate on $G$.  We then have
 $\OG=\fps{\OX}{x}$, and it is easy to check that $I=(x)$ so
 $I^2=(x^2)$ so $\om_{G/X}$ is freely generated over $\OX$ by
 $d_0(x)$.

 Now let $K$ be the ideal in $\OO_{G\tm_XG}$ of functions $k$ such that
 $k(0,0)=0$.  In terms of the usual description
 $\OO_{G\tm_XG}=\fps{\OX}{x',x''}$, this is just the ideal generated by
 $x'$ and $x''$.  Given $g\in I$, we define
 $\dl(g)(u,v)=g(u+v)-g(u)-g(v)$.  We claim that $\dl(g)\in K^2$.
 Indeed, we clearly have $\dl(g)(0,v)=0$, so $\dl(g)$ is divisible by
 $x'$.  We also have $\dl(g)(u,0)=0$, so $\dl(g)$ is divisible by
 $x''$.  It follows easily that $\dl(g)$ is divisible by $x'x''$ and
 thus that it lies in $K^2$ as claimed.  

 Next, let $J$ be the ideal of functions on $G\tm_XG$ that vanish on
 the diagonal (so we have $\Om^1_{G/X}=J/J^2$).  For any function
 $g\in I$ we define $\lm(g)\in J$ by $\lm(g)(u,v)=g(u-v)$.  As
 $g(0)=0$ we see that $\lm(g)\in J$, so $\lm$ induces a map
 $\om_{G/X}\xra{}\Om^1_{G/X}$.  We claim that
 $\lm(g)\in\Prim(\Om^1_{G/X})$.  To make this more explicit, let $L$
 be the ideal of functions $l$ on $G^4_X$ such that $l(s,s,u,u)=0$.
 The claim is that $\sg^*\lm(g)-\pi_0^*\lm(g)-\pi_1^*\lm(g)=0$ in
 $L/L^2$, or equivalently that the function
 \[ k\:(s,t,u,v)\mapsto\lm(g)(s+u,t+v)-\lm(g)(s,t)-\lm(g)(u,v) \]
 lies in $L^2$.  To see this, note that $k=\dl(g)\circ\tht$, where
 $\tht(s,t,u,v)=(s-t,u-v)$.  It is clear that $\tht^*K\subset L$ and
 thus that $\tht^*K^2\subset L^2$, and we have seen that
 $\dl(g)\in K^2$ so $k\in L^2$ as claimed.  Thus, we have a map
 $\lm\:\om_{G/X}\xra{}\Prim(\Om^1_{G/X})$.  

 Next, given a function $h(u,v)$ in $J$, we have a function
 $\mu(h)(u)=h(u,0)$ in $I$.  It is clear that $\mu$ induces a map
 $\Om^1_{G/X}\xra{}\om_{G/X}$ with $\mu\circ\lm=1$.  Now suppose that
 $h$ gives an element of $\Prim(\Om^1_{G/X})$ and that $\mu(h)\in I^2$.
 Define $k(s,t,u,v)=h(s+u,t+v)-h(s,t)-h(u,v)$.  The primitivity of $h$
 means that $k\in L^2$.  Define
 $\phi\:G\tm_XG\xra{}G\tm_XG\tm_XG\tm_XG$ by $\phi(s,t)=(t,t,s-t,0)$.
 One checks that $\phi^*L\sse J$ and that 
 \[ h(s,t) = k(t,t,s-t,0) + h(t,t) + h(s-t,0). \]
 Noting that $h(t,t)=0$, we see that $h=\phi^*k+\psi^*\mu(h)$, where
 $\psi(u,v)=u-v$.  As $\mu(h)\in I^2$ and $k\in L^2$ we conclude that
 $h\in J^2$.  This means that $\mu$ is injective on
 $\Prim(\Om^1_{G/X})$.  As $\mu\lm=1$, we conclude that $\lm$ and $\mu$
 are isomorphisms.  

 Finally, we need to show that the map $f\ot\al\mapsto f\lm(\al)$
 gives an isomorphism $\OG\ot_{\OX}\om_{G/X}\xra{}\Om^1_{G/X}$.  As
 $\Om^1_{G/X}$ is freely generated over $\OG$ by $d(x)$, we must have
 $\lm(d_0(x))=u(x)d(x)$ for some power series $u$.   As $\om_{G/X}$ is
 freely generated over $\OX$ by $d_0(x)$, it will suffice to check
 that $u$ is invertible, or equivalently that $u(0)$ is a unit in
 $\OX$.  To see this, observe that $\mu(f\,d(g))=f(0)d_0(g)$, so that
 $d_0(x)=\mu\lm(d_0(x))=\mu(u(x)d(x))=u(0)d_0(x)$, so $u(0)=1$.
\end{proof}

More explicitly, let $F$ be the formal group law such that
$x(a+b)=F(x(a),x(b))$, and define $H(s)=D_2F(s,0)$, where $D_2F$ is
the partial derivative with respect to the second variable.  We
observe that $H(0)=1$, so $H$ is invertible in $\fps{R}{s}$.  We then
define $\al=H(x)^{-1}dx\in\Om^1_{G/X}$.  One can check that, in the
notation of the above proof, we have $\al=\lm(d_0(x))$, and thus that
$\al$ generates $\Prim(\Om^1_{G/X})$.  \index{O1@$\Om^1$}
\index{PrimOGX@$\Prim(\Om^1_{G/X})$}

\subsection{Heights}

\begin{proposition}\label{prop-hom-zero-diff}
 Let $G$ and $H$ be ordinary formal groups over an affine scheme $X$,
 and let $s\:G\xra{}H$ be a homomorphism.  Suppose that the induced
 map $s^*\:\om_H\xra{}\om_G$\index{oGX@$\om_{G/X}$} is zero.
 \begin{itemize}
  \item[(a)] If $X$ is a scheme over $\spec(\rat)$, then $s=0$.
  \item[(b)] If $X$ is a scheme over $\spec(\Fp)$ for some prime $p$
   then there is a unique homomorphism $s'\:F_X^*G\xra{}H$ of formal
   groups over $X$ such that $s=s'\circ F_{G/X}$.
   \index{FWX@$F_{W/X}$}
 \end{itemize}
\end{proposition}
\begin{proof}
 It follows from the definitions that our identification of
 $\om_{G/X}$ with $\Prim(\Om_{G/X})$ is natural for homomorphisms.
 Thus, if $\al\in\Prim(\Om_{H/X})$ then $s^*\al=0$.  We also know that
 $\Om_{H/X}=\OO_H\ot_{\OX}\om_{H/X}$, so any element of $\Om_{H/X}$
 can be written as $f\al$ with $f\in\Prim(\Om_{H/X})$.  Thus
 $s^*(f\al)=(f\circ s).s^*\al=0$.  Thus,
 Proposition~\ref{prop-zero-differential} applies to $s$.  If $X$ lies
 over $\spec(\rat)$ then we conclude that $s$ is constant on each
 fibre.  As it is a homomorphism, it must be the zero map.  Suppose
 instead that $X$ lies over $\spec( Fp)$.  In that case we know that
 there is a unique map $s'\:G'= F_X^*G\xra{}H$ such that
 $s=s'\circ F_{G/X}$, and we need only check that this is a
 homomorphism.  In other words, we need to check that the map
 $t'(u,v)=s'(u+v)-s'(u)-s'(v)$ (from $G'\tm_XG'$ to $H$) is zero.
 Because $s$ and $ F_{G/X}$ are homomorphisms, we see that
 $t'\circ F_{G\tm_XG/X}=0\:G\tm_XG\xra{}H$.  Using the uniqueness
 clause in Proposition~\ref{prop-zero-differential}, we conclude that
 $t'=0$ as required.
\end{proof}

\begin{corollary}\label{cor-height}
 Let $G$ and $H$ be ordinary formal groups over an affine scheme $X$,
 which lies over $\spec(\Fp)$.  Let $s\:G\xra{}H$ be a homomorphism.
 Then either $s=0$ or there is an integer $n\ge 0$ and a homomorphism
 $s'\:( F^n_X)^*G\xra{}H$ such that $s=s'\circ F^n_{G/X}$ and $(s')^*$
 is nonzero on $\om_{H/X}$. \index{oGX@$\om_{G/X}$}
 \index{FWX@$F_{W/X}$}
\end{corollary}
\begin{proof}
 Suppose that there is a largest integer $n$ (possibly $0$) such that
 $s$ can be factored in the form $s=s'\circ F^n_{G/X}$.  Write
 $G'=(F^n_X)^*G$, so that $s'\:G'\xra{}H$.  If $(s')^*=0$ on
 $\om_{H/X}$ then the proposition gives a factorisation
 $s'=s''\circ F_{G'/X}$ and thus $s=s''\circ F^{n+1}_{G/X}$
 contradicting maximality.  Thus $(s')^*\neq 0$ as claimed.  On the
 other hand, suppose that there is no largest $n$.  Choose coordinates
 $x$ and $y$ on $G$ and $H$, so there is a series $g$ such that
 $y(s(u))=g(x(u))$ for all points $u$ of $G$.  As $s$ is a
 homomorphism we have $g(0)=0$.  If $s$ factors through $F^n_{G/X}$ we
 see that $g(x)=h(x^{p^n})$ for some series $h$.  As this happens for
 arbitrarily large $n$, we see that $g$ is constant.  As $g(0)=0$ we
 conclude that $g=0$ and thus $s=0$.
\end{proof}

\begin{definition}\label{defn-height}
 Let $G$ and $H$ be ordinary formal groups over an affine scheme $X$,
 which lies over $\spec(\Fp)$.  Let $s\:G\xra{}H$ be a homomorphism.
 If $s=0$, we say that $s$ has \emph{infinite height}.  Otherwise, the
 \dfn{height} of $s$ is defined to be the integer $n$ occurring in
 Corollary~\ref{cor-height}.  The height of the group $G$ is defined
 to be the height of the endomorphism $p_G\:G\xra{}G$ (which is just
 $p$ times the identity map).
\end{definition}

\begin{definition}\label{defn-Landweber-exact}
 Let $G$ be an ordinary formal group over an affine scheme $X$.  Let
 $X_n$ be the largest closed subscheme of $X$ on which $G$ has height
 at least $n$, and write $G_n=G\tm_X X_n$.  We then have a map
 $s_n\:H_n=(F^n_X)^*G_n\xra{}G_n$ such that
 $p_{G_n}=s_n\circ F^n_{G/X}$, and thus a map
 $s_n^*\:\om_{G_n}\xra{}\om_{H_n}$ of trivialisable line bundles over
 $X_n$.  If we trivialise these line bundles then $s_n^*$ becomes an
 element $u_n\in\OO_{X_n}$, which is well-defined up to multiplication
 by a unit, and $X_{n+1}=V(u_n)=\spec(\OO_{X_n}/u_n)$.  Note also that
 $u_0=p$.  

 We say that $G$ is \dfn{Landweber exact} if for all $p$ and $n$, the
 element $u_n$ is not a zero-divisor in $\OO_{X_n}$.  Because $X_0=X$
 and $u_0=p$, this implies in particular that $\OX$ is torsion-free.
\end{definition}

\subsection{Logarithms}

\begin{definition}\label{defn-log}
 A \dfn{logarithm} for an ordinary formal group $G$ is a map of formal
 schemes $u\:G\xra{}\haf^1$ satisfying $u(g+h)=u(g)+u(h)$, or in other
 words a homomorphism $G\xra{}\AFG$.  A logarithm for a formal group
 law $F$ over a ring $R$ is a power series $f(s)\in\fps{R}{s}$ such
 that $f(F(s,t))=f(s)+f(t)\in\fps{R}{s,t}$.  Clearly, if $x$ is a
 coordinate on $G$ and $F$ is the associated formal group law then
 logarithms for $F$ biject with logarithms for $G$ by $u(g)=f(x(g))$.
 It is also clear that when $u$ is a logarithm, the differential $du$
 lies in $\om_G$.  We thus have a map
 $d\:\Hom(G,\AFG)\xra{}\aff(\om_G)$.
\end{definition}

\begin{proposition}\label{prop-rat-log}
 If $\OX$ is a $\rat$-algebra then the map
 $d\:\Hom(G,\AFG)\xra{}\aff(\om_G)$ is an isomorphism.
\end{proposition}
\begin{proof}
 If $u=f(x)$ is a logarithm and $du=f'(x)dx=0$ then $f$ is constant
 (because $\OX$ is rational so we can integrate) but $f(0)=0$ (because
 $u(0)=u(0+0)=u(0)+u(0)$) so $f=0$ so $u=0$.  Thus $d$ is injective.
 Conversely, suppose that $\al=g(x)dx\in\om_G$.  Let $f$ be the
 integral of $g$ with $f(0)=0$, so $u=f(x)\:G\xra{}\haf^1$ and
 $du=\al$.  Consider the function $w(g,h)=u(g+h)-u(g)-u(h)$, so
 $w\:G\tX G\xra{}\haf^1$ and $dw=\sg^*\al-\pi_1^*\al-\pi_2^*\al=0$.
 Thus $w$ is constant and $w(0,0)=0$ so $u(g+h)=u(g)+u(h)$ as
 required.  
\end{proof}

\begin{corollary}\label{cor-fg-rational}
 Any ordinary formal group over a scheme $X$ over $\spec(\rat)$ is
 isomorphic to the additive group $\aff^1\tm X$. \qed
\end{corollary}

\subsection{Divisors}

An ordinary formal group $G$ over $X$ is in particular a pointed
formal curve over $X$, so it makes sense to consider the schemes
$\Div_n^+(G)=G^n_X/\Sg_n$\index{Divpn@$\Div^+_n(C)$} and so on.
Moreover, Proposition~\ref{prop-Mplus-div} tells us that
$\Div^+(G)=M^+(G)$ and so on. \index{MUa@$M^+(U)$}
\index{Diva@$\Div^+(C)$} \index{divisor}

\begin{proposition}
 The formal scheme $\Div^+(G)$ has a natural structure as a
 commutative semiring object in the category $\FX_X$.
\end{proposition}
\begin{proof}
 Everywhere in this proof, products really mean fibre products over
 $X$.  

 We define a map $\nu_{i,j}\:G^i\tm G^j\xra{}G^{ij}$ by
 \[ \nu_{i,j}(a_1,\ldots,a_i,b_1,\ldots,b_j)=
      (a_1+b_1,\ldots,a_i+b_j).
 \]
 Using the fact that the colimits involved are strong, we see that
 there is a unique map
 $\mu_{i,j}\:G^i/\Sg_i\tm G^j/\Sg_j\xra{}G^{ij}/\Sg_{ij}$ that is
 compatible with the maps $\nu_{i,j}$ in the evident sense.  We can
 use the isomorphisms $\Div_i^+(G)=G^i/\Sg_i$ and
 $\Div^+(G)=\coprod_i\Div_i^+(G)$ to piece these maps together, giving
 a map $\mu\:\Div^+(G)\tm\Div^+(G)\xra{}\Div^+(G)$.  Given two
 divisors $D$ and $E$ we write $D*E=\mu(D,E)$.  The above discussion
 really just shows that the definition
 $(\sum_i[a_i])*(\sum_j[b_j])=\sum_{i,j}[a_i+b_j]$ makes sense.  It is
 easy to check (although tedious to write out in detail) that the
 operation $*$ is associative and commutative, and that the divisor
 $[0]$ is a unit for it, and that it distributes over addition.  Thus,
 $\Div^+(G)$ is a semiring object in $\FX_X$ as claimed.
\end{proof}

\begin{remark}\label{rem-lambda-semiring}
 One can also interpret and prove the statement that $\Div^+(G)$
 \index{Diva@$\Div^+(C)$} is a graded $\lm$-semiring object in
 $\FX_X$, with
 \[ \lm^k(\sum_{i=1}^n[a_i]) = 
    \sum_{i_1<\ldots<i_k} [a_{i_1}+\ldots+a_{i_k}].
 \]
\end{remark}

\begin{proposition}\label{prop-Div-ring}
 The formal scheme $\Div(G)$\index{Div@$\Div(C)$} has a natural
 structure as a commutative ring object in the category $\FX_X$.
\end{proposition}
\begin{proof}
 We know that $\Div(G)=M(G)$ is a group under addition.  It thus makes
 sense to define a map
 $\mu(n,m)\:\Div^+(G)\tm_X\Div^+(G)\xra{}\Div(G)$ by 
 \[ \mu(n,m)(D,E) = D*E - mE - nD + nm[0]. \]
 It is easy to check that 
 \[ \mu(n+i,m+j)(D+i[0],E+j[0])=\mu(n,m)(D,E). \]
 Recall that $\Div(G)=\colim_n\Div^+(G)$, where the maps in the
 diagram are of the form $D\mapsto D+i[0]$.  This is a filtered
 colimit and thus a strong one, so
 $\Div(G)\tm_X\Div(G)=\colim_{m,n}\Div^+(G)\tm_X\Div^+(G)$, where the
 maps have the form $(D,E)\mapsto(D+i[0],E+j[0])$.  It follows that
 the maps $\mu(n,m)$ fit together to give a map
 $\mu\:\Div(G)\tm_X\Div(G)\xra{}\Div(G)$.  We leave it to the reader
 to check that this product makes $\Div(G)$ into a ring object.
\end{proof}

\section{Formal schemes in algebraic topology}
\label{sec-topology}

In this section, we show how suitable cohomology theories give rise to
functors from suitable categories of spaces to formal schemes.  In
particular, the space $\cpi$ gives rise to a formal group $G$.  We
show how vector bundles over spaces give rise to divisors on $G$ over
the corresponding formal schemes, and we investigate the schemes
arising from classifying spaces of Abelian Lie groups.  We then give a
related construction that associates informal schemes to ring
spectra.  Using this we relate the Thom isomorphism to the theory of
torsors, and maps of ring spectra to homomorphisms of formal groups.

\subsection{Even periodic ring spectra}
\label{subsec-even-rings}

In this section, we define the class of cohomology theories that we
wish to study.  We would like to restrict attention to commutative
ring spectra, but unfortunately that would exclude some examples that
we really need to consider.  We therefore make the following \emph{ad
  hoc} definition, which should be ignored at first reading.
\begin{definition}\label{defn-quasi-commutative}
 Let $E$ be an associative ring spectrum, with multiplication
 $\mu\:E\Smash E\xra{}E$.  A map $Q\:E\xra{}\Sg^dE$ is a
 \dfn{derivation} if we have 
 \[ Q\circ \mu = \mu\circ(1\Smash Q + Q\Smash 1). \]
 A ring spectrum $E$ is \dfn{quasi-commutative} if there is a
 derivation $Q$ of odd degree $d$ and a central element
 $v\in\pi_{2d}E$ such that $2v=0$ and
 \[ \mu - \mu\circ\tau = v \mu \circ (Q\Smash Q). \]
 Note that if $2$ is invertible in $\pi_*E$ then $v=0$ and $E$ is
 actually commutative.
\end{definition}

\begin{definition}\label{defn-even-periodic}
 An \dfn{even periodic ring spectrum} is a quasi-commutative ring
 spectrum $E$ such that
 \begin{enumerate}
 \item $\pi_1E=0$
 \item $\pi_2E$ contains a unit.
 \end{enumerate}
 This implies that $\pi_{\odd}(E)=0$.  Thus, the derivation $Q$ in
 Definition~\ref{defn-quasi-commutative} acts trivially on $E_*$, so
 $E_*$ is a commutative ring.  Similarly, if $X$ is any space such
 that $E^1X=0$ then $E^0X$ is commutative.
\end{definition}

\begin{example}\label{eg-HP}
 The easiest example is $E^*X=H^*(X;\Zh[u^{\pm 1}])$, where we give
 $u$ degree $2$.  This is represented by the even periodic ring
 spectrum
 \[ HP = \bigWedge_{k\in\Zh} \Sg^{2k}H. \]
 \index{HP@$HP$}
\end{example}

\begin{example}\label{eg-KU}
 The next most elementary example is the complex $K$-theory spectrum
 $KU$.  This is an even periodic ring spectrum, by the Bott
 periodicity theorem.  If $p$ is a prime then we can smash this with
 the mod $p$ Moore space to get a spectrum $KU/p$.  It is true but not
 obvious that this is a ring spectrum.  It is commutative when $p>2$,
 but only quasi-commutative when $p=2$.  The derivation $Q$ in
 Definition~\ref{defn-quasi-commutative} is just the Bockstein map
 $\bt\:KU/2\xra{}\Sg KU/2$.
\end{example}

\begin{example}\label{eg-MP}
 Let $MP$ \index{MP@$MP$} be the Thom spectrum associated to the
 tautological virtual bundle over $\ZtBU$.  It is more usual to
 consider the connected component $BU=0\tm BU$ of $\ZtBU$, giving the
 Thom spectrum $MU$.  It turns out that
 $MP=\bigvee_{k\in\Zh}\Sg^{2k}MU$, and that this is an even periodic
 ring spectrum.  Moreover, a fundamental theorem of Quillen tells us
 that $MP_0=L=\OO_\FGL$\index{FGL@$\FGL$}.
\end{example}

\begin{example}\label{eg-EKMM}
 It turns out~\cite{ekmm:rma,st:pmm} that given any ring $E_0$ that
 can be obtained from $MP_0[\half]$ by inverting some elements and
 killing a regular sequence, there is a canonical even periodic ring
 spectrum $E$ with $\pi_0E=E_0$.  If we work over $MP_0$ rather than
 $MP_0[\half]$ then things are more complicated, but typically not too
 different in cases of interest, except that we only have
 quasi-commutativity rather than commutativity.  Because
 $MP_0=\OO_\FGL$, the theory of formal group laws provides us with
 many naturally defined rings $E_0$ to which we can apply this result.
\end{example}

\subsection{Schemes associated to spaces}
\label{subsec-spaces}

Let $E$ be an even periodic ring spectrum.  We write $S_E=\spec(E^0)$.
\index{SE@$S_E$}

\begin{example}\label{eg-Quillen}
 As mentioned above, Quillen's theorem tells us that
 $S_{MP}=\FGL$\index{FGL@$\FGL$}.  Less interestingly, we have
 $S_{HP}=S_K=1=\spec(\Zh)$, the terminal scheme.
\end{example}

If $Z$ is a finite complex, we write $Z_E=\spec(E^0Z)\in\CX_{S_E}$.
\index{ZE@$Z_E$} This is a covariant functor of $Z$.  If $Z$ is an
arbitrary space, we write $\Lm(Z)$\index{LZ@$\Lm(Z)$} for the category
of pairs $(W,w)$, where $W$ is a finite complex and $w$ is a homotopy
class of maps $W\xra{}Z$.
\begin{lemma}
 The category $\Lm(Z)$\index{LZ@$\Lm(Z)$} is filtered and essentially
 small.
\end{lemma}
\begin{proof}
 It is well-known that every finite CW complex is homotopy equivalent
 to a finite simplicial complex, and that there are only countably
 many isomorphism types of finite simplicial complexes.  It follows
 easily that $\Lm(Z)$ is essentially small.  If $(W,w)$ and $(V,v)$
 are objects of $\Lm(Z)$ then there is an evident map
 $u\:U=V\amalg W\xra{}Z$ whose restrictions to $V$ and $W$ are $v$ and
 $w$.  Thus $(U,u)\in\Lm(Z)$, and there are maps
 $(V,v)\xra{}(U,u)\xla{}(W,w)$ in $\Lm(Z)$.

 On the other hand, suppose we have a parallel pair of maps
 $f_0,f_1\:(V,v)\xra{}(W,w)$ in $\Lm(Z)$.  Let $U$ be the space
 $(W\amalg V\tm I)/\sim$, where $(x,t)\sim f_t(x)$ whenever $x\in V$
 and $t\in\{0,1\}$.  Let $g\:W\xra{}U$ be the evident inclusion, so
 clearly $gf_0\simeq gf_1$.  We are given that $wf_0$ and $wf_1$ are
 homotopic to $v$.  A choice of homotopy between $wf_0$ and $wf_1$
 gives a map $u\:U\xra{}X$ with $ug=w$.  Thus $g$ is a map
 $(W,w)\xra{}(U,u)$ in $\Lm(Z)$ with $gf_0=gf_1$.  This proves that
 $\Lm(Z)$ is filtered.
\end{proof}
\begin{remark}\label{rem-subcomplexes}
 Let $Z$ be a space with a given CW structure, and let $\Lm_\CW(Z)$
 be the directed set of finite subcomplexes of $Z$.  Then there is an
 evident functor $\Lm_\CW(Z)\xra{}\Lm(Z)$, which is easily seen to be
 cofinal.  We can also define $\Lm_{\text{stable}}(Z)$ to be the
 filtered category of finite spectra $W$ equipped with a map
 $w\:W\xra{}\Sgi Z_+$.  There is an evident stabilisation functor
 $\Lm(Z)\xra{}\Lm_{\text{stable}}(Z)$, and one checks that this is
 also cofinal.
\end{remark}
\begin{remark}\label{rem-Lambda-product}
 Given two spaces $Y$ and $Z$, there is a functor
 $\Lm(Y)\tm\Lm(Z)\xra{}\Lm(Y\tm Z)$ given by
 $((V,v),(W,w))\mapsto(V\tm W,v\tm w)$.  This is always cofinal, as
 one can see easily from the previous remark (for example).  
\end{remark}

\begin{definition}\label{defn-ZE}
 For any space $Z$, we write \index{ZE@$Z_E$}
 \[ Z_E=\colim_{(W,w)\in\Lm(Z)}\spec(E^0W)\in \FX_{S_E}. \]
 We also give $E^0Z$ the linear topology defined by the ideals
 $I_{(W,w)}=\ker(E^0Z\xra{w^*}E^0W)$.  Thus 
 \[ \spf(E^0Z) = \colim_{\Lm(Z)} \spec(\image(E^0Z\xra{}E^0W)). 
 \]
 We write $\hE^0Z$ for the completion of $E^0Z$.  There is an evident
 map $Z_E\xra{}\spf(E^0Z)$.  Also, if $Y$ is another space then the
 projection maps $Y\xla{}Y\tm Z\xra{}Z$ give rise to a map
 $(Y\tm Z)_E\xra{}Y_E\tm_{S_E}Z_E$.
\end{definition}
\begin{remark}\label{rem-phantoms}
 We know from~\cite{ad:vbr} that the map
 $E^0Z\xra{}\invlim_{\Lm(Z)}E^0W$ is surjective; the kernel is the
 ideal of phantom maps.  It is clear that the map
 $E^0(Z)/I_{(W,w)}\xra{}E^0W$ is injective, so the same is true of the
 map 
 \[ \invlim E^0(Z)/I_{(W,w)}\xra{}\invlim E^0W. \]
 It follows by diagram chasing that
 $\hE^0Z=\invlim E^0(Z)/I_{(W,w)}=\invlim E^0W$, and that this is a
 quotient of $E^0Z$.  {}From this we see that $E^0Z$ is complete if and
 only if there are no phantom maps $Z\xra{}E$.
\end{remark}

\begin{definition}\label{defn-tolerable}
 We say that $Z$ is \dfn{tolerable} (relative to $E$) if
 $Z_E=\spf(E^0Z)$ and $(Y\tm Z)_E=Y_E\tm_{S_E}Z_E$ for all finite
 complexes $Y$.
\end{definition}

\begin{proposition}\label{prop-tolerable-product}
 If $Z$ is \idx{tolerable} and $Y$ is arbitrary then
 \[ (Y\tm Z)_E=Y_E\tm_{S_E}Z_E. \]
 If $Y$ is also tolerable then so is $Y\tm Z$, and
 $\hE^0(Y\tm Z)=\hE^0(Y)\hot_{E^0}\hE^0(Z)$.  Of course if $E^0Y$,
 $E^0Z$ and $E^0(Y\tm Z)$ are complete this means that
 $E^0(Y\tm Z)=E^0Y\hot_{E^0}E^0Z$.
\end{proposition}
\begin{proof}
 If we fix $V\in\Lm(Y)$ then the functor from $\Lm(Z)$ to
 $\Lm(V\tm Z)$ given by $W\mapsto V\tm W$ is clearly cofinal, so
 $\colim_W(V\tm W)_E=(V\tm Z)_E$, and this is the same as
 $V_E\tm_{S_E}Z_E$ because $Z$ is tolerable and $V$ is finite.  If we
 now take the colimit over $V$ and use the fact that filtered colimits
 of formal schemes commute with finite limits, we find that
 $\colim_{V,W}(V\tm W)_E=Y_E\tm_{S_E}Z_E$.  It follows from
 Remark~\ref{rem-Lambda-product} that
 $(Y\tm Z)_E=\colim_{V,W}(V\tm W)_E$, so the first claim follows.

 Now suppose that $Y$ is tolerable.  Then
 \begin{align*}
  (Y\tm Z)_E      &= Y_E\tm_{S_E}Z_E                    \\
                  &= \spf(E^0Y)\tm_{S_E}\spf(E^0Z)      \\
                  &= \spf(\hE^0Y)\tm_{S_E}\spf(\hE^0Z)  \\
                  &= \spf(\hE^0Y\hot_{E^0}\hE^0Z).
 \end{align*}
 It follows that
 $\hE^0(Y\tm Z)=\OO_{(Y\tm Z)_E}=\hE^0Y\hot_{E^0}\hE^0Z$ as claimed.
 It also follows that $(Y\tm Z)_E$ is solid, and thus that
 $(Y\tm Z)_E=\spf(E^0(Y\tm Z))$.

 Now let $X$ be a finite complex.  We need to show that
 $(X\tm Y\tm Z)_E=X_E\tm_{S_E}(Y\tm Z)_E=X_E\tm_{S_E}Y_E\tm_{S_E}Z_E$.
 In fact, we have $(X\tm Y)_E=X_E\tm_{S_E}Y_E$ because $Y$ is
 tolerable, and $((X\tm Y)\tm Z)_E=(X\tm Y)_E\tm_{S_E}Z_E$ because $Z$
 is tolerable, and the claim follows.
\end{proof}

\begin{definition}\label{defn-decent}
 A space $Z$ is \dfn{decent} if $H_*Z$ is a free Abelian group,
 concentrated in even degrees.
\end{definition}
\begin{example}
 The spaces $\cpi$, $BU(n)$, $\ZtBU$, $BSU$ and $\Om S^{2n+1}$ are all
 decent.
\end{example}

\begin{proposition}\label{prop-even-cells}
 Let $Z$ be a decent space.  Then $Z$ is \idx{tolerable} for any $E$,
 and $Z_E$ is coalgebraic over $S_E$.  Moreover, for any map
 $E\xra{}E'$ of even periodic ring spectra, the resulting diagram
 \begin{diag}
  \node{Z_{E'}} \arrow{s} \arrow{e} \node{Z_E} \arrow{s} \\
  \node{S_{E'}} \arrow{e} \node{S_E}
 \end{diag}
 is a pullback.
\end{proposition}
\begin{proof}
 We may assume that $Z$ is connected (otherwise treat each component
 separately).  As $H_1Z=0$ we see that $\pi_1Z$ is perfect, so we can
 use Quillen's plus construction to get a homology equivalence
 $Z\xra{}Z^+$ such that $\pi_1(Z^+)=0$.  By the stable Whitehead
 theorem, this map is a stable equivalence, so
 $E^0(Y\tm Z^+)=E^0(Y\tm Z)$ for all $Y$.  We may thus replace $Z$ by
 $Z^+$ and assume that $\pi_1Z=0$.  This step is not strictly
 necessary, but it seems the cleanest way to avoid trouble from the
 fundamental group.  Given this, it is well-known that $Z$ has a CW
 structure in which all the cells have even dimension.  It follows
 that the Atiyah-Hirzebruch spectral sequence collapses and that
 $E_*Z$ is a free module over $E_*$, with one generator $e_i$ for each
 cell.  As $E_*$ is two-periodic, we can choose these generators in
 degree zero.  Similarly, $E_*(Z\tm Z)$ is free on generators
 $e_i\ot e_j$ and thus is isomorphic to $E_*(Z)\ot_{E_*}E_*(Z)$, so we
 can use the diagonal map to make $E_*Z$ into a coalgebra over $E_*$.
 By periodicity, $E_0(Z\tm Z)=E_0(Z)\ot_{E_0}E_0(Z)$ and $E_0Z$ is a
 coalgebra over $E_0$, and is freely generated as an $E_0$-module by
 the $e_i$.

 If $W$ is a finite subcomplex of $Z$, it is easy to see that $E_0W$
 is a standard subcoalgebra of $E_0Z$ (in the language of
 Definition~\ref{defn-good-basis}).  Moreover, any finite collection
 of cells lies in a finite subcomplex, so it follows that any finitely
 generated submodule of $E_0Z$ lies in a standard subcoalgebra.  It
 follows that $\{e_i\}$ is a good basis for $E_0Z$, so that
 $E_0Z\in\CC'_{S_E}$.  

 It follows from the above in the usual way that $E\Smash Z_+$ is
 equivalent as an $E$-module spectrum to a wedge of copies of $E$ (one
 for each cell), and thus that $E^*Z=\Hom_{E^*}(E_*Z,E_*)$.  Using the
 periodicity we conclude that $E^0Z=\Hom_{E^0}(E_0Z,E^0)$.  It follows
 that $\spf(E^0Z)=\sch_{S_E}(E_0Z)$ is a solid formal scheme, which is
 coalgebraic over $S_E$.  It is also easy to check that $\spf(E^0Z)$
 is the colimit of the schemes $\spec(E^0W)$ as $W$ runs over the
 finite subcomplexes.  It follows from Remark~\ref{rem-subcomplexes}
 that $\spf(E^0Z)=Z_E$.

 Now let $Y$ be another space.  Let $W$ be a finite subcomplex of $Z$,
 and let $(V,v)$ be an object of $\Lm(Y)$.  The usual K\"unneth
 arguments show that $E^0(W\tm V)=E^0W\ot_{E^0}E^0V$, and thus that
 $(W\tm V)_E=W_E\tm_{S_E}V_E$.  Using Remark~\ref{rem-Lambda-product}
 we conclude that 
 \[ (Z\tm Y)_E=\colim_{W,V}W_E\tm_{S_E}V_E=
    (\colim_W W_E)\tm_{S_E}(\colim_V V_E) = Z_E\tm_{S_E} Y_E.
 \]
 This proves that $Z$ is tolerable.  We leave it to the reader to
 check that a map $E\xra{}E'$ gives an isomorphism
 $Z_{E'}=Z_E\tm_{S_E}S_{E'}$.  
\end{proof}
\begin{example}
 It follows from the proposition that the spaces $\cpi$, $BU(n)$,
 $\ZtBU$, $BSU$ and $\Om S^{2n+1}$ are all \idx{tolerable}, and the
 corresponding schemes are \idx{coalgebraic} over $S_E$.  The case of
 $\cpi$ is particularly important.  We note that $\cpi=BS^1=K(\Zh,2)$
 is an Abelian group object in the homotopy category, so $G_E=\cpi_E$
 \index{GE@$G_E$} is an Abelian formal group over $S_E$.  Because
 $H^*\cpi=\fps{\Zh}{x}$, the Atiyah-Hirzebruch spectral sequence tells
 us that $E^0\cpi=\fps{E^0}{x}$ (although this does not give a
 canonical choice of generator $x$).  This means that
 $G_E\simeq\haf^1\tm S_E$ in $\Based\FX_{S_E}$, so that $G_E$ is an
 ordinary formal group.\index{formal group!ordinary}
\end{example}

We next recall that for $n>0$ there is a quasicommutative rings
spectrum $P(n)=BP/I_n$ with $P(n)^*=\Fp[v_k\st k\ge 0]$, where $v_k$
has degree $-2(p^k-1)$.  The cleanest construction now available is
given in~\cite{ekmm:rma,st:pmm}, although of course there are much
older constructions using Baas-Sullivan theory.  We also have
$P(0)=BP$, with $P(0)^*=\Zh_{(p)}[v_k\st k>0]$.

\begin{definition}\label{defn-exact-module}
 Let $E$ be an even periodic ring spectrum.  We say that $E$ is an
 \emph{exact $P(n)$-module}
 \index{exact P(n) module@exact $P(n)$-module} (for some $n\ge 0$) if
 it is a module-spectrum over $P(n)$, and the sequence
 $(v_n,v_{n+1},\ldots)$ is regular on $E_*$.  
\end{definition}
\begin{proposition}\label{prop-rawiya-tolerable}
 Let $E$ be an exact $P(n)$-module.  Let $Z$ be a CW complex of finite
 type such that $K(m)_*Z$ is concentrated in even degrees for
 infinitely many $m$.  If $n=0$, assume that $H^s(Z;\rat)=0$ for
 $s\gg 0$.  Then $Z$ is tolerable for $E$.
\end{proposition}
\begin{remark}
 When combined with Proposition~\ref{prop-tolerable-product} this
 gives a useful K\"{u}nneth theorem.
\end{remark}
The proof will follow after Corollary~\ref{cor-E-kunneth}.  Many
spaces are known to which this applies: simply connected finite
Postnikov towers of finite type, classifying spaces of many finite
groups and compact Lie groups, the spaces $QS^{2m}$, $BO$, $\ImJ$ and
$BU\langle 2m\rangle$ for example.  See~\cite{rawiya:bpc} for more
details.  The proof of our proposition will also rely heavily on the
results of that paper.

We next need some results involving the pro-completion of the category
of graded Abelian groups, which we denote by $\Pro(\Ab_*)$.  It is
necessary to distinguish this carefully from the category
$\Pro(\Ab)_*$ of graded systems of pro-groups.  A tower of graded
groups can be regarded as an object in either category, but the
morphisms are different.  A tower $\{A_{0*}\xla{}A_{1*}\xla{}\cdots\}$
in $\Pro(\Ab_*)$ is pro-trivial if for all $j$, there exists $k>j$
such that the map $A_{k*}\xra{}A_{j*}$ is zero.  It is pro-trivial in
$\Pro(\Ab)_*$ if for all $j$ and $d$ there exists $k$ such that the
map $A_{kd}\xra{}A_{jd}$ is zero.  Because $k$ is allowed to depend on
$d$, this is a much weaker condition than triviality in $\Pro(\Ab_*)$.
Note also that if $R_*\xra{}R'_*$ is a map of graded rings, and
$\{M_{\al*}\}$ is a pro-system of $R_*$-modules that is trivial in
$\Pro(\Ab_*)$, then the same is true of $R'_*\ot_{R_*}M_*$.  However,
the corresponding statement for $\Pro(\Ab)_*$ is false.

\begin{remark}\label{rem-yagita}
 If $E$ is an exact $P(n)$-module, we know from work~\cite{ya:efb} of
 Yagita that the functor $M\mapsto E^*\ot_{P(n)^*}M$ is
 an exact functor on the category of $P(n)^*P(n)$-modules that are
 finitely presented as modules over $P(n)^*$.  (This category is
 Abelian, because the ring $P(n)^*$ is coherent.)  It follows that
 $E^*Z=E^*\ot_{P(n)^*}P(n)^*Z$ for all finite complexes $Z$.
\end{remark}

The following lemma is largely a paraphrase of results
in~\cite{rawiya:bpc}.  
\begin{lemma}\label{lem-rawiya}
 Fix $n\ge 0$.  Suppose that $Z$ is a CW complex of finite type, and
 write $Z^r$ for the $r$-skeleton of $Z$.  If $n=0$ we also assume
 that $H^s(Z;\rat)=0$ for $s\gg 0$.  Let
 $F^{r+1}=\ker(P(n)^*Z\xra{}P(n)^*Z^r)$ denote the $(r+1)$'st
 Atiyah-Hirzebruch filtration in $P(n)^*Z$.  Then the tower
 $\{P(n)^*Z^r\}_{r\ge 0}$ is isomorphic to
 $\{P(n)^*(Z)/F^{r+1}\}_{r\ge 0}$ in $\Pro(\Ab_*)$, and thus is
 Mittag-Leffler.  Moreover, the groups $P(n)^*(Z)/F^{r+1}$ are
 finitely presented modules over $P(n)^*$, and their inverse limit is
 $P(n)^*Z$.
\end{lemma}
\begin{proof}
 Write $P=P(n)$ for brevity.  Write $A_r=P^*Z^r$ and
 \[ B_r=P^*(Z)/F^{r+1}=\image(P^*Z\xra{}A_r). \]
 We then have an inclusion of towers $\{B_r\}\xra{}\{A_r\}$, for which
 we need to provide an inverse in the $\Pro$-category.  We claim that
 for each $r$, there exists $m(r)>r$ such that the image of the map
 $A_{m(r)}\xra{}A_r$ is precisely $B_r$.  We will deduce the lemma
 from this before proving it.  Define $m_0=0$ and
 $m_{k+1}=m(m_k)>m_k$.  By construction, the map
 $A_{m_{k+1}}\xra{}A_{m_k}$ factors through $B_{m_k}\sse A_{m_k}$.
 One checks that the resulting maps $A_{m_{k+1}}\xra{}B_{m_k}$ are
 compatible as $k$ varies, and that they provide the required inverse.
 We also know that $P^*$ is a coherent ring, so the category of
 finitely presented modules is Abelian and closed under extensions.
 It follows in the usual way that $A_r$ is finitely presented for all
 $r$, and thus that $B_r=\image(A_{m(r)}\xra{}A_r)$ is finitely
 presented.

 We now need to show that $m(r)$ exists.  By the basic setup of the
 Atiyah-Hirzebruch spectral sequence, it suffices to show that for
 large $m$, the first $r+1$ columns in the spectral sequence for
 $P^*Z^m$ are the same as in the spectral sequence for $P^*Z$.
 This is Lemma~4.4 of~\cite{rawiya:bpc}.  (When $n=0$, we need to
 check that we are in the case $P(0)=BP$ of their Definition~1.5.
 This follows from our assumption that $H^s(Z;\rat)=0$ for $s\gg 0$.)

 Finally, we need to show that $P^*Z=\invlim_rP^*(Z)/F^{r+1}$.
 This is essentially~\cite[Corollary 4.8]{rawiya:bpc}.
\end{proof}

\begin{corollary}\label{cor-rawiya}
 Let $Z$ and $n$ be as in the Lemma, and let $E$ be an exact
 $P(n)$-module.  Then $E^0Z$ is complete, and $Z_E=\spf(E^0Z)$, and
 $E^*Z=E^*\hot_{P(n)^*}P(n)^*Z$.  Moreover we have isomorphisms 
 \begin{align*}
    \{E^*(Z^r)\}\simeq\{E^*(Z)/F^{r+1}\} &\simeq
    \{E^*\ot_{P(n)^*}P(n)^*Z^r\} \\
    & \simeq
    \{E^*\ot_{P(n)^*}(P(n)^*(Z)/F^{r+1})\}
 \end{align*}
 in $\Pro(\Ab_*)$.
\end{corollary}
\begin{proof}
 We reuse the notation of the previous proof.  We also define
 $A'_r=E^*\ot_{P^*}A_r$ and $B'_r=E^*\ot_{P^*}B_r$.  As $Z^r$ is
 finite we see that $A'_r=E^*Z^r$.  Next recall that for any $r$ we
 can choose $m>r$ such that $B_r=\image(A_m\xra{}A_r)$.  As the
 functor $E^*\ot_{P^*}(-)$ is exact on finitely presented comodules,
 we see that $B'_r$ is the image of the map $A'_m\xra{}A'_r$ and in
 particular that the map $B'_r\xra{}A'_r$ is injective.  Next, the map
 $E^*\ot_{P^*}P^*Z\xra{}E^*\ot_{P^*}P^*Z^m=E^*Z^m=A'_m$ clearly
 factors through $E^*Z$, so our epimorphism $P^*Z\xra{}A_m\xra{}B_r$
 gives an epimorphism $E^*\ot_{P^*}P^*Z\xra{}A'_m\xra{}B'_r$ which
 factors through $E^*Z$, so the map $E^*Z\xra{}B'_r$ is surjective.
 Thus $B'_r=\image(E^*Z\xra{}E^*Z^r)=E^*(Z)/F^{r+1}$.  We can now
 apply the functor $E^*\ot_{P^*}(-)$ to the pro-isomorphisms in the
 Lemma to get the pro-isomorphisms in the present corollary.  This
 makes it clear that the tower $\{E^*Z^r\}$ is Mittag-Leffler so the
 Milnor sequence tells us that 
 \[ E^*Z=\invlim_r E^*Z^r=\invlim_rE^*(Z)/F^{r+1}. \]
 This means in particular that $E^0Z$ is complete with respect to the
 linear topology generated by the ideals $F^{r+1}$, which is easily
 seen to be the same as the topology in Definition~\ref{defn-ZE}.
 Moreover, we have an isomorphism $\{A'_r\}\simeq\{B'_r\}$
 in the Pro category of groups, and it is easy to see from the
 construction that this is actually an isomorphism in the Pro category
 of rings as well, so by applying $\spec(-)$ we get an isomorphism in
 the Ind category of schemes, which is just the category of formal
 schemes.  {}From the definitions we have $Z_E=\colim_r\spec(A'_r)$ and
 $\spf(E^0Z)=\colim_r\spec(B'_r)$, so we conclude that
 $Z_E=\spf(E^0Z)$. 
\end{proof}
 
\begin{lemma}\label{lem-rawiya-FZE}
 Let $E$ and $Z$ be as in Corollary~\ref{cor-rawiya}, and suppose that
 $K(m)^*Z$ is concentrated in even degrees for infinitely many $m$.
 Then the ring $E^*Z$ is Landweber exact over $P(n)^*$, so the
 function spectrum $F(Z_+,E)$ is an exact $P(n)$-module.
\end{lemma}
\begin{proof}
 We know from~\cite[Lemma 5.3]{rawiya:bpc} that $P(m)^*Z$ is
 concentrated in even degrees for all $m$, and
 from~\cite[Corollary 4.6]{rawiya:bpc} that the tower $\{P(m)^*Z^r\}$
 has the Mittag-Leffler property.  It follows that the tower
 $\{P(m)^\odd Z^r\}$ is pro-trivial.  Next, consider the cofibration
 $\Sg^{2(p^m-1)}P(m)\xra{v_m}P(m)\xra{}P(m+1)\xra{}\Sg^{2p^m-1}P(m)$.
 This gives a pro-exact sequence of towers
 \[ 0 \xra{} \{P(m)^\ev Z^r\} \xra{v_m} \{P(m)^\ev Z^r\} \xra{} 
    \{P(m+1)^\ev Z^r\} \xra{} 0.
 \]
 It follows that the sequence $(v_n,v_{n+1},\ldots)$ acts regularly on
 the tower $\{P(n)^*Z^r\}$.  Next, for any spectrum $X$ we have a
 map $P(m)\Smash X\xra{}P(m)\Smash BP\Smash X$ which makes $P(m)_*X$ a
 comodule over $P(m)_*BP=BP_*BP/I_n$.  Moreover, we have
 $P(m)_*X\ot_{P(m)_*}P(m)_*BP=P(m)_*X\ot_{BP_*}BP_*BP$ so this
 actually makes $P(m)_*X$ into a comodule over $BP_*BP$.  One can
 check from this construction that the maps
 $\Sg^{2(p^m-1)}P(m)\xra{v_m}P(m)\xra{}P(m+1)\xra{}\Sg^{2p^m-1}P(m)$
 give rise to maps of comodules, so our whole diagram of towers is a
 diagram of finitely-presented comodules over $P(n)_*BP$.  The functor
 $E^*\ot_{P(n)^*}(-)$ is exact on this category.  It is easy to
 conclude by induction that
 $\{E^*\ot_{P(n)^*}P(m)^*Z^r\}\simeq\{E^*(Z^r)/I_m\}$, that the odd
 dimensional part of these towers is pro-trivial, the towers are
 Mittag-Leffler, and the sequence $(v_n,v_{n+1},\ldots)$ is regular on
 the tower $\{E^*(Z^r)\}$.  We can now pass to the inverse limit
 (using the Mittag-Leffler property to show that the $\invlim^1$ terms
 vanish) to see that the sequence $(v_n,v_{n+1},\ldots)$ is regular on
 $E^*(Z)$.
\end{proof}

Our next few results are closely related to those 
of~\cite[Section 9]{rawiya:bpc}, although a precise statement of the
relationship would be technical.
\begin{lemma}\label{lem-Pn-kunneth}
 Let $Z$ be a CW complex of finite type such that $K(m)^*Z$ is
 concentrated in even degrees for infinitely many $m$.  If $n=0$ we
 also assume that $H^s(Z;\rat)=0$ for $s\gg 0$.  Then for any finite
 spectrum $W$ we have pro-isomorphisms
 \begin{align*}
    \{P(n)^*(Z^r\tm W)\} & \simeq
    \{P(n)^*Z^r\ot_{P(n)^*}P(n)^*W\} \\
    &\simeq
    \{P(n)^*(Z)/F^{r+1}\ot_{P(n)^*}P(n)^*W\},
 \end{align*}
 and these towers are Mittag-Leffler.  Moreover, we have isomorphisms
 \[ P(n)^*(Z\tm W) =
    P(n)^*Z\ot_{P(n)^*}P(n)^*W = 
    P(n)^*Z\hot_{P(n)^*}P(n)^*W.
 \]
\end{lemma}
\begin{proof}
 Write $P=P(n)$ for brevity.  The usual Landweber exactness argument
 shows that $W\mapsto P^*Z\ot_{P^*}P^*W$ is a cohomology theory and
 thus that it coincides with $P^*(Z\tm W)$.  We can also do the same
 argument with pro-groups.  We saw in the proof of the previous lemma
 that the sequence $(v_n,v_{n+1},\ldots)$ acts regularly on the
 pro-group $\{P^*Z^r\}$, so the pro-group
 $\{\Tor^{P^*}_1(P(m)^*,P^*Z^r)\}$ is trivial for all $m\ge n$.  Any
 finitely presented comodule $M^*$ has a finite Landweber filtration
 whose quotients have the form $P(m)^*$ for $m\ge n$, and we see by
 induction on the length of the filtration that
 $\{\Tor^{P^*}_1(M^*,P^*Z^r)\}$ is trivial.  This implies that the
 construction $M^*\mapsto\{M^*\ot_{P^*}P^*Z^r\}$ gives an exact
 functor from finitely presented comodules to $\Pro(\Ab_*)$, so that
 $W\mapsto\{P^*W\ot_{P^*}P^*Z^r\}$ is a $\Pro(\Ab_*)$-valued cohomology
 theory on finite complexes.  The construction
 $W\mapsto\{P^*(W\tm Z^r)\}$ gives another such cohomology theory, and
 we have a natural transformation from the first to the second that is
 an isomorphism when $W$ is a sphere, so it is an isomorphism in
 general.  Thus $\{P^*(Z^r\tm W)\}=\{P^*Z^r\ot_{P^*}P^*W\}$, as
 claimed.  We have seen that the tower $\{P^*Z^r\}$ is pro-isomorphic
 to $\{P^*(Z)/F^{r+1}\}$, so it follows that
 $\{P^*Z^r\ot_{P^*}P^*W\}\simeq\{P^*(Z)/F^{r+1}\ot_{P^*}P^*W\}$.  The
 second of these is a tower of isomorphisms, so all three of our
 towers are Mittag-Leffler as claimed.  As $Z\tm W$ is the homotopy
 colimit of the spaces $Z^r\tm W$, the Milnor sequence gives an
 isomorphism $P^*(Z\tm W)=\invlim_rP^*(Z)/F^{r+1}\ot_{P^*}P^*W$, and
 the right hand side is by definition $P^*(Z)\hot_{P^*}P^*W$, which
 completes the proof.
\end{proof}

\begin{corollary}\label{cor-E-kunneth}
 Let $E$ be an exact $P(n)$-module.  Let $Z$ be a CW complex of finite
 type such that $K(m)^*Z$ is concentrated in even degrees for
 infinitely many $m$.  If $n=0$ we also assume that $H^s(Z;\rat)=0$
 for $s\gg 0$.  Then for any finite spectrum $W$ we have
 pro-isomorphisms
 \[ \{E^0(Z^r\tm W)\}\simeq
    \{E^0Z^r\ot_{E^0}E^0W\}\simeq
    \{E^0(Z)/F^{r+1}\ot_{E^0}E^0W\},
 \]
 and these towers are Mittag-Leffler.  Moreover, we have isomorphisms
 \[ E^0(Z\tm W) =
    E^0Z\ot_{E^0}E^0W = 
    E^0Z\hot_{E^0}E^0W.
 \]
\end{corollary}
\begin{proof}
 If we apply the functor $E^*\ot_{P(n)^*}(-)$ to the pro-isomorphisms
 in the lemma, we get the pro-isomorphisms in the corollary.  We
 deduce in the same way as in the lemma that
 $E^0(Z\tm W)=E^0Z\hot_{E^0}E^0W$.  On the other hand, we see from
 Lemma~\ref{lem-rawiya-FZE} that 
 \begin{multline*}
   E^*(Z\tm W)=[W_+,F(Z_+,E)]^*=E^*Z\ot_{P(n)^*}P(n)^*W= \\
    E^*Z\ot_{E^*}(E^*\ot_{P(n)^*}P(n)^*W)=E^*Z\ot_{E^*}E^*W.  
 \end{multline*}
 Thus $E^0(Z\tm W)=E^0Z\ot_{E^0}E^0W$ as claimed.
\end{proof}

\begin{proof}[Proof of Proposition~\ref{prop-rawiya-tolerable}]
 Corollary~\ref{cor-rawiya} shows that $Z_E=\spf(E^0Z)$.  Write
 \[ F^{r+1}=\ker(E^0Z\xra{}E^0Z^r), \]
 so $Z_E=\colim_rV(F^r)$.  Let $W$ be a finite complex.  We then have
 \begin{align*}
  Z_E \tm_{S_E} W_E
    &= \colim_r V(F^r)\tm_{S_E} W_E                     \\
    &= \colim_r \spec(E^0(Z)/F^r\ot_{E^0}E^0(W))        \\
    &= \spf(E^0(Z)\hot_{E^0}E^0(W))                     \\
    &= \spf(E^0(Z\tm W)),
 \end{align*}
 where we have used Corollary~\ref{cor-E-kunneth}.  We can apply
 Lemma~\ref{lem-rawiya} to $Y\tm Z$ and conclude that
 $\spf(E^0(Y\tm Z))=(Y\tm Z)_E$, giving the required isomorphism
 $(Y\tm Z)_E=Y_E\tm_{S_E}Z_E$.
\end{proof}

\subsection{Vector bundles and divisors}
\label{subsec-bundles-divisors}

Let $V$ be a complex vector bundle of rank $n$ over a \idx{tolerable}
space $Z$.  We write $P(V)$\index{PV@$P(V)$} for the space of pairs
$(z,W)$, where $z\in Z$ and $W$ is a line (i.e.\ a one-dimensional
subspace) in $V_z$.  This is clearly a fibre bundle over $Z$ with
fibres $\cp^{n-1}$.  We write $\BD(V)=P(V)_E$\index{DV@$\BD(V)$}.
There is a tautological line bundle $L$ over $P(V)$, whose fibre over
a pair $(z,W)$ is $W$.  This is classified by a map $P(V)\xra{}\cpi$.
By combining this with the projection to $Z$, we get a map
$P(V)\xra{}\cpi\tm Z$ and thus a map $\BD(V)\xra{}G\tm_SZ_E$.  The
well-known theorem on projective bundles translates into our language
as follows.

\begin{proposition}\label{prop-DV-divisor}
 The above map is a closed inclusion, making
 $\BD(V)$\index{DV@$\BD(V)$} into an effective \idx{divisor} of degree
 $n$ on $G$.
\end{proposition}
\begin{proof}
 Choose an orientation $x$ of $E$, so $x\in\tilde{E}^0\cpi$.  We also
 write $x$ for the image of $x$ under the map $P(V)\xra{}\cpi$, which
 is just the Euler class of $L$.  We claim that $E^*P(V)$ is freely
 generated over $E^*Z$ by $\{1,x,\ldots,x^{n-1}\}$, which will prove
 the claim.  This is clear when $V$ is trivialisable.  For the general
 case, we may assume that $Z$ is a regular CW complex.  The claim
 holds when $Z$ is a finite union of subcomplexes on which $V$ is
 trivialisable, by a well-known Mayer-Vietoris argument.  It thus
 holds when $Z$ is a finite complex, and the general case follows by
 passing to colimits.  
\end{proof}

\begin{proposition}\label{prop-whitney-sum}
 If $V$ and $W$ are two vector bundles over a \idx{tolerable} space
 $Z$ then $\BD(V\oplus W)=\BD(V)+\BD(W)$.\index{DV@$\BD(V)$}
\end{proposition}
\begin{proof}
 Choose an orientation, and let $x$ be the Euler class of the
 usual line bundle over $P(V\oplus W)$.  The polynomial
 $f_{\BD(V\oplus W)}(t)$ is the unique one of degree
 $\dim(V\oplus W)$ of which $x$ is a root, so it suffices to check
 that $f_{\BD(V)}(x)f_{\BD(W)}(x)=0$.  There are evident inclusions
 $P(V)\xra{}P(V\oplus W)\xla{}P(W)$ with $P(V)\cap P(W)=\emptyset$.
 Write $A=P(V\oplus W)\setminus P(V)$ and
 $B=P(V\oplus W)\setminus P(W)$, so that $A\cup B=P(V\oplus W)$.  By a
 well-known argument, if $a,b\in E^0P(V\cup W)$ and $a|_A=0$ and
 $b|_B=0$ then $ab=0$, so it suffices to check that
 $f_{\BD(V)}(x)|_B=0$ and $f_{\BD(W)}(x)|_A=0$.  It is not hard to see
 that the inclusions $P(V)\xra{}B$ is a homotopy equivalence and thus
 that $f_{\BD(V)}(x)|_B=0$, and the other equation is proved
 similarly. 
\end{proof}

\begin{proposition}\label{prop-div-line-bundle}
 If $M$ is a complex line bundle over a tolerable space $Z$, which is
 classified by a map $u\:Z\xra{}\cpi$, then $\BD(M)$ is the image of
 the map $(u,1)_E\:Z_E\xra{}(\cpi\tm Z)_E=G\tm_SZ_E$.
\end{proposition}
\begin{proof}
 This follows from the definitions, using the obvious fact that
 $P(M)=Z$. 
\end{proof}

\begin{proposition}\label{prop-BUn}
 There is a natural isomorphism
 $BU(n)_E=\Div_n^+(G)$. \index{Divpn@$\Div^+_n(C)$} 
\end{proposition}
\begin{proof}
 This is essentially well-known, but we give some details to
 illustrate how everything fits together.  Let $T(n)$ be the maximal
 torus in $U(n)$, so that $BT(n)\simeq(\cpi)^n$ and $BT(n)_E=G^n_S$.
 Thus, the inclusion $i\:T(n)\xra{}U(n)$ gives a map
 $G^n_S\xra{}BU(n)_E$.  If $\sg\in\Sg_n$ is a permutation, then the
 evident action of $\sg$ on $T(n)$ is compatible with the action on
 $U(n)$ given by conjugating with the associated permutation matrix.
 This matrix can be joined to the identity matrix by a path in $U(n)$,
 so the conjugation is homotopic to the identity.  Thus, our map
 $G^n_S\xra{}BU(n)_E$ factors through a map
 $\Div_n^+(G)=G^n_S/\Sg_n\xra{}BU(n)_E$.  On the other hand, the
 tautological bundle $V_n$ over $BU(n)$ gives rise to a divisor
 $\BD(V_n)$ over $BU(n)_E$ and thus a map $BU(n)_E\xra{}\Div_n^+(G)$.
 The composite
 $G^n_S=BT(n)_E\xra{}BU(n)_E\xra{}\Div_n^+(G)=G^n_S/\Sg_n$ classifies
 the divisor $\BD(i^*V_n)$.  Let $M_1,\ldots,M_n$ be the evident line
 bundles over $BT(n)$, so that $i^*V_n=M_1\oplus\ldots\oplus M_n$.
 One checks from this and Propositions~\ref{prop-whitney-sum}
 and~\ref{prop-div-line-bundle} that the composite is just the usual
 quotient map $G^n_S\xra{}G^n_S/\Sg_n$, and thus the composite
 $\Div_n^+(G)\xra{}BU(n)_E\xra{}\Div_n^+(G)$ is the identity.

 Next, we take the space of $n$-frames in $\cplx^\infty$ as our model
 for $EU(n)$.  There is then a homeomorphism
 $EU(n)/(S^1\tm U(n-1))\xra{}P(V_n)$ (sending $(w_1,\ldots,w_n)$ to
 the pair $(L,W)$, where $W$ is the span of $\{w_1,\ldots,w_n\}$ and
 $L$ is the span of $w_1$).  The left hand side is a model for
 $\cpi\tm BU(n-1)$.  By induction on $n$, we may assume that
 $BU(n-1)_E=G^{n-1}/\Sg_{n-1}$.  This gives a commutative diagram as
 follows. 
 \begin{diag}
  \node{G\tm G^{n-1}/\Sg_{n-1}} \arrow{e,t}{\simeq} \arrow{s}
  \node{P(V_n)_E} \arrow{s,A} \\
  \node{G^n/\Sg_n} \arrow{e,V} \node{BU(n)_E}
 \end{diag}
 The top horizontal is an isomorphism by induction and the right hand
 vertical is faithfully flat, and thus a categorical epimorphism.  It
 follows that the bottom map is an epimorphism, but we have already
 seen that it is a split monomorphism, so it is an isomorphism as
 required. 
\end{proof}

\begin{definition}\label{defn-chern}
 Let $x$ be a coordinate on $G$.  If $V$ is a vector bundle of rank
 $n$ over a tolerable space $Z$, then we have
 $\BD(V)=\spf(\fps{E^0Z}{x}/f(x))$ \index{DV@$\BD(V)$} for a unique
 monic polynomial $f(x)=\sum_{i=0}^nc_i(V)x^{n-i}$, with
 $c_i(V)\in E^0Z$.  We call $c_i(V)$ the $i$'th \dfn{Chern class} of
 $V$.
\end{definition}

\begin{definition}\label{defn-thom-sheaf-top}
 We write $\BL(V)$\index{LV@$\BL(V)$} for $L(\BD(V))$, the
 \idx{Thom sheaf} of $\BD(V)$, which is a line bundle over $Z_E$.  It
 is easy to see that $\BL(V)=\tilde{E}^0Z^V$, where
 $Z^V=P(\cplx\oplus V)/P(V)$ is the \idx{Thom space} of $V$.
\end{definition}

\begin{remark}\label{rem-Lin-residue}
 Let $E$ be an even periodic ring spectrum and put $G=G_E=(\cpi)_E$
 and $S=S_E=\spec(E^0)$ as usual.  Then the Thom spectra
 $\cp^\infty_{-n}$ form a tower, and there is a natural identification
 $\CM_{G/S}=\colim_nE^0(\cp^\infty_{-n})$.  We also have
 $\om_{G/S}=\widetilde{E}^0\cp^1=\widetilde{E}^0S^2=\pi_2E$.  The
 theory of invariant differentials identifies $\CM\Om^1_{G/S}$ with
 $\CM_{G/S}\ot_{E^0}\om_{G/S}=\colim_nE^0(\Sg^2\cp^\infty_{-n})$.  The
 $S^1$-equivariant Segal conjecture gives an equivalence between
 $\holim_n\Sg^2\cp^\infty_{-n}$ and the profinite completion of $S^0$,
 and one can show that the resulting map
 $\CM\Om^1_{G/S}=\colim_nE^0(\Sg^2\cp^\infty_{-n})\xra{}E^0$ is just
 $\res_{G/S}$.
\end{remark}

\begin{proposition}
 There are natural isomorphisms
 \begin{align*}
  (\coprod_n BU(n))_E &= M^+(G) = \Div^+(G)         \\
  BU_E                &= N^+(G) = N(G) = \Div_0(G)  \\
  (\ZtBU)_E           &= M(G)   = \Div(G)           \\
  (\ZtBU)^E           &= \Map_S(G,\MG).
 \end{align*}
 \index{MUa@$M^+(U)$} \index{NVa@$N^+(V)$}
 \index{MU@$M(U)$} \index{NV@$N(V)$}
 \index{Diva@$\Div^+(C)$} \index{Div@$\Div(C)$}
\end{proposition}
\begin{proof}
 This is well-known, and follows easily from
 Proposition~\ref{prop-BUn} and the remarks following
 Definition~\ref{defn-Div}.  The fourth statement follows from the
 third one by \idx{Cartier dual}ity.
\end{proof}

Next, recall that there is a ``complex reflection map''
$r\:S^1\tm\cp^{n-1}_+\xra{}U(n)$, where $r(z,L)$ has eigenvalue $z$ on
the line $L<\cplx^n$ and eigenvalue $1$ on $L^\perp$.  This gives an
unbased map $\cp^{n-1}\xra{}\Om U(n)$.  We can also fix a line
$L_0<\cplx^n$ and define $\ov{r}(z,L)=r(z,L)r(z,L_0)^{-1}$, giving a
map $\ov{r}\:\cp^{n-1}\xra{}\Om SU(n)$.  Moreover, the Bott
periodicity isomorphisms $\Om U=\Zh\tm BU$ and $\Om SU=BU$ give us
maps $\Om U(n)\xra{}\Zh\tm BU$ and $\Om SU(n)\xra{}BU$.  It is easy to
see that $(\cp^{n-1})_E$ is the divisor
$D_n=n[0]=\spec(\fps{E^0}{x}/x^n)$ on $G_E$ over $S_E$.
\begin{proposition}
 There are natural isomorphisms
 \begin{align*}
  (\Om U(n))_E  &= M(D_n)               \\
  (\Om SU(n))_E &= N(D_n)               \\
  (\Om U(n))^E  &= \Map_S(D_n,\MG)      \\
  (\Om SU(n))^E &= \text{BasedMap}_S(D_n,\MG).
 \end{align*}
 Under these identifications, the map $\Om U(n)\xra{}\Zh\tm BU$ gives
 the obvious map $M(D_n)\xra{}M(G_E)$ and so on.
\end{proposition}
\begin{proof}
 For the second statement, it is enough (by
 Remark~\ref{rem-symmetric-algebra}) to check that $E_*(\Om SU(n))$ is
 the symmetric algebra generated by the reduced $E$-homology of
 $\cp^{n-1}$.  This is well-known for ordinary homology, and it
 follows for all $E$ by a collapsing Atiyah-Hirzebruch spectral
 sequence.  See~\cite{ra:ccs,ra:nps} for more details.  The inclusion
 $S^1=U(1)\xra{}U(n)$ and the determinant map $\det\:U(n)\xra{}S^1$
 give a splitting $U(n)=S^1\tm SU(n)$ of spaces and thus
 $\Om U(n)=\Zh\tm\Om SU(n)$ of $H$-spaces and the first claim follows
 in turn using this.  The last two statements follow by Cartier
 duality. 
\end{proof}

\subsection{Cohomology of Abelian groups}
\label{subsec-abelian-groups}

Let $A$ be a compact Abelian Lie group, and write $A^*$ for the
character group $\Hom(A,S^1)$, which is a finitely generated discrete
Abelian group.  Let $G$ be an ordinary formal group over a base $S$.
For any point $s\in S(R)$ we write
$\Gm(G_s)=\FX_{\spec(R)}(\spec(R),G_s)$ for the associated group of
sections.  A coordinate gives a bijection between $\Gm(G_s)$ and
$\Nil(R)$, which becomes a homomorphism if we use an appropriate
formal group law to make $\Nil(R)$ a group.  We define a formal scheme
$\Hom(A^*,G)$ by
\[ \Hom(A^*,G)(R)=
    \{(s,\phi)\st s\in S(R) \text{ and } \phi\:A^*\xra{}\Gm(G_s)\}.
\]
(If $A^*$ is a direct sum of $r$ cyclic groups then this can be
identified with a closed formal subscheme of $G^r_S$ in an evident
way, which shows that it really is a scheme.)

\begin{proposition}\label{prop-BA}
 For any finite Abelian group $A$, there is a natural map
 $BA_E\xra{}\Hom(A^*,G)$.  \index{HomAG@$\Hom(A^*,G)$}  This is an
 isomorphism if $E$ is an exact $P(n)$-module for some $n$.
\end{proposition}
\begin{proof}
 An element $\al\in A^*=\Hom(A,S^1)$ gives a map $BA\xra{}BS^1$ of
 spaces and thus a map $BA_E\xra{}(BS^1)_E=G$ of formal groups over
 $S$.  One checks that the resulting map $A^*\xra{}\Ab\FX_S(BA_E,G)$
 is a homomorphism, so by adjointing things around we get a map
 $BA_E\xra{}\Hom(A^*,G)$.  If $A$ is a torus then $A^*\simeq\Zh^r$ and
 $BA_E=G^r=\Hom(A^*,G)$, so our map is an isomorphism.  Moreover, in
 this case $BA\simeq(\cpi)^r$ which is decent and thus tolerable for
 any $E$.  If $A=\Zh/m$ then there is a well-known way to identify
 $BA$ with the circle bundle in the line bundle $L^m$, where $L$ is
 the tautological bundle over $\cpi$.  This gives a long exact Gysin
 sequence
 \[ E^*BA \xla{} E^*\cpi \xla{[m](x)} E^*\cpi. \]
 The second map here is multiplication by $[m](x)$, which is the image
 of $x$ under the map $G\xra{\tm m}G$.  If this map is injective then
 the Gysin sequence is a short exact sequence and we have
 $E^0BA=E^0\cpi/[m](x)$, and we conclude easily that
 $\spf(E^0BA)=\ker(G\xra{m}G)=\Hom(A^*,G)$.  One can apply similar
 arguments to the skeleta $S^{2k+1}/(\Zh/m)$ of $BA$ and find that
 $\spf(E^0BA)=BA_E$.  

 In the case of two-periodic Morava $K$-theory we recover the
 well-known calculation showing that $K(n)^*BA$ is concentrated in
 even degrees for all $n$.  We also have $H^s(BA,\rat)=0$ for $s>0$ so
 Proposition~\ref{prop-rawiya-tolerable} tells us that $BA_E$ is
 tolerable for any $E$ that is an exact $P(n)$-module for any $n$.
 Moreover, it is easy to see that $[m](x)$ is not a zero-divisor in
 this case so $BA_E=\Hom(A^*,G)$.  We have just shown this when $A^*$
 is cyclic, but it follows easily for all $A$ by
 Proposition~\ref{prop-tolerable-product}.
\end{proof}

\subsection{Schemes associated to ring spectra}
\label{subsec-ring-spectra}

If $R$ is a commutative ring spectrum with a ring map $E\xra{}R$, we
have a scheme $\spec(\pi_0R)$ over $S_E$.  If $Z$ is a finite complex
we can take $R=F(Z_+,E)$ and we recover the case
$Z_E=\spec(E^0Z)=\spec(\pi_0R)$.  If $M$ is an arbitrary commutative
ring spectrum, we can take $R=E\Smash M$.  In this case we write
$M^E=\spec(E_0M)$\index{ME@$M^E$} for the resulting scheme.  If $Y$ is
a commutative $H$-space we can take $M=\Sgi Y_+$, and we write
$Y^E$\index{YE@$Y^E$} for $M^E=\spec(E_0Y)$ in this case.  If we have
a K\"unneth isomorphism $E_0Y^k=(E_0Y)^{\ot k}$ then $E_0Y$ is a Hopf
algebra, so $Y^E$ is a group scheme over $S$.  If $Y$ is \idx{decent}
then $E_0Y$ is a coalgebra with good basis.  In this case
Proposition~\ref{prop-cartier} applies, and we have a Cartier duality
$Y^E=D(Y_E)=\Hom_S(Y_E,\MG)$ and $Y_E=D(Y^E)=\Hom_S(Y^E,\MG)$.
\index{Cartier dual}\index{DG@$DG$}

If $\{R_\al\}$ is an inverse system of ring spectra as above, we have
a formal scheme $\colim_\al\spec(\pi_0R_\al)$.  If $Z_\al$ runs over
the finite subcomplexes of a CW complex $Z$, then the rings
$F(Z_{\al+},E)$ give an example of this, and the associated formal
scheme is just $Z_E$.  Another good example is to take the tower of
spectra $E/p^k$, where $E$ is an even periodic ring spectrum such that
$E^0$ is torsion-free.  More generally, if $E$ has suitable
Landweber exactness properties then we can smash $E$ with a
generalised Moore spectrum $S/I$ (see~\cite[Section 4]{host:mkl}, for
example) and get a new even periodic ring spectrum $E/I$, and then we
can consider a tower of these.  There are technicalities about the
existence of products on the spectra $E/I$, which we omit here.

\subsection{Homology of Thom spectra}
\label{subsec-thom-spectra}

Let $Z$ be a space equipped with a map $Z\xra{z}\ZtBU$, and let
$T(Z,z)$\index{TZz@$T(Z,z)$} be the associated \idx{Thom spectrum}.
It is well-known that $T$ is a functor from spaces over $\ZtBU$ to
spectra, which preserves homotopy pushouts.  Moreover, if $(Y,y)$ is
another space over $\ZtBU$ then we can use the addition on $\ZtBU$ to
make $(Y\tm Z,(y,z))$ into a space over $\ZtBU$ and we find that
$T(Y\tm Z,(y,z))=T(Y,y)\Smash T(Z,z)$.

The above construction really needs an actual map $Z\xra{z}\ZtBU$ and
not just a homotopy class.  However, we do have the following result.
\begin{lemma}
 If $Z$ is a decent space then the spectrum $T(Z,z)$ depends only on
 the homotopy class of $z$, up to canonical homotopy equivalence.
 Thus $T$ can be regarded as a functor from the homotopy category of
 decent spaces over $\ZtBU$ to spectra.  In particular, we can define
 $T(Z,V)$ when $V$ is a virtual bundle over $Z$.
\end{lemma}
\begin{proof}
 Suppose we have two homotopic maps $z_0,z_1\:Z\xra{}\ZtBU$.  We can
 then choose a map $w\:Z\tm I\xra{}\ZtBU$ such that $wj_0=z_0$ and
 $wj_1=z_1$, where $j_t(a)=(a,t)$.  The maps $j_t$ induce maps of
 spectra $T(Z,z_t)\xra{f_t}T(Z\tm I,w)$, and the Thom isomorphism
 theorem implies that these give equivalences in homology so they are
 weak equivalences.  We thus have a weak equivalence
 $f_1^{-1}f_0\:T(Z,z_0)\xra{}T(Z,z_1)$.  This much is true even when
 $Z$ is not decent.

 To see that our map is canonical when $Z$ is decent, note that
 $KU^*Z$ is concentrated in even degrees, so the space $F$ of
 unpointed maps from $Z$ to $\ZtBU$ has trivial odd-dimensional
 homotopy groups with respect to any basepoint.  We can think of $z_0$
 and $z_1$ as points of $F$, and $w$ as a path between them.  If $w'$
 is another path then then we can glue $w$ and $w'$ to get a map of
 $S^1$ to $F$, which can be extended to give a map $u\:D^2\xra{}F$
 because $\pi_1F=0$.  It follows that we have a commutative diagram as
 follows:
 \begin{diag}
  \node{T(Z,z_0)} 
  \arrow[2]{r,t}{f_0}
  \arrow{se}
  \arrow[2]{s,l}{f'_0}
  \node[2]{T(Z\tm I,w)}
  \arrow{sw}                    \\
  \node[2]{T(Z\tm D^2,u)}       \\
  \node{T(Z\tm I,w')}
  \arrow{ne}
  \node[2]{T(Z,z_1)}
  \arrow[2]{w,b}{f'_1}
  \arrow{nw}
  \arrow[2]{n,r}{f_1}
 \end{diag}
 It follows easily that $f_1^{-1}\circ f_0=(f'_1)^{-1}\circ f'_0$, as
 required. 
\end{proof}

A coordinate on $G_E$ is the same as a degree zero complex orientation
of $E$, which gives a multiplicative system of Thom classes for all
virtual complex bundles.  In particular, this gives isomorphisms
$E_*T(Y,y)\simeq E_*Y$, which are compatible in the evident way with
the isomorphisms $T(Y\tm Z,(y,z))=T(Y,y)\Smash T(Z,z)$.

If $Z\xra{z}\{n\}\tm BU(n)$ classifies an honest $n$-dimensional
bundle $V$ over $Z$ then we have $T(Z,z)=\Sgi Z^V$.  In particular,
the inclusion $\cpi=BU(1)\xra{}\{1\}\tm BU$ just gives the Thom
spectrum $\Sgi(\cpi)^L$, which is well-known to be the same as
$\Sgi\cpi$ (without a disjoint basepoint).

Now let $Z$ be a decent commutative $H$-space.  Let $z\:Z\xra{}\ZtBU$
be an $H$-map, and write $M=T(Z,z)$.  We note that addition gives a
map $(Z\tm Z,(z,z))\xra{}(Z,z)$ of spaces over $\ZtBU$ and thus a map
of spectra $M\Smash M\xra{}M$, which makes $M$ into a commutative ring
spectrum.  Similarly, the diagonal gives a map 
$(Z,z)\xra{}(Z\tm Z,(0,z))$ and thus a map 
$M\xra{\dl}\Sgi Z_+\Smash M$.  Finally, we consider the shearing map
$(a,b)\mapsto(a,a+b)$.  This is an isomorphism 
$(Z\tm Z,(z,z))\xra{}(Z\tm Z,(0,z))$ over $\ZtBU$, which gives an
isomorphism $M\Smash M\xra{}\Sgi Z_+\Smash M$ of spectra.

A choice of coordinate gives a Thom isomorphism $E_*M\simeq E_*Z$,
which shows that $E_*M$ is free and in even degrees.  For the moment
we just use this to show that we have K\"unneth isomorphisms, from
which we will recover a more natural statement about the relationship
between $E_*Z$ and $E_*M$.

Recall that we defined define $Z^E=\spec(E_0Z)=\spec(E_0\Sgi Z_+)$
(which is a commutative group scheme over $S=S_E$) and
$M^E=\spec(E_0M)$.  Our diagonal map $\dl$ gives an action of $Z^E$ on
$M^E$.  The shearing isomorphism $M\Smash M=\Sgi Z_+\Smash M$ shows
that the action and projection maps give an isomorphism
$Z^E\tm_SM^E\xra{}M^E\tm_SM^E$.

A choice of coordinate on $G$ gives an isomorphism $E_0M\simeq E_0Y$.
One can check (using the multiplicative properties of Thom classes)
that this is an isomorphism of $E_0Y$-comodule algebras, so it gives
an isomorphism $Y^E\simeq M^E$ of schemes, compatible with the action
of $Y^E$.  This means that $M^E$ is a trivialisable \idx{torsor} for
$Y^E$.

In the universal case $Y=\ZtBU$, this works out as follows.  As
mentioned previously, we have a map $\cpi=\{1\}\tm BU(1)\xra{}\ZtBU$,
and the Thom functor gives a map $\Sgi\cpi\xra{}MP$.  In particular,
the bottom cell gives a map $S^2=\cp^1\xra{}MP$, or an element
$u\in\pi_2MP$.  The inclusion $\{-1\}\xra{}\ZtBU$ also gives an
element of $\pi_{-2}MP$, which one checks is inverse to $u$.  Thus, a
ring map $E_0MP\xra{}R$ gives an $E_0$-algebra structure on $R$, and
an $E_0$-module map $\tilde{E}_0\cpi\xra{}R$, which sends
$\tilde{E}_0S^2$ into $R^\tm$.  In other words, it gives a point 
$s\in S_E(R)$ together with an element $y\in
R\hot_{E^0}\tilde{E}^0\cpi$.  We can identify
$R\hot_{E^0}\tilde{E}^0\cpi$ with the ideal of functions on $G_s$ that
vanish at zero, and the extra condition on the restriction to $S^2$
says that $y$ is a coordinate.  This gives a natural map
$MP^E\xra{}\Coord(G)$.  Well-known calculations show that $E_0MP$ is
the symmetric algebra over $E_0$ on $\tilde{E}_0\cpi$, with the bottom
class inverted.  This implies easily that the map
$MP^E\xra{}\Coord(G)$\index{CoordC@$\Coord(C)$} is an isomorphism.
Recall also that $(\ZtBU)^E=\Map(G,\MG)$.  Clearly, if $u\:G\xra{}\MG$
and $x$ is a coordinate on $G$, then the product $ux$ is again a
coordinate.  This gives an action of $\Map(G,\MG)$ on $\Coord(G)$,
which makes $\Coord(G)$ into a torsor over $\Map(G,\MG)$.  One can
check that this structure arises from our geometric coaction of
$\ZtBU$ on $MP$.

\subsection{Homology operations}
\label{subsec-operations}

Let $G$ be an ordinary formal group over $S$, and let $H$ be an
ordinary formal group over $T$.  Let $\pi_S$ and $\pi_T$ be the
projections from $S\tm T$ to $S$ and $T$.  We write $\Hom(G,H)$
\index{HomGH@$\Hom(G,H)$} for $\Hom_{S\tm T}(\pi_S^*G,\pi_T^*H)$,
which is a scheme over $S\tm T$ by Proposition~\ref{prop-exists-Hom}.
Recall that $\Hom(G,H)(R)$ is the set of triples $(s,t,u)$ where 
$s\in S(R)$ and $t\in T(R)$ and $u\:G_s\xra{}H_t$ is a map of formal
groups over $\spec(R)$.  We write $\Iso(G,H)(R)$
\index{IsoGH@$\Iso(G,H)$} for the subset of triples for which $u$ is
an isomorphism.  If we choose coordinates $x$ and $y$ on $G$ and $H$,
then for any $u$ we have $y(u(g))=\phi(x(g))$ for some power series
$\phi\in\fps{R}{t}$ with $\phi(0)=0$, and $u$ is an isomorphism if and
only if $\phi'(0)$ is invertible.  It follows that $\Iso(G,H)$ is an
open subscheme of $\Hom(G,H)$.

\begin{proposition}\label{prop-stable-operations}
 Let $E$ and $E'$ be even periodic ring spectra.  Then there is a
 natural map $S_{E\Smash E'}\xra{}\Iso(G_E,G_{E'})$ of schemes over
 $S_E\tm S_{E'}$.  This is an isomorphism if $E$ or $E'$ is 
 \idx{Landweber exact} over $MP$.
\end{proposition}
\begin{proof}
 We write $S'=S_{E'}$ and $G'=G_{E'}$.  The evident ring maps
 $E\xra{}E\Smash E'\xla{}E'$ give maps
 $S\xla{q}S_{E\Smash E'}\xra{q'}S'$, and pullback squares
 \begin{diag}
  \node{G} \arrow{s}
  \node{G_{E\Smash E'}} \arrow{w} \arrow{s} \arrow{e}
  \node{G'} \arrow{s} \\
  \node{S}
  \node{S_{E\Smash E'}} \arrow{w,b}{q} \arrow{e,b}{q'}
  \node{S'}
 \end{diag}
 This gives an isomorphism $v\:q^*G\xra{}(q')^*G'$.  Using
 this, we easily construct the required map.

 Now consider the case $E'=MP$, so that $S'=\FGL$.  Then
 $\Iso(G,G')(R)$ is the set of triples $(s,F,x)$, where $s\in S(R)$
 and $F$ is a formal group law over $R$ and
 $x\:G_s\xra{}\spec(R)\tm\haf^1$ is an isomorphism over $\spec(R)$
 such that $x(g+h)=F(x(g),x(h))$.  In other words, $x$ is a coordinate
 on $G_s$ and $F$ is the unique formal group law such that
 $x(g+h)=F(x(g),x(h))$.  Thus, we find that
 $\Iso(G,G')=\Coord(G)=MP^E=\spec(\pi_0MP)$ (see
 Section~\ref{subsec-thom-spectra}).  It follows after a comparison of
 definitions that our map $S_{E\Smash E'}\xra{}\Iso(G,G')$ is an
 isomorphism.  

 Now suppose that $E''$ is Landweber exact over $E'$, in the sense
 that there is a ring map $E'\xra{}E''$ which induces an isomorphism
 $E''_0\ot_{E'_0}E'_0Z=E''_0Z$ for all spectra $Z$.  We then find that
 $G''=G'\tm_{S'}S''$ and that 
 \[ S_{E\Smash E''} = S_{E\Smash E'}\tm_{S'}S''
    = \Iso(G,G')\tm_{S'}S'' = \Iso(G,G''), 
 \]
 as required.
\end{proof}

\begin{remark}
 If there are enough K\"unneth isomorphisms, then $E_0\Omi E'$ will be
 a Hopf ring over $E_0$ and thus the $*$-indecomposables
 $\Ind(E_0\Omi E')$ will be an algebra over $E_0$ using the circle
 product.  The procedure described in~\cite{kasttu:mkh} will then give
 a map $\spec(\Ind(E_0\Omi E'))\xra{}\Hom(G,G')$, which is an
 isomorphism in good cases.\index{HomGH@$\Hom(G,H)$}
\end{remark}

\begin{definition}
 Let $G$ and $G'$ be formal groups over $S$ and $S'$, respectively.  A
 \dfn{fibrewise isomorphism} from $G$ to $G'$ is a square of the form
 \begin{diag}
  \node{G} \arrow{e,t}{f} \arrow{s} \node{G'} \arrow{s} \\
  \node{S} \arrow{e,b}{g}           \node{S'}
 \end{diag}
 such that the induced map $G\xra{}f^*G'$ is an isomorphism of formal
 groups over $S$.
\end{definition}

\begin{definition}
 We write $\OFG$\index{OFG@$\OFG$} for the category of ordinary formal
 groups over affine schemes and fibrewise isomorphisms, and
 $\EPR$\index{EPR@$\EPR$} for the category of even periodic ring
 spectra.  We thus have a functor $\EPR^{\text{op}}\xra{}\OFG$ sending
 $E$ to $G_E$.  We write $\LOFG$\index{LOFG@$\LOFG$} for the
 subcategory of $\OFG$ consisting of Landweber exact formal groups,
 and $\LEPR$\index{LEPR@$\LEPR$} for the category of those $E$ for
 which $G_E$ is \idx{Landweber exact}.
\end{definition}

\begin{proposition}\label{prop-landweber-maps}
 If $E\in\EPR$ and $E'\in\LEPR$ then the natural map 
 \[ \EPR(E',E) \xra{} \OFG(G_E,G_{E'}) \]
 is an isomorphism.  Moreover, the functor
 $\LEPR^{\text{op}}\xra{}\LOFG$ is an equivalence of categories.
\end{proposition}
\begin{proof}
 Using~\cite[Proposition 2.12 and Corollary 2.14]{host:mkl}, we see
 that there is a cofibration $P\xra{}Q\xra{}E'\xra{}\Sg P$, in which
 $P$ and $Q$ are retracts of wedges of finite spectra with only even
 cells, and the connecting map $E'\xra{}\Sg P$ is phantom.  If $W$ is
 an even finite spectrum then we see from the Atiyah-Hirzebruch
 spectral sequence that $E_1W=0$ and $E_0W$ is projective over $E_0$
 and $[W,E]=\Hom(E_0W,E_0)$ and $[\Sg W,E]=0$.  It follows that all
 these things hold with $W$ replaced by $P$ or $Q$.  Using the
 cofibration we see that $E_1E'=0$, and there is a short exact
 sequence
 \[ E_0P \mra E_0Q \era E_0E'. \]
 Now consider the diagram
 \begin{diag}
  \node{0}                                      \arrow{e}
  \node{[E',E]}           \arrow{s,l}{\al_{E'}} \arrow{e}
  \node{[Q ,E]}           \arrow{s,l}{\al_Q}    \arrow{e}
  \node{[P ,E]}           \arrow{s,l}{\al_P}            \\
  \node{0}                                      \arrow{e}
  \node{\Hom(E_0E',E_0)}                        \arrow{e}
  \node{\Hom(E_0Q ,E_0)}                        \arrow{e}
  \node{\Hom(E_0P ,E_0).}
 \end{diag}
 The short exact sequence above implies that the bottom row is exact.
 The top row is exact because of our cofibration and the fact that
 $[\Sg P,E]=0$.  We have seen that $\al_P$ and $\al_Q$ are
 isomorphisms, and it follows that $\al_{E'}$ is an isomorphism.
 Thus, $[E',E]$ is the set of maps of $E_0$-modules from $E_0E'$ to
 $E_0$.  One can check that the ring maps $E'\xra{}E$ biject with the
 maps of $E_0$-algebras from $E_0E'$ to $E_0$ (using~\cite[Proposition
 2.19]{host:mkl}).  We see from
 Proposition~\ref{prop-stable-operations} that these maps biject with
 sections of $S_{E\Smash E'}=\Iso(G_E,G_{E'})$ over $S_E$, and these
 are easily seen to be the same as fibrewise isomorphisms from $G_E$
 to $G_{E'}$.  Thus $\EPR(E',E)=\OFG(G_E,G_{E'})$, as claimed.  This
 implies that the functor $\LEPR^{\text{op}}\xra{}\LOFG$ is full and
 faithful, so we need only check that it is essentially surjective.
 Suppose that $G$ is a Landweber exact ordinary formal group over an
 affine scheme $S$.  Define a graded ring $E_*$ by putting
 $E_{2k+1}=0$ and $E_{2k}=\om_{G/S}^{\ot k}$ for all $k\in\Zh$, so in
 particular $E_0=\OO_S$.  A choice of coordinate on $G$ gives a formal
 group law $F$ over $\OO_S=E_0$ and thus a map $S\xra{}\FGL$ or
 equivalently a map $u\:MP_0=\OO_{\FGL}\xra{}E_0$.  If $G_0=G_{MP}$ is
 the evident formal group over $\FGL$ then one sees from the
 construction that $S\tm_{\FGL}G_0=G$.  Given this, we see that our
 map $u$ extends to give a map $MP_*\xra{}E_*$.  We define a functor
 from spectra to graded Abelian groups by 
 \[ E_*Z = E_*\ot_{MP_*}MP_*Z = E_*\ot_{MU_*} MU_*Z, \]
 where we have used the map $MU\xra{}MP$ of ring spectra to regard
 $E_*$ as a module over $MU_*$.  One can also check that
 $E_0Z=E_0\ot_{MP_0}MP_0Z$.  The classical Landweber exact functor
 theorem implies that this is a homology theory, represented by a
 spectrum $E$.  The refinements given in~\cite[Section 2.1]{host:mkl}
 show that $E$ is unique up to canonical isomorphism, and that it
 admits a canonical commutative ring structure, making it an even
 periodic ring spectrum.  It is easy to check that
 $E^0\cpi=E_0\hot_{MP_0}MP^0\cpi$ and thus that
 $G_E=S\tm_{\FGL}G_0=G$, as required.
\end{proof}


\end{document}



Local Variables:
max-specpdl-size: 1500
max-lisp-eval-depth: 500
End: